\documentclass[11pt, reqno]{amsart}

\usepackage[T1]{fontenc}
\usepackage{lmodern}
\usepackage{microtype}

\usepackage{amsmath, amsfonts,amssymb, latexsym, mathrsfs}
\usepackage{amstext, amsthm, mathtools}
\usepackage[all]{xy}
\usepackage{fullpage} 
\usepackage{hyperref,pb-diagram,epic,eepic,verbatim,graphicx,graphics,epsfig,psfrag}
\usepackage{tikz-cd}
\usepackage{float}
\hypersetup{colorlinks}
\usepackage{color}
\usepackage{resizegather}
\definecolor{darkred}{rgb}{0.5,0,0}
\definecolor{darkgreen}{rgb}{0,0.5,0}
\definecolor{darkblue}{rgb}{0,0,0.5}
\hypersetup{ colorlinks,
linkcolor=darkblue,
filecolor=darkgreen,
urlcolor=darkred,
citecolor=darkblue }
\usepackage{cleveref}
\usepackage{ulem}
\usepackage{extarrows}

\usepackage{enotez}
\setenotez{backref=true}

\usepackage{pbox}

\makeatletter 
\def\makeCal#1{%
\expandafter\newcommand\csname c#1\endcsname{\mathcal{#1}}}

\def\makeBB#1{%
\expandafter\newcommand\csname b#1\endcsname{\mathbb{#1}}}

\def\makeFrak#1{%
\expandafter\newcommand\csname f#1\endcsname{\mathfrak{#1}}}

\def\makeScr#1{%
\expandafter\newcommand\csname s#1\endcsname{\mathscr{#1}}}

\count@=0
\loop
\advance\count@ 1
\edef\y{\@Alph\count@} 
\expandafter\makeCal\y
\expandafter\makeBB\y
\expandafter\makeFrak\y
\expandafter\makeScr\y
\ifnum\count@<26
\repeat

\newcommand{\bseries}[1]{ [\hspace{-0,5mm}[ {#1} ]\hspace{-0,5mm}] }

\def\makelowercaseFrak#1{%
\expandafter\newcommand\csname mf#1\endcsname{\mathfrak{#1}}}
\count@=0
\loop
\advance\count@ 1
\edef\y{\@alph\count@} 
\expandafter\makelowercaseFrak\y
\ifnum\count@<26
\repeat

\theoremstyle{plain}
\newtheorem{thm}{Theorem}[section]
\newtheorem{cor}[thm]{Corollary}
\newtheorem{lem}[thm]{Lemma}

\newtheorem{prop}[thm]{Proposition}

\newtheorem*{thm*}{Theorem}

\theoremstyle{definition}
\newtheorem{rem}[thm]{Remark}
\newtheorem{remark}[thm]{Remark}

\newtheorem{defn}[thm]{Definition}

\newtheorem{context}[thm]{Context}
\newtheorem{notn}[thm]{Notation}

\newtheorem{ex}[thm]{Example}
\newtheorem{example}[thm]{Example}

\newtheorem*{claim*}{Claim}

\newcommand{\op}[1]{\!\!\mathop{\rm ~#1}\nolimits}

\DeclareMathOperator{\Span}{Span}
\DeclareMathOperator{\GL}{GL}
\DeclareMathOperator{\Bun}{Bun}
\newcommand{\gen}{\mathrm{gen}}
\newcommand{\mrk}{\mathrm{mrk}}
\newcommand{\aux}{\mathrm{aux}}
\newcommand{\amp}{\mathrm{amp}}
\newcommand{\std}{\mathrm{std}}
\DeclareMathOperator{\ev}{ev}
\DeclareMathOperator{\id}{id}
\DeclareMathOperator{\triv}{triv}
\DeclareMathOperator{\Frac}{Frac}
\DeclareMathOperator{\Div}{Div}

\DeclareMathOperator{\wt}{wt}
\DeclareMathOperator{\maxwt}{maxWt}
\DeclareMathOperator{\Forget}{Forget}
\DeclareMathOperator{\Sect}{Sect}

\DeclareMathOperator{\Set}{Set}

\DeclareMathOperator{\DCoh}{D^b_{coh}}

\DeclareMathOperator{\Sym}{Sym}
\DeclareMathOperator{\Dqc}{D_{qc}}

\DeclareMathOperator{\ch}{ch}

\newcommand{\pt}{\op{pt}}

\DeclareMathOperator{\Perf}{Perf}

\DeclareMathOperator{\Filt}{Filt}
\DeclareMathOperator{\uMap}{\underline{Map}}
\newcommand{\MpC}{\mathcal{M}}
\DeclareMathOperator{\MpCrat}{\mathcal{M}^{rat}}
\newcommand{\uSpec}{\underline{\rm Spec}}
\DeclareMathOperator{\uAut}{\underline{\on{Aut}}}
\DeclareMathOperator{\uIso}{\underline{\on{Iso}}}
\DeclareMathOperator{\Grad}{Grad}
\DeclareMathOperator{\Tr}{Tr}

\newcommand{\iDeg}{{\mathscr Deg}}
\DeclareMathOperator{\Flag}{Flag}

\DeclareMathOperator{\colim}{colim}
\DeclareMathOperator{\cofib}{cofib}

\DeclareMathOperator{\QCoh}{QCoh}
\DeclareMathOperator{\rank}{rank}
\DeclareMathOperator{\gr}{gr}
\DeclareMathOperator{\Map}{\operatorname{Map}}
\DeclareMathOperator{\DF}{\mathbf{DF}}

\DeclareMathOperator{\Gr}{Gr}
\DeclareMathOperator{\GGW}{I}
\DeclareMathOperator{\marked}{\underline{p}}
\DeclareMathOperator{\proj}{proj}

\newcommand{\closedpt}{\mathfrak{o}}

\newcommand\D{\mathcal{D}}
\newcommand{\XX}{\mathfrak{X}}
\newcommand{\dash}{{\operatorname{-}}}

\newcommand{\on}{\operatorname}

\newcommand{\Quot}{\on{Quot}}

\newcommand{\rk}{\on{rk}}

\newcommand{\dual}{\vee}

\newcommand{\Proj}{\on{Proj}}

\newcommand{\Aut}{ \on{Aut} } 

\newcommand{\Iso}{ \on{Iso} }

\newcommand{\Ad}{ \on{Ad} }

\newcommand{\Rep}{\on{Rep}}
\newcommand{\Hom}{ \on{Hom}}
\newcommand{\Mor}{\on{Mor}}

\renewcommand{\ker}{ \on{ker}}

\newcommand{\im}{ \on{Im}}

\newcommand{\ob}{\on{ob}}

\newcommand{\Spec}{\on{Spec}}
\newcommand{\Sch}{\on{Sch}}
\newcommand{\Aff}{\on{Aff}}

\newcommand{\Pic}{\on{Pic}}

\newcommand\dirac{/\kern-1.2ex\partial} 
\newcommand\qu{/\kern-.7ex/} 
\newcommand\lqu{\backslash \kern-.7ex \backslash} 

\newcommand\dr{r_+ \kern-.7ex - \kern-.7ex r_-}
 
\tikzset{
  symbol/.style={
    draw=none,
    every to/.append style={
      edge node={node [sloped, allow upside down, auto=false]{$#1$}}}
  }
}

    \newtheoremstyle{TheoremNum}
        {\topsep}{\topsep}              
        {\itshape}                      
        {}                              
        {\bfseries}                     
        {.}                             
        { }                             
        {\thmname{#1}\thmnote{ \bfseries #3}}
    \theoremstyle{TheoremNum}

    \usepackage{trimclip}
    \makeatletter
\newcommand{\joinsubset}{%
  \mathrel{\hbox{\clipbox{{0.45\width} 0 0 0}{$\m@th\subset$}}}%
}
\newcommand{\xsubset}[2][]{%
  \ext@arrow 3095 \subsetfill@{#1}{#2}%
}
\newcommand{\subsetfill@}{\arrowfill@\subset\joinsubset\joinsubset}

\author{Daniel Halpern-Leistner}
\address{Mathematics Department,
Cornell University, 310 Malott Hall, Cornell University, Ithaca, New York 14853,
U.S.A.}
\email{daniel.hl@cornell.edu}

\author{Andres Fernandez Herrero}
\address{Department of Mathematics,
Columbia University, 3990 Broadway, New York, NY 10027,
U.S.A.} 
\email{af3358@columbia.edu}

\begin{document}
\title{The structure of the moduli of gauged maps from a smooth curve}

\begin{abstract}
For a reductive group $G$, Harder-Narasimhan theory gives a structure theorem for principal $G$ bundles on a smooth projective curve $C$. A bundle is either semistable, or it admits a canonical parabolic reduction whose associated Levi bundle is semistable. We extend this structure theorem by constructing a $\Theta$-stratification of the moduli stack of gauged maps from $C$ to a projective-over-affine $G$-variety $X$. The open stratum coincides with the previously studied moduli of Mundet semistable maps, and in special cases coincides with the moduli of stable quasi-maps. As an application of the stratification, we provide a formula for K-theoretic gauged Gromov-Witten invariants when $X$ is an arbitrary linear representation of $G$. This can be viewed as a generalization of the Verlinde formula for moduli spaces of decorated principal bundles. We establish our main technical results for smooth families of curves over an arbitrary Noetherian base. Our proof develops an infinite-dimensional analog of geometric invariant theory and applies the theory of optimization on degeneration fans.
\end{abstract}

\maketitle

\setcounter{secnumdepth}{2}
\setcounter{tocdepth}{2}

\tableofcontents

\section{Introduction}
Since the work of Weil \cite{weil-fibre-spaces}, the moduli of vector bundles on an algebraic curve has been one of the most important and most studied moduli problems in mathematics. It has connections to number theory, representation theory, algebraic geometry, differential geometry, and mathematical physics. One of the early achievements of geometric invariant theory (GIT) was the introduction of a notion of ``semistability'' for vector bundles, and the construction of a projective moduli space of semistable bundles \cite{mumford_bundles, seshadri-unitary}. This moduli problem has been generalized in innumerable ways, but one of the most basic is the \emph{moduli of gauged maps}.

For simplicity, we will restrict our discussion here to a fixed smooth projective curve $C$ over an algebraically closed field of characteristic $0$, a connected reductive group $G$, and an affine $G$-scheme $X$ of finite type. At the end of the introduction, we will state our main results in full generality, which apply over arbitrary Noetherian base schemes and to $X$ that are projective over an affine $G$-scheme.

The moduli of gauged maps $\MpC(X/G)$ parameterizes a principal $G$-bundle $E$ over $C$ along with a section of the associated fiber bundle $E(X) \to C$. Succinctly, $\MpC(X/G)$ is the mapping stack $\Map(C,X/G)$. The moduli of vector bundles corresponds to $G=\GL_n$ and $X=\pt$, and the most studied situation beyond that is where $X$ is a linear representation of $G$. Specific examples, such as the moduli of (twisted) Higgs bundles on $C$ or the moduli of Bradlow pairs, have been of considerable interest in the fields mentioned above.



\begin{ex}[Noncommutative Hurwitz spaces]\label{E:locally_free_algebras}
Consider the moduli problem of a locally free sheaf of finite rank $\cA$ on $C$ together with the structure of a unital associative $\cO_C$-algebra. If $\cA$ is commutative, this is equivalent to a pure one dimensional scheme $C' = \Spec_C(\cA)$ with a finite morphism $C' \to C$, so in general we can think of $\cA$ as a ``noncommutative'' curve with a finite map to $C$. For any finite dimensional vector space $V$, there is an affine subvariety $\rm{Assoc}(V) \hookrightarrow \Hom(V \otimes V, V) \times V$ parameterizing points $(m,e)$ such that $m$ is an associative product on $V$ and $e$ is a unit. Then the moduli of locally free algebras over $C$ is equivalent to $\MpC(\rm{Assoc}(V)/\GL(V))$. There is a stability condition that, for commutative $\cA$, amounts to $\cA$ being reduced, and the resulting moduli space is a compactification of the Hurwitz space of $C$. The analogous stability condition for general associative $\cA$ is more subtle, and the resulting moduli space can be regarded as a noncommutative Hurwitz space of $C$.
\end{ex}

The moduli of gauged maps first arose in the study of the 2D generalized vortex equations associated to a hamiltonian group action on a K\"{a}hler manifold. Mundet i Riera introduced an algebraic notion of semistable points of $\MpC(X/G)$ and established a Kobayashi-Hitchin correspondence identifying them with solutions of the vortex equation \cite{mundet-semistability}, which have moduli spaces in the sense of differential geometry. Later, these moduli spaces were constructed as algebraic varieties by Schmitt using GIT \cite{schmitt-decorated-principal-bundles}.

The original motivation was to define \emph{gauged Gromov-Witten invariants} for hamiltonian group actions \cite{cgs-j-holomorphic, cgms-gw-hamiltonian}, which are integrals (in cohomology or $K$-theory) of tautological classes on these moduli spaces. The hope is that these invariants can be used to compute the (genus zero) Gromov-Witten invariants of a GIT quotient $X^{\rm{ss}}/\!/G$. Below, we recall two well-developed techniques to recover Gromov-Witten invariants from gauged Gromov-Witten invariants. These techniques are only effective, however, if one can actually compute the gauged Gromov-Witten invariants. In this paper, we give an explicit formula for the $K$-theoretic gauged Gromov-Witten invariants whenever $\pi_1(G)$ is free, and $G$ acts linearly on $X=\bA^n$. (See \Cref{T:index_formula}.)


\subsubsection{The stratification}

Returning to vector bundles on a curve, a result of Harder and Narasimhan \cite{harder-narasimhan} states that every bundle $E$ has a unique filtration $E_p \subset \cdots \subset E_0 = E$ whose associated subquotients $E_i / E_{i+1}$ are semistable bundles with increasing slopes. (The filtration is trivial if $E$ is semistable.) This stratifies the moduli problem into locally closed pieces \cite{shatz-decomposition-specialization, nitsure-schematichn}, where each stratum parameterizes bundles such that $E_i/E_{i+1}$ has fixed rank and degree for $i=0,\ldots,p$. In addition to its many applications, we view this stratification as a structure of fundamental interest that essentially completes the classification of vector bundles on $C$: Every bundle is uniquely classified by a point on a moduli space of semistable (graded) bundles, along with the extension data describing how to construct the Harder-Narasimhan filtration from its graded pieces.


The theory of $\Theta$-stratifications \cite{hl_instability} provides a framework for generalizing the Harder-Narasimhan stratification to potentially any moduli problem that can be represented by an algebraic stack $\cX$ (satisfying mild hypotheses). For any field $k$, let $\Theta_k$ be the quotient stack $\bA^1_k / \bG_m$. The key idea is that a ($\bZ$-weighted) \emph{filtration} of a point $x\in \cX$ is a morphism $f : \Theta_k \to \cX$ along with an isomorphism $f(1) \cong x$, and the \emph{associated graded point} is $f(0)$, which by construction has a canonical $\bG_m$-action $(\bG_m)_k \to \Aut_\cX(f(0))$. A \emph{numerical invariant} $\mu$ is a function on filtrations satisfying certain conditions (\Cref{defn: numerical invariant}). Given a numerical invariant, we define a point $x \in \cX$ to be \emph{unstable} if there is a filtration $f$ with $f(1) \cong x$ and $\mu(f)>0$, and semistable otherwise. An HN filtration of an unstable point $x$ is a maximizer of $\mu(f)$ subject to the constraint $f(1) \cong x$.

A $\Theta$-stratification of $\cX$ is a stratification by locally closed substacks, where each stratum parameterizes a point along with an HN filtration. More precisely, the strata are canonically identified with open substacks of the stack $\Filt(\cX):=\Map(\Theta,\cX)$ of filtered points of $\cX$. In addition, each stratum deformation retracts onto it \emph{center}, which is an open substack of $\Grad(\cX):=\Map(B\bG_m,\cX)$ parameterizing semistable graded points.

In \Cref{subsection: numerical invariants} and \Cref{defn: l-chi} below, we introduce a line bundle $\cL(\chi)$ on $\MpC(X/G)$, which depends on a choice of representation and character $\chi$ of $G$. We also introduce a ``norm'' on filtrations, which depends on a Weyl-invariant norm $b$ on the co-weight lattice of $G$. On any filtration $f$, the value of the numerical invariant is $\mu(f) = -\wt(\cL(\chi)_{f(0)}) / \|f\|_b \in \bR$, where $\wt$ denotes the weight of the fiber at $f(0)$ with respect to the canonical $\bG_m$-automorphism. This is a direct generalization of the normalized Hilbert-Mumford numerical invariant studied in geometric invariant theory \cite{mumford-git, dolgachev-hu}. For the special case $\Bun_G(C):=\MpC(\pt/G)$, the line bundle $\cL(0)$ formed using the adjoint representation of $G$ and the trivial character $\chi=0$ recovers the usual notion of semistability of $G$-bundles \cite{ramanathan-stable}.


\begin{thm}[cf.\Cref{thm: theta stratification in general}+\Cref{thm: moduli space in the case of centrally contractible}] \label{thm: main thm intro}
The numerical invariant $\mu$ defines a $\Theta$-stratification of $\MpC(X/G)$. Each connected component of the semistable locus admits a good moduli space that is proper over $X/\!/G := \Spec_S(\cO_X^G)$, via a generalized Hitchin morphism. 

Furthermore, for each stratum there is a one-parameter subgroup $\lambda$ of $G$ and a character $\chi'$ of the associated Levi subgroup $L_{\lambda}\subset G$ such that the center of the stratum is a connected component of the $\cL(\chi')$-semistable locus of $\MpC(X^{\lambda =0}/L_{\lambda})$, where $X^{\lambda=0}$ is the fixed locus of $\lambda$.
\end{thm}

The main contribution of this paper is to construct the $\Theta$-stratification, but our results also give a new construction of the moduli space that does not use geometric invariant theory. The proof of \Cref{thm: main thm intro} uses the collection of criteria on the numerical invariant $\mu$ known as ``intrinsic GIT,'' developed in \cite[Sect. 5]{hl_instability} and recalled in \Cref{thm: theta stability paper theorem}. They are:
\begin{enumerate}
    \item Strict $\Theta$-monotonicity and strict S-monotonicity (\Cref{defn: strictly theta monotone and STR monotone}).
    \item HN boundedness (\Cref{defn: HN boundedness}).
\end{enumerate}
The methods we introduce to verify strict monotonicity in \Cref{section: monotonicity} and HN-boundedness in \Cref{section: hn boundedness} constitute the main innovations of this paper.

\begin{ex}[Moduli of semistable noncommutative Hurwitz covers] \label{ex: hurwitz 2}
Continuing \Cref{E:locally_free_algebras}, one can identify a morphism $f : \Theta \to \MpC(\rm{Assoc}(V)/\GL(V))$ with a $\bZ$-weighted filtration of vector bundles $\cdots \subset \cA_{w+1} \subset \cA_w \subset \cdots \subset \cA$ such that $\cA_w = 0$ for $w \gg 0$, $\cA_w = \cA$ for $w\ll 0$, $\cA_i \cdot \cA_j \subset \cA_{i+j}$ for all $i,j$, and the unit section satisfies $e_{\cA} \in \Gamma(C,\cA_0)$. The Harder-Narasimhan filtration of $\cA$ as a vector bundle need not lift to a $\bZ$-weighted filtration of $\cA$ as an algebra. Furthermore, filtrations of $\cA$ as an algebra can not be constructed iteratively from $2$-step filtrations, as in the case of vector bundles, so without the general machinery that we use below, it is not clear how to begin searching for a maximizer of $\mu$.\endnote{The usual method of constructing the HN filtration is to show the existence of a maximal subbundle of maximal slope $F \subset E$, and then to combine this inductively with the HN filtration for $E/F$.} Nevertheless, \Cref{thm: main thm intro} says that every unstable locally free algebra over $C$ admits a canonical $\mathbb{G}_m$-equivariant degeneration that maximizes $\mu$, and the limiting algebra will be semistable as a graded algebra for a twisted numerical invariant $\mu'$. This stratifies the moduli stack of \emph{all} noncommutative Hurwitz covers into pieces that retract to stacks with good moduli spaces.

\end{ex}

Let us explain the recursive nature of the stratification of \Cref{thm: main thm intro} in more detail. Fix a choice $T \subset B \subset G$ of maximal split torus $T$ and Borel subgroup $B\supset T$. The strata in $\MpC(X/G)$ are indexed by certain pairs $\nu= (d, \lambda)$ (see \Cref{defn: chi indexing datum}) consisting of a $B$-dominant cocharacter $\lambda \in X_*(T)$ and $d \in X_*(T)_{\mathbb{Q}}$. The only $d$ that appear lie in $\frac{1}{H} X_\ast(T) \subset X_\ast(T)_\bQ$ for some integer $H$, and there are finitely many $\lambda$ for each $d$. When $G=\GL_n$, $d$ encodes the slopes of the associated graded pieces of the filtration of the underlying vector bundle, and $\lambda$ encodes a canonical choice of weights of this filtration. A posteriori, when $X=\pt$, the canonical weights of the Harder-Narasimhan filtration are just the slopes of the associated graded bundles, but this does not hold in general.

The stratum $\cS_{\nu}$ corresponding to $\nu = (d,\lambda)$ is isomorphic to an open substack of $\MpC(X^{\lambda \geq 0}/P_{\lambda})_d$, where $P_{\lambda} \subset G$ is the associated parabolic subgroup, $X^{\lambda \geq0}\subset X$ is the corresponding concentrator subscheme, and the subscript $d$ indicates that we restrict to the substack of $P_{\lambda}$-bundles of degree $d$. The stratum admits a morphism to its center $\cZ_{\nu} \subset \MpC(X^{\lambda=0}/L_{\lambda})_d$, which is the $\mu'$-semistable locus in $\MpC(X^{\lambda=0}/L_{\lambda})_d$, where $\mu'$ is the numerical invariant for a shifted stability character of $L_{\lambda}$ (cf. \Cref{subsection: centers of unstable strata}). We summarize the previous discussion with the following diagram, which is Cartesian with respect to the downward-pointing arrows:
\[ \begin{tikzcd}
  & \cS_{\nu} \ar[d] \ar[r, symbol = \xsubset{\;\;\;\; \; \; \; } ] & \MpC(X^{\lambda \geq 0}/P_{\lambda})_d \ar[d] \ar[r] & \MpC(X/G) \\ \cZ_{\nu} \ar[ur] \ar[r, symbol = \xlongequal{\;\;\;\;}] & \MpC(X^{\lambda =0}/L_{\lambda})_d^{\mu'\dash {\rm ss}} \ar[r, symbol = \xsubset{\;\;\;\; \;}] & \MpC(X^{\lambda =0}/L_{\lambda})_d \ar[u, bend right] &
\end{tikzcd}\]

In the process of establishing this description of the stratification, we prove certain bounds on the degrees of semistable maps $C \to X/G$ which are analogous to Bogomolov-Gieseker inequalities, and we show that for $\chi$ sufficiently large, any $\cL(\chi)$-semistable map $C \to X/G$ must map the generic point of $C$ to the $\chi$-semistable locus of $X$ in the sense of GIT (\Cref{prop: constraints degree semistable locus}).


\subsection{Gauged Gromov-Witten theory}

Suppose that $X$ is a linear representation of $G$, and we consider maps $C \to X/G$ of a fixed degree $d$. For a certain kind of tautological perfect complex $F$ on $\MpC(X/G)$, called an \emph{Atiyah-Bott} complex (\Cref{defn: atiyah-bott complexes}), the associated $K$-theoretic gauged Gromov-Witten invariant is defined as the $K$-theoretic integral
\[
I_d^\chi(X/G, F) := \sum_i (-1)^i \dim \mathbb{H}^i(\MpC(X/G)^{\chi\dash\rm{ss}}_d, \cL \otimes F),
\]
where $\chi\dash\rm{ss}$ denotes the $\cL(\chi)$-semistable locus, and we assume that $X/\!/G = \pt$ so that the Euler characteristic is well-defined. More generally, when $X/\!/G$ is non-compact, \Cref{thm: main thm intro} implies that $\mathbb{H}^i(\MpC(X/G)^{\chi\dash\rm{ss}}_d, \cL \otimes F)$ is a coherent sheaf on $X/\!/G$ that is equivariant for the scaling action of $\bG_m$, and we consider instead the graded Euler characteristic of $R\Gamma(\MpC(X/G)^{\chi\dash\rm{ss}}_d, \cL \otimes F)$.

\begin{rem}
An important technical note is that $\MpC(X/G)$ is not smooth, but it has a canonical enhancement to a quasi-smooth derived algebraic stack. Whenever we discuss cohomology of complexes on $\MpC(X/G)$, we are using $\MpC(X/G)$ to denote this derived enhancement. Our formulas compute the index of $\cL \otimes F$ over this derived stack, which is equal to the index of $\cL \otimes F \otimes \cO^{vir}$ on the classical stack, where $\cO^{vir} \in \DCoh(\MpC(X/G))$ is a virtual structure sheaf. $\cO^{\rm vir}$ is the analogue in $K$-theory of the virtual fundamental class in cohomology.
\end{rem}

Our approach to computing these invariants extends earlier work on the Verlinde formula, which computes the dimension of the space of sections of $\cL$ on $\Bun_G(C)$. That this space (and in fact the total cohomology of $\cL$) is finite dimensional at all is somewhat remarkable, and can be seen by identifying it with a space of conformal blocks \cite{beauville-conformal}. A different approach was developed in \cite{teleman-quantization, teleman-woodward}, using the Harder-Narasimhan stratification of $\Bun_G(C)$. In fact, \cite{teleman-woodward} gives a vast generalization of the Verlinde formula, allowing one to compute the ``index'' $\chi(\Bun_G(C),\cL \otimes F)$ for any Atiyah-Bott complex $F$. Furthermore, a non-abelian localization theorem associated to the Harder-Narasimhan stratification provides an explicit comparison between the index over $\Bun_G(C)$ and over $\Bun_G(C)^{\rm ss}$.



In order to generalize this approach from $\Bun_G(C)=\MpC(BG)$ to $\MpC(X/G)$ for an arbitrary linear representation $X$, we combine the $\Theta$-stratification of \Cref{thm: main thm intro} with the \emph{virtual non-abelian localization formula} of \cite{verlinde, HL-D-equivalence} to show:
\begin{thm}[\Cref{T:main_index}]\label{thm: main thm index formula}
After replacing $\cL$ with $\cL^q$ for some integer $q>0$, we have that for any Atiyah-Bott complex $F$:
\begin{enumerate}
    \item For any $m$, $H^i R\Gamma(\MpC(X/G)_d, \cL(\chi)^m \otimes F)$ is a coherent graded $\cO_X^G$-module for all $i$ and vanishes for $|i|\gg 0$, and the resulting graded Euler characteristic admits an explicit formula in terms of the Teleman-Woodward index theorem. (See \Cref{T:index_formula}.) \\
    
    \item For $m \gg 0$, the following restriction map is an isomorphism,
    \[
    R\Gamma(\MpC(X/G)_d, \cL(\chi)^m \otimes F) \to R\Gamma(\MpC(X/G)_d^{\chi\dash\rm{ss}},\cL(\chi)^m \otimes F).\\
    \]
\end{enumerate}
\end{thm}
Furthermore, in \Cref{P:recursive_index_formula} we use virtual non-abelian localization to prove a recursive formula for the gauged Gromov-Witten invariants $I_d^\chi(X/G,F)$ in terms of the index over the whole stack $I_d(X/G,F) := \sum_i (-1)^i \dim \mathbb{H}^i(\MpC(X/G)_d,\cL \otimes F)$ and correction terms that are gauged Gromov-Witten invariants for subspaces of $X$ with respect to Levi subgroups of $G$. One might question the value of this somewhat complicated recursive formula in light of the simpler formula for $I_d(X/G,F)$ in \Cref{T:main_index}(1). It is $I_d^\chi(X/G,F)$, however, that can be used to compute the Gromov-Witten invariants of $X^{\chi \dash \rm{ss}}/\!/G$. 


When $G$ acts with finite stabilizers on the $\chi$-semistable locus $X^{\chi\dash \rm{ss}} \subset X$, and when $\|\chi\|$ is large with respect to $\|d\|$, the stack $\MpC(X/G)_d^{\chi\dash\rm{ss}}$ is often referred to as the moduli of stable quasi-maps. These moduli problems have played an important role in enumerative geometry and geometric representation theory \cite{braverman-quasi-maps}. In this special case, the stable maps $C \to X/G$ are precisely those that take the generic point of $C$ to $X^{\chi\dash\rm{ss}}$. There are at least two well-developed methods for relating quasi-maps invariants, which are integrals of tautological classes over the stacks of stable quasi-maps, with the Gromov-Witten invariants of the GIT quotient $X^{\chi\dash\rm{ss}} /\!/G$:
\begin{itemize}

\item The ``adiabatic limit theorem'' \cite{gs-adiabatic, woodwardIII, gonzalez2020wallcrossing} relates a certain generating function for gauged Gromov-Witten invariants with an analogous generating function for graph Gromov-Witten invariants, which count stable maps $\tilde{C} \to C \times (X^{\chi\dash\rm{ss}}/\!/G)$ whose composition $\tilde{C} \to C$ has degree $1$. Although $C$ has arbitrary genus, it is a fixed smooth curve, so these graph Gromov-Witten invariants morally only capture genus $0$ Gromov-Witten theory.\\

\item Givental's Mirror theorem \cite{givental-equivariant, Givental1998} and its generalizations \cite{lian-liu-yau-mirror,iritani-quantum-D-modules,clct-quantum-weighted-projective}, relate the $I$-function, which is a generating function for ($\bG_m$-localized) quasi-maps invariants on $\bP^1$, with the $J$-function, which is a generating function for Gromov-Witten invariants of $X^{\chi\dash\rm{ss}}$. In this setting, quasi-maps invariants \cite{quasimap-to-git} and the wall-crossing formulas used to relate them to Gromov-Witten invariants \cite{ciocan-fontanine-kim-wall-crossing, ciocan-fontanine-kim-higher-genus, ciocan-fontanine-kim-ihes,clader-janda-ruan} have been established for moduli problems that allow the source curve to vary in arbitrary nodal families, so in principle it can be applied in arbitrary genus, whereas gauged Gromov-Witten invariants have only been defined for a fixed smooth curve.\\
\end{itemize}



We hope that this work stimulates further investigation in enumerative geometry. For one thing, the generating functions in our formulas differ from the generating functions typically studied in enumerative geometry -- the $I$ function and the gauged Gromov-Witten potential. It is an interesting problem to develop a good conceptual framework for packaging our recursive formulas and relating our generating functions to both of these more commonly studied generating functions.

\vspace{10pt}
In the remainder of the introduction, we discuss the proof of \Cref{thm: main thm intro} and its generalization.

\subsection{Monotonicity via infinite dimensional GIT}

Codimension-2 filling conditions for an algebraic stack $\cX$, known as $\Theta$-reductivity and $S$-completeness, have played an important role in recent developments in algebraic moduli theory \cite{hl_instability, AHLH}. More precisely, there are certain regular affine $2$-dimensional schemes $Y$ with a $\bG_m$-action such that $\{\closedpt\} = Y^{\bG_m}$ is the unique closed orbit. $S$-completeness and $\Theta$-reductivity require that any morphism $f : (Y \setminus \closedpt)/\bG_m \to \cX$ extends uniquely to a morphism $\tilde{f} : Y/\bG_m \to \cX$. Unfortunately, the stacks $\MpC(X/G)$ are not $\Theta$-reductive or $S$-complete in general (cf. \cite{torsion-freepaper}).

The notions of strict monotonicity introduced in \cite[Sect. 5]{hl_instability}, which we recall in \Cref{subsection: theta-sratifications}, are relaxations of these filling conditions that still suffice to construct $\Theta$-stratifications and moduli spaces. They require an extension $\tilde{f} : \Sigma/\bG_m \to \cX$ to exist for some equivariant blowup $\Sigma \to Y$, such that the value of the numerical invariant satisfies certain mononicity conditions on the $\bG_m$-fixed points in the special fiber $\Sigma_\closedpt$. (See \Cref{defn: strictly theta monotone and STR monotone}.)

\begin{example}[Why monotonicity holds in projective GIT]\label{E:infinite_dimensional_GIT}
Let $\cX$ be a quotient stack $Z/G$ for a $G$-projective scheme $Z$ over the ground field, and consider the Hilbert-Mumford numerical invariant determined by a $G$-equivariant ample bundle $\cL$ on $Z$. Suppose we are given a morphism $f : (Y \setminus \closedpt) / \bG_m \to Z/G$. $Z/G$ is not $\Theta$-reductive or $S$-complete in general, but the stack $BG$ is, so the morphism $(Y \setminus \closedpt)/\bG_m \to BG$ extends uniquely to a morphism defined on $Y/\bG_m$. In other words, there is a unique dashed arrow filling the diagram
\[
\begin{tikzcd}
 (Y \setminus \closedpt) / \bG_m \ar[d, symbol= \hookrightarrow] \ar[r] & Z/G \ar[d, "\varphi"] \\  Y/\bG_m  \ar[r, dashed] &  BG
\end{tikzcd}.
\]
The original morphism $(Y \setminus \closedpt)/\bG_m \to Z/G$ defines a section of the projective morphism $Z/G \times_{BG} Y \to Y$ over $Y \setminus \closedpt$. The closure of this section is a surface $\Sigma$ with a projective $\bG_m$-equivariant morphism to $Y$ that fits into a commutative diagram
\[
\begin{tikzcd}
 (Y\setminus \closedpt)/\bG_m \ar[r, symbol=\subset] & \Sigma/ \mathbb{G}_m \ar[d] \ar[r] & Z/G \ar[d, "\varphi"] \\  & Y/\bG_m \ar[r]  &  BG
\end{tikzcd}
\]
Because $\cL$ is relatively ample for the projective morphism $\varphi : Z/G \to BG$, its pullback to the special fiber $\Sigma_{\closedpt}$ will be ample. This implies that the weights of $\cL$ at the fixed points of $\Sigma_\closedpt$ are strictly monotone increasing, which is the key to establishing strict monotonicity for the Hilbert-Mumford numerical invariant.
\end{example}

Our strategy for establishing strict monotonicity is analogous to that in \Cref{E:infinite_dimensional_GIT}, but it is ``infinite dimensional'' in the sense that it leaves the world of algebraic stacks. Although we do not give a general definition of infinite dimensional GIT, the approach applies to a numerical invariant $\mu$ on an algebraic stack $\cX$ equipped with the following data (roughly speaking):
\begin{enumerate}
    \item a potentially non-algebraic stack $\cG$ that satisfies the codimension-$2$ filling conditions above;
    \item a morphism $\varphi : \cX \to \cG$ that is representable by ind-projective ind-schemes;
    \item a $\varphi$-ample (formal) invertible sheaf $\cL$ on $\cX$ and a norm on filtrations in $\cG$ such that \[\mu(f) = \wt(\cL_{f(0)})/ \|\varphi \circ f\|.\]
\end{enumerate}
This set up is inspired by ideas from geometric representation theory \cite{barlev, lurie-gaitsgory-tamagawabook} as well as the technique of infinite dimensional symplectic reduction in differential geometry \cite{atiyah-bott-yangmills, donaldson-surfaces}. For technical reasons, we develop infinite dimensional GIT for $\MpC(X/G)$ only when $G=\GL_N \times H$, where $H$ is an extension of a finite \'etale group by a torus, and we reduce monotonicity in general to this special case.

\subsubsection{Infinite dimensional GIT for gauged maps} In our setting, $\cG$ is the non-algebraic stack $\MpCrat(X/G)$ of ``rational maps'' from $C$ to $X/G$. Its points consist of a $G$-bundle $E$ over $C$, an effective divisor $D \subset C$, and a section $s$ over $C\setminus D$ of the associated fiber bundle $E(X) \to C$. An isomorphism between two points $(E,D,s)$ and $(E',D',s')$ in $\MpCrat(X/G)$ is an equality $D =D'$ and an isomorphism of restrictions $E|_{C\setminus D} \simeq E'|_{C \setminus D'}$ that maps the section $s$ to $s'$. For any base scheme $T$ and relative effective Cartier divisor $D \subset C_T$, we define a morphism
\[ \varphi_D: \MpC(X/G)_T \to \MpCrat(X/G)_T, \; \; \; (E,s) \mapsto (E, D, s|_{C \setminus D}), \]
that plays the role of the morphism $Z/G \to BG$ in projective GIT. The fiber of $\varphi_D$ over a $T$-point $(E,D,s)$, which we denote $\Gr_{X/G}^{E,D,s}$, is an ind-projective-ind-scheme over $T$, a closed sub-ind-scheme of the Beilinson-Drinfeld grassmannian for $G$.

The strategy of infinite dimensional GIT is to show that for any morphism $f : (Y \setminus \closedpt) / \bG_m \to \MpC(X/G)$, one can choose a $\bG_m$-equivariant relative effective Cartier divisor $D \hookrightarrow C_Y$ such that the morphism $\varphi_D \circ f$ extends to a morphism $(E,D,s) : Y/\bG_m \to \MpCrat(X/G)$. In other words, we can fill the lower horizontal arrow in the diagram below, and then fill the rest with the dotted arrows so that the square is Cartesian,
\[
\begin{tikzcd}
  (Y \setminus \closedpt)/ \bG_m \ar[r, symbol=\dashrightarrow] \ar[dr] \ar[rr, bend left, "f"] & (\Gr^{E, D, s}_{X/G})/\bG_m \ar[d, symbol=\dashrightarrow] \ar[r, symbol=\dashrightarrow] & \MpC(X/G)_{Y/ \bG_m} \ar[d, "\varphi_D"]  \\  & Y / \bG_m  \ar[r, symbol=\dashrightarrow, "{(E,D,s)}"] & \MpCrat(X/G)_{Y/ \bG_m}
 \end{tikzcd}.
\]
We now let $\Sigma \subset \Gr^{E, D,s}_{X/G}$ be the scheme-theoretic image of $Y \setminus \closedpt$. As in the finite dimensional setting, the key to establishing monotonicity for $\MpC(X/G)$ is that the line bundle $\cL(\chi)^\dual$ is relatively ample for the morphism $\varphi_D$, which guarantees that it is relatively ample for the projection $\Sigma \to Y$.

\subsection{HN boundedness}

A numerical invariant $\mu$ on $\cX$ satisfies \emph{HN boundedness} (\Cref{defn: HN boundedness}) if for any quasi-compact open substack $\cU \subset \cX$, there is another quasi-compact open substack $\cU' \subset \cX$ such that for any (finite type) point $x \in \cU$, in order to maximize $\mu(f)$ among all filtrations with $f(1) \cong x$, it suffices to consider only those with $f(0) \in \cU'$. It is a non-trivial fact that, under mild hypotheses, this condition implies the existence of a filtration $f$ that maximizes $\mu(f)$ subject to $f(1)\cong x$ \cite[Prop.~4.4.2]{hl_instability}, i.e., an HN filtration. Note, in particular, that if $\cX$ is quasi-compact, then HN boundedness holds trivially with $\cU'=\cX$ for any $\cU$, so HN filtrations always exist.

In order to study the optimization problem $\sup \{ \mu(f) | f(1) \cong x\}$, we use the theory of degeneration fans developed in \cite[Sect.~3]{hl_instability}. Given a $k$-point $x \in \MpC(X/G)(k)$, the degeneration fan $|\DF(\MpC(X/G),x)|$ is a topological space that parameterizes all filtrations of $x$. It is constructed as a union of simplicial Euclidean cones, and we recall its properties in the beginning of \Cref{section: hn boundedness}. The numerical invariant defines a function on $|\DF(\MpC(X/G),x)|$, and the goal is to show that the maximum of $\mu$ must occur on a subspace of $|\DF(\MpC(X/G),x)|$ that parameterizes filtrations in a quasi-compact substack.

While the degeneration fan of a point in $\Bun_G(C)$ has a simple Lie-theoretic description amenable to convex optimization techniques, the degeneration fan of a point $x \in |\MpC(X/G)|$ in the case of an arbitrary affine scheme $X$ is complicated. The cones involved depend on the choice of the scheme $X$ and the point $x$. As a consequence, it becomes difficult to control the solutions of the convex optimization problem systematically (cf. \Cref{ex: hurwitz 2}). Instead of approaching the problem directly, we use a compactification and wall-crossing argument.  Our strategy proceeds in three stages:

\subsubsection*{Stage 1} In \Cref{subsection: hn boundedness for bung} we focus on the case when $X=\Spec(k)$, so $\MpC(X/G) = \Bun_{G}(C)$ is the stack of $G$-bundles on $C$. We are able to solve the optimization problem explicitly in this case, and we obtain a new proof of HN boundedness (and hence the existence of Harder-Narasimhan parabolic reductions) that differs from the classical arguments in the case of vector bundles.

\subsubsection*{Stage 2} We compactify $X \subset \overline{X}$ relative to $X/\!/G$, i.e., $\overline{X} \to X/\!/G$ is a projective scheme over $X/\!/G$ equipped with a linearizable $G$-action. Using Kontsevich stable maps, we define a stack of gauged maps $\cM^{G}_0(\overline{X})$ (cf. \Cref{defn: gauged maps}) that partially ``compactifies'' the mapping stack $\MpC(\overline{X}/G) \supset \MpC(X/G)$ over the base $X/\!/G$. Several variations of this stack were already considered in \cite{woodward-quantum-quotients, gonzalez-solis-woodward-stable-gauged, gonzalez-woodward-quantum}. A key fact is that the degeneration fan of a point in the compactified stack $\cM_0^{G}(\overline{X})$ is identified with the degeneration fan of the underlying $G$-bundle in $\Bun_{G}(C)$, so we can directly compare these two optimization problems. Using convex optimization techniques, we establish HN boundedness for $\cM^{G}_0(\overline{X})$ by controlling the distance between the optimizing filtration in $\cM^{G}_0(\overline{X})$ and the HN filtration of the underlying $G$-bundle in $\Bun_{G}(C)$.

\subsubsection*{Stage 3} To return to our original problem for the stack $\MpC(X/G)$ we use a wall-crossing strategy. The boundary divisor $\overline{X} \setminus X$ provides an additional stability parameter $\gamma > 0$ on the stack $\cM^{G}_0(\overline{X})$. We show that for points in any quasi-compact open substack $\cU \subset \MpC(X/G) \subset \cM_0^G(\overline{X})$, the optimizing filtrations in $\cM^{G}_0(\overline{X})$ stabilize uniformly for a large enough value of $\gamma$, and moreover they are forced to factor through the open substack $\MpC(X/G) \subset \cM^{G}_0(\overline{X})$. This yields a bounded set of optimizing filtrations in $\MpC(X/G)$ for the points of $\cU$, thus concluding the proof of HN boundedness.

\subsection{Full generality} \label{subsection: intro full generality}
In this subsection, we explain the general version of our results. We work over a Noetherian base scheme $S$, and fix a smooth projective family of curves $C \to \mathcal{Y}$ over a Noetherian $S$-stack $\cY$. Our group $G$ will be an affine smooth geometrically reductive group scheme over $S$ (the $S$-fibers need not be connected). We fix a finite type projective-over-affine scheme $X \to S$ equipped with an action of $G$; this is a general context in which the results of relative GIT apply to the stack $X/G$. Finally, we fix a number $n$; the number of marked points. We refer the reader to \Cref{subsection: notation} for more details on our setup.

In \Cref{defn: gauged maps} we introduce an algebraic stack of gauged maps $\cM_n^G(X) \to \cY$, which parameterizes Kontsevich stable maps $u : \tilde{C} \to E(X)$ with marked points $\marked=(p_1,\ldots,p_n)$ such that the composition $\tilde{C} \to E(X) \to C$ has degree $1$. There is a generalized Hitchin morphism $\cM_n(X) \to (A_{X} /\!/G) \times_{S} \cY$, which satisfies the existence part of the valuative criterion for properness (\Cref{prop: valuative criterion over affine GIT quotient}). Here we denote by $A_X/\!/G$ the $S$-affine GIT quotient of the affinization $A_X := \uSpec_S(\Gamma(\cO_X))$ of $X$ relative to $S$. 

In \Cref{subsection: numerical invariants} we define a broad class of numerical invariants on $\cM_n^G(X)$. On every cone of the degeneration fan $|\DF(\cM_n^G(X),x)|$ for a point $x = (E,u,\marked)$ in $\cM_n^G(X)$, each numerical invariant has the form
\[ \mu = \frac{\ell_{V} + \ell_{\mrk} + \ell_{\gen} + \ell_{\epsilon}}{\sqrt{b}},\]
where $b$ is a positive definite quadratic form, and each $\ell$ is a linear form. More specifically:
\begin{enumerate}
    \item $\ell_V$ depends on a choice of representation $G \to \GL(V)$ for some vector bundle $V$ over $S$. It measures the instability of the underlying $G$-bundle $E$.
    \item $\ell_{\mrk}$ depends on a choice of an $n$-tuple of stability parameters for $X/G$ (one for each marked point). It measures the instability of the image $u(p_i)$ of each marked point $p_i$.
    \item $\ell_{\gen}$ depends on a choice of a stability parameter on $X/G$. It measures the instability of the image of the generic point of $C$ under the gauged map $u$.
    \item $b$ depends on a choice of rational quadratic norm on graded points of the classifying stack $BG$. It is a generalization of the notion of Weyl invariant norm of graded points to the case of a non-split reductive group scheme $G$ over $S$.
    \item $\ell_{\epsilon}$ is a formally infinitesimal term coming from an ample line bundle in the stack of Kontsevich stable maps. Its inclusion in the numerical invariant $\mu$ is a necessary technicality that comes up in the total order on $\Theta$-strata, and further cuts down the semistable locus with respect to $\ell_V + \ell_\mrk+\ell_\gen$ if the parameters are non-generic.
\end{enumerate}

Our main general result in this paper is the following.
\begin{thm}[\Cref{thm: theta stratification in general} + \Cref{thm: moduli space in the case of centrally contractible}]
The numerical invariant $\mu$ defines a weak $\Theta$-stratification on $\cM_{n}^{G}(X)$. Moreover, if the base scheme $S$ is defined over $\bQ$, then $\mu_{\delta}$ defines a $\Theta$-stratification, and the semistable locus $\cM_{n}^{G}(X)^{\mu \dash \mathrm{ss}}$ admits a relative good moduli space that is separated and locally of finite type over $\cY$. Each connected component of the moduli space is proper over the generalized Hitchin base $(A_{X} /\!/G) \times_{S} \cY$.
\end{thm}

The same arguments should go through for a version of the stack $\cM_n^G(X)$ where we replace $C$ with a smooth projective orbifold curve and we allow the over curve $\widetilde{C}$ to acquire orbifold points as well. We have chosen not to pursue that direction in this paper, and have contented ourselves with showing how the argument can be generalized to the orbifold setting in the case of an affine target $X$ in \Cref{section: orbifold curves}.

\subsection{Author's note}

Preliminary results on this project (for projective $X$ over a field, and $G=\GL_n$) were first presented in June 2019, at the Vector Bundles over Algebraic Curves conference in Sønderborg, Denmark. At that time the project was a collaboration between the first author, Eduardo Gonzalez, and Pablo Solis. The second author joined the project as he became involved in developing the foundations of infinite dimensional GIT. At the same time, Pablo withdrew due to a career change. The paper was delayed by the discovery that we could extend our techniques to the more general setting where $X$ is projective-over-affine. This significantly improved our results, but it meant that other projects building on the ideas of infinite dimensional GIT \cite{torsion-freepaper, rho-sheaves-paper} were released prior to this paper. It also lengthened the manuscript, so we have split the project into this more foundational paper, and a planned followup paper with Eduardo applying these results to wall-crossing problems. We thank Eduardo and Pablo both for their early contributions to this paper, and many inspiring conversations on these topics. We also thank interested readers for their patience with the haphazard roll-out of these results.



\subsubsection{On endnotes} We are experimenting with a new expository technique in this paper. A typical math paper, in order to conform with stylistic and length norms, leaves many small claims unjustified, and many arguments abridged. Rather than depriving the reader of the details in our notebooks, we have included these small arguments and elaborations as endnotes. They will be included in the preprint, but not in the published manuscript. We warn the reader that they are informal, and not as carefully edited as the body of the paper.

\subsubsection{Related work}

\noindent \emph{Harder-Narasimhan theory.} As far as we are aware, this is the first algebro-geometric development of Harder-Narasimhan theory for a general stack of the form $\MpC(X/G)$. Previous instances of Harder-Narasimhan parabolic reductions were known for stacks of $G$-bundles \cite{behrend-thesis, biswas-holla-hnreduction}, Higgs bundles \cite{higgs-hn-reduction} and Bradlow pairs \cite{thaddeus-stable-pairs, mochizuki-donaldson-invariants}. Over $\bC$, our $\Theta$-stratification is analogous to the Morse stratification for the Yang-Mills-Higgs functional arising in the study of the vortex equations. This was studied in \cite{venugopalan-yang-milss,trautwein-hn-symplectic} under the assumption that the moment map is proper, so the affine GIT quotient $X/\!/G$ is a point.\\

\noindent \emph{Index formulas.} Beyond the case of the stack of $G$-bundles $\Bun_{G}(C)$ \cite{teleman-quantization,teleman-woodward}, a similar Verlinde formula for the stack of semistable $G$-Higgs bundles over a curve $C$ was obtained in \cite{verlinde, andersen2017verlinde}. Other instances of index computations for gauged maps appeared in \cite{ruan-zhange-verlinde} and \cite{oprea-sinha}. We are not aware of other general index computations of K-theory classes on stacks of bundles with affine decorations. We note that our setting directly applies to the stack of Higgs bundles only when $C$ is an elliptic curve. The authors are currently investigating possible generalizations of the setting of this paper that would in particular allow twisting the target $X$ by a fixed line bundle on $C$, and hence would also recover the Verlinde formulas for stacks of Higgs bundles.

\subsubsection{Acknowledgements}

In addition to Pablo and Eduardo, the first author would like to thank Jochen Heinloth for first suggesting the usefulness of affine Grassmannians for checking the intrinsic GIT criteria. We would also like to thank the following people for helpful discussions on these topics: Konstantin Aleshkin, Ron Donagi, Giovanni Inchiostro, Melissa Liu, Davesh Maulik, Andrei Okounkov, Shubham Sinha, and Chris Woodward.

The first author was supported by NSF CAREER grant DMS-1945478, NSF FRG grant DMS-2052936, and a Sloan Foundation Research Fellowship. The second author would like to thank the mathematics departments at Cornell University and Columbia University for supporting him while working on this project.

\section{Preliminaries}
\subsection{Notation} \label{subsection: notation}
\noindent{\textit{Bases.}} We work over a Noetherian scheme $S$. We fix a finite type algebraic stack $\mathcal{Y}$ over $S$ that will parametrize our ``reference'' family of smooth curves. 

\medskip

\noindent{\textit{Reference curve.}} Fix once and for all a smooth projective morphism $C \to \mathcal{Y}$ such that the geometric fibers are connected curves. We also fix the choice of a relatively very ample line bundle $\mathcal{O}_{C}(1)$ on $C$. We denote by $K_{C} \vcentcolon = \Omega^1_{C/S}$ the relative sheaf of Kahler differentials, which is a line bundle on $C$ in this case. We will denote by $\sqrt{K_{C}}$ the class in rational K-theory of $C$ given by $\frac{1}{2}[\cO_{C}] + \frac{1}{2}[K_{C}] \in K_0(C)_{\mathbb{Q}}$\endnote{We think of $\sqrt{K_{C}}$ as a square root of $K_{C}$. In fact the restriction at each geometric fiber of $C \to \mathcal{Y}$ is indeed represented by a line bundle that squares to $K_{C}$.}.

\medskip

\noindent{\textit{Torsors and mapping stacks.}} For any flat group scheme $G$ locally of finite type over $S$, we denote the quotient stack $S/G$ as $BG$. For any $S$-stack $\cX$, we denote by $\MpC(\cX)$ the mapping stack $\uMap_{\cY}(C, \cX\times_{S} \cY)$ which sends a $\cY$-scheme $T$ to the groupoid of morphisms $f: C_{T} \to \cX_{T}$. We denote by $\Bun_G(C)$ the moduli stack of $G$-bundles over $C$. By definition, $\Bun_G(C)$ is the mapping stack $\MpC(BG)$.

\medskip

\noindent{\textit{Group schemes.}} In this article, we work with a fixed smooth affine group scheme $G$ over $S$ with (possibly disconnected) reductive fibers. This means that there is a short exact sequence of group schemes
\[ 1 \to G_0 \to G \to  \pi_0(G) \to 1\]
where $G_0$ is a smooth group scheme with connected reductive $S$-fibers, and $\pi_0(G)$ is a separated \'etale group scheme of finite type over $S$ \cite[Prop. 3.1.3]{conrad_reductive}. We impose the additional condition that $\pi_0(G)$ is finite over $S$; this is equivalent to $G$ being geometrically reductive by \cite[Thm. 9.7.6]{alper-adequate-moduli}. We denote by $\D(G)$ the derived closed subgroup scheme of $G_0$, as defined in \cite[Thm. 5.3.1]{conrad_reductive}. We will use the following notion.
\begin{defn}
Let $V$ be a vector bundle on $S$. We say that a homomorphism $G \to \GL(V)$ is \emph{almost faithful} if the kernel of the composition $\D(G) \to G \to \GL(V)$ is finite and flat over $S$.
\end{defn}

\noindent{\textit{Target.}} Let $X$ be a projective-over-affine $S$-scheme equipped with a left action of $G$. We fix the choice of an $S$-very ample $G$-equivariant polarization $L_{\amp} \in \Pic^{G}(X)$. We shall also fix a $G$-equivariant $S$-semi-ample line bundle $L_{\gen} \in \Pic^{G}(X)_{\bQ}$. Set $A_{X} \vcentcolon = \uSpec_{S}(\mathcal{O}_{X})$ to be the relative affinization of $X$. Note that $A_{X}$ admits an induced left $G$-action. By assumption $A_{X}$ is of finite type over $S$, and the canonical morphism $X \to A_{X}$ is projective with $A_X$-very ample $G$-polarization $L_{\amp}$.

\medskip

\noindent{\textit{Associated fiber bundles.}} For any $G$-bundle $E$ on an $S$-scheme $Y$, we denote by $E(X)$ the associated $X$-bundle over $Y$. In other words, $E(X) =(E \times X)/_{\rm diag} G= (E\times X)/(eg,x)\sim(e,gx)$. \'Etale descent for polarized morphisms implies that $E(X)$ is represented by a projective-over-affine scheme over $Y$.

\medskip

\noindent{\textit{The stack $\Theta$.}} We denote by $\Theta$ the quotient stack $\mathbb{A}^1 / \mathbb{G}_m$. By convention, we equip $\mathbb{A}^1= \Spec(\mathbb{Z}[t])$ with the $\mathbb{G}_m$-action that gives $t$ weight $-1$.

\medskip

\noindent\textit{Convention on line bundles associated to characters.} For any character $\chi: G \to \mathbb{G}_m$ of $G$, we can think of $\chi$ as a one-dimensional representation of $G$, and its associated $\cO_G$-comodule corresponds to a line bundle on the stack $BG$. We denote by $\cO_{X}(\chi) = \cO_X \otimes_{\cO_S} (-\chi)$ the equivariant line bundle in $\Pic^G(X)$ obtained by pulling back the line bundle associated to $-\chi$ under the morphism $X/G \to BG$. This convention is chosen so that the subspace of $\Gamma(X,\cO_X)$ of sections $s$ such that $g \cdot s = \chi(g) s$ is canonically identified with $\Gamma(X,\cO_X(\chi))^G$.

Under this convention, let $k$ be a field and choose a morphism $\mathbb{A}^1_k/\mathbb{G}_m \to X/G$. The induced homomorphism of automorphism groups at the origin $0 \in \mathbb{A}^1_k(k)$ induces a one-parameter subgroup $\lambda: \mathbb{G}_m \to G_k$. Then the pullback of the line bundle $\cO_X(\chi)$ to $\mathbb{A}^1_k/\mathbb{G}_m$ has a non-vanishing invariant global section if and only if $\langle \lambda, \chi \rangle \leq 0$.

\medskip

\noindent{\textit{Categories of schemes.}} For any algebraic stack $T$ over $S$, we let $\Sch_T$ be the category of schemes over $T$. We write $\Aff_T$ to denote the category of (absolutely) affine schemes equipped with a morphism to $T$.

\medskip

\noindent{\textit{Rational line bundles.}} We will generally work with rational line bundles on schemes and stacks, meaning an element of the Picard group tensored by $\mathbb{Q}$. Sometimes we might omit the word ``rational'' and refer to these simply as line bundles whenever the distinction is harmless.

\subsection{The moduli of gauged maps} \label{subsection: gauged maps introduction}
\begin{defn}[Moduli of gauged maps] \label{defn: gauged maps}
For any integer $n \geq 0$, the stack $\cM_{n}^{G}(X)$ of $n$ marked \emph{gauged maps} from $C$ to $X$ sends any $\mathcal{Y}$-scheme $T$ to the groupoid $\cM_{n}^{G}(X)(T)$ of tuples $(E, u,\marked)$, where 
\begin{enumerate}
    \item $E$ is a $G$-bundle on $C_{T}$.
    \item $u \colon \tilde{C} \to E(X)$ is a $T$-family of Kontsevich stable maps into $E(X)$ with markings $\marked=(p_i)_{i=1}^n : T \to \tilde{C}^n$. 
    \item The map $\tilde{C}_t \to E_{t}(X) \to C_{k(t)}$ has degree $1$ for every $t \in T$.\endnote{This condition means that for every $t$ there is a unique component $\tilde{C}_{t,0} \subset \tilde{C}_t$ such that $\pi \colon \tilde{C}_{t,0} \to C_{t}$ is an isomorphism and moreover any other subcurve of $D \subset \tilde{C}_t$ that does not contain $\tilde{C}_{t,0}$ is collapsed by $\pi$ to a divisor in $C_t$.}
\end{enumerate}
\end{defn}
Let $N^1_S(X/G)$ denote the quotient of the group $\Pic(X/G)$ by the subgroup generated by line bundles $\cL$ such that for any proper irreducible curve $C$ over a field $k$ and any morphism $C \to X/G$ such that the composition $C \to S$ factors through $\Spec(k)$, we have $\deg(\cL|_C)=0$. We will use the notation $H_2(X/G) \vcentcolon = \Hom(N^1_S(X/G),\bZ)$. Let $(E,u, \marked)$ be a $T$-point of $\cM_{n}^{G}(X)$. The pair $(E|_{\tilde{C}}, u)$ defines a morphism $\tilde{C} \to X/G$. Choose $t \in T$. For each $[\mathcal{L}] \in N^1_S(X/G)$, we pull back $\mathcal{L}$ to the reduced curve $\widetilde{C}_{t}$ and take the degree to get a well-defined homomorphism
\[\delta_{t}: N^1_{S}(X/G) \to \mathbb{Z}, \; \; \; \; \; \;  [\mathcal{L}] \mapsto  \deg(\mathcal{L}|_{\tilde{C}_{t}})\]
 We say that the morphism $\tilde{C} \to X/G$ has degree $d \in H_2(X/G)$ if $\delta_{t}$ is equal to $d$ for all $t \in T$.
\begin{defn}
We define $\cM_{n}^{G}(X)_d$ to be the full subfunctor of $\cM_{n}^{G}(X)$ consisting of $T$-points such that the induced morphism $(E|_{\tilde{C}}, u): \tilde{C} \to X/G$ has degree $d$. If $d \in H_2(BG)$, then we define by $\cM_{n}^{G}(X)_d: = \sqcup_{d' \in \psi^{-1}(d)} \cM_{n}^{G}(X)_{d'}$, where $\psi: H_2(X/G) \to H_2(BG)$ is the group homomorphism induced by $X/G \to BG$.
\end{defn}
For each $d \in H_2(X/G)$, $\cM_{n}^{G}(X)_d$ is an open and closed substack of $\cM_{n}^{G}(X)$, and we have $\cM_{n}^{G}(X) =  \bigsqcup_{d \in H_2(X/G)} \cM_{n}^{G}(X)_d$.
\begin{example} \label{example: degree of a G-bundle}
If $X = S$, then we have $\cM_{0}^G(X) = \Bun_{G}(C)$. The relative Neron-Severi group $N^1_S(BG)$ is the group $X^*(G)$ consisisting characters $G \to (\bG_m)_S$. For every element of the dual group $d \in H_2(BG) =  X^*(G)^{\vee}$, we denote by $\Bun_{G}(C)_d$ the corresponding open and closed substack of $\Bun_{G}(C)$. We say that a $G$-bundle corresponding to a point in $\Bun_{G}(C)_d$ has degree $d \in H_2(BG)$.
\end{example}
\subsubsection{The morphism $p$}
For any $T$-point of $\cM_{n}^{G}(X)$, as in \Cref{defn: gauged maps}, the fact that $\tilde{C} \to C_T$ has degree $1$ and that $E(A_X) \to C_T$ is affine implies that the composition $\tilde{C} \to E(X) \to E(A_{X})$ factors uniquely through a section $s: C_{T} \to E(A_{X})$.\endnote{Let $\pi: \tilde{C} \to E(X) \to C_{T}$ denote the composition. The canonical morphism $\mathcal{O}_{C_{T}} \to \pi_{*}(\mathcal{O}_{\tilde{C}})$ is an isomorphism, because the fibers $\pi_{t}$ have degree 1 for all $t \in T$. Since the associated fiber bundle $E(A_{X})$ is affine over $C$, we have the following chain of equalities
\[  \Mor(\widetilde{C},\,  E(A_{X})) = \Hom_{\mathcal{O}_{C_{T}}}(\mathcal{O}_{E(A_{X})}, \, \pi_*(\mathcal{O}_{\tilde{C}})) = \Hom_{\mathcal{O}_{C_{T}}}(\mathcal{O}_{E(A_{X})}, \, \mathcal{O}_{C_{T}}) = \Mor(C_{T},\,  E(A_{X})).\]} This defines a morphism of stacks
\begin{equation}\label{E:morphism_p}
p:  \cM_{n}^{G}(X)_d \to \MpC(A_{X}/ G) , \; \; \; \; \; (E, u,\marked) \mapsto (E, s)
\end{equation}
\begin{example}
If $X = A_{X}$ is affine over $S$, then $p$ induces an isomorphism $\cM_{0}^G(X) \cong \MpC(A_{X}/G)$.
\end{example}

\subsubsection{Determinant and marking line bundles on $\cM_{n}^{G}(X)$.}
\begin{defn} \label{defn: line bundle almost faithful reps}
Let $V$ be an almost faithful representation $G \to \GL(V)$. We denote by $\cD(V)$ the rational line bundle on $\Bun_{G}(C)$ given by \[\cD(V) \vcentcolon = \det\left(R\pi_\ast \left(E_{univ}(V) \otimes \sqrt{K_{C}}|_{\Bun_{G}(C) \times C} \right) \right)^{\vee},\]where $E_{univ}$ is the universal $G$-bundle on $\Bun_{G}(C) \times C$ and $\pi: \Bun_{G}(C) \times C \to \Bun_{G}(C)$ is the projection onto the first factor. 
We also denote by $\cD(V)$ the corresponding pullback line bundles on $\cM_{n}^G(X)$ and $\MpC( A_{X}/G)$ under the morphisms $\cM_{n}^G(X) \to \MpC( A_{X}/G) \to \Bun_{G}(C)$.
\end{defn}

\begin{defn} \label{defn: line bundles pullbacked markings}
Fix an $n$-tuple $(L_i^{\mrk})_{i=1}^n$ of $G$-equivariant $S$-semi-ample line bundles on $X$. We denote by $L^{\mrk}: = \boxtimes_{i=1}^n L_i^{\mrk}$ the tensor line bundle on the stack $(X/G)^n$. We shall consider the line bundle $\cL_{\mrk} := \ev^*(L^{\mrk}) \in \Pic(\cM_{n}^{G}(X))$, where 
\[\ev: \cM_{n}^{G}(X) \to (X/G)^n, \; \; \; \; (E,u,\marked) \mapsto (u\circ p_1, \ldots, u \circ p_n)\] 
is the evaluation morphism induced by the markings.
\end{defn}

\subsubsection{The Cornalba line bundle on $\cM_{n}^{G}(X)$}
We recall a construction due to Cornalba \cite{cornalba_projectivity}. 
Let $\pi \colon Y \to B$ be a projective morphism of schemes. The stack of relative Kontsevich stable maps $\cM_{g,n}(Y/B)$ over $B$ is defined by
\[
\cM_{g,n}(Y/B) \colon T \mapsto \left\{ \vcenter{\xymatrix{\cC \ar[r]^u \ar[d] & Y \ar[d] \\ T \ar@/_/[u]_{p_i} \ar[r] & B}} \begin{array}{c} \text{such that } \forall t \in T, \cC_t \to Y_t \\ \text{is Kontsevich stable} \end{array} \right\},
\]
Here $\cC \to T$ is a flat family of connected genus $g$ nodal curves and $p_i : T \to \cC$ for $i=1,\ldots,n$ are sections landing in the smooth locus of $\cC / T$. 

\begin{notn}
    For any such family of nodal curves $\pi: \cC \to T$, we will make use of the Deligne pairing $\langle -,- \rangle \colon \Pic(\cC)^{\otimes 2} \to \Pic(T)$ \cite[Chapt.XIII, Sect. 5]{acg-algebraic-curves-ii}. By \cite[Chapt. XIII, Thm. 5.8]{acg-algebraic-curves-ii}, the Deligne pairing of two line bundles $\cH, \cQ \in \Pic(\cC)$ admits the following expression
\[\langle \cH, \cQ \rangle := \det \left( R\pi_\ast \left( ([\cH]-[\cO_{\cC}]) \otimes ([\cQ]-[\cO_{\cC}]) \right) \right), \]
which we take as a definition for our purposes. We will use the notation $\langle \cH \rangle_2 := \langle \cH, \cH \rangle$.
\end{notn}

In the following, we let $\pi : \cC \to \cM_{g,n}(Y/B)$ denote the universal curve, $u : \cC \to Y$ the universal stable map, and $D_i \hookrightarrow \cC$ the image of the universal marking $p_i : \cM_{g,n}(Y/B) \to \cC$ viewed as a Cartier divisor.

\begin{prop} \label{P:ample_kontsevich}
$\cM_{g,n}(Y/B)$ is a DM stack locally of finite type over $B$ whose connected components are proper over $B$. If $\dim(Y)>0$ and $M$ is a line bundle on $Y$ that is very ample relative to $B$, then on each component of $\cM_{g,n}(Y/B)$ the line bundle $\left\langle \omega_{\pi}(\sum_i D_i) \otimes u^\ast(M^3)\right\rangle_2$
is very ample relative to $B$.
\end{prop}

\begin{proof}
Using local projective embeddings of $Y$, one can reduce the claim to the special case where $B = \Spec(\mathbb{Z})$ and $Y = \bP^n$.\endnote{This reduction uses the following three properties of the functor $\cM_{g,n}(Y/B)$:
\begin{enumerate}
    \item Given a map $B' \to B$, the diagram
\[
\xymatrix{\cM_{g,n}(Y \times_B B' / B') \ar[r] \ar[d] & \cM_{g,n}(Y/B) \ar[d] \\ B' \ar[r] & B}
\]
is Cartesian.
\item The formation of $\cL_M$ is compatible with base change along $B' \to B$.
\item If $Z \to Y$ is a closed immersion of projective $B$-schemes then the induced map $\cM_{g,n}(Z/B) \to \cM_{g,n}(Y/B)$ is also a closed immersion.
\end{enumerate} } 
The fact that $\cM_{g,n}(\bP^n/\Spec(\mathbb{Z}))$ is DM with proper connected components is standard \cite[Sect. 2]{ao-stable-mixedchar}. The ampleness of the line bundle on each component is the main theorem of \cite{cornalba_projectivity}.
\end{proof}

\begin{defn}By a formal line bundle on a stack $\cZ$, we mean a linear combination $\cL(\epsilon) = \sum_{j \geq 0} \cL_j \epsilon^j$, where $\cL_j \in \Pic(\cZ)$ and $\cL_j = 0$ for $j \gg 0$.
\end{defn}

We shall mimic the construction of the Cornalba line bundle for the stack of gauged maps $\cM_{n}^{G}(X)$. For the following definition we denote again by $\pi: \cC \to \cM_{n}^{G}(X)$ the universal curve, $u: \cC \to X/G$ the universal morphism and $D_i \subset \cC$ the divisors corresponding to the universal markings. For simplicity of notation, we also denote by $\cO_{C}(1)$ the pullback of the very ample line bundle on $C$ to $\cC$. We use the formal line bundle $\cO_C(1) + \epsilon \cdot u^*(L_{\amp})$ on $\cC$.
\begin{defn}
We denote by $\cL_{Cor}(\epsilon)$ the formal line bundle on $\cM_{n}^{G}(X)$ defined by
\begin{align*}
    \cL_{Cor}(\epsilon)  = & \epsilon^4 \left\langle \omega_{\pi} + \sum_i D_i + \frac{3}{\epsilon^2} \left(\cO_C(1) + \epsilon \cdot u^*(L_{\amp})\right) \right \rangle_2 \\
     := & 9 \left\langle \cO_{C}(1) \right\rangle_2 + 18\epsilon \left\langle \cO_{C}(1), u^*(L_{amp})\right\rangle + 6\epsilon^2 \left\langle \omega_\pi(\sum_i D_i), \cO_{C}(1)\right\rangle + \\
     & +9 \epsilon^2 \langle u^*(L_{amp}) \rangle_2 + 6 \epsilon^3 \left\langle \omega_\pi(\sum_i D_i), u^*(L_{amp})\right\rangle + \epsilon^4 \left\langle \omega_\pi(\sum_i D_i)\right\rangle_2
\end{align*}
\end{defn}

We regard $\bR[\epsilon]$ as a totally ordered vector space, where $f \geq g$ if $f(\epsilon) \geq g(\epsilon)$ for all $0 < \epsilon \ll 1$.

\begin{defn}
If $\cZ$ is separated with finite inertia, we say that a formal line bundle $\cL(\epsilon)$ is \emph{positive} if for any proper reduced curve $\Sigma$ with a quasi-finite morphism $f : \Sigma \to \cZ$, we have $\sum_{j \geq 0} \deg(\cL_j|_{\Sigma}) \epsilon^j > 0$. Given a morphism $q : \cX \to \cY$ that is separated with finite relative inertia, we say that a formal line bundle $\cL(\epsilon)$ is positive relative to $q$ if it is positive after base change along any map $\Spec(A) \to \cY$.
\end{defn}

\begin{prop} \label{prop: formal line bundle is ample on p fibers}
The morphism $p : \cM_{n}^{G}(X) \to \MpC(A_X/G)$ is relatively representable by a disjoint union of projective DM stacks. The formal line bundle $\cL_{Cor}(\epsilon)$ is positive relative to $p$.
\end{prop}
\begin{proof}
Choose a morphism from an affine scheme $T \to \MpC(A_X/G)$, which corresponds to a $G$-bundle $E$ over $C_{T}$ and a section $s: C_{T} \to E(A_{X})$. Let $E^s(X)$ be the fiber product
\[
\begin{tikzcd}
 E^s(X) \ar[d] \ar[r] & C_{T} \ar[d, "s"]\\  E(X)  \ar[r] &  E(A_{X})
\end{tikzcd}.
\]
By descent for polarized morphisms, $E^s(X) \to C_{T}$ is projective, and hence $E^s(X) \to T$ is projective. The fiber product $T \times_{\MpC(A_X/G)} \cM_{n}^{G}(X)$ is the subfunctor of $\cM_{g,n}(E^s(X)/T)$ parameterizing families such that the composition $\tilde{C} \to E^s(X) \to C$ has degree 1 on every fiber. Since this is an open and closed condition, the connected components of $T \times_{\MpC(A_X/G)} \cM_{n}^{G}(X)$ are projective DM stacks.

The $G$-equivariant $S$-ample line bundle $L_{\amp}$ on $X$ induces a line bundle on $E^s(X)$, which we also denote by $L_{\amp}$. Set $\epsilon = \frac{1}{k}$ where $k$ is a nonzero integer. Using the bilinearity properties of the Deligne pairing (which in particular allow us to formally extend it to rational line bundles), we see that $\cL_{Cor}\left(\frac{1}{k}\right)$ on the fiber product $T \times_{\MpC(A_X/G)} \cM_{n}^{G}(X) \subset \cM_{g,n}(E^s(X)/T)$ can be rewritten as
\[\cL_{Cor}\left(\frac{1}{k} \right) = \frac{1}{k^4} \left\langle \omega_\pi(\sum_i D_i) \otimes u^*\left((\cO_C(1)|_{E^s(X)})^{k^2} \otimes L_{\amp}^k\right)^3 \right\rangle_2. \]%
Since $(\cO_C(1)|_{E^s(X)})^{k^2} \otimes L_{\amp}^k$ is very ample on $E^s(X) \to T$ for all sufficiently large positive integers $k$, \Cref{P:ample_kontsevich} implies that $\cL_{Cor}\left(\frac{1}{k}\right)$ is $T$-ample on $T \times_{\MpC(A_X/G)} \cM_{n}^{G}(X) \subset \cM_{g,n}(E^s(X)/T)$ for any $k\gg0$. Hence, for every $0<\epsilon \ll 1$ the line bundle $\cL_{Cor}(\epsilon)$ pairs positively with every curve in a $T$-fiber of $T \times_{\MpC(A_X/G)} \cM_{n}^{G}(X)$.
\end{proof}

\subsubsection{A completeness result}

Consider the affine GIT quotient $A_{X}/ \! / G= \uSpec_{S}(\cO_{A_{X}}^{G})$, where $\cO_{A_{X}}^{G}$ denotes the sheaf of rings of $G$-invariants.

\begin{defn}[Generalized Hitchin morphism]
For any $T$-point $(E,u,\marked)$ of $\cM_n^G(X)$, the composition $\tilde{C} \to E(X) \to (A_X/\!/G) \times_S T$ factors through $T$, because $\tilde{C} \to T$ is a flat projective family with geometrically integral fibers. This defines a generalized Hitchin morphism
\[
h : \cM_n^G(X) \to (A_X/\!/G) \times_S \cY.
\]
\end{defn}



\begin{prop} \label{prop: valuative criterion over affine GIT quotient}
Suppose that $S$ is a scheme over $\mathbb{Q}$ (i.e., in characteristic $0$). Then the morphism $h$ satisfies the existence part of the valuative criterion of properness (for stacks) with respect to discrete valuation rings.
\end{prop}
\begin{proof}
The morphism $h$ is the composition of $p : \cM_{n}^{G}(X) \to \MpC(A_X/G)$ and the canonical morphism of mapping stacks $\MpC(A_X/G) \to \MpC(A_X/\!/G)=(A_X /\!/G) \times_S \cY$ induced by the projection $A_X/G \to A_X/\!/G$. By \Cref{prop: formal line bundle is ample on p fibers}, it suffices to show that $\MpC( A_{X}/G) \to (A_{X}/\!/G)\times_{S} \cY$ satisfies the existence part of the valuative criterion. Let $R$ be a discrete valuation ring over $(A_{X}/\!/G)\times_{S} \cY$. By base-changing to $\Spec(R)$, we can assume without loss of generality that $S = \cY = A_{X}/\!/G = \Spec(R)$. In this case $G$ is embeddable, so we can choose an isomorphism $A_X / G \cong \Spec(B)/\GL_N$ for some $R$-algebra $B$ and assume that $G=\GL_N$. Let $\eta$ and $s$ denote the generic point and special point of $\Spec(R)$ respectively. Choose an $R$-morphism $\eta \to \MpC( A_{X}/G)$, corresponding to a morphism $f: C_{\eta} \to A_{X}/G$ consisting of a pair $(E, \psi)$, where $E$ is a $G$-bundle on $C_{\eta}$ and $\psi$ is a section $\psi: C_{\eta} \to E(A_{X})$. Since $A_{X}/G \to\Spec(R)$ is a good moduli space morphism, it is universally closed \cite[Thm. 4.16(ii)]{alper-good-moduli}. By \cite[Th. 3.9]{dilorenzo2023degenerations} we can extend the morphism $f: C_{\eta} \to A_{X}/G$ over the generic point of the special fiber $C_{s}$ after perhaps passing to an extension of $R$. By spreading out, we extend the morphism $f: C_{\eta} \to A_{X}/G$ further to $\widetilde{f}: U \to A_{X}/G$, where $U \subset C_{R}$ is the open complement of finitely many codimension $2$ points in the special fiber $C_{s}$ of the two-dimensional regular scheme $C_{R}$. Now the proof of \cite[Prop.~3.21(2)]{AHLH} applies verbatim to show that morphisms to $\Spec(B)/\GL_N$ extend uniquely over regular codimension-2 closed points.\endnote{We recall the argument for the benefit of the reader. By using Hartogs's theorem for maps into affine schemes, we are reduced to extending the composition $U \to \Spec(B)/\GL_N \to B\GL_N$ to the total space $C_R \supset U$. In this case $U \to B\GL_N$ corresponds to a rank $N$ vector bundle $\mathcal{E}$ on $U$, and we extend to a vector bundle on $C_R$ by taking the unique reflexive extension $(j_* \mathcal{E})^{\vee \vee}$, where $j: U \hookrightarrow C_R$ denotes the open immersion \cite[\href{https://stacks.math.columbia.edu/tag/0B3N}{Tag 0B3N}]{stacks-project}.}

\end{proof}

\subsection{Numerical invariants} \label{subsection: numerical invariants}
Recall that we equip the polynomial ring $\bR[\epsilon]$ with the total order defined by $f \leq g$ if and only if there exists some real number $\rho>0$ such that $f(\epsilon) \leq g(\epsilon)$ for all $0 \leq \epsilon< \rho$.
\begin{defn}\label{D:numerical_invariant}
A \emph{numerical invariant} on the stack $\cM_{n}^{G}(X)$ with values in $\bR[\epsilon]$ is a rule that assigns to any field $k$ over $\mathcal{Y}$ and any non-degenerate map $g : (B\bG_m^h)_{k} \to \cM_{n}^{G}(X)$ a scale invariant function $\bR^h \setminus 0 \to \bR[\epsilon]$. This assignment is required to be
\begin{enumerate}
\item invariant under field extensions $k \subset k'$,
\item compatible with restriction along maps $(B\bG_m)^m \to (B\bG_m)^h$ with finite kernel, and
\item locally constant in algebraic families of maps $g$.
\end{enumerate}
\end{defn}
See \cite[Defn.~0.0.3]{hl_instability} or \cite[2.5]{torsion-freepaper} for more details. We will consider numerical invariants of the form $\ell / \sqrt{b}$, where $\ell$ assigns a linear form $\bR^h \to \bR[\epsilon]$ and $b$ assigns an $\mathbb{R}$-valued rational quadratic form on $\bR^h$ to any such map $g : (B\bG_m)^h \to \cM_{n}^{G}(X)$.

One should think of $\bR^h$ via its canonical isomorphism with the group of real cocharacters $X_{*}((\bG_m^h)_{k}) \otimes \bR$. We can view any invertible sheaf $\cL$ on $(B\bG_m^h)_k$ as a character in $X^*((\bG_m^h)_{k})$, and we denote by $\wt(\cL)$ the unique real-valued linear functional on $\bR^h \cong X_{*}((\bG_m^h)_{k}) \otimes \bR$ determined by the natural pairing $X_{*}((\bG_m^h)_{k}) \otimes X^{*}((\bG_m^h)_{k}) \to \bZ$. Similarly, for any formal line bundle $\cL = \sum_i \cL_i \, \epsilon^i$ on $(B\bG_m^h)_{k}$, we define the $\bR[\epsilon]$-valued linear functional $\wt(\cL)$ on $\bR^h$ by $\sum_{i} \wt(\cL_i) \epsilon^i$. In the case when $h=1$, we will often abuse notation and also denote by $\wt(\mathcal{L})$ the image of $1 \in \mathbb{R}$ under the linear functional $\wt(\mathcal{L})$.

\begin{notn}
For any line bundle $\mathcal{L}$ on a stack $\mathcal{M}$, we denote by $\wt(\mathcal{L})$ the assignment that takes a map $g: (\mathbb{G}_m^h)_k \to \mathcal{M}$ to the linear functional $\wt(g^*(\mathcal{L}))$ defined as above. We similarly define the analogous assignment $\wt(\mathcal{L})$ as $g \mapsto \wt(g^*(\mathcal{L}))$  for formal line bundles $\mathcal{L}$.
\end{notn}

In our context $\ell$ will be a combination of weights of line bundles along with the following ``generic Hilbert-Mumford'' linear form.

\begin{defn} \label{defn: linear functionals on the stack}
For any field $k$ over $\mathcal{Y}$ and any multigraded point $g : B(\bG_m^{h})_{k} \to \cM_{n}^{G}(X)$, we define $\ell_{\gen}(g)$ to be the linear form on $\bR^h$ corresponding to the weight of the line bundle $-L_{\gen}$ at the generic point of the curve $C_k$. More precisely, the morphism $g$ consists of a $\mathbb{G}_m^h$-equivariant tuple $(E, u, p_1, \cdots, \, p_n)$, where $E$ is a $G$-bundle on $C_{k}$ and $u: \tilde{C} \to E(X)$. Since $\tilde{C} \to C_{k}$ is an isomorphism over the generic point $\eta \in C_{k}$, we can view the restriction $u_{\eta}$ as a $\bG_m^h$-equivariant section of the restriction $E_{\eta}(X)$ of $E(X)$ to $\eta$. The data $(E_{\eta}, u_{\eta})$ amounts to a morphism $g_{res}: (B\bG_m)^h_{k(\eta)} \to X/G$. We define the linear functional $\ell_{\gen}(g)$ by $-\wt(g_{\eta}^* L_{\gen})$.\endnote{Note that this makes sense as a linear functional on the cocharacters of $(\bG_m^h)_{k}^h$, because there is a canonical isomorphism between the group of cocharacters $X_{*}\left((\bG_m^h)_{k}\right)$ and the group of cocharacters $X_{*}\left((\bG_m^h)_{k(\eta)}\right)$.}
\end{defn}

A rational quadratic norm $b$ on graded points of a stack $\cM$ is an assignment to any non-degenerate graded point $g : (B\bG_m)_k^h \to \cM$ a positive definite rational quadratic form $g^\ast(b)$ on $X_\ast((\bG_m)_k^h) \cong \bR^h$ that satisfies conditions (1), (2), and (3) of \Cref{D:numerical_invariant}. When the group $G$ is split, there is a simple description of quadratic norms $b$ on $BG$. For this we need to introduce a notion of Weyl group for $G$.
\begin{defn} \label{defn: weyl group}
Suppose that $G$ is split, so that we can choose a maximal split torus $T \subset G$. We define the Weyl group $W$ to be the separated \'etale $S$-group scheme $W_{G}(T) = N_{G}(T)/Z_{G}(T)$ defined in \cite[Thm. 2.3.1]{conrad_reductive}.
\end{defn}
Our assumption that $\pi_0(G)$ is finite implies that $W$ is finite \'etale over $S$. Indeed, the short exact sequence
\[ 1 \to G_0 \to G \to \pi_0(G) \to 1\]
induces a short exact sequence of smooth group schemes
\[ 1 \to N_{G_0}(T) \to N_{G}(T) \to F \to 1\]
Here $F$ is a closed \'etale subgroup scheme of $\pi_0(G)$, and hence $F$ is also finite \'etale. Quotienting by the closed smooth subgroup scheme $Z_{G_0}(T)$ we get
\[ 1 \to W_{G_0}(T) \to N_{G}(T)/Z_{G_0}(T) \to F \to 1\]
By \cite[Pro. 3.2.8]{conrad_reductive}, the group $W_{G_0}(T)$ is finite \'etale. It follows that $N_{G}(T)/Z_{G_0}(T)$ is finite \'etale. $W$ is the quotient of $N_{G}(T)/Z_{G_0}(T)$ by the smooth group scheme $Z_{G}(T)/Z_{G_0}(T)$, and hence we have a smooth surjection $N_{G}(T)/Z_{G_0}(T) \to W$. The finiteness of $N_{G}(T)/Z_{G_0}(T)$ implies that the \'etale group scheme $W$ is also finite.

\begin{lem} \label{lemma: description of norm in the split case}
Suppose that $S$ is connected and $G$ admits a split maximal torus $T \cong (\bG_m^n)_S \subset G$. Then restriction to cocharacters in $T$ gives a bijection between the set of rational quadratic norms on graded points of $BG$ and the set of Weyl group invariant positive definite rational quadratic forms on the vector space of real cocharacters $X_{*}(T)_{\bR}$.
\end{lem}
\begin{proof}
By \cite[Thm. 1.4.8]{hl_instability}, a graded point of $BG$ corresponds to a point $x \in S(k)$ and a homomorphism $\gamma_x : (\bG_m)_{k} \to G_x$ up to conjugation. Given a norm on graded points of $BG$, one obtains a $W$-invariant rational quadratic norm on $X_\ast(T)_\bR$ by fixing a point $s \in S(k)$ and defining the norm of a cocharacter of $T$ to be the norm of the graded point at $s$ induced by this cocharacter -- this is independent of the choice of $s$ because by definition a norm on graded points is locally constant in families. Conversely, if one fixes a $W$-invariant rational quadratic norm on $X_\ast(T)_\bR$, then for any graded point as above one can conjugate $\gamma_x$ to lie in $T_x$, and define the norm of the graded point $(x,\gamma_x)$ to be the norm of the corresponding global cocharacter $(\bG_m)_S \to T$. This is well-defined because the norm on $X_\ast(T)_\bR$ is $W$-invariant. (note that the Weyl group of $G_x$ and its action on the coweight lattice is independent of $x$.)
\end{proof}

A rational quadratic norm $b$ on graded points of $BG$ induces a rational quadratic norm on graded points of $\cM_n^G(X)$, which we also denote by $b$, as follows:

\begin{defn} \label{defn: quadratic norm graded points}
For every field $k$ over $S$ and every nondegenerate graded point $g : (B\bG_{m}^h)_{k} \to \cM_{n}^{G}(X)$, we let $g_\eta : (B\bG_m^h)_{k(\eta)} \to X/G$ be the restriction to the generic point discussed in \Cref{defn: linear functionals on the stack}. Then $g^\ast(b)$ is the rational quadratic norm on $X_\ast((\bG_m^h)_k)_\bR$ identified with $(g_\eta)^\ast(b)$ under the canonical isomorphism $X_{*}((\bG_m^h)_{k})_{\bR} \cong X_{*}((\bG_m^h)_{k(\eta)})_{\bR}$. We will often denote the norm associated to $g^\ast(b)$ by $\|w\|_b := \sqrt{b(w)}$.
\end{defn}

\begin{defn}[Numerical invariant for gauged maps] \label{defn: numerical invariant}
Let $\delta = (\delta_{\gen}, \delta_{\mrk}) \in (\bR_{\geq 0})^{2}$ be a pair of nonnegative real numbers. The numerical invariant $\mu_{\delta}$ on $\cM_{n}^{G}(X)$ is defined by
\[\mu_{\delta} = \frac{-\wt(\cD(V)) - \delta_{\mrk} \wt(\cL_{\mrk}) +\delta_{\gen} \ell_{\gen} -  \epsilon \wt(\cL_{Cor}(\epsilon))}{\sqrt{b}}.\]
We shall denote by $\mu_{\delta}^0 := (-\wt(\cD(V))- \delta_{\mrk} \wt(\cL_{\mrk})+\delta_{\gen} \ell_{\gen})/\sqrt{b}$ the constant term of $\mu_{\delta}$ with respect to the formal variable $\epsilon$.
\end{defn}

\subsection{Generalities on \texorpdfstring{\(\Theta\)}{Theta}-stratifications}\label{subsection: theta-sratifications}
Let $\cM$ be an algebraic stack over $\cY$ that is locally of finite type, quasi-separated, and has separated inertia with affine relative automorphism groups. The stack of graded objects $\Grad(\cM)$ (resp. the stack of filtrations $\Filt(\cM)$) is defined to be the mapping stack $\Map_{\cY}((B\bG_m)_{\cY}, \cM)$ (resp. $\Map_{\cY}(\Theta_{\cY}, \cM)$. The restriction via the inclusion $0/\bG_m \to \Theta$ induces a morphism $\gr:\Filt(\cM) \to \Grad(\cM)$.

There is also an evaluation map $\ev_1: \Filt(\cM)\to \cM$ given by restricting a morphism to the open substack $(\bA^1\setminus \{0\})/\bG_m$ inside $\Theta$. By \cite[Lemma 1.1.13]{hl_instability} this morphism is representable by algebraic spaces. For any $\cY$-scheme $T$ and any $T$-point $\xi: T \to \cM$, the flag space $\Flag(\xi)$ is the fiber over $\xi$ of the evaluation map $\ev_1$. In other words, $\Flag(\xi)$ parametrizes maps $f:\Theta_{T'} \to \cM$ along with an isomorphism
$f(1)\simeq \xi_{T'}$ for any scheme $T'$ over $T$.

Let $\cU$ be an open substack of $\cM$. A weak $\Theta$-stratum of $\cU$ is a union of connected components of $\Filt(\cU)$ such that the restriction of $\ev_1$ is finite and radicial. A $\Theta$-stratum is a weak $\Theta$-stratum such that $\ev_1$ is a closed immersion.
\begin{defn} \label{defn: theta stratification}
A (weak) $\Theta$-stratification of $\cM$ consists of a collection of open substacks $(\cM_{\leq c})_{c \in \Gamma}$ indexed by a totally ordered set $\Gamma$. We require the following conditions to be satisfied
\begin{enumerate}
\item $\cM_{\leq c} \subset \cM_{\leq c'}$ for all $c< c'$.
\item $\cM = \bigcup_{c \in \Gamma} \cM_{\leq c}$.
\item For all $c$, there exists a (weak) $\Theta$-stratum $\cS_c \subset \Filt(\cM_{\leq c})$ of $\cM_{\leq c}$ such that
\[ \cM_{\leq c} \setminus \ev_1(\cS_c) = \bigcup_{c' < c} \cM_{\leq c'}\]
\item For every point $x \in \cM$, the set $\left\{ c \in \Gamma \, \mid \, x \in \cM_{\leq c}\right\}$ has a minimal element.
\end{enumerate}
\end{defn}

One way to construct (weak) $\Theta$-stratifications on a given stack is via the use of numerical invariants. For a given $\cY$-field $k$, a numerical invariant $\mu$ on $\cM$ assigns to any nondegenerate filtration $f:\Theta_{k} \to \cM$ a polynomial $\mu(f) \vcentcolon = \mu(\gr(f))(1)$ in $\mathbb{R}[\epsilon]$. $\mu(f)$ is invariant under field extension, deformation in an algebraic family, and under rescaling, that is, under pre-compositition with the
ramified covering $\Theta_{k}\to \Theta_{k}; z\mapsto z^n$. A point $x\in |\cM|$ is said to be \emph{$\mu$-unstable} if there
is a nondegenerate filtration $f:\Theta_{k}\to \cM$ with $f(1) = x \in |\cM|$ such that $\mu(f)>0$. If $x$ is not unstable, then it is called $\mu$-semistable.

Let $\Gamma$ denote the Dedekind–MacNeille completion of the totally ordered set $\bR[\epsilon]$, i.e., the smallest totally ordered set that contains $\bR[\epsilon]$ and such that any subset of $\Gamma$ has an infimum and a supremum. We define the associated \emph{stability function} $M^\mu : |\cM| \to \Gamma$ by
\[
  M^\mu(x) = \sup \left\{\mu(f) | f\in
	|\Filt(\cM) | , \ f(1) \simeq x \right\}.
\]
A filtration $f$ such that $\mu(f) = M^\mu(f(1))$ is called an \emph{HN filtration}. We say that $\mu$ \emph{defines a (weak) $\Theta$-stratification of $\cM$} if:
\begin{enumerate}
    \item Every unstable point of $\cM$ has a unique HN filtration up to scaling;
    \item $\forall c \in \bR[\epsilon]$, the subset $\{s \in |\cM| \mid M^\mu(x) \leq c\}$ is open, and thus define open substacks $\cM_{\leq c} \subset \cM$; and
    \item $\forall c\in\bR[\epsilon]$, the subset $\{f \in |\Filt(\cM_{\leq c})| \mid \mu(f)=c\}$ is open, and if $\cS_c \subset \Filt(\cM_{\leq c})$ denotes a complete set of orbit representatives of the scaling action of $\bN^\times$ on the set of connected components of this open substack of $\Filt(\cM_{\leq c})$, then $\cS_c$ is a (weak) $\Theta$-stratum in $\cM_{\leq c}$.
\end{enumerate}  
The following theorem provides a criterion for when a numerical invariant defines a (weak) $\Theta$-stratification. It yields a convenient procedure for constructing $\Theta$-stratifications and good moduli spaces.

\begin{thm}[{\cite[Theorem B]{hl_instability}}] \label{thm: theta stability paper theorem}
Let $\mu$ be a numerical invariant on $\cM$, and defined by $\ell / \sqrt{b}$ for an assignment $\ell$ of a rational $\bR[\epsilon]$-valued linear functional and a rational quadratic norm $b$ on graded points.
\begin{enumerate}
    \item If $\mu$ is \emph{strictly $\Theta$-monotone} (\Cref{defn: strictly theta monotone and STR monotone}), then it defines a weak $\Theta$-stratification of $\cM$ if and only if it satisfies the \emph{HN-boundedness} condition (\Cref{defn: HN boundedness}). If moreover $\cY$ is defined over $\bQ$ (i.e. in characteristic $0$), then $\mu$ defines a $\Theta$-stratification.\\
    \item Suppose that all the conditions in (1) above are satisfied. If $\mu$ is \emph{strictly $S$-monotone} (\Cref{defn: strictly theta monotone and STR monotone}),
    then every quasi-compact closed substack of the semistable locus $\cM^{\mu \dash \rm{ss}}$, which is an open substack of $\cM$ by (1), has a relative good moduli space. This relative good moduli space is separated and locally of finite type over $\cY$, and it is proper over $\cY$ if $\cM$ satisfies the existence part of the valuative criterion for properness for complete discrete valuation rings relative to $\cY$ \cite[\href{https://stacks.math.columbia.edu/tag/0CLK}{Tag 0CLK}]{stacks-project}.
\end{enumerate}
\end{thm}

\begin{remark}
\cite[Theorem B]{hl_instability} applies to numerical invariants of the form $\ell/\sqrt{b}$ as described above because they are standard and satisfy property (R).\endnote{To check property (R), choose a field $k$ and a nondegenerate morphism $g: (B\mathbb{G}_m^h)_k \to \cM$. By assumption there are finitely many rational linear functionals $\ell_0, \ell_1, \ell_2, \ldots, \ell_n$ on $\mathbb{R}^h$ such that the associated scale invariant function $\mu: \mathbb{R}^h \setminus 0 \to \mathbb{R}[\epsilon]$ is of the form $\mu = \frac{\ell_0  + \epsilon \ell_1 + \ldots + \epsilon^n \ell_n}{\sqrt{b}}$. We need to check that $\mu$ admits a maximum at a rational point in $\mathbb{Q}^n \subset \mathbb{R}^n$. We proceed by induction on the number of terms $n$. If there are no terms, then the function $\mu$ is identically $0$ and there is nothing to show. For the induction step, we first consider $\frac{\ell_0}{\sqrt{b}}$. A Lagrange multiplier computation (cf. \Cref{lemma:lagrange}) shows that one of the following holds:
\begin{enumerate}
    \item $\frac{\ell_0}{\sqrt{b}}$ admits a unique maximum up to scaling which can be chosen to be rational.
    \item $\ell_0$ is identically $0$.
\end{enumerate}
If (1) holds then we are done, as the unique rational maximum of $\frac{\ell_0}{\sqrt{b}}$ up to scaling must also be a maximum of $\frac{\ell}{\sqrt{b}}$ in our lexicographic ordering of $\mathbb{R}[\epsilon]$. Otherwise, maximizing $\frac{\ell_0 + \epsilon \ell_1 + \ldots + \epsilon^n \ell_n}{\sqrt{b}}$ is equivalent to maximizing $\frac{\ell_1 + \epsilon \ell_2 + \ldots + \epsilon^{n-1} \ell_n}{\sqrt{b}}$, and we conclude by induction.
}
\end{remark}

We end this section by explaining each of the hypotheses that one needs to check in \Cref{thm: theta stability paper theorem}. For this we will need some notation.

\begin{notn} \label{notn: notation for monotonicity}
For any discrete valuation ring $R$ with uniformizer $\varpi$, we define \[Y_{\Theta_R} \vcentcolon = \Spec(R[t]) \text{     and     } Y_{\overline{ST}_R} := \Spec(R[s,t]/(st-\varpi)),\] equipped with the $\bG_m$-action that assigns $t$ weight $-1$ and $s$ weight $1$. The isomorphism class of the $\mathbb{G}_m$-scheme $Y_{\overline{ST}_{R}}$ is independent of the choice of uniformizer $\varpi$. We denote $\Theta_R = Y_{\Theta_R} / \bG_m$ and $\overline{ST}_R = Y_{\overline{ST}_R}/\bG_m$. Note that $Y_{\Theta_R}$ and $Y_{\overline{ST}_R}$ each contain a unique $\bG_m$-invariant closed point cut out by the ideals $(t, \varpi)$ and $(s,t)$ respectively. We denote this closed point by $\closedpt$ in both cases.
\end{notn}

\begin{notn}Let $\kappa$ be a field and let $a \geq 1$ be an integer. We denote by $\mathbb{P}^1_{\kappa}[a]$ the $\mathbb{G}_m$-scheme $\mathbb{P}^{1}_{\kappa}$ equipped with the $\mathbb{G}_m$-action determined by the equation $t \cdot [x:y] = [t^{-a}x : y]$. We set $0 = [0: 1]$ and $\infty = [1:0]$.
\end{notn}

\begin{defn}[Strict monotonicity] \label{defn: strictly theta monotone and STR monotone} A numerical invariant $\mu$ on $\mathcal{M}$ is strictly $\Theta$-monotone (resp. strictly $S$-monotone) relative to $\cY$ if the following condition holds using the stack $\XX = \Theta_R$ (resp. $\XX=\overline{ST}_R$) for any complete discrete valuation ring $R$ over $\cY$:

For any $\cY$-map $\varphi: \XX \setminus \closedpt \rightarrow \mathcal{M}$, after replacing $R$ with a finite DVR extension, there exists a reduced and irreducible algebraic stack $\cW$ with $\cY$-maps $f: \cW \rightarrow Y_{\XX}/\bG_m$ and $\widetilde{\varphi}: \cW \rightarrow  \mathcal{M}$ such that:
\begin{enumerate}
\item The map $f$ is proper with finite relative inertia, and its restriction induces an isomorphism $f : \, \cW_{Y_{\XX} \setminus \closedpt} \xrightarrow{\sim} \XX \setminus \closedpt$;\\

\item The morphism $\widetilde{\varphi}$ has quasi-finite relative inertia and makes the following diagram commute
\[
\begin{tikzcd}
  \cW_{Y_{\XX} \setminus \closedpt} \ar[rd, "\widetilde{\varphi}"] \ar[d, "f"'] & \\   \XX \setminus \closedpt \ar[r, "\varphi"'] &  \mathcal{M}
\end{tikzcd}; \text{ and}
\]\\

\item Let $\kappa$ denote a finite extension of the residue field of $R$. For any $a \geq 1$ and any finite morphism $\mathbb{P}^1_{\kappa}[a]/\bG_m \to \cW_{\closedpt/\bG_m}$ fitting into a commutative diagram
\[
\centering
\begin{tikzcd}
  \mathbb{P}^1_{\kappa}[a]/\bG_m \ar[dr] \ar[r] & \cW_{\closedpt/\bG_m} \ar[d]\\  &  \closedpt/\bG_m
\end{tikzcd}
\]
we have $\mu\left( \;\widetilde{\varphi}|_{\infty / \mathbb{G}_m} \;\right) >  \mu\left(\; \widetilde{\varphi}|_{0 / \mathbb{G}_m} \;\right)$.\\
\end{enumerate}
\end{defn}

\begin{defn}[HN Boundedness] \label{defn: HN boundedness}
We say that a numerical invariant $\mu$ satisfies the \emph{HN boundedness} condition if for any morphism $g: T \rightarrow \mathcal{M}$ from an affine Noetherian scheme $T$ of finite type over $\cY$, there exists a quasi-compact open substack $\mathcal{U} \subset \mathcal{M}$ such that the following holds:
For all geometric points $t : \Spec(k) \to T$ and all nondegenerate $\cY$-filtrations $f: \Theta_{k} \rightarrow \mathcal{M}$ of the point $g(t)$ with $\mu(f)>0$, there exists another filtration $f'$ of $g(t)$ satisfying $\mu(f') \geq \mu(f)$ and $f'(0) \in \mathcal{U}$.
\end{defn}

\begin{rem}
The HN boundedness condition defined in \cite[Thm.~2.2.2(4')]{hl_instability} only requires this condition for \emph{finite type} points of $T$. In this paper, it will be convenient to work with geometric points, and we will end up verifying the a priori stronger condition in \Cref{defn: HN boundedness}.
\end{rem}

Finally, we explain an additional version of monotonicity that we will need.
\begin{notn}
For any field $k$, we define $Y_{\Theta^2_k} = \Spec(k[t,s])$. We equip $Y_{\Theta^2_k}$ with the action of $\mathbb{G}_m^2$ that acts on $t$ via the tuple of weights $(-1,0)$ and on $s$ via $(0,-1)$. There is a unique point $\mathbb{G}^2_m$-invariant closed point in $Y_{\Theta^2_k}$ given by the vanishing of the ideal $(t.s)$, we shall also denote this closed point by $\closedpt$. Note that $Y_{\Theta_k^2} = Y_{\Theta_k^2}/\mathbb{G}_m^2$.
\end{notn}

\begin{defn}[$\Theta^2$-monotonicity] \label{defn: theta square monotonicity}
A numerical invariant $\mu$ on $\mathcal{M}$ is strictly $\Theta^2$-monotone relative to $\cY$ if for any field $k$ over $\cY$, and any $\cY$-map $\varphi: \Theta^2_k \setminus \closedpt \rightarrow \mathcal{M}$, after replacing $k$ with a finite extension, there exists a reduced and irreducible algebraic stack $\cW'$ with $\cY$-maps $f': \cW' \rightarrow Y_{\Theta^2_k}/\bG_m$ and $\widetilde{\varphi}': \cW' \rightarrow  \mathcal{M}$ such that the base change $f: \cW = (\Spec(k\bseries{t}) \times_k \Theta_k) \to \Spec(k\bseries{t}) \times_k \Theta_k$ and the composition $\widetilde{\varphi}: \cW \to \cW' \xrightarrow{\widetilde{\varphi}'} \mathcal{M}$ satisfy the properties (1), (2) and (3) in \Cref{defn: strictly theta monotone and STR monotone}. 

We say that $\nu$ is $\Theta^2$-monotone if we can always find $f'$ and $\widetilde{\varphi}'$ satisfying (1) and (2), and also satisfying the weakening of condition (3) where we replace the strict inequality $>$ with $\geq$.
\end{defn}

We refer the reader to \cite{hl_instability} and \cite[Sect. 2.5]{torsion-freepaper} for more detailed discussion on the material above.



\section{Infinite dimensional GIT for the moduli of gauged maps} \label{section: monotonicity}

For the rest of this section $H$ will denote a choice of smooth relatively affine $S$-group scheme. We assume that $H$ is an extension
\[ 1 \to T \to H \to F \to 1\]
where $T$ is a torus and $F$ is a finite \'etale group scheme over $S$. We shall show that the numerical invariants $\mu_{\delta}$ are strictly monotone on the stack $\cM_{n}^{G}(X)$ by reducing to the case when $G= \GL_N \times H$.

\subsection{Infinite dimensional GIT for \texorpdfstring{$\cM^{\GL_N \times H}_{n}(X)$}{M{GLN x H}n(X)}} \label{subsection: affine grassmannian}

The following stack of ``rational maps" $C \to A_X/(\GL_N \times H)$ is inspired by analogous notions defined in \cite{barlev}.
\begin{defn}
The \emph{stack of rational maps} $\MpCrat(A_X/(\GL_N\times H))$ is the pseudofunctor that sends any $S$-scheme $T \in \Aff_{S}$ to the groupoid of triples $(\mathcal{E} \times \cH, D, s)$, where
\begin{enumerate}
    \item $\mathcal{E} \times \cH$ is $\GL_N \times H$-bundle on $C_T$ consisting of a vector bundle $\mathcal{E}$ of rank $N$ and an $H$-bundle $\cH$.
    \item $D$ is a $T$-flat effective Cartier divisor on $C_T$.
    \item $s: C_{T} \setminus D \to (\mathcal{E} \times \cH)(A_{X})$ is a section of $(\mathcal{E} \times \cH)(A_X) \to C_{T}$ defined away from $D$.
\end{enumerate}
A morphism from $(\mathcal{E}_1\times \cH_1, D_1, s_1)$ into $(\mathcal{E}_2\times \cH_2, D_2, s_2)$ is the data of equality $D_1 = D_2$ and a pair $(\phi, \psi)$ where
\begin{enumerate}
\item $\psi$ is an isomorphism of $H$-bundles $\psi: \cH_2 \xrightarrow{\sim} \cH_1$.
\item $\phi$ is an isomorphism of vector bundles $\phi: \mathcal{E}_1|_{C_{T} \setminus D_1} \to \mathcal{E}_2|_{C_{T} \setminus D_2}$ that is compatible with the sections $s_1$ and $s_2$.
\end{enumerate}
\end{defn}

In the introduction, we discussed a uniformization morphism $\varphi_D : \MpC(A_{X}/(\GL_N \times H))_T \to \MpCrat(A_X/(\GL_N\times H))_T$ for any family of relative Cartier divisors in $C_T$. In fact, all of this data is encoded by a single uniformization morphism
\begin{align*}
q: \MpC(A_{X}/(\GL_N \times H)) \times \Div(C) &\to \MpCrat(A_X/(\GL_N\times H)).\\
(\mathcal{E} \times \cH, s, D) &\mapsto (\mathcal{E}\times \cH, D, s|_{C \setminus D}),
\end{align*}
defined relative to the stack $\Div(C)$ that parametrizes flat families of effective Cartier divisors in $C$. $\Div(C)$ is equipped with a schematic morphism $\Div(C) \to \cY$ represented by unions of quasi-projective schemes. This allows us to define a uniformization morphism
\begin{equation} \label{E:uniformization_morphism}
q \circ (p \times \id) : \cM^{\GL_N \times H}_{n}(X) \times \Div(C) \to \MpCrat(A_X/(\GL_N\times H)),
\end{equation}
where $p : \cM_{n}^{G}(X) \to \MpC(A_X/(\GL_N \times H))$ is the morphism from \eqref{E:morphism_p}. We will ultimately show that
\begin{enumerate}
    \item The fibers of \eqref{E:uniformization_morphism} are ind-projective ind-Deligne-Mumford stacks, and the numerical invariant $\mu$ is constructed using formal line bundles that are relatively ample on these fibers; and
    \item The stack $\MpCrat(A_X/(\GL_N\times H))$ satisfies filling conditions analogous to $\Theta$-reductivity and $S$-completeness.
\end{enumerate}
This will allow us to carry out the strategy for establishing strict monotonicity described in \Cref{E:infinite_dimensional_GIT}.


\subsubsection{Generalized affine Grassmannians}
We shall prove that the fibers of $q$ are closed sub-ind-schemes of the Beilinson-Drinfeld Grassmannian for $\GL_N$. We will use the notation $A = A_X$ to stress that in the following discussion, $A$ can be any choice of relatively affine scheme of finite type over $S$, equipped with an action of $\GL_N \times H$.
\begin{defn}
Let $T$ be an $S$-scheme of finite type. Choose a morphism $g: T \to \MpCrat(A/(\GL_N \times H))$ represented by a triple $(\mathcal{E}\times \cH, D, s)$. We define the generalized BD Grassmanian $\Gr_{A/(\GL_N \times H)}^{(\mathcal{E}\times \cH, D, s)}$ to be the fiber product $(\MpC(A/(\GL_N \times H)) \times \Div(D)) \times_{\MpCrat(A/(\GL_N\times H))} T$.
\end{defn}
 By definition, $\Gr_{A/(\GL_N \times H)}^{(\mathcal{E}\times \cH, D, s)}$ is the pseudofunctor that sends $Q \in \Aff_{T}$ to the set of isomorphism classes of triples $(\mathcal{F}, \phi, s)$ where
\begin{enumerate}
    \item $\mathcal{F}$ is a vector bundle of rank $N$ on $C_{Q}$.
    \item $\phi$ is an isomorphism $\varphi: \mathcal{E}|_{C_{Q} \setminus D_{Q}}\xrightarrow{\sim} \mathcal{F}|_{C_{Q} \setminus D_{Q}}$.
    \item $\widetilde{s}: C_{Q} \to (\mathcal{F} \times \cH|_{C_{Q}})(A)$ is an extension of the section
    \[ C_{Q} \setminus D_{Q} \xrightarrow{s_{Q}} (\mathcal{E} \times \cH)|_{C_{Q} \setminus D_{Q}}(A) \xrightarrow{\sim} (\mathcal{F} \times \cH)|_{C_{Q} \setminus D_{Q}}(A)\]
with the last isomorphism induced by the identification $\phi: \cE|_{C_{Q} \setminus D_{Q}} \xrightarrow{\sim} \cF|_{C_{Q} \setminus D_{Q}}$.
\end{enumerate}
There is a projection morphism
\[\Gr_{A/(\GL_N \times H)}^{(\mathcal{E}\times \cH, D, s)} \to \MpC(A/(\GL_N \times H)), \; \; \; \; \; \; (\mathcal{F}, \phi,\widetilde{s}) \mapsto (\mathcal{F}\times (\cH|_{C_{Q}}), \widetilde{s}) \]

\begin{defn} \label{defn: V standard} Let $V_{\std}$ the $N$-dimensional almost faithful representation of $\GL_N \times H$ that is trivial on $H$ and restricts to the standard representation on the $\GL_N$ factor. We equip $\Gr_{A/(\GL_N \times H)}^{(\mathcal{E}\times \cH, D, s)}$ with the pullback of the line bundle $\cD(V_{\rm{std}})$ under the projection above.
\end{defn}

We study first the special case when $H$ is trivial and $A =S$. In this case the $T$-point $g \in \MpCrat(A/(\GL_N \times H))(T)$ amounts to a pair $(\mathcal{E}, D)$ of a vector bundle on $C_{T}$ and a relative Cartier divisor $D \hookrightarrow C_T$.
\begin{lem} \label{lemma: affine grassmannian for GLn}
$\Gr_{B\GL_N}^{(\mathcal{E}, D)}$ is represented by a strict-ind-scheme that is ind-projective over $T$. The line bundle $\cD(V_{\rm{std}})$ is $T$-relatively ample on each projective stratum of $\Gr_{B\GL_N}^{(\mathcal{E}, D)}$.
\end{lem}
\begin{proof}
$D$ is a relative effective Cartier divisor, defined by a section $\sigma \colon \cO_{C_{T}} \to \cO_{C_{T}}(D)$. For any $M > 0$ we define the subfunctor
\[\Gr_{B\GL_N}^{(D,\cE), \leq M} (T'/T) := \left\{ \begin{array}{c} (\cF,\phi) \in \Gr_{B\GL_N}^{(\mathcal{E}, D)}(T') \text{ such that $\cF \otimes \cO_{C_{T'}}[\sigma^{\pm 1}] \xrightarrow{\phi} \cE|_{C_{T'}} \otimes \cO_{C_{T'}}[\sigma^{\pm 1}]$} \\ \text{satisfies} \; (\cE\cdot \sigma^{M})|_{C_{T'}} \subset \phi(\cF) \subset (\cE \cdot \sigma^{-M})|_{C_{T'}}. \end{array} \right\}
\]
For any point $(\cF, \phi)$ in $\Gr_{B\GL_N}^{(D, \cE)}$, the morphism $\phi$ is encoded by the embedding $\cF \hookrightarrow \cE|_{C_{T'}} \otimes \cO_{C_{T'}}[\sigma^{\pm 1}]$. For any finitely presented sheaf $(\cE\cdot \sigma^{M})|_{C_{T'}} \subset \cF \subset (\cE\cdot \sigma^{-M})|_{C_{T'}}$, the short exact sequence
\[
0 \to \cF / (\cE\cdot \sigma^{M})|_{C_{T'}}  \to \cE[\sigma^{\pm 1}]|_{C_{T'}} / (\cE\cdot \sigma^{M})|_{C_{T'}}  \to (\cF[\sigma^{\pm 1}]) / \cF\to 0
\]
implies that $\cF / \cE|_{C_{T'}}$ is flat over $T'$, and that $\cF$ is locally free (because its restriction to every fiber is torsion free). This shows that the data of a vector bundle $\cF$ which fits into a sequence of inclusions $(\cE\cdot \sigma^{M})|_{C_{T'}} \subset \cF \subset (\cE \cdot \sigma^{-M})|_{C_{T'}}$ is equivalent to specifying a flat family of quotients of $(\cE\cdot \sigma^{-M} / \cE \cdot \sigma^{M} ) |_{C_{T'}}$. Hence, we have
\[
\Gr_{B\GL_N}^{(D,\cE), \leq M} \simeq \Quot_{C_{T}/T}(\cE\cdot \sigma^{-M} / \cE \cdot \sigma^{M} ).
\]
This is represented by a projective scheme over $T$. To simplify notation, we will denote this Quot scheme by $Q$. Let 
$$(\cE\cdot \sigma^{-M} / \cE \cdot \sigma^{M} )|_{C_{T} \times_T Q} \to \cG$$
be the universal quotient sheaf on $C_{T} \times_T Q$. Set $\pi_{Q}: C_{T} \times_T Q \to Q$. The line bundle $\det\left((\pi_{Q})_\ast\left(\cG \otimes \sqrt{K_{C_{Q}}}\right)\right)$ is relatively ample for the morphism $Q \to T$.\endnote{In order to see this it suffices to check this on every fiber of $T$ \cite[Cor. 9.6.5]{egaivIII}, so we can assume that $T$ is $\Spec(k)$ for some field $k$ over $S$. In that case we have that $2MD$ is a divisor on the one-dimensional curve $C_{k}$, and hence $2MD$ is Artinian. The restrictions of the line bundles $K_{C_k}$ and $\mathcal{O}(1)$ to $2MD$ are therefore trivial. This shows that the restrictions of $\mathcal{O}(1)$ and $K_{C_{Q}}$ to $Q \times_{k} 2MD$ are trivial. Since $\mathcal{G}$ is supported on $Q \times_{k} 2MD$, this implies that
 \[ (\pi_{Q})_\ast\left(\cG \otimes \sqrt{K_{C_{Q}}}\right) \cong (\pi_{Q})_\ast(\cG) \cong  (\pi_{Q})_\ast(\cG(n)) \]
    for any $n$-twist. By \cite[Prop. 2.2.5]{huybrechts.lehn}, $\det (\pi_{Q})_\ast(\cG(n)))$ is a relatively ample line bundle for $n \gg 0$, thus concluding the proof of the claim that $\det\left((\pi_{Q})_\ast\left(\cG \otimes \sqrt{K_{C_{Q}}}\right)\right)$ is a relatively ample for $Q \to T$.} For every finitely presented sheaf $\mathcal{F} \subset \cE[\sigma^{\pm}]|_{C_{T'}}$ there exists some $M$ such that $(\cE\cdot \sigma^{M})|_{C_{T'}} \subset \cF \subset (\cE \cdot \sigma^{-M})|_{C_{T'}}$, and so we have $\Gr_{B\GL_N}^{(\mathcal{E}, D)} \simeq \colim_{M} \Gr_{B\GL_N}^{(D, \mathcal{E}), \leq M}$. The connecting maps between different $\Gr_{B\GL_N}^{(\mathcal{E}, D), \leq M}$ are forced to be closed immersions, since they are proper monomorphisms. There is a canonical short exact sequence
\[
0 \to \cF_{univ} \to \cE \otimes \cO_{C_{T}}(-M D) |_{C_{T} \times_T Q} \to \cG \to 0
\]
Tensoring with $\sqrt{K_{C_{Q}}}$, applying pushforward and determinant, and using the base change formula, we see that
\begin{gather*}
\det\left((\pi_{Q})_\ast\left(\cF_{univ} \otimes \sqrt{K_{C_{Q}}}\right)\right)^\dual \simeq \det\left((\pi_{Q})_\ast\left(\cG \otimes \sqrt{K_{C_{Q}}}\right)\right) \otimes \left(\det\left(R\pi_\ast\left(\cE \otimes \cO_{C_{T}}(-M D) \otimes \sqrt{K_{C_{T}}}\right)\right)\right)|_{C_{T} \times_T Q}.
\end{gather*}
In particular, $\cD(V_{\rm{std}})$ restricted to $\Gr_{B\GL_N}^{(D,\cE), \leq M} \simeq Q$ is isomorphic to a twist of the relatively ample bundle $\det\left((\pi_{Q})_\ast\left(\cG \otimes \sqrt{K_{C_{Q}}}\right)\right)$ by a line bundle which is pulled back from $T$. Thus $\cD(V_{\rm{std}})^{\vee}$ is relatively ample on $\Gr_{B\GL_N}^{(D,\cE), \leq M}$.
\end{proof}
For any finite type $\cY$-scheme $T$, and a $T$-point $g = (\cE \times \cH, D, s) \in \MpCrat(A/(\GL_N\times\cH))(T)$, there is a morphism of BD Grassmannians
\[ \Forget:  \Gr_{A/(\GL_N \times H)}^{(\mathcal{E}\times \cH, D, s)} \to \Gr_{B\GL_N}^{(\cE, D)}, \; \; \; \; \; \; (\mathcal{F}, \phi,s) \mapsto (\mathcal{F}, \phi) \]
The pullback of the line bundle $\cD(V_{\rm{std}})$ on $\Gr_{B\GL_N}^{(\mathcal{E}, D)}$ coincides with $\cD(V_{\rm{std}})$ on $\Gr_{A/(\GL_N \times H)}^{(\mathcal{E}\times \cH, D, s)}$.
\begin{lem} \label{lemma: forgetful morphism of affine grassmannians is closed}
The morphism $Forget$ is a closed immersion of strict ind-schemes. In particular, $\Gr_{A/(\GL_N \times H)}^{(\cE \times \cH, D, s)}$ is represented by an ind-scheme that is ind-projective over $T$, and $\cD(V_{\rm{std}})$ is $T$-relatively ample on each projective stratum.
\end{lem}
\begin{proof}
Let $Q \in \Aff_{T}$. Choose a morphism $f: Q \to \Gr_{B\GL_N}^{(\mathcal{E}, D)}$ corresponding to a pair $(\mathcal{F}, \phi)$. The fiber product $\Gr_{A/(\GL_N \times H)}^{(\mathcal{E}\times \cH, D, s)} \times_{\Gr_{B\GL_N}^{(\mathcal{E}, D)}} Q$ is the functor classifying extensions $\widetilde{s}: C_{Q} \to (\mathcal{F} \times \cH|_{C_{Q}})(A)$ of the section
\[ C_{Q} \setminus D_{Q} \xrightarrow{s_{Q}} (\mathcal{E} \times \cH)|_{C_{Q} \setminus D_{Q}}(A) \xrightarrow{\sim} (\mathcal{F} \times \cH)|_{C_{Q} \setminus D_{Q}}(A)\]
By \cite[Lemma 4.12]{rho-sheaves-paper} this functor is represented by a closed subscheme of $Q$. The ampleness of $\cD(V_{\rm{std}})$ now follows because $\cD(V_{\rm{std}})$ is $T$-ample on each $\Gr_{B\GL_N}^{(\mathcal{E}, D), \leq M}$ by Lemma \ref{lemma: affine grassmannian for GLn}.
\end{proof}

\subsection{Monotonicity for \texorpdfstring{$\MpC(A/ G)$}{MC(A/G)}}
For this subsection, $A= A_X$ will denote a $G$-scheme that is affine and of finite type over $S$. We fix an almost faithful representation $G \to \GL(V)$, and a rational quadratic norm $b$ pulled back from a quadratic norm on $BG$ by restricting to the generic point of $C$ (Definition \ref{defn: quadratic norm graded points}).
\begin{prop} \label{prop: monotonicity of mapping stack general group}
The numerical invariant $\mu := -\wt(\cD(V)) / \sqrt{b}$ on $\MpC( A/G)$ is strictly $\Theta$-monotone, strictly $S$-monotone, and strictly $\Theta^2$-monotone.
\end{prop}
We first prove a special case of the proposition, and we use this to prove the general statement at the end of this subsection.
\begin{proof}[Proof of \Cref{prop: monotonicity of mapping stack general group} for $G = \GL_N \times H$ and $V = V_{\rm std}$]
We first prove strict $\Theta$-monotonicity and strict $S$-monotonicity. Let $\XX = Y_\XX / \bG_m$ denote either $\Theta_{R}$ or $\overline{ST}_{R}$ for some complete discrete valuation ring $R$ over $\cY$, and set $\cM = \MpC( A/(\GL_N \times H))$. Choose a $\cY$-morphism $\varphi: \XX \setminus \closedpt \to \mathcal{M}$ represented by a pair $(\mathcal{F} \times \cH, s)$, where $\mathcal{F} \times \cH$ is a $\mathbb{G}_m$-equivariant $\GL_N \times H$-bundle on $C_{Y_{\XX} \setminus \closedpt}$ and $s$ is a $\mathbb{G}_m$-equivariant section $s: C_{Y_{\XX} \setminus \closedpt} \to (\mathcal{F} \times \cH)(A)$. The composition 
\[\XX \setminus \closedpt \xrightarrow{\varphi} \mathcal{M} \to \Bun_{\GL_N}(C) \times \Bun_{H}(C)\]
amounts to the data of two morphisms $\varphi_1: \XX \setminus \closedpt \to \Bun_{\GL_N}(C)$ and $\varphi_2: \XX \setminus \closedpt \to \Bun_{H}(C)$.

By the rational filling condition in \cite[Lemma 4.7]{torsion-freepaper} there is a $\bG_m$-equivariant locally free sheaf $\cE$ on $C_{Y_\XX}$, a $\bG_m$-invariant $Y_{\XX}$-flat effective Cartier divisor $D \hookrightarrow C_{Y_\XX}$ with complement $U = C_{Y_\XX} \setminus D$, and an isomorphism $\psi : \cF|_{U \setminus U_{\closedpt}} \xrightarrow{\sim} \cE|_{U \setminus U_{\closedpt}}$. By \Cref{lemma: stack of h-bundles is theta red and s complete}, $\Bun_{H}(C) \to \cY$ is $\Theta$-reductive and S-complete, so $\varphi_2$ extends uniquely to morphism $\XX \to \Bun_{H}(C)$, which represents a $\mathbb{G}_m$-equivariant $H$-bundle $\widetilde{\cH}$ on $C_{Y_{\XX}}$ extending $\cH$.

The isomorphism $\psi$ allows us to identify $s|_{U \setminus U_{\closedpt}}$ as a section $s|_{U \setminus U_{\closedpt}} : U \setminus U_{\closedpt} \to (\mathcal{E} \times \widetilde{\cH})(A)_{U}$. The projection $(\mathcal{E} \times \widetilde{\cH})(A)_{U} \to U$ is affine and of finite type. All points contained in the codimension $2$ subscheme $U_{\closedpt}$ have depth at least $2$ in $U$, because $U$ is regular. It follows that $s|_{U \setminus U_{\closedpt}}$ extends to a unique $\mathbb{G}_m$-equivariant section $\tilde{s}: U \to (\mathcal{E} \times \widetilde{\cH})(A)_U$ (see, e.g., the proof of \cite[Lemma 2.2 (c)]{rho-sheaves-paper}). We regard the triple $(\cE \times \tilde{\cH}, D, \tilde{s})$ as a $\bG_m$-equivariant morphism $g : Y_{\XX} \to \MpCrat(A/(\GL_N \times H))$. The isomorphism $\psi$ identifies the composition 
\[Y_{\XX} \setminus \closedpt \to \MpC(A/(GL_N \times H)) \times \Div(C) \to \MpCrat(A/(\GL_N \times H))\] with $g|_{Y_\XX \setminus \closedpt}$, so we get a commutative diagram:
\[
\begin{tikzcd}
  Y_{\XX} \setminus \closedpt \ar[ddr, bend right=20] \ar[rrd, bend left=20, "{(\varphi, D)}"]  \ar[dr, "{(\cF, \psi)}"] &  & \\ & \Gr_{A/(\GL_N\times H)}^{(\cE \times \tilde{\cH}, D, \tilde{s})} \ar[d] \ar[r] & \MpC( A/(GL_N \times H)) \times \Div(C) \ar[d] \\
  & Y_{\XX} \ar[r, "g"] & \MpCrat(A/(\GL_N \times H))
\end{tikzcd}
\]

By \Cref{lemma: forgetful morphism of affine grassmannians is closed}, $\Gr_{A/(\GL_N \times H)}^{(\cE \times \tilde{\cH}, D, \tilde{s})}$ is an ind-scheme that is ind-projective over $Y_{\XX}$. It admits a natural $\mathbb{G}_m$-action such that the projection to $Y_{\XX}$ is $\bG_m$-equivariant, because the data used to define it is $\mathbb{G}_m$-equivariant. Let $\Sigma$ denote the scheme theoretic image of the $\mathbb{G}_m$-equivariant morphism labeled $(\cF,\psi)$ in the diagram above. $f: \Sigma \to Y_{\XX}$ is a reduced projective $\mathbb{G}_m$-scheme over $Y_{\XX}$. By construction, the restriction $f|_{Y_{\XX} \setminus \closedpt}$ is an isomorphism. We set $\cW = \Sigma / \bG_m$ and define $\widetilde{\varphi} : \cW \to \MpC( A/(\GL_N \times H))$ to be induced by the $\mathbb{G}_m$-equivariant composition 
\[\Sigma \hookrightarrow \Gr_{A/(\GL_N \times H)}^{(\cE \times \tilde{\cH},D,\tilde{s})} \to \MpC( A/(\GL_N \times H)).\]
To complete the proof of monotonicity, we are only left to check condition (3) in \Cref{defn: strictly theta monotone and STR monotone}. Let $\kappa$ be a finite extension of the residue field, and let $a$ be a positive integer. Choose a finite $\bG_m$-equivariant morphism $\bP^1_{\kappa}[a] \to \Sigma_{\closedpt}$. We shall show that $\nu\left(\widetilde{\varphi}|_{\infty  / \mathbb{G}_m}\right) > \nu\left(\widetilde{\varphi}|_{0  / \mathbb{G}_m}\right)$, where
\[ \nu\left(\widetilde{\varphi}|_{\infty  / \mathbb{G}_m}\right) = \frac{-\wt( \cD(V_{\rm{std}})|_{\infty})}{\sqrt{b\left(\widetilde{\varphi}|_{ \infty  / \mathbb{G}_m} \right)}} \qquad \text{and} \qquad \nu\left(\widetilde{\varphi}|_{ 0 / \mathbb{G}_m}\right) = \frac{-\wt(\cD(V_{\rm{std}})|_{0})}{\sqrt{b\left(\widetilde{\varphi}|_{  0  / \mathbb{G}_m} \right)}}.\]

By \Cref{lemma: affine grassmannian for GLn} and \Cref{lemma: forgetful morphism of affine grassmannians is closed}, the line bundle $\cD(V_{\rm{std}})$ is $Y_{\XX}$-relatively ample on each projective stratum of $\Gr_{A/(\GL_N \times H)}^{(\cE \times \tilde{\cH},D,\tilde{s})}$, and hence $Y_{\XX}$-relatively ample on the closed projective subscheme $\Sigma \subset \Gr_{A/(\GL_N \times H)}^{(\cE \times \tilde{\cH},D,\tilde{s})}$. It follows that $\cD(V_{\rm{std}})|_{\mathbb{P}^1_{\kappa}}$ is ample, and so $\cD(V_{\rm{std}})|_{\mathbb{P}^{1}_{\kappa}} \cong \mathcal{O}_{\mathbb{P}^1_{\kappa}}(n)$ for some $n >0$. Note that $\mathcal{O}_{\mathbb{P}_{\kappa}^1}(n)$ admits a unique $\mathbb{G}_m$-equivariant structure up to twisting by a character, and we have $-\wt \,\mathcal{O}_{\mathbb{P}^1_{\kappa}}(n)|_{\infty} > -\wt \,\mathcal{O}_{\mathbb{P}^1_{\kappa}}(n)|_{0}$. Hence, we will be done if we can prove $b\left(\widetilde{\varphi}|_{\infty  / \mathbb{G}_m} \right) =b\left(\widetilde{\varphi}|_{ 0  / \mathbb{G}_m} \right)$. By construction, both morphisms $\widetilde{\varphi}|_{\infty / \mathbb{G}_m}$ and $\widetilde{\varphi}|_{0  / \mathbb{G}_m} : (B\bG_m)_\kappa \to \MpC(A/(\GL_N \times H)) \times \Div(C)$ agree after composition with the morphism $\MpC(A/(\GL_N \times H)) \times \Div(C) \to \MpCrat(A/(\GL_N \times H))$. This means that, when viewed as morphisms $C_{\kappa} \times (B\bG_m) \to A/(\GL_n \times H)$, both coincide away from the divisor $D_{\kappa} \hookrightarrow C_{\kappa}$. In particular, the restrictions to the generic point of $C_{\kappa}$ agree, and therefore  $b\left(\widetilde{\varphi}|_{ \infty  / \mathbb{G}_m} \right) =b\left(\widetilde{\varphi}|_{ 0  / \mathbb{G}_m} \right)$ by definition.

The proof of strict $\Theta^2$-montonicity follows from the same argument, replacing the scheme by $Y= \mathbb{A}^2_k$ and the group by $\mathbb{G}_m^2$.
\end{proof}

We now proceed to prove Proposition \ref{prop: monotonicity of mapping stack general group} in the case of an arbitrary geometrically reductive group scheme $G$. We will need the following lemma.
\begin{lem} \label{lemma: zariski's main theorem stacks}
 Let $U$ and $\cM$ be Noetherian algebraic stacks. Let $f: U \to \cM$ be a representable morphism that is quasi-finite and separated. Then, there exists a stack $\cM'$ finite over $\cM$ and a factorization
\[
\centering
\begin{tikzcd}
  & \cM'  \ar[rd] & \\   U \ar[ru, dashrightarrow] \ar[rr, "f"] &  & \cM
 \end{tikzcd}
\]
 such that the morphism $U \to \cM'$ is an open immersion.
 \end{lem}
 \begin{proof}
 This is a special case of \cite[Thm. 8.6 (ii)]{rydh-noetherian-approximation}.
 \end{proof}For the purpose of monotonicity, we can restrict to the case when $\cY = S$ is the spectrum of a complete discrete valuation ring. After passing to a finite extension of this DVR, we assume without loss of generality that $G$ admits a maximal split torus over $S$. In particular, the norm $b$ on $BG$ comes from a Weyl invariant quadratic norm on the group of real cocharacters of a maximal split torus inside $G$ (Lemma \ref{lemma: description of norm in the split case}).

The vector bundle $V$ on $S$ is trivializable, and so we can identify $\GL(V) \cong \GL_N$ for some $N$. We denote by $K$ the kernel of the composition $\D(G) \to G \to \GL(V)$. The fppf quotient $\overline{G} = G / K$ is a geometrically reductive group scheme over $S$. The morphism $G \to \overline{G}$ induces an isomorphism on real cocharacters of maximal split tori, and therefore our choice of quadratic norm $b$ of graded points of $BG$ comes from a quadratic norm $b$ on $\overline{G}$. The group scheme $\overline{G}$ acts naturally on the affine GIT quotient $A/\!/K := \Spec(\cO_A^K)$, and the morphism $A \to A/\!/K$ is equivariant with respect to the quotient homomorphism $G \to G/K$.\endnote{ By \cite[\href{https://stacks.math.columbia.edu/tag/03BM}{Tag 03BM}]{stacks-project} we have $\overline{G} = \uSpec_{S}(\mathcal{O}_{G}^K)$, where $\cO_{G}^K \subset \cO_{G}$ denotes the sub $\cO_{S}$-algebra of $K$-invariants. We have that $\overline{G}$ is a geometrically reductive group scheme over $S$ by the criterion in \cite[Thm. 9.7.6]{alper-adequate-moduli}. Indeed smoothness can be checked after passing to the fppf cover $G \to G/K$. We have that $(G/K)_0 = G_0/K$ has connected reductive fibers, because this is true for quotients over fields and forming fppf quotient commutes with passing to fibers. Moreover, $\pi_0(G/K) = \pi_0(G)$ is finite over $S$.

We define $A / \! / K \vcentcolon = \uSpec_{S}(\cO_{A}^K)$, where $\cO_{A}^K \subset \cO_{A}$ denotes the sub $\cO_{S}$-algebra of $K$-invariants of $\cO_{A}$. The $\cO_{G}$-comodule structure on $\cO_{A}$ restricts to a $\cO_{\overline{G}}$-comodule structure on $\cO_{A}^K$, thus inducing an action on $A / \! / K$. By construction, the natural morphism $A \to A / \! / K$ induced by the inclusion $\cO_{A}^K \hookrightarrow \cO_{A}$ is equivariant for the actions of $G$ and $\overline{G}$ respectively via the quotient morphism $G \to \overline{G}$.} This induces a morphism of quotient stacks $A/ G \to (A/\!/K)/\overline{G}$, which in turn induces a morphism $\MpC( A/G) \to \MpC( (A/\!/K)/\overline{G})$.
\begin{lem} \label{lemma: finiteness of map to isogenous group}
The morphism $\MpC( A/G) \to \MpC( (A/\!/K)/ \overline{G})$ is quasi-finite and proper.
\end{lem}
\begin{proof}
The statement is \'etale local on the base $S$, hence we can assume that $G$, $\cD(G)$ and $\overline{G}$ are split \cite[Lemma 5.1.3]{conrad_reductive}. $K$ is a $R$-flat closed subgroup scheme of any maximal split torus $\cT \subset \D(G)$, and hence $K$ is of multiplicative type by \cite[B.3.3]{conrad_reductive}. We can further assume $K$ is diagonalizable, and so we can write $K \cong \prod_{i} \mu_{n_i}$ for some finite tuple of positive integers $(n_i)$.

Factor the morphism $\MpC( A/G) \to \MpC( (A/\!/K)/\overline{G})$ as
\[ \MpC( A/G) \; \to \; \MpC( (A/\!/K)/ \overline{G}) \times_{\Bun_{\overline{G}}(C)} \Bun_{G}(C) \to \;  \MpC( (A/\!/K)/ \overline{G})\]
The right-most morphism is quasi-finite and proper by \Cref{lemma: finiteness stack of bundles under isogeny} in the appendix. We shall conclude by proving that $\MpC( A/G) \; \to \; \MpC( (A/\!/K)/ \overline{G}) \times_{\Bun_{\overline{G}}(C)} \Bun_{G}(C)$ is schematic and finite. Choose a point $T \to \MpC( (A/\!/K)/ \overline{G}) \times_{\Bun_{\overline{G}}(C)}$ consisting of a $G$-bundle $\cP$ on $C_{T}$ and a section $s$ of the associated $A/\!/K$-fiber bundle morphism $(A/ \!/ K)(\cP \times^{\overline{G}} \overline{G}) \to C_{T}$. Form the fiber product
\[
\begin{tikzcd}
  Y \ar[d] \ar[r] & \cP(A) \ar[d] & \\ C_{T} \ar[r, "s"]& (\cP \times^{\overline{G}} \overline{G})(A/ \!/ K)
\end{tikzcd}
\]
The fiber of $\MpC( A/G) \to \MpC( (A/\!/K)/ \overline{G}) \times_{\Bun_{\overline{G}}(C)} \Bun_{G}(C)$ over $T$ parametrizes sections of the morphism $Y \to S$, in other words it is given by the functor $\Gamma_{C_{T}}(Y)$ as in Lemma \ref{lemma: sections of finite morphism finite}. By Lemma \ref{lemma: sections of finite morphism finite}, in order to prove that this $T$-fiber is finite it is sufficient to show that $Y \to C_{T}$ is a finite morphism. By the fiber product diagram defining $Y$, we can reduce to showing that $\cP(A) \to (\cP \times^{\overline{G}} \overline{G})(A/ \!/ K)$ is a finite morphism. Since finiteness can be checked \'etale locally, it suffices to prove that $A \to A/\!/K$ is finite. This follows because $A \to A/\!/K$ is integral \cite[\href{https://stacks.math.columbia.edu/tag/03BJ}{Tag 03BJ}]{stacks-project} and of finite type \endnote{Because $A$ is of finite type over $S$ and $A/\!/K$ is separated over $S$.}.
\end{proof}

\begin{proof}[Proof of Proposition \ref{prop: monotonicity of mapping stack general group}]
We  first study $\MpC( (A/\!/K)/ \overline{G})$. Set $H \vcentcolon = G/\D(G)$, which fits into an extension diagram
\[ 1 \to G_0/\D(G) \to H \to \pi_0(G) \to 1\]
The natural homomorphism $\overline{G} \to \GL(V) \times H$ is a monomorphism; it induces a schematically free action of $\overline{G}$ on $\GL(V) \times H$ via left multiplication. The quotient $Z = (A/\!/K) \times^{G} (\GL(V) \times H)$ is a relatively affine scheme over $S$ \cite[Thm. 9.4.1]{alper-adequate-moduli}, and is equipped with an action of $\GL(V) \times H$ via right multiplication. Under the natural identification of stacks $(A / \!/ K) / \overline{G} \cong Z / ( \GL(V) \times H)$, the line bundle $\cD(V_{\rm{std}})$ on $\MpC( Z/ ( \GL(V) \times H))$ (Definition \ref{defn: V standard}) corresponds to the line bundle $\cD(V)$. 

Let us prove strict $\Theta$-monotonicity and strict $S$-monotonicity. Let $\XX$ denote either $\Theta_{R}$ or $\overline{ST}_{R}$ for some complete discrete valuation ring $R$ over $\cY$. Choose a $\cY$-morphism $\varphi: Y_{\XX} \setminus \closedpt \to \MpC( A/G)$. Consider the composition $\psi: Y_{\XX} \setminus \closedpt \to \MpC( A/G) \to \MpC( (A/\!/K )/\overline{G})$. Using the identification $(A / \!/ K) / \overline{G} \cong Z / ( \GL(V) \times H)$, the proof of monotonicity in the special case $G = \GL_N \times H$ yields an irreducible and reduced algebraic stack $\cW$ with morphisms $p: \cW \to \XX$ and $\widetilde{\psi}: \cW \to \MpC( (A/\!/K)/\overline{G})$ satisfying conditions (1)-(3) in the definition of monotonicity (Definition \ref{defn: strictly theta monotone and STR monotone}). Consider the naturally induced morphism $h: \XX \setminus \closedpt \to \MpC( A/G) \times_{\MpC( (A/\!/K)/\overline{G})} \cW$. Since the morphism $\MpC( A/G) \to \MpC( (A/\!/K)/\overline{G})$ is quasi-finite and proper, it follows that $\MpC( A/G) \times_{\MpC( (A/\!/K)/\overline{G})} \cW$ is quasi-finite and proper over $\cW$, and so in particular it has finite relative diagonal. Note that $h$ factors as a composition
\begin{gather*} \XX \setminus \closedpt \; \to \; \MpC( A/G) \times_{\MpC( (A/\!/K)/\overline{G})} (\XX \setminus \closedpt) \; \hookrightarrow \; \MpC( A/G) \times_{\MpC( (A/\!/K)/\overline{G})} \cW 
\end{gather*}
The second morphism is an open immersion. The first morphism is representable and finite, since it is a section of a morphism of stacks with finite relative diagonal. In particular, $h$ is representable, quasi-finite and separated. Lemma \ref{lemma: zariski's main theorem stacks} yields a diagram
\[
\begin{tikzcd}
  & \cM'  \ar[dr] & \\ \XX \setminus \closedpt  \ar[ru, "g"] \ar[rr, "h"] & & \MpC( A/G) \times_{\MpC( (A/\!/K)/\overline{G})} \cW
\end{tikzcd}
\]
where $\cM' \to \MpC( A/G) \times_{\MpC( (A/\!/K)/\overline{G})} \cW$ is schematic and finite and $g$ is an open immersion. Define $\cW'$ to be the stack theoretic image of $g$. By construction, $\cW'$ admits a morphism $\widetilde{\varphi}$ to $\MpC( A/G)$ and the natural projection $\cW' \to \cW$ is proper and quasi-finite. Moreover, since $g$ is an open immersion, the natural projection $\cW' \to \XX$ restricts to an isomorphism over $\XX \setminus \closedpt$. We are left to check condition (3) in Definition \ref{defn: strictly theta monotone and STR monotone}. By construction, we have a commutative diagram
\[
\begin{tikzcd}
  \cW'  \ar[d] \ar[r, "\widetilde{\varphi}"] & \MpC( A/G) \ar[d] \\ \cW  \ar[r] & \MpC( (A/\!/K)/\overline{G})
\end{tikzcd}
\]
Recall that the numerical invariant $\mu = -\wt(\cD(V))/ \sqrt{b}$ on $\MpC( A/G)$ is pulled back from the corresponding invariant on $\MpC( (A/\!/K)/\overline{G})$, because both the line bundle $\cD(V)$ and the norm $b$ (coming from $B\overline{G}$) are pulled back from $\MpC( (A/\!/K)/\overline{G})$.

Let $\kappa$ be an extension of the residue field of $R$, and choose a finite morphism $\bP^1_{\kappa}[a]/\bG_m \to \cW'_{\closedpt/\bG_m}$ fitting into a commutative diagram
\[
\begin{tikzcd}
  \mathbb{P}^1_{\kappa}[a]/\bG_m \ar[drr] \ar[r] & \cW'_{\closedpt/\bG_m} \ar[r] &\cW_{\closedpt/\bG_m} \ar[d] \\ &  &  \closedpt/\bG_m
\end{tikzcd}
\]
Since $\mathbb{P}^1_{\kappa}[a]/\mathbb{G}_m \to \closedpt/\mathbb{G}_m$ is representable, it follows that $\mathbb{P}^1_{\kappa}[a]/\mathbb{G}_m \to \cW_{\closedpt/\bG_m}$ is representable. Since the morphism $\cW_{\closedpt/\bG_m}' \to \cW_{\closedpt/\bG_m}$ is proper and quasi-finite, the morphism $\mathbb{P}^1_{\kappa}[a]/\mathbb{G}_m \to \cW_{\closedpt/\bG_m}$ is finite. We can therefore conclude monotonicity by using strict monotonicity of the numerical invariant on $\cW_{\closedpt/\bG_m}$, since $\mu$ is pulled back via the composition $\cW' \to \cW \to \MpC( (A/\!/K)/\overline{G})$.

The proof of strict $\Theta^2$-monotonicity is exactly the same, by reducing to the case when $G = \GL_N \times H$.
\end{proof}

\subsection{Monotonicity for \texorpdfstring{$\cM^G_{n}(X)$}{MGn(X)}}
\begin{thm} \label{thm: monotonicity of gauged maps}
For any pair $\delta = (\delta_{\mrk}, \delta_{\gen})$ of nonnegative real numbers, the numerical invariant $\mu_{\delta}$ on $\cM_{n}^{G}(X)$ defined in \Cref{subsection: numerical invariants} is strictly $\Theta$-monotone and strictly S-monotone. Furthermore, the constant term $\mu_{\delta}^0$ is $\Theta^2$-monotone.
\end{thm}
\begin{proof}
Let $\XX$ be either $\Theta_{R}$ or $\overline{ST}_{R}$ for some discrete valuation ring $R$ over $\cY$. Choose a $\cY$-morphism $\varphi: \XX \setminus \closedpt \to \cM_{n}^{G}(X)$. Denote by $\psi$ the composition $\XX \setminus \closedpt \xrightarrow{\varphi} \cM_{n}^{G}(X) \xrightarrow{p} \MpC( A_{X}/G)$. By the proof of Proposition \ref{prop: monotonicity of mapping stack general group}, there exists a reduced and irreducible algebraic stack $\cW$ with morphisms $f: \cW \to Y_{\XX}/\bG_m$ and $\widetilde{\psi}: \cW \to \MpC( A_{X}/G)$ satisfying the following conditions
\begin{enumerate}
\item The map $f$ is proper with finite relative diagonal and its restriction induces an isomorphism $f : \, \cW_{Y_{\XX} \setminus \closedpt} \xrightarrow{\sim} \XX \setminus \closedpt$.
\item The morphism $\widetilde{\psi}$ has quasi-finite relative inertia, and under the identification induced by $f$ we have $\widetilde{\psi}|_{\XX \setminus \closedpt} \simeq \psi$.
\item The pullback $\widetilde{\psi}^*(\cD(V)^{\vee})$ is positive relative to $\cW \to \XX$.
\item Suppose that $\kappa$ is a finite extension of the residue field of $R$. Let $g_1, g_2: (B\bG_m)_{k} \to \cW_{\closedpt/\bG_m}$ be two morphisms such that the compositions $g_1, g_2: (B\bG_m)_{k} \to \cW_{\closedpt/\bG_m} \to \closedpt/\bG_m$ coincide. Then, we have $b(\widetilde{\psi} \circ g_1) = b(\widetilde{\psi} \circ g_2)$.
\end{enumerate}
The morphism $\varphi$ induces $h: \XX \setminus \closedpt \to \cM_{n}^{G}(X) \times_{\MpC( A_{X}/G)} \cW$. Since $\cM_{n}^{G}(X) \xrightarrow{p} \MpC( A_{X}/G)$ is proper with finite relative diagonal, it follows that $\cM_{n}^{G}(X) \times_{\MpC( A_{X}/G)} \cW \to \cW$ also has finite relative diagonal. Note that $h$ factors as a composition
\begin{gather*} \XX \setminus \closedpt \; \to \;  \cM_{n}^{G}(X) \times_{\MpC( A_{X}/G)} \XX \setminus \closedpt \; \hookrightarrow \; \cM_{n}^{G}(X) \times_{\MpC( A_{X}/G)} \cW
\end{gather*}
The second morphism is an open immersion. The first morphism is representable and finite, since it is a section of a morphism of stacks with finite relative diagonal. In particular, $h$ is representable, quasi-finite and separated. Lemma \ref{lemma: zariski's main theorem stacks} yields a diagram
\[
\begin{tikzcd}
  & \cM'  \ar[dr] & \\ \XX \setminus \closedpt  \ar[ru, "g"] \ar[rr, "h"] & & \cM_{n}^{G}(X) \times_{\MpC( A_{X}/G)} \cW
\end{tikzcd}
\]
where $\cM' \to \cM_{n}^{G}(X) \times_{\MpC( A_{X}/G)} \cW$ is schematic and finite and $g$ is an open immersion. Define $\cW'$ to be the scheme theoretic image of $g$. By construction, $\cW'$ admits a morphism $\widetilde{\varphi}$ to $\cM_{n}^{G}(X)$ and the natural projection $\cW' \to \cW$ is proper with finite relative diagonal. Since $g$ is an open immersion, the natural projection $\cW' \to \cW$ restricts to an isomorphism over $\XX \setminus \closedpt$. We have a commutative diagram
\[
\begin{tikzcd}
  \cW'  \ar[d] \ar[r, "\widetilde{\varphi}"] & \cM_{n}^{G}(X) \ar[d] \\ \cW  \ar[r] & \MpC( A_{X}/G)
\end{tikzcd}
\]
where the vertical morphisms are proper with finite diagonal. By Proposition \ref{prop: formal line bundle is ample on p fibers}, the formal line bundle $\cL_{Cor}(\epsilon)$ is positive relative to the vertical morphisms. Since $\cD(V)^{\vee}|_{\cW}$ positive relative to $\cW \to \XX$, it follows that the formal line bundle $(\cD(V)^{\vee}+ \epsilon \cL_{Cor}(\epsilon))|_{\cW'}$ is positive relative to $\cW' \to \XX$. 

Let $\kappa$ be an extension of the residue field of $R$, and choose a finite morphism $\bP^1_{\kappa}[a]/\bG_m \to \cW'_{\closedpt/\bG_m}$ fitting into a commutative diagram
\[
\begin{tikzcd}
  \mathbb{P}^1_{\kappa}[a]/\bG_m \ar[dr] \ar[r] & \cW'_{\closedpt/\bG_m} \ar[d]\\  &  \closedpt/\bG_m
\end{tikzcd}
\]
We want to show that $\mu_{\delta}(\widetilde{\varphi}|_{\infty/\bG_m}) > \mu_{\delta}(\widetilde{\varphi}|_{0/\bG_m})$. Recall that
\begin{gather*} \mu_{\delta}(\widetilde{\varphi}|_{\infty/\bG_m}) = \frac{ -\wt((\cD(V) + \epsilon \cL_{Cor}(\epsilon))|_{\infty}) - \delta_{\mrk} \wt(\cL_{\mrk}|_{\infty}) + \delta_{\gen} \ell_{\gen}(\widetilde{\varphi}|_{\infty/\bG_m})}{\sqrt{b(\widetilde{\varphi}|_{\infty/\bG_m})}}\\
\mu_{\delta}(\widetilde{\varphi}|_{0/\bG_m}) = \frac{ -\wt((\cD(V) + \epsilon \cL_{Cor}(\epsilon))|_{0}) - \delta_{\mrk} \wt(\cL_{\mrk}|_0) + \delta_{\gen} \ell_{\gen}(\widetilde{\varphi}|_{0/\bG_m})}{\sqrt{b(\widetilde{\varphi}|_{0/\bG_m})}}
\end{gather*}
Since the formal line bundle $(\cD(V) + \epsilon \cL_{Cor}(\epsilon))|_{\cW'}$ is positive relative to $\cW' \to \XX$, we can prove that $-\wt((\cD(V) + \epsilon \cL_{Cor}(\epsilon))|_{\infty}) > -\wt((\cD(V) + \epsilon \cL_{Cor}(\epsilon))|_{0})$ exactly as in the proof of Proposition \ref{prop: monotonicity of mapping stack general group} in the case $G = \GL_N \times H$. On the other hand, property (4) above shows that $b(\widetilde{\varphi}|_{\infty/\bG_m})) = b(\widetilde{\varphi}|_{0/\bG_m})$, because the norm $b$ comes from the stack $\MpC( A_X/G) \to BG$. Therefore, we can conclude our proof by showing the following two claims.

\noindent (a) We have $\wt(\cL_{\mrk}|_{\infty}) \leq \wt(\cL_{\mrk}|_{0})$.\\
(b) We have $\ell_{\gen}(\widetilde{\varphi}|_{\infty/\bG_m})\geq \ell_{\gen}(\widetilde{\varphi}|_{0/\bG_m})$.

For (a) it suffices to show that $\deg(\ev^*(L^{\mrk})|_{\bP^1_{\kappa}[a]}) \geq 0$. This follows because we are pulling back the $S$-semi-ample line bundle $L^{\mrk}$, which restricts to a numerically effective line bundle on $S$-fibers. For (b), denote by $\kappa(C_{\kappa})$ the fraction field of $C_{\kappa}$. The restriction $\varphi|_{\bP^1_{\kappa}[a]}$ yields a $\bG_m$-equivariant $\bP^1_{\kappa}$-point $(E, u, \marked)$ of $\cM_{n}^{G}(X)$. Consider the base change
\[
\begin{tikzcd}
 B\ar[d] \ar[r] & \widetilde{C} \ar[d]\\ \bP^1_{\kappa(C_{\kappa})} \ar[r] & C_{\kappa} \times_{\kappa} \bP^1_{\kappa}
\end{tikzcd}
\]
The map $B \to \bP^1_{\kappa(C_{\kappa})}$ is an isomorphism,\endnote{We claim that $B \to \bP^1_{\kappa(C_{\kappa})}$ is an isomorphism. Indeed, consider the $\bG_m$-equivariant flat locus $U \subset C_{\kappa} \times \bP^1_{\kappa}$ of the morphism $\widetilde{C} \to C_{\kappa} \times \bP^1_{\kappa}$ of flat $\bP^1_{\kappa}$-families of curves. This coincides with the locus where $\widetilde{C} \to C_{\kappa} \times \bP^1_{\kappa}$ is an isomorphism. The intersection $U \cap \bP^1_{\kappa(C_{\kappa})}$ yields a $\bG_m$-equivariant open subset of $\bP^1_{\kappa(C_{\kappa})}$. In order to show that this is the whole of $\bP^1_{\kappa(C_{\kappa})}$, it suffices to prove that $U \cap \bP^1_{\kappa(C_{\kappa})}$ contains both of the points $0_{\kappa(C_{\kappa})}$ and $\infty_{\kappa(C_{\kappa})}$. In other words, we need to show that $\widetilde{C} \to C_{\kappa} \times \bP^1_{\kappa}$ is flat at the generic points of the $0$ (resp. $\infty$) fiber. By the fiberwise criterion for flatness \cite[\href{https://stacks.math.columbia.edu/tag/05VK}{Tag 05VK}]{stacks-project}, this is can be checked after restricting to $\widetilde{C}_0 \to C_0$ (resp. $\widetilde{C}_{\infty} \to C_{\infty}$), and then the flatness at the generic point is obvious.} hence we can use the restrictions $E|_{\bP^1_{\kappa(C_{\kappa})}}$ and $u|_{B}: B \to E(X)$ to define a morphism $\alpha \vcentcolon = (E|_{\bP^1_{\kappa(C_{\kappa})}}, u|_{B}): \bP^1_{\kappa(C_{\kappa})}[a]/\bG_m \to X/G$. Since the line bundle $L_{\gen}$ on $X$ is numerically effective on each $S$-fiber, it follows that the pullback $\alpha^*(L_{\gen})$ on $\bP^1_{\kappa(C_{\kappa})}[a]$ has nonnegative degree. This implies the desired inequality
\[ \ell_{\gen}(\widetilde{\varphi}|_{\infty/\bG_m})= \wt(-\alpha^*(L_{\gen})|_{\infty_{\kappa(C_{\kappa})}}) \geq \wt(-\alpha^*(L_{\gen})|_{0_{\kappa(C_{\kappa})}}) = \ell_{\gen}(\widetilde{\varphi}|_{0/\bG_m}).\]

The proof of $\Theta^2$-monotonicity of $\mu^0_{\delta}$ follows by the same argument. Note that in this case the monotonicity is not strict, since we are not incorporating the component $\cL_{Cor}(\epsilon)$ that is ample along the fibers of $\cM_{n}^G(X) \to \cM_{C}(A_X/G)$.
\end{proof}

\section{HN-boundedness} \label{section: hn boundedness}

The main result of this section is the following theorem, whose proof appears at the end of the section. It applies in the general context of \Cref{subsection: notation}.
\begin{thm} \label{thm: HN-boundedness nonconnected}
HN-boundedness (\Cref{defn: HN boundedness}) holds for the numerical invariant $\mu_{\delta}$ on $\cM_{n}^{G}(X)$ introduced in \Cref{defn: numerical invariant}.
\end{thm}



Before we begin, we will reduce \Cref{thm: HN-boundedness nonconnected} to a simpler context.

\begin{context} \label{context:boundedness_simplification}
$G$ admits a split maximal torus $T \subset G$ defined over $S$, $S$ is affine and connected, $\cY = S = A_X /\!/ G$, and $G = G_0$ has connected fibers. It follows from \Cref{lemma: description of norm in the split case} that the quadratic norm $b$ on $BG$ comes from a bilinear form $(-,-)_{b}$ on real cocharacters $X_{*}(T)_{\bR}$ of $T$. Finally, we shall fix $d \in H_2(X/G)$ and restrict ourselves to the open and closed stack $\cM_{n}^{G}(X)_d$.
\end{context}


\begin{prop} \label{prop: reduction to connected groups}
It suffices to prove \Cref{thm: HN-boundedness nonconnected} in \Cref{context:boundedness_simplification}.
\end{prop}
\begin{proof}
It suffices to check HN-boundedness after base change to a surjective finite type cover of $S$. Therefore, it suffices to consider only the case when $G$ is split and $S$ is affine and connected by \cite[Lemma 5.1.3]{conrad_reductive}. Similarly, we can replace $S$ with $\Spec_{S}(\cO_X^{G})$, and hence we may assume without loss of generality that $A_X /\!/ G = S$.\endnote{Indeed, we are interested in the behavior of bounded families $T \to \cM_{n}^{G}(X)$ for some Noetherian $S$-scheme $T$. Such family is represented by a tuple $(E, u: \widetilde{C} \to E(X), \marked)$. The morphism $u: \widetilde{C} \to E(X) \to E(A_{X})$ factors uniquely through $C_{T}$ (as in \Cref{subsection: gauged maps introduction}). So we have a diagram
\[\xymatrix{\widetilde{C} \ar[d] \ar[r] & E(X) \ar[d]\\  C_{T}  \ar[r] \ar[rd]  &  E(A_{X}) \ar[d]\\ & E(\Spec_{T}(\cO_{X_{T}}^G)) }.\]
Note that $E(\Spec_{T}(\cO_{X_{T}}^{G})) = C_{T} \times_{T} \Spec_{T}(\cO_{X_{T}}^G)$. Therefore the section $C_{T} \to E(\Spec_{T}(\cO_{X_{T}}^{G}))$ is uniquely determined by a morphism $C_{T} \to \Spec_{T}(\cO_{X_{T}}^G)$. This is in turn equivalent to a morphism $\cO_{X_{T}}^G \to (\pi_{C_{T}})_* \cO_{C_{T}}$ of sheaves of $\cO_{T}$-algebras on $T$, where $\pi_{C_{T}}: C_{T} \to T$ is the structure morphism. Since $\pi_{C_{T}}$ is proper, smooth and the fibers are geometrically connected, it follows from cohomology and base-change that $(\pi_{C_{T}})_*\cO_{C_{T}} = \cO_{T}$. Therefore $C_{T} \to \Spec_{T}(\cO_{X_{T}}^G)$ factors through a unique section $T \to \Spec_{T}(\cO_{X_{T}}^G)$. Replacing $S$ with $\Spec_{T}(\cO_{X_{T}}^G)$, replacing $G$ with $G \times_{S} \Spec_{T}(\cO_{X_{T}}^G)$, and changing $X$ to $X_{T}$, we may assume without loss of generality that $A_X /\!/ G = S$.} Also, any filtration of a point in a particular open and closed substack $\cM^G_n(X)_d$ must have associated graded point in the same substack, so it suffices to prove HN boundedness separately for each $\cM^G_n(X)_d$.

Finally, we must show that it suffices to assume $G$ is connected. Because we can check HN boundedness \'etale locally on the base, we can assume that $\pi_0(G) = \underline{\Gamma}$ is a constant group scheme \cite[\href{https://stacks.math.columbia.edu/tag/04HN}{Tag 04HN}]{stacks-project}. Let $\Bun_\Gamma(C)^{o} \subset \Bun_{\Gamma}(C)$ denote the open and closed substack parameterizing torsors whose total space has connected geometric fibers over the base, and let $\cM_{n}^G(X)^{o}$ denote the preimage of $\Bun_{\Gamma}(C)^{o}$ under the canonical map $\cM_{n}^G(X) \to \Bun_{\Gamma}(C)$. There is a finite surjective morphism $\bigsqcup_{H \subset \Gamma} \Bun_H(C)^{o} \to \Bun_{\Gamma}(C)$ induced by the inclusions $H \subset \Gamma$, where the disjoint union is over conjugacy classes of finite subgroups $H \subset \Gamma$.\endnote{For finiteness $\bigsqcup_{H \subset \Gamma} \Bun_H(C)^{o} \to \Bun_{\Gamma}(C)$, it suffices to prove finiteness of $\Bun_H(C) \to \Bun_{\Gamma}(C)$, since each $\Bun_{H}(C)^{o} \hookrightarrow \Bun_{H}(C)$ is a closed immersion. Let $T$ be a $\mathcal{Y}$-scheme and choose a point $T \to \Bun_{\Gamma}(C)$ corresponding to a $\Gamma$-torsor $P$ on $C_{T}$. The fiber product $T \times_{\Bun_{\Gamma}(C)} \Bun_{H}(C)$ classifies sections of the finite morphism $P \times^{\Gamma} (\Gamma/H) \to C_T$. By \Cref{lemma: sections of finite morphism finite}, it follows that of $T \times_{\Bun_{\Gamma}(C)} \Bun_{H}(C)$ is finite, as desired.

For surjectivity, we need to prove that $\bigsqcup_{H \subset \Gamma} \Bun_H(C)^{o} \to \Bun_{\Gamma}(C)$ is surjective on geometric points, where $H \subset \Gamma$ runs over a set of representatives of conjugacy classes of subgroups of $\Gamma$. Let $k$ be an algebraically closed field over $\mathcal{Y}$ and choose a $k$-point of $\Bun_{\Gamma}(C)$ corresponding to a $\Gamma$-torsor $P$ over the smooth projective curve $C_{k}$. Write $P$ as a disjoint union of it finitely many connected components $P = \sqcup_i P_i$. Each $P_i$ is geometrically connected over $k$, since $k$ is algebraically closed. The group $\Gamma$ acts on $P$, and it permutes transitively the connected components $P_i$. Fix a connected component $P_0$ and let $H\subset\Gamma$ denote the subgroup that stabilizes the connected component $P_0$. By possibly changing the choice of connected component $P_0$, we can replace $H$ with a conjugate and assume that $H$ lies in one of the representatives we fixed to form the union $\bigsqcup_{H \subset \Gamma} \Bun_H(C)^{o}$. Now the group $H$ acts on $P_0$, and under this action $P_0$ becomes an $H$-torsor over $C_k$. Since $P_0$ is geometrically connected, it corresponds to a $k$-point in $\Bun_{H}(C)^{o}$. By transitivity of the $G$-action on the set of connected components $P_i$, can write $P = \sqcup_{\sigma \in H \backslash \Gamma} \sigma \cdot P_0$. This shows that $P$ is the $\Gamma$-torsor obtained by extension of structure group from $P_0$, and so it is in the image of the morphism $\Bun_{H}(C)^{o} \to \Bun_{\Gamma}(C)$.} The fiber product of $\Bun_{H}(C)^{o}$ with $\cM_{n}^G(X)$ is identified with $\cM_{n}^{G \times_{\underline{\Gamma}} \underline{H}}(X)^{o}$, and the pullback of the numerical invariant on $\cM_{n}^G(X)$ coincides with the corresponding numerical invariant on $\cM_{n}^{G \times_{\underline{\Gamma}}\underline{H}}(X)$. Since finite morphisms induce bijections on filtrations \cite[Prop. 3.2.12 (1)]{hl_instability}, it suffices to prove HN boundedness for the component $\cM_{n}^{G}(X)^{o}$ only.

Let $\widetilde{C}$ denote the universal torsor over $C \times \Bun_{\Gamma}(C)^{o}$, regarded as a smooth family of curves over $\Bun_{\Gamma}(C)^{o}$. Then $\widetilde{C}$ has a free $\Gamma$-action, and $\widetilde{\cM} := \cM_{\widetilde{C},n\cdot m}^{G_0}(X)$ admits a natural $\Gamma$-action as well, where $m = |\Gamma|$ and $G_0$ is the kernel of $G \to \underline{\Gamma}$, and we use the subscript $\tilde{C}$ to indicate the stack of gauged maps from $\tilde{C}$, rather than $C$. In fact, \'etale descent provides an isomorphism
\begin{equation}\label{E:fixed_points}
\cM_{n}^G(X)^{o} \cong \widetilde{\cM}^{\Gamma},
\end{equation}
where the right-hand side denotes the stack of $\Gamma$-invariant points (see \cite[Def.~2.3]{romagny-actions-stacks}) of $\widetilde{\cM}$.\endnote{Fix a choice of indexing $\Gamma = \{ \gamma_1, \gamma_2, \ldots , \gamma_m\}$, where $m = |\Gamma|$. The stack $\cM_{n \cdot m}^{G_0}(X)$ of gauged maps from $\widetilde{C}$ into $X$ admits a natural $\Gamma$-action induced by its action on $\widetilde{C}$ (here the action on the markings is done by grouping them into sets of $m$ consecutive elements, and acting on the indexes of each set). Let $\widetilde{\cM}^{\Gamma}$ denote the stack of $\Gamma$-equivariant points in $\cM_{h \cdot n}^{G_0}(X)$. By definition $\widetilde{\cM}$ sends a $\Bun_{\Gamma}(C)^{o}$-scheme $T$ to the group of tuples $(\widetilde{E}, u, \cH_1, \cH_2, \ldots \cH_n, \widetilde{p}_1, \widetilde{p}_2, \ldots, \widetilde{p_n})$, where
\begin{enumerate}
    \item $\widetilde{E}$ is a $\Gamma$-equivariant $G_0$-bundle on $\widetilde{C}_{T}$.
    \item $\mathcal{C} \to T$ is $T$-family of nodal curves equipped with an action of $\Gamma$.
    \item $u: \mathcal{C} \to \widetilde{E}(X)$ is a $\Gamma$-equivariant $T$-morphism. 
    \item Each $\cH_i$ is a $\Gamma$-bundle over $T$, and $\widetilde{p_i}: \cW_i \to \mathcal{C}$ are $\Gamma$-equivariant $T$-morphisms.
    \item For each geometric point $\overline{t}$ and a choice of trivializations $\cW_i|_{\overline{t}} \cong \Gamma \times \overline{t}$, we can use our fixed choice of indexing $\Gamma = \{\gamma_l\}_{l=1}^m$ to interpret the restriction $(\mathcal{C}_{\overline{t}}, (\widetilde{p}_1)_{\overline{t}}, \ldots, (\widetilde{p}_n)_{\overline{t}})$ as a nodal curve with $n \cdot m$ markings. We then require that $(\mathcal{C}_{\overline{t}}, u_{\overline{t}}, (\widetilde{p}_1)_{\overline{t}}, \ldots, (\widetilde{p}_n)_{\overline{t}})$ is a Kontsevich stable morphism (Note that this does not depend on the choice of trivialization or indexing of $\Gamma$).
\end{enumerate}
Let $\cP$ denote the universal $\Gamma$-bundle on $C \times \Bun_{\Gamma}(C)^{o}$. The pullback $\Gamma$-bundle $\cP_{\widetilde{C}}$ on its total space $\widetilde{C}$ admits a canonical trivialization (given by the diagonal). Consider the morphism $\cM_{n}^{G}(X) \to C\times \Bun_{\Gamma}(\Gamma)^{o}$. For any $T$-point $(E, s: \cC \to E(X), \marked) \in \cM_{n}^{G}(X)$, the pullback of $E|_{\widetilde{C}_{T}}$ acquires a canonical reduction of structure group to a $G_0$-bundle $E_{\widetilde{C}_{T}}^0$, coming from the canonical trivialization of the associated $\Gamma$-bundle on $\widetilde{C}_{T}$. By construction this reduced pullback $E_{\widetilde{C}_T}^0$ acquires a canonical $\Gamma$-equivariant structure with respect to the Galois cover action of $\widetilde{C}_T \to C_T$. We have an induced fiber product section $s_{\widetilde{C_{T}}}: \cC\times_{C_T} \widetilde{C_T} \to (E(X)) \times_{C_T} \widetilde{C_{T}} \xrightarrow{\sim} E_{\widetilde{C}_T}^0(X)$, which is naturally $\Gamma$-equivariant. Moreover, each fiber product
\[\xymatrix{\cH_i \ar[r] \ar[d] & \cC \times_{C_T} \widetilde{C}_T \ar[d]\\ T \ar[r]^{p_i} &  \cC}\]
is naturally a $\Gamma$-bundle on $T$, equipped with a $\Gamma$-equivariant morphism $\widetilde{p}_i: \cH_i \to \cC \times_{C_T} \widetilde{C}_T$. It follows from construction that the geometric fibers of the tuple $(E_{\widetilde{C}_T}, s: \cC\times_{C_T} \widetilde{C}_T \to E_{\widetilde{C}_T}(X), \cH_1, \ldots, \cH_n, \widetilde{p}_1, \ldots, \widetilde{p}_n)$ satisfy the required Kontsevich stability condition in (5) above. Since this ``pullback'' assignment $(E, s: \cC \to E(X), \marked) \mapsto ((E_{\widetilde{C}_T}, s: \cC\times_{C_T} \widetilde{C}_T \to E_{\widetilde{C}_T}(X), \cH_1, \ldots, \cH_n, \widetilde{p}_1, \ldots, \widetilde{p}_n))$ is compatible with base-change on $T$, it induces a morphism of stacks $\Bun_{\Gamma}(C)^{o}$-stacks $\cM_{n}^{G}(X) \to \widetilde{\cM}^{\Gamma}$. Conversely, for any tuple $(\widetilde{E}, u, \cH_1, \cH_2, \ldots \cH_n, \widetilde{p}_1, \widetilde{p}_2, \ldots, \widetilde{p_n})$ in $\widetilde{\cM}$, taking the quotient of everything by the action of $\Gamma$ yields a point of $\cM_{n}^{G}(C)$. By \'etale descent for the Galois covers $\widetilde{C}_T \to C_T$, these operations are inverses of each other, and hence we have established a natural isomorphism $\cM_{n}^{G}(C) \xrightarrow{\sim} \widetilde{\cM}^{\Gamma}$ induced by pulling back to $\widetilde{C}$.} The natural forgetful morphism $\varphi : \widetilde{\cM}^{\Gamma} \to \widetilde{\cM}$ induces a bijection between the set of filtrations of any point $x$ of $\widetilde{\cM}^{\Gamma}$ and the set of $\Gamma$-invariant filtrations of $\varphi(x)$.\endnote{Since the morphism $\varphi: \widetilde{\cM}^{\Gamma} \to \widetilde{\cM}$ is representable and separated \cite[Thm. 3.3]{romagny-actions-stacks}, it follows that for every point $x$ in $\widetilde{\cM}^{\Gamma}$ the induced morphism $\Flag_{\mathcal{Y}}(x) \to \Flag_{\mathcal{Y}}(\varphi(x))$ on sets of filtrations relative to $\mathcal{Y}$ is an injection \cite[Lemma 3.7(iii)]{AHLH}. More concretely, suppose that $x$ is a geometric point defined over an algebraically closed field $k$ over $\mathcal{Y}$. Let $E$ denote the underlying $G_k$-bundle on $C_k$ corresponding to $x$. We denote by $\widetilde{E}$ the pullback $G_k$-bundle to the cover $\widetilde{C}_k$; this is the underlying $G_k$-bundle for $\varphi(x)$. We can think of filtrations of $x$ as tuples $(\lambda, E_{P_{\lambda}})$, where $\lambda: (\mathbb{G}_m)_k \to G_k$ and $E_{P_{\lambda}}$ is a reduction of structure group of $E$ to the parabolic subgroup $P_{\lambda}$ corresponding to the one-parameter subgroup $\lambda$ (note that not all such tuples arise from filtrations if $X \to S$ is not projective, and we should identify tuples up to conjugation by elements of $G_k$). Filtrations of $\varphi(x)$ admit a similar description as tuples $(\lambda, \widetilde{E}_{P_{\lambda}})$ for the $G_k$-bundle $\widetilde{E}_{P_{\lambda}}$. Since $\widetilde{E}$ has a canonical $\Gamma$-equivariant structure, the group $\Gamma$ acts on the set of parabolic reductions $\widetilde{E}_{P_{\lambda}}$, and hence on the set of tuples $(\lambda, \widetilde{E}_{P_{\lambda}})$. This action is compatible with the action of $\Gamma$ on the set of filtrations $\Flag_{\mathcal{Y}}(\varphi(x))$. The morphism $\Flag_{\mathcal{Y}}(x) \to \Flag_{\mathcal{Y}}(\varphi(x))$ is induced by the morphism on tuples $(\lambda, E_{P_{\lambda}})  \mapsto (\lambda, \widetilde{E}_{P_{\lambda}})$ defined by pulling back the parabolic reduction. So a point $\Flag_{\mathcal{Y}}(\varphi(x))$ corresponding to $(\lambda, \widetilde{E}_{P_{\lambda}})$ comes from a filtration in $\Flag_{\mathcal{Y}}(x)$ if and only if the parabolic reduction $\widetilde{E}_{P_{\lambda}}$ descends to a parabolic reduction of $E$. By \'etale descent, this happens exactly when $\widetilde{E}_{P_{\lambda}}$ is fixed by the action of $\Gamma$, as desired.}

Under the isomorphism \eqref{E:fixed_points}, the numerical invariant on $\cM_{n}^G(X)^{o}$ corresponds to the restriction along $\varphi$ of a numerical invariant on $\widetilde{\cM}$, which by hypothesis satisfies HN boundedness and thus defines a weak $\Theta$-stratification by \Cref{thm: monotonicity of gauged maps} and \Cref{thm: theta stability paper theorem}. Now consider a morphism $\xi : T \to \widetilde{\cM}^\Gamma$ from a finite type affine scheme $T$. Every geometric point of $T$ whose image under $\varphi \circ \xi$ is unstable admits a maximizing filtration in $\widetilde{\cM}$. The uniqueness (up to scale) of this filtration implies that it is $\Gamma$-invariant and thus defines a filtration of $\xi(t)$ in $\widetilde{\cM}^\Gamma$, which must maximize the restricted numerical invariant. Furthermore, as $t \in T$ varies, the associated graded points of these optimal filtrations in $\widetilde{\cM}$ must lie in a quasi-compact open substack $U_\xi \subset \widetilde{\cM}$, so the optimal filtrations in $\widetilde{\cM}^\Gamma$ must lie in $\varphi^{-1}(U_\xi)$, which is quasi-compact by \cite[Prop.~3.7]{romagny-actions-stacks}.

\end{proof}

For the rest of this section we will assume without loss of generality that we are in \Cref{context:boundedness_simplification}. We will use the following:

\begin{notn} \label{notn:root data torus}
Let $Z \subset T$ denote the maximal central torus of $G$. We denote by $T'$ the intersection of $T \cap \D(G)$. 
The natural inclusions induce a direct sum decomposition $X_{*}(T)_{\mathbb{R}}  = X_{*}(Z)_{\mathbb{R}} \oplus X_*(T')_{\mathbb{R}}$. We will denote by $(-)^{\dagger}: X^*(T)_{\mathbb{R}} \xrightarrow{\sim} X_*(T)_{\mathbb{R}}$ the identification induced by the nondegenerate quadratic norm $b$. Since $b$ is invariant under the Weyl group, the isomorphism $(-)^{\dagger}$ preserves the direct sum decomposition above. We shall denote by $(-,-)_b$ the bilinear form on $X_*(T)_{\mathbb{R}}$ induced by $b$.

We fix once and for all the choice of a Borel subgroup $B \supset T$. This yields a set of simple roots $\Delta = \{ \alpha_1, \alpha_2, \ldots, \alpha_l \} \subset X_{*}(T')_{\mathbb{R}}^{\vee}$, which form a basis for $X^*(T')_{\mathbb{R}}$. Similary we get a set of simple coroots $\Delta^{\vee} = \{ \alpha_1^{\vee}, \alpha_2^{\vee}, \ldots, \alpha_l^{\vee} \} \subset X_{*}(T')_{\mathbb{R}}$. We denote by $\omega_{i}^{\vee} \in X_{*}(T')_{\mathbb{R}}$ the fundamental coweight dual to the simple root $\alpha_i$, and similarly we denote by $\omega_i$ the fundamental weight dual to $\alpha_i^{\vee}$. By classical theory of root spaces applied to the dual root datum, we have that $(\omega_i^{\vee}, \alpha_j^{\vee})_b = \frac{1}{2} \|\alpha_j^{\vee}\|_b^2 \cdot \delta_{i,j} \geq 0$ for all $i,j$ (cf. the classical definition of fundamental weights in \cite[Sect. 13.1]{humphreys-intro-lie}).
\end{notn}

\subsubsection{Degeneration fans and spaces} \label{subsection: degeneration fans}

Let $\cX$ be a stack locally of finite type over $\cY$ with quasi-affine diagonal. Let $k$ be a field over $\cY$. For a $k$-point $x: \Spec(k) \to \cX$, the \emph{degeneration
fan} $\mathbf{DF}(\cX,x)$ (relative to $\cY$) is a formal fan whose $n$-cones are non-degenerate $k$-points of the $n$-flag space of $x$ relative to $\cY$, i.e., maps $f : \Theta_k^n \to \cX \times_{\cY} \Spec(k)$ equipped with an isomorphism $f(1) \cong x$. The geometric realization $|\mathbf{F}_\bullet|$ of a
formal fan $\mathbf{F}_\bullet$ is defined as a colimit of a copy of $\bR^n_{\geq 0}$ for each $n$-cone $\sigma \in F_n$, similar to 
the construction of the geometric realization
of a semi-simplicial set.
There is a natural scaling action of $\bR^{\times}$ on the geometric realization $\lvert \DF(\cX, x) \rvert$, with a canonical fixed cone point $0$. The \emph{degeneration space} $\iDeg(\cX, x)$ is defined as the quotient
\[
\iDeg(\cX, x) := ( \lvert\DF(\cX,x)\rvert \setminus 0 ) / \bR_{>0}^\times.
\]
We direct the reader to \cite[Sect. 3]{hl_instability} for a more detailed discussion of this general construction.

Turning to the specific example of $\Bun_G(C)$, a $k$-point for a field $k$ over $S$ is a $G$-bundle $E$ on $C_{k}$. By \cite{heinloth2018hilbertmumford}, a nondegenerate filtration consists of a tuple $(\lambda, E_{\lambda})$ where $\lambda$ is a dominant coweight and $E_{\lambda}$ is a reduction of structure group of $E$ to the parabolic subgroup $P_{\lambda}$ associated to $\lambda$.


\begin{notn} \label{notation: beginning proof hn-boundedness bun_g}
Fix a parabolic $P \supseteq B_k$ and a reduction $E_{P}$ to $P$. $P$ is determined by the subset $I_{P} \subset \{1, \ldots, l\}$ of indices of simple roots $\alpha_i \in \Delta$ such that the root group $U_{-\alpha_i} \not\subset P$ (e.g., $I_{B} = \{1, \ldots, l\}$ and $I_{G} = \emptyset$). Then we define the cones
\begin{gather*}
    \sigma_{E_{P}} = X_{*}(Z)_{\bR} + \sum_{ j \in I_{P}} \bR_{> 0}\omega^{\vee}_j, \quad \text{and} \quad
    \overline{\sigma}_{E_{P}} = X_{*}(Z)_{\bR} + \sum_{ j \in I_{P}} \bR_{\geq 0}\omega^{\vee}_j.
\end{gather*}
The vector subspace $\Span(\sigma_{E_P}) \subset X_{*}(T)_{\bR}$ spanned by $\sigma_{E_{P}}$ is equal to $X_\ast(T)_\bR^{W_P}$, where $W_P \subset W$ is the subgroup generated by the simple reflections indexed by $I_P$. Similarly, we define the cone $\overline{\tau}_{E_{P}} = X_{*}(Z)_{\bR} + \sum_{ j \in I_{P}} \bR_{\geq 0}\alpha^{\vee}_j$. Note that $\overline{\tau}_{E_{P}} \supset \overline{\sigma}_{E_P}$. We set $\rho_P^{\vee} \vcentcolon = \sum_{j \in I_{P}} \omega_j^{\vee}$.
\end{notn} 

The cone $\overline{\sigma}_{E_{P}}$ can be identified with the closed cone of $P$-dominant real cocharacters in $X_{*}(T)_{\bR}$, and it parameterizes all filtrations $(\lambda,E_\lambda)$ of $E$ with $P = P_\lambda$ and $E_P = E_\lambda$. The elements $w \in \overline{\sigma}_{E_{P}}$ in the boundary correspond to extensions $E_{P'}$ of the reduction $E_{P}$, where $P' \supset P$ is the parabolic subgroup associated to the facet of the Weyl chamber of $X_{*}(T)_{\bR}$ containing $w$ in its relative interior. In this way, we realize $\lvert \DF(\Bun_{G}(C), E)\rvert$ as the union of the cones $\overline{\sigma}_{E_{P}}$ where $E_{P}$ ranges over all reductions of structure group to a parabolic $P \supseteq B_k$, equipped with the colimit topology.

\subsection{HN-boundedness for \texorpdfstring{$\Bun_{G}(C)$}{BunG(C)}} \label{subsection: hn boundedness for bung}

\begin{prop} \label{prop: hn boundedness bung}
In \Cref{context:boundedness_simplification}, the numerical invariant $\mu \vcentcolon = -\wt(\cD(V)) / \sqrt{b}$ on $\Bun_{G}(C)$ satisfies HN-boundedness.
\end{prop}

We will prove this at the end of this subsection, after analyzing $\mu$ on $\overline{\sigma}_{E_{P}}$ in more detail.

\begin{defn} \label{defn: ch_2(V) and adjoint}

We define a Weyl-invariant symmetric bilinear form $(-,-)_V$ on $X_{*}(T)_{\bR}$ given by
\[
(w_1, w_2)_{V} = \sum_{\chi \in X^\ast(T)}\rk(V_{\chi}) \chi(w_1) \chi(w_2),
\]
where $V_\chi$ denotes the $T$-eigenspace of $V$ of weight $\chi$. Let $\phi_{V}: X_{*}(T)_{\mathbb{R}} \to X_{*}(T)_{\mathbb{R}}$ be the linear endomorphism determined by $(\phi_{V}(w), v)_b = (w, v)_{V}$ for all $w, v \in X_{*}(T)_{\mathbb{R}}$. Note that $\phi_{V}$ preserves each subspace $\Span(\sigma_{E_P}) \subset X_{*}(T)_{\mathbb{R}}$. Since $b$ and $(-,-)_{V}$ are Weyl-invariant, the summands in the direct sum decomposition $\Span(\sigma_{E_P})= X_{*}(Z)_{\bR} \oplus (X_{*}(T')_{\bR} \cap \Span(\sigma_{E_P}))$ are both $b$- and $(-,-)_{V}$-orthogonal. Therefore, $\phi_{V}$ respects the direct sum decomposition.
\end{defn}

Any parabolic reduction $E_{P}$ of a $G$-bundle $E$ on $C_{k}$ has an associated \emph{degree} $d_{E_P} \in X^\ast(P)^\dual_\bR$, defined by $d_{E_P}(\chi) = \deg_C(E_P(\chi))$. Under the canonical identification $X^\ast(P)^\dual_\bR \cong X_\ast(T)_\bR^{W_P}$, we regard $d_{E_P}$ as an element of $\Span(\sigma_{E_P}) = X_\ast(T)_\bR^{W_P}$.\endnote{One way to think of the identification $X^\ast(P)^\dual_\bR \cong X_\ast(T)_\bR^{W_P}$ is as follows. Let $M_{P}$ denote the unique Levi subroup of $P$ containing $T_{k}$. For any parabolic reduction $E_{P}$ of a $G$-bundle $E$ on $C_{k}$, we can form the associated $M_{P}$-bundle $E_{M_{P}}$ using the natural quotient homomorphism $P \to M_{P}$. The subspace $\Span(\sigma_{E_P})$ is canonically identified with the space of real cocharacters $X_*(Z_{M_{P}})_{\bR}$ of the center $Z_{M_{P}}$ of $M_{P}$.  The degree $d_{E_{P}}$ of the $M_{P}$-bundle $E_{M_{P}}$ can be viewed as a linear functional in $X^*(M_{P})_{\mathbb{R}}^{\vee}$ (see Example \ref{example: degree of a G-bundle}). Since the composition $Z_{M_{P}} \to M_{P} \to M_{P}/\D(M_{P})$ is an isogeny of tori, it induces a canonical identification $X^{*}(M_{P})_{\bR} \xrightarrow{\sim} X^{*}(Z_{M_{P}})_{\bR} =\Span(\sigma_{E_P})$.}

\begin{lem} \label{lemma: riemann-roch expression line bundle}
Given a parabolic reduction $E_{P}$, the linear functional $-\wt(\cD(V))$ on the cone $\overline{\sigma}_{E_P}$ is given by 
\[-\wt(\cD(V))(w) = (d_{E_P}, w)_{V} = (\phi_{V}(d_{E_P}), w)_b.\]
\end{lem}
\begin{proof}
This follows from a Riemann-Roch computation, cf. \cite[1.F]{heinloth2018hilbertmumford}.
\end{proof}

\begin{lem}\label{lemma:lagrange}
Let $W$ be a finite dimensional real vector space, $\ell^\dagger \in W$ a non-zero element, and $b$ a positive definite quadratic form on $W$. Then $\ell^\dagger$ is the unique local maximum of the function $\mu(w) = (\ell^\dagger,w)_b / \|w\|_b$ on $(W \setminus 0) / \bR_{>0}^\times$.

\end{lem}
\begin{proof}
Locally maximizing $\mu$ is equivalent to locally maximizing $(\ell^\dagger,w)_b$ subject to the constraint $\|w\|_b^2=1$. The Lagrange multiplier equations have exactly two solutions $w=\pm \ell^\dagger/\lVert \ell^\dagger \rVert_b$. By compactness of the set $\{w \, | \, \|w\|=1\}$ one must be a local minimum and the other a local maximum. It is clear that $w = \ell^\dagger$ is the local maximum.
\end{proof}

The key geometric input to our proof of HN boundedness here and in later subsections is the following.
\begin{lem} \label{lemma: boundedness of filtrations vs cones}
Let $\mathfrak{S} \subset \lvert \Bun_{G}(C)\rvert$ be a bounded collection of $G$-bundles on the geometric fibers of $C \to S$. Then the following hold.
\begin{enumerate} 
    \item There exists a uniform constant $c_{\mathfrak{S}}$ such that for all $E \in \mathfrak{S}$ and for all parabolic reductions $E_{P}$, we have $\phi_{V}(d_{E_P}) \in c_{\mathfrak{S}} \cdot \rho_P^{\vee} - \overline{\tau}_{E_{P}}$.
    \item For each fixed constant $c \in \bR$, the collection of parabolic reductions
    \[ \left\{ (E, E_{P}) \, \vert \, E \in \mathfrak{S} \text{ and} \; E_{P} \text{ is a parabolic reduction of $E$ such that $\phi_{V}(d_{E_P}) \in c \cdot \rho_P^{\vee} + \overline{\sigma}_{E_{P}}$} \right\} \]
    is bounded.
\end{enumerate}
\end{lem}
\begin{proof}
Choose a field $k$, a $k$-point of $S$, and a parabolic subgroup $P \subset G_{k}$ containing $B_{k}$.  We first claim that, in order to prove the lemma, we can replace $\phi_{V}(d_{E_P})$ with $d_{E_{P}}$ in both (1) and (2). 

For any constant $c$, $\phi_{V}(d_{E_P})$ lies in $c \cdot \rho_P^{\vee} - \overline{\tau}_{E_{P}}$ (resp. $c \cdot \rho_P^{\vee} + \overline{\sigma}_{E_{P}}$) if and only if $\langle \phi_{V}(d_{E_P}), \omega_j \rangle \leq c \cdot \langle \rho_P^{\vee}, \omega_j\rangle$ (resp. $\langle \phi_{V}(d_{E_P}), \alpha_j \rangle \geq c$) for all $j \in I_P$.\endnote{Indeed, by definition we have $\overline{\tau}_{E_{P}} = \left\{ w \in \Span(\sigma_{E_{P}}) \; \mid \; \langle w, \omega_i \rangle \geq 0 \; \text{for all $i \in I_P$} \right\}$. This means that for $c \cdot \rho_P^{\vee} \in \Span(\sigma_{E_{P}})$ we have
 \begin{align*} c \cdot \rho_P^{\vee} - \overline{\tau}_{E_{P}} &= \left\{ w \in \Span(\sigma_{E_{P}}) \; \mid \; \langle w, \omega_i \rangle \leq \langle c \cdot \rho_P^{\vee}, \omega_i \rangle \; \text{for all $i \in I_P$} \right\} \\ &=  \left\{ w \in \Span(\sigma_{E_{P}}) \; \mid \; \langle w, \omega_i \rangle \leq c \langle \rho_P^{\vee}, \omega_i \rangle \; \text{for all $i \in I_P$} \right\} \end{align*}
 Similarly, using $\overline{\sigma}_{E_{P}} = \left\{ w \in \Span(\sigma_{E_{P}}) \; \mid \; \langle w, \alpha_i \rangle \geq 0 \; \text{for all $i \in I_P$} \right\}$ we get that
 \[ c \cdot \rho_P^{\vee} + \overline{\sigma}_{E_{P}} = \left\{ w \in \Span(\sigma_{E_{P}}) \; \mid \; \langle w, \alpha_i \rangle \geq c\langle \cdot \rho_P^{\vee}, \alpha_i \rangle \; \text{for all $i \in I_P$} \right\} =  \left\{ w \in \Span(\sigma_{E_{P}}) \; \mid \; \langle w, \alpha_i \rangle \geq c \; \text{for all $i \in I_P$} \right\},\]
 where we are using that $\langle \rho_P^{\vee}, \alpha_i \rangle=1$ for all $i \in I_P$ to derive the last equality of sets.
 } This is completely determined by the component of $\phi_{V}(d_{E_P})$ in $X_{*}(T')_{\mathbb{Q}}$. Since $(-,-)_V$ is Weyl invariant, it must agree with $(-,-)_b$ up to scaling on each component of $X_*(T')_{\mathbb{Q}}$ corresponding to the torus of an almost simple factor. Therefore, the component of $\phi_{V}(d_{E_{P}})$ on $X_*(T')_{\mathbb{Q}}$ agrees with the component of $d_{E_{P}}$ up to certain fixed scaling on the cocharacter space of each almost simple component. Hence, a uniform bound on the quantities $\langle \phi_{V}(d_{E_P}), \alpha_j \rangle$ is equivalent to a corresponding uniform bound on $\langle d_{E_{P}}, \alpha_j \rangle$. We therefore proceed with the proof for $d_{E_{P}}$.

\medskip
\noindent \textit{Proof of (1):} For any field $k$ over $S$ and every fundamental weight $\omega_j$, there is a representation of $G_k$ with highest weight $\omega_j$ \cite[Chpt. 2]{jantzen-repn}. By spreading out and using Noetherian induction, we get a finite stratification of $S$ by locally closed subschemes such that the restriction of $G$ admits a representation $V_j$ with highest weight $\omega_j$ for all $j$. Therefore, without loss of generality, we may assume that such representations $G \to \GL(V_j)$ exist over $S$ for all $j$. For each $j$, we denote by $\mu_{max}(E(V_j))$ the maximum among the slopes of subbundles of the associated vector bundle $E(V_j)$.

For each $j$, the set $\{ E(V_j) \, \mid \, E \in \mathfrak{S}\}$ of vector bundles on fibers of $C \to S$ is bounded, so by \cite[Lemma 1.7.6+Lemma 1.7.9]{huybrechts.lehn} there exists some $c_{\mathfrak{S}}$ such that $\mu_{max}(E(V_j)) \leq c_{\mathfrak{S}}\cdot \langle \rho_P^{\vee}, \omega_j\rangle$ for all $E \in \mathfrak{S}$. For each parabolic subgroup $P \supset B$ in $G$, choose a coweight $\lambda$ such that $P = P_\lambda$. For all $E \in \mathfrak{S}$ and all parabolic reductions $E_{P_{\lambda}}$, $\lambda$ induces a $\bZ$-weighted filtration of $E(V_j) = E_{P_\lambda}(V_j)$ by subbundles, where the highest weight subbundle of $E(V_j)$ corresponds to the highest weight $\lambda$-eigenspace of $V_j$. The pairing $\langle d_{E_{P}}, \omega_j \rangle$ is the slope of this highest weight subbundle, because $\omega_j$ is a highest weight of $V_j$.\endnote{Note that there can be several $T$-weights $w$ of $V_j$ such that $\langle \lambda, w\rangle = \langle \lambda, \omega_j \rangle$ is maximal. However, any such $w$ satisfies $\langle d_{E_P} , w\rangle = \langle d_{E_P}, \omega_j \rangle$. Indeed, since $\omega_j$ is the highest $T$-weight of $V_j$, we must have $w = \omega_j - \sum_{i \in I} c_i \alpha_i$ for some nonnegative numbers $c_j\geq 0$. We know that $\langle \lambda, \alpha_i \rangle\geq 0$ for any $i \in I$, with strict inequality if $i \in I_P$. Since $\langle \lambda, w\rangle = \langle \lambda, \omega_j\rangle$, we must have $c_i =0$ for all $i \in I_P$. Using that $\langle d_{E_P}, \alpha_i \rangle =0$ for all $i \notin I_P$, we conclude that \[\langle d_{E_P}, w\rangle = \langle d_{E_P}, \omega_j - \sum_{i \notin I_P} c_i \alpha_i \rangle = \langle d_{E_P}, \omega_j\rangle.\]} Therefore $\langle d_{E_{P}}, \omega_j \rangle \leq \mu_{max}(E(V_j)) \leq c_{\mathfrak{S}}\cdot \langle \rho_P^{\vee}, \omega_j\rangle $, as desired.

\medskip
\noindent \textit{Proof of (2):} Let $\Ad: G \to \GL(\mathfrak{g})$ be the adjoint representation. For any parabolic $P \supseteq B$, fix a parabolic subgroup $a(P) \subset \GL(\mathfrak{g})$ such that $\rm{Ad}^{-1}(a(P))=P$. In order to bound the set of parabolic reductions $(E, E_{P})$ it is sufficient to bound the set of associated $a(P)$-reductions of the adjoint vector bundles $\Ad(E)$.\endnote{Indeed, since there are finitely many parabolic types for $P$, we can focus on a single $P$. Suppose that we are given bounded compatible families of parabolic reductions $\Ad(E)_{a(P)}$ and $G$-bundles parametrized by a morphism $D \to \Bun_{a(P)}(C) \times_{\Bun_{\GL(\mathfrak{g})}(C)} \Bun_{G}(C)$, where $D$ is a Noetherian scheme. We want to bound the fiber product
\[\xymatrix{Y \ar[d] \ar[r] & \Bun_{P}(C) \ar[d]\\ D \ar[r] & \Bun_{a(P)}(C)\times_{\Bun_{\GL(\mathfrak{g})}(C)} \Bun_{G}(C)}\]
parametrizing all possible $P$-parabolic reductions that yield elements in the compatible family. We can think of the family $D \to \Bun_{a(P)}(C) \times_{\Bun_{\GL(\mathfrak{g})}(C)} \Bun_{G}(C)$ as representing a $G$-bundle $R$ on $C_{D}$ and a section $s$ of the $\GL(\mathfrak{g})/a(P)$-fiber bundle $\Ad(R)/a(P) \to C_{D}$. There is a morphism of $C_{D}$-schemes $R/P \to \Ad(R)/a(P)$, inducing the fiber product:
\[\xymatrix{R/P \times_{C_{D}} \Ad(R)/a(P)\ar[d] \ar[r] & R/P \ar[d]\\ C_{D} \ar[r]^{s} & \Ad(R)/a(P)}\]
$Y$ parametrizes sections of $R/P \times_{C_{D}} \Ad(R)/a(P) \to C_{D}$. By Lemma \ref{lemma: sections of finite morphism finite}, the boundedness of $Y$ would follow if we prove that $R/P \times_{C_{D}} \Ad(R)/a(P) \to C_{D}$ is finite. For this it suffices to show that $R/P \to \Ad(R)/a(P)$ is finite. Since this can be checked \'etale locally on $C_{D}$, we can further reduce to showing that $G/P \to \GL(\mathfrak{g})/a(P)$ is finite. In fact $G/P \to \GL(\mathfrak{g})/a(P)$ is a closed immersion: it is proper because both the target and the source are $S$-proper, and it is a monomorphism because $(\Ad)^{-1}(a(P)) = P$.} Choose any pair $(E, E_{P})$ in the subset described in (2). The slopes of the associated graded pieces of the filtration of $\Ad(E)$ induced by $E_P$ are of the form $\langle d_{E_{P}}, \alpha \rangle$ for some root $\alpha$ of $(G,T)$. The assumption on $(E,E_{P})$ implies that $\langle d_{E_{P}}, \alpha_i \rangle \geq c$ for all simple roots $\alpha_i$. In combination with the upper bound in (1), this also implies a uniform upper bound $\langle d_{E_{P}}, \alpha_i\rangle \leq h$ for all simple roots $\alpha_i$.\endnote{By (1) we know that there is a uniform upper bound $\langle d_{E_{P}}, \omega_i\rangle \leq H$ for all fundamental weights $\omega_i$. Every fundamental weight can be expressed as a nonnegative linear combination of simple roots
\[ \omega_i = \sum_{j \in I} c_{i,j} \alpha_j\]
where $c_{i,i} >0$ and $c_{i,j} \geq 0$ for all $i, j$. This implies
\[ \alpha_i = \frac{1}{c_{i,i}} \left( \omega_i - \sum_{j \neq i} c_{i,j} \alpha_j\right)\]
By using the upper bound $\langle d_{E_{P}}, \omega_i\rangle \leq H$ and the lower bound $\langle d_{E_{P}}, \alpha_i \rangle \geq c$, we get
\[ \langle d_{E_{P}}, \alpha_i \rangle  = \frac{1}{c_{i,i}} \left( \langle d_{E_{P}}, \omega_i\rangle - \sum_{j \neq i} c_{i,j} \langle d_{E_{P}}, \alpha_j \rangle \right) \leq \frac{1}{c_{i,i}}\left(H - \sum_{j \neq i} c_{i,j}\cdot c\right)\]
This yields the desired uniform upper bound on $\langle d_{E_{P}}, \alpha_i \rangle$ depending only on the root data, and the constants $c$ and $H$.} Since every root is a combination of simple roots, this yields uniform uppper and lower bounds on the slopes of the associated graded pieces of the filtration of $\Ad(E)$ induced by any $E_{P}$ in the set. By \cite[Lemma 6.8]{torsion-freepaper}, this implies that the set of all filtrations of elements in $\{\Ad(E)\, \vert \, E \in \mathcal{S}\}$ arising from such $E_{P}$ forms a bounded family, as desired.

\end{proof}

\begin{proof}[Proof of \Cref{prop: hn boundedness bung}]
Let $\mathfrak{S} \subset \lvert \Bun_{G}(C)\rvert$ be a bounded set of $G$-bundles on the geometric fibers of $C \to S$. Consider a filtration $f$ of an $E \in \mathfrak{S}$ with $\mu(f)>0$, which corresponds to a parabolic $P \supset B_k$, a parabolic reduction $E_{P}$, and a $\lambda \in \sigma_{E_{P}} \cap X_{*}(T)$. Let $w^*$ denote a maximizer of $\mu$ on $\overline{\sigma}_{E_P} \setminus 0$, which exists because $\mu$ is scale invariant and continuous and $(\overline{\sigma}_{E_{P}} \setminus 0)/\bR^{\times}_{>0}$ is compact. Then either $w^\ast \in \sigma_{E_P}$, or $w^\ast$ lies in the relative interior of a boundary facet $\sigma_{E_{P'}} \subset \overline{\sigma}_{E_P}$ corresponding to the extension $E_{P'}$ of $E_P$ to a larger parabolic $P' \supset P$. $w^\ast$ is a local maximizer of $\mu$ on $(\Span(\sigma_{E_{P'}}) \setminus 0) / \bR^{\times}_{>0}$, so \Cref{lemma: riemann-roch expression line bundle} and \Cref{lemma:lagrange} imply that $w^\ast = \phi_V(d_{E_{P'}})$ up to scaling by a positive constant. Because $\mu(w^\ast) \geq \mu(f)$ by construction, we have shown that for any filtration of a point $E \in \mathfrak{S}$, there is another filtration of $E$ whose underlying parabolic reduction $(E,E_{P'})$ lies in 
\[ \mathfrak{S}_{par} \vcentcolon = \left\{ (E, E_{P}) \, \vert \, E \in \mathfrak{S} \text{ and} \; E_{P} \text{ is a parabolic reduction of $E$ such that $\phi_{V}(d_{E_P}) \in \sigma_{E_{P}}$} \right\}. \]
This set is bounded by \Cref{lemma: boundedness of filtrations vs cones}(2), and hence so is the set of associated graded points of such filtrations.

\end{proof}

\subsection{HN-boundedness for \texorpdfstring{$\cM_{n}^{G}(X)$}{MGn(X)}: projective \texorpdfstring{$X$}{X}} \label{subsection: hn boundedness projective}


We keep the assumptions of \Cref{context:boundedness_simplification}. The fibers of the forgetful morphism $\cM_{n}^{G}(X) \to \Bun_{G}(C)$ have finite automorphism groups, so composition filtrations $\Theta^n \to \cM_{n}^{G}(X)$ with the forgetful morphism induces a map of degeneration fans \cite[Sect.~3.2.2]{hl_instability} for any geometric point $(E,u,\marked)$
\[\DF(\cM_{n}^{G}(X),(E, u,\marked)) \to \DF(\Bun_{G}(C),E).\]
\begin{lem}\label{lemma:filt}
If $X \to S$ is projective, then the resulting map on geometric realizations is a homeomorphism
\[
  |\DF(\cM_{n}^{G}(X),(E, u,\marked))| \xrightarrow{\cong}
|\DF(\Bun_G(C),E)|.
\] 
\end{lem}
\begin{proof}
When $X \to S$ is projective, \Cref{prop: formal line bundle is ample on p fibers} implies that $p : \cM_{n}^{G}(X) \to \Bun_{G}(C)$ is relatively representable by a disjoint union of projective DM stacks, so the claim follows from \Cref{T:degeneration_fan_DM_morphism}.
\end{proof}

Using \Cref{lemma:filt} we will regard $|\DF(\cM_{n}^{G}(X),(E, u,\marked))|$, the space of filtrations of a point in the stack
$\cM_{n}^{G}(X)$, as a union of the cones $\overline{\sigma}_{E_{P}}$. We first analyze the constant term $\mu^0_{\delta} := (-\wt(\cD(V))- \delta_{\mrk} \wt(\cL_{\mrk}) + \delta_{\gen} \ell_{\gen}) / \sqrt{b}$ (see \Cref{defn: numerical invariant}) on these cones.

\begin{lem} \label{lemma:mumford_weight} There is a finite subset $\Gamma \subset X^\ast(T)_\bQ$ such that for any given geometric point $x = (E, u,\marked)$ of $\cM_{n}^{G}(X)$, there are subsets $\Sigma_{x,\gen}, \Sigma_{x,\mrk} \subset \Gamma$ such that for every cone $\overline{\sigma}_{E_{P}} \subset |\DF(\cM_{n}^{G}(X),x)|$ we have $\mu^0_{\delta}(w) = \ell_{E_P}^{\delta,x}(w) / \| w\|_b$, where
\begin{equation}\label{E:numerical_invariant_linear_part} \ell_{E_P}^{\delta,x}(w) \vcentcolon = (\phi_V(d_{E_P}),w)_b +\delta_{\mrk} \min\{\ell(w)| \ell \in \Sigma_{x,\mrk}\} + \delta_{\gen} \min\{\ell(w) | \ell \in \Sigma_{x,\gen}\}.
\end{equation} 
\end{lem}
\begin{proof}
Let $\pi: X/G \to BG$ be the natural schematic and projective morphism. Since the line bundle $L_{\gen}$ is $S$-semiample, there is a $G$-equivariant coherent sheaf $\cW$ and an equivariant morphism $X \to \bP_S(\cW):= \Proj_{S}(\Sym^{\bullet}(\cW))$ such that $\cO_{\bP_{S}(\cW)}(1)$ pulls back to $L_{\gen}^{\otimes m}$ on $X$ for some $m>0$. By restricting the $G$-comodule structure to the maximal split torus $T$, we obtain a weight decomposition $\cW = \bigoplus_{\xi\in X^*(T)} \cW_{\xi}$. Since the sheaf $\cW$ is coherent, there are only finitely many $\xi \in X^*(T)$ such that $\cW_{\xi} \neq 0$, let us denote by $\Psi \subset X^*(T)$ this finite subset of weights. Set $\Gamma := \frac{1}{m} \Psi$, viewed as a set of linear functionals on $X_{*}(T)_{\bQ}$. By the Hilbert-Mumford criterion applied to the $S$-fibers of $X$, it follows that for all algebraically closed $\cY$-fields $k$ and all $k$-points $x = (E, u,\marked) \in \cM_{n}^{G}(X)$ there is a subset $\Sigma_{x,\gen} \subset \Gamma$ such that for all cones $\overline{\sigma}_{E_{P}} \subset |\DF(\cM_{n}^{G}(X),x)|$ the function  $\ell_{\gen}$ on $\overline{\sigma}_{E_{P}}$ is given by $\ell_{\gen}(w) = \min \{ \ell(w) \, \mid \, \ell \in \Sigma_{x,\gen}\}$.\endnote{In order to see this, we can restrict to the images of $x$ in $\cY$ and $S$, and so we can assume without loss of generality that $\cY = S = \Spec(k)$. The corresponding restriction $\cW|_{k}$ is a linear $G_{k}$-representation, with $T_{k}$-weights given by a subset of $\Psi$. Therefore, after passing to the $k$-fibers, we can assume that $\bP(\cW)$ is a projective space over a field equipped with a linear representation of $G$. The data $(E,u)$ coming from the point $x$ can be viewed as a morphism $\widetilde{C} \to X/G$. By restricting to the generic point $\eta_{C} \in \widetilde{C}$, we obtain a morphism $p:\Spec(k(C)) \to X_{k(C)}/G_{k(C)}$. The cone $\overline{\sigma}_{E_{P}} \subset |\DF(\cM_{n}^{G}(X),x)|$ can be identified with a cone $\overline{\sigma}_{E_{P_{k(C)}}}$ in the degeneration fan $|\DF(X_{k(\eta_{C}}/G_{k(C)},p)|$. The function $\ell_{\gen}$ on $\overline{\sigma}_{E_{P}}$ is by definition the Hilbert-Mumford weight of the $G_{k(C)}$-equivariant line bundle $L_{\gen}$ on $X_{k(C)}$. Therefore $m \cdot \ell_{\gen}$ is given by the Hilbert-Mumford weight of $L_{\gen}^{\otimes m}$. The morphism $f: X_{k(C)}/ G_{k(C)} \to \bP(\cW_{k(C)}) / G_{k(C)}$ induces a morphism of degeneration fans $|\DF(X_{k(C)}/G_{k(C)},p)| \to |\DF(\bP(\cW_{k(C)})/G_{k(C)},f(p))|$, and identifies the cone $\overline{\sigma}_{E_P}$ with $\overline{\sigma}_{E_{P_{k(C)}}} \subset |\DF(\bP(\cW_{k(C)})/G_{k(C)},f(p))|$. Since $L_{\gen}^{\otimes m}$ is the pullback of $\cO_{\bP_{S}(\cW)}(1)$, the function $m \cdot \ell_{\gen}$ on $\overline{\sigma}_{E_{P_{k(C)}}} \subset |\DF(\bP(\cW_{k(C)})/G_{k(C)},f(p))|$ is given by the Hilbert-Mumford weight of the $G_{k(C)}$-equivariant line bundle $\cO_{\bP_{S}(\cW)}(1)$. Now $f(p)$ is a point of the projective space $\cO_{\bP_{S}(\cW_{k(C)})}$ equipped with a linear $G$-action. If we denote by $m \cdot \Sigma_{x,\gen} \subset \Psi$ the subset of $T$-weights appearing on the nonzero homogeneous coordinates of any lift of $f(p)$ to the total space of $\cW_{k(C)}$, then the Hilbert-Mumford weight $m\cdot \ell_{\gen}$ of $\cO_{\bP_{S}(\cW)}(1)$ is given by $m \cdot \ell_{\gen}(w) = \min \{ \xi(w) \, \mid \, \xi \in m \cdot \Sigma_{x,\gen}\}$ \cite[Chpt. 2.2, Prop. 2.3]{mumford-git}. It follows that $\ell_{\gen}(w) = \min \{ \ell(w) \, \mid \, \ell \in \Sigma_{x,\gen}\}$.}

A similar reasoning yields finitely many linear functionals $\Gamma'$ for $-\wt(\cL_{\mrk})$, by mapping the projective $S$-scheme $X^n$ equivariantly to some $\bP_{S}(\cW')$. We conclude by taking our set of linear functionals to be the union $\Gamma \cup \Gamma'$, and applying the expression for $-\wt(\cD(V))$ from \Cref{lemma: riemann-roch expression line bundle}.
\end{proof}

Note that $\ell_{E_P}^{\delta,x}$ is a continuous, concave, piecewise linear function on $\overline{\sigma}_{E_P}$. This allows us to reformulate the problem of maximizing $\mu_{\delta}^0$ on $(|\DF(\cM_{n}^{G}(X),(E, u,\marked))| \setminus 0) / \bR^\times_{>0}$ as a convex optimization problem on $|\DF(\cM_{n}^{G}(X),(E, u,\marked))|$ itself. 

\begin{lem} \label{L:optimization}
Let $W$ be a finite vector space over $\bQ$, let $\overline{\sigma} \subset W_\bR$ be a closed polyhedral cone, let $\ell(w)$ be a continuous, concave, piecewise linear function on $W_\bR$, and let $b$ be a positive definite quadratic form on $W_\bR$. There is a unique local maximum $w^\ast$ for 
\begin{equation} \label{eqn: quadratic}
\ell(w) - \frac{1}{2} \|w\|_b^2
\end{equation}
on $\overline{\sigma}$, which is also a global maximum. Furthermore:
\begin{enumerate}
    \item $w^\ast$ lies in the interior $\sigma^\circ \subset \overline{\sigma}$ if and only if $w^\ast$ is also the global maximum of \eqref{eqn: quadratic} on $W$.
    \item If $\ell$, $b$, and $\overline{\sigma}$ are rational, then $w^\ast$ is rational.
    \item If $\ell(w)>0$ for some $w \in \overline{\sigma}$, then $w^\ast$ is the unique local maximum of $\mu(w)=\ell(w)/\|w\|_b$ on $(\overline{\sigma} \setminus 0) / \bR^\times_{>0}$ with $\ell(w)>0$, and it is therefore the unique global maximum of $\mu$.
    \item If $\ell(w) = (\ell^\dagger, w)$ for some $\ell^\dagger \in W$, then the unique maximum of \eqref{eqn: quadratic} on $W$ occurs at $w^\ast = \ell^\dagger$.
\end{enumerate} 

\end{lem}

\begin{proof}

The objective function \eqref{eqn: quadratic} is strictly concave, and it tends to $-\infty$ as $\lVert w\rVert_b \to \infty$, so standard results in convex optimization imply the existence and uniqueness of a global and local maximum on a closed convex domain. Strict concavity also implies that if the unconstrained maximizer $u$ of \eqref{eqn: quadratic} lies outside of $\sigma^\circ$, then $w^\ast$ must lie on the boundary $\partial \overline{\sigma}$, because the objective function is strictly increasing along the line segment from $w^\ast$ to $u$. This establishes (1).

For (2)-(4), note that $\overline{\sigma}$ is partitioned into cones on which $\ell(w)$ is linear, and a point $w$ is a local maximum for either \eqref{eqn: quadratic} or $\mu(w)$ if and only if it is a local maximum for the function on each cone whose closure contains $w$. It suffices to prove these claims separately on each sub-cone, so we may assume $\ell(w) = (\ell^\dagger,w)_b$ is linear for the remainder of the proof.

Assume that $\overline{\sigma}$ is defined by inequalities $(g_i,w)_b \leq 0$ for $i=1,\ldots,n$. The Karush–Kuhn–Tucker conditions \cite[7.2]{nonlinear-programmming} imply that $w^\ast$ is a local maximum of $f(w) = \ell(w) - \frac{1}{2} \lVert w \rVert^2_b$ subject to these constraints if and only if there are constants $c_i \geq 0$ such that
\[
w^\ast = \ell^\dagger - \sum c_i g_i,
\]
satisfies $(w^\ast, g_i)_b \leq 0$, and $c_i=0$ if $(w^\ast, g_i)_b<0$.

Note at this point that if $\overline{\sigma}=W$, then $w^\ast = \ell^\dagger$, which establishes (4). More generally, if $W'_\bR \subset W_\bR$ is the linear span of the unique facet of $\overline{\sigma}$ for which $w^\ast$ lies in the relative interior, then $w^\ast$ is the unique unconstrained maximizer of $(\ell^\dagger,w)_b - \frac{1}{2} \lVert w \rVert^2$ on $W'_\bR$, which by (4) is the orthogonal projection of $\ell^\dagger$ onto $W'_\bR$. One can check that this projection is rational if $b$, $\ell$, and $\overline{\sigma}$ are rational, which implies (2).

To show (3), we reformulate the optimization of $\mu$ on $(\overline{\sigma} \setminus 0) / \bR^\times_{>0}$ as the optimization of $f(w) = \ell(w)$ subject to the previous constraints $(w,g_i)_b \leq 0$ and the additional constraint $(w,w)_b=1$. Now the Karush-Kuhn-Tucker conditions state that
\[
2 \lambda w^\ast = \ell^\dagger - \sum c_i g_i,
\]
where $\lambda$ is an arbitrary constant, and $c_i$ satisfy the same conditions as before. The constraint $(w^\ast,w^\ast)_b=1$ then implies that $\lambda = (w^\ast,\ell^\dagger)/2 = \ell(w^\ast)/2$. If $\ell(w^\ast)>0$, then $\lambda>0$ and the conditions for $w^\ast$ to be a local maximum of $\mu$ on $\overline{\sigma} \cap \{\lVert w \rVert_b=1\}$ imply that $w^\ast$ is the normalization of the (unique) global maximum of the previous problem, and vice versa.
\end{proof}

The discussion above allows us to reformulate the condition of HN boundedness in our context.

\begin{defn}
For any $x = (E, u, \marked) \in |\cM^G_n(X)|$ and any parabolic reduction $E_{P}$ of $E$, let $w^\ast_{\delta,x} \in \Span(\sigma_{E_P})$ denote the unique unconstrained maximum of $\ell_{E_P}^{\delta,x}(w)-\frac{1}{2} \lVert w \rVert^2_b$ on $\Span(\sigma_{E_P})$, with $\ell_{E_P}^{\delta,x}(w)$ as in \eqref{E:numerical_invariant_linear_part}.
\end{defn}

\begin{cor} \label{C:hn_boundedness_criterion}
In \Cref{context:boundedness_simplification}, assume $X \to S$ is projective. For any $\mathfrak{S} \subset |\cM^G_n(X)|$, let
\[
\mathfrak{S}_{par,\delta} := \left\{ E_P \left| \exists x = (E,u,\marked) \in \mathfrak{S} \text{ s.t. } E_P \text{ is a }P \text{-reduction of }E, \text{ and } w^\ast_{\delta,x} \in \sigma_{E_P} \right. \right\}.
\]
Then for any $x \in \mathfrak{S}$ and filtration $f$ of $x$ with $\mu^0_{\delta}(f)>0$, there is another filtration $f'$ of $x$ such that $\mu_{\delta}(f') \geq \mu_{\delta}(f)$ and the underlying parabolic bundle $E_{P'}$ of $f'$ lies in $\mathfrak{S}_{par,\delta}$.
\end{cor}
\begin{proof}
By \Cref{lemma:filt}, $f$ amounts to a parabolic reduction $E_P$ of the $G$-bundle underlying $x$ and an integral point of the cone $\sigma_{E_P}$. By \Cref{lemma:mumford_weight}, \Cref{L:optimization} applies to the function $\mu_{\delta}^0$. It implies that there is an integral maximizer $w$ of $\mu_{\delta}^0$ on $\overline{\sigma}_{E_P}$, and if $P' \supseteq P$ is the unique parabolic such that $w$ lies in the relative interior of the facet $\sigma_{E_{P'}} \subset \overline{\sigma}_{E_P}$, where $E_{P'}$ is the extension of $E_P$ to $P'$, then $w$ is a positive multiple of the vector $w^\ast_{\delta,x}$ for the reduction $E_{P'}$.

So either $\mu^0_{\delta}(w)>\mu^0_{\delta}(f)$, in which case $w$ corresponds to the sought-after filtration $f'$, or $\mu^0_{\delta}(w)=\mu^0_{\delta}(f)$, in which case the uniqueness of $w$ up to positive scale implies that $E_P \in \mathfrak{S}_{par,\delta}$ already, so we can take $f'=f$.
\end{proof}

The basic idea behind the proof of \Cref{prop: hn boundedness bung} above was to observe that if $\mathfrak{S}$ is bounded, then so is $\mathfrak{S}_{par,\delta}$, and thus \Cref{C:hn_boundedness_criterion} implies HN boundedness in the case where $X=S$. This strategy will also work to verify HN boundedness in the following.

\begin{prop} \label{prop: hn boundedness gauged maps projective case}
Under the assumptions that $X$ is projective over $S$ and \Cref{context:boundedness_simplification}, the numerical invariant $\mu_{\delta}$ on $\cM_{n}^{G}(X)_d$ satisfies HN-boundedness.
\end{prop}

We will prove this after establishing three technical lemmas.
\begin{lem} \label{lemma:bounded_movement}
There is a constant $c>0$ such that for any point $x= (E, u, \marked)$, any parabolic reduction $E_{P}$, and any $\delta, \gamma \in \bR^{2}_{\geq 0}$, we have $\|w^\ast_{\delta,x} - w^\ast_{\gamma,x}\|_b \leq c \lvert \delta-\gamma\rvert$, where $\lvert - \rvert$ denotes the $L^1$ norm on $\bR^2$.
\end{lem}
\begin{proof}
Let $\delta_t = (\delta^{\mrk}_{t},\delta^{\gen}_{t})$ for $t \in [0,1]$ denote the linearly parameterized line segment with $\delta_0 = \delta$ and $\delta_1 = \gamma$. One can partition $\Span(\sigma_{E_P}) = \bigcup_{i\in I} \cC_i$ into a finite disjoint union of relatively open convex polyhedral cones $\cC_i$ on which $\ell_t(w) := \ell_{E_P}^{\delta_t,x}(w)$ is linear. More specifically,
\[
\ell_t(w) = (w,\phi_V(d_{E_P}) + \delta^{\gen}_{t} \ell_{i,\gen}^\dagger + \delta^{\mrk}_{t} \ell_{i,\mrk}^\dagger)_b \text{ for } w \in \cC_i,
\]
where $\ell_{i,\mrk}\in \Sigma_{x,\mrk}$ and $\ell_{i,\gen} \in \Sigma_{x,\gen}$ are weights realizing the respective minima in \eqref{E:numerical_invariant_linear_part} on the cone $\cC_i$.

The unconstrained maximizer $w_t^\ast \in \Span(\sigma_{E_P})$ must occur on some $\cC_i$, in which case it is the unique maximizer of $\ell_t(w)-\frac{1}{2} \lVert w \rVert^2_b$ on $W_i := \Span(\cC_i)$. By \Cref{L:optimization}(4), this means
\[
w_t^\ast = u_{i,t} := \proj^\perp_{W_i} \left(\phi_V(d_{E_P}) + \delta^{\gen}_{t} \ell_{i,\gen}^\dagger + \delta^{\mrk}_{t} \ell_{i,\mrk}^\dagger \right),
\]
for some $i$ for which $u_{i,t} \in \cC_i$, and the maximal value is $\frac{1}{2} \lVert u_{i,t} \rVert^2_b$. One can write $[0,1]$ as a union of closed intervals such that on each interval $w_t^\ast = u_{i,t}$ for a fixed $i$,\endnote{$w_t^\ast$ is the minimizer of $\frac{1}{2} \lVert u_{i,t} \rVert^2_b$ among all $i$ such that $u_{i,t} \in \overline{\cC_i}$. For each $i$, the set of $t \in [0,1]$ for which $u_{i,t} \in \overline{\cC_i}$ is a closed interval, so we first partition $[0,1]$ into intervals on which the set of such $i$ is constant. On each of these intervals, $\lVert u_{i,t} \rVert^2_b$ is a real quadratic function, and one can further partition the interval into closed subintervals on which the minimizing index $i$ is constant.} so by the triangle inequality it suffices to prove the inequality assuming that $w_t^\ast = u_{i,t}$ for a single $i$ and for all $t \in [0,1]$. In this case, the triangle inequality gives
\begin{align*}
\lVert w_1^\ast - w_0^\ast \rVert_b &= \lVert \proj_{W_i}^\perp \left( (\delta^{\gen}_1-\delta^{\gen}_0) \ell_{i,\gen}^\dagger + (\delta^{\mrk}_1-\delta^{\mrk}_0) \ell_{i,\mrk}^\dagger \right) \rVert_b \\
&\leq |\delta^{\gen}_1-\delta^{\gen}_0| \lVert \ell_{i,\gen}^\dagger \rVert_b + |\delta^{\mrk}_1-\delta^{\mrk}_0| \lVert \ell_{i,\mrk}^\dagger \rVert_b
\end{align*}
The claim now follows from the observation that by \Cref{lemma:mumford_weight}, as $x$ and $E_P$ vary, only finitely many characters appear as $\ell_{i,\mrk}$ and $\ell_{i,\gen}$ in this inequality.
 \end{proof}
 \begin{lem} \label{lemma: distance vs distance to cone}
There exists a constant $c_1>0$ such that for any $G$-bundle $E$, parabolic reduction $E_P$, $h \in \bR$, and $w_1, w_2 \in \Span(\sigma_{E_P})$, $w_1 \in h\cdot \rho_P^{\vee} + \overline{\sigma}_{E_{P}}$ implies that $w_2 \in (h - c_1 \|w_2 - w_1\|_b) \cdot \rho_P^{\vee} + \overline{\sigma}_{E_{P}}$.
 \end{lem}
 \begin{proof}
 Let $c_1$ denote the maximum of the operator norms with respect to $b$ of the simple roots $\alpha_i$, which we view as linear functionals on $X_{*}(T)_{\bR}$. If $w_1 \in h\cdot \rho_P^{\vee} + \overline{\sigma}_{E_{P}}$, then $\alpha_i(w_1) \geq h$ for all $i \in I_P$. Hence, for all $i \in I_P$ we have 
 \[\alpha_i(w_2) = \alpha_i(w_1) + \alpha_i(w_2 -w_1)  \geq \alpha_i(w_1) -c_1\| w_2 - w_1\|_b \geq h- c_1\| w_2 - w_1\|_b\]
 These inequalities imply that $w_2 \in (h - c_1 \|w_1 - w_2\|_b) \cdot \rho_P^{\vee} + \overline{\sigma}_{E_{P}}$, as desired.
\end{proof}

\begin{lem} \label{lemma: boundedness of constant term semsitable stack}
Suppose that $X \to S$ is projective and we are in \Cref{context:boundedness_simplification}. The set $\mathcal{H}$ of $\mu_\delta^0$-semistable points, i.e., $x\in  |\cM_{n}^{G}(X)_d|$ such that $\mu_{\delta}^0(f) \leq 0$ for every filtration $f$ of $x$, is bounded.
\end{lem}
\begin{proof}

If we let $d' \in H_2(BG)$ be the image of $d$ under the homomorphism $H_2(X/G) \to H_2(BG)$, then the forgetful morphism $\cM_n^G(X)_d \to \Bun_G(C)_{d'}$ is quasi-compact by the proof of \Cref{prop: formal line bundle is ample on p fibers}. It therefore suffices to show that the image of $\cH$ in $\Bun_G(C)_{d'}$ is bounded. By \Cref{lemma: boundedness of HN stratification relative case}, it suffices to show that the image of $\cH$ is contained in $\Bun_G(C)^{\mu \leq \gamma}$ for some $\gamma$. We will therefore complete the proof by showing that for any $x \in \cH$ the underlying $G$-bundle $E$ of $x$ lies in $\Bun_G(C)^{\mu \leq h |\delta|}$, where $h := \max_{\ell \in \Gamma} \lVert \ell \rVert_b$ is the maximum operator norm on the finite subset $\Gamma \subset X^{*}(T)_{\mathbb{R}}$ from \Cref{lemma:mumford_weight}.

Consider a parabolic reduction $E_P$ of $E$. \Cref{lemma:mumford_weight} states that on $\sigma_{E_P}$, there are subsets $\Sigma_{p,\gen}, \Sigma_{p,\mrk} \subset \Gamma$ such that the constant term $\mu_\delta^0$ is related to the numerical invariant $\mu = -\wt(\cD(V)) / \sqrt{b}$ for $E$ by the formula
\[ \mu_{\delta}^0(w) = \mu(w) + \frac{ \delta_{\mrk} \min_{\ell \in \Sigma_{p,\mrk}}\{\ell(w)\} +\delta_{\gen} \min_{\ell\in \Sigma_{p,\gen}}\{\ell(w)\}}{\|w\|_b}.\]
By the definition of $h$, we have $\min_{\ell \in \Sigma_{p,\gen}}\{\ell(w)\} \geq - h \cdot \|w\|_b$, and similarly for $\Sigma_{p,\mrk}$. Therefore, we can bound the expression above by
\begin{equation} \label{equation: last inequality boundedness of strata general}
    \mu_{\delta}^0(w) \geq \mu(w) + \frac{ -\delta_{\mrk} \cdot h\cdot \|w\|_b -\delta_{\gen} \cdot h \cdot \|w\|_b}{\|w\|_b} \geq \mu(w) - h \cdot \lvert\delta\rvert
\end{equation}
The assumption $x \in \cH$ implies that $\mu^0_{\delta}(w)\leq 0$, hence $\mu(w) \leq h\lvert\delta\rvert$. This applies to any filtration of $E$ in $\Bun_G(C)$ by \Cref{lemma:filt}, so $E \in |\Bun_G(C)^{\mu \leq h|\delta|}|$ as desired.
\end{proof}

\begin{proof}[Proof of \Cref{prop: hn boundedness gauged maps projective case}]

Let $\mathfrak{S} \subset \lvert \cM_{n}^{G}(X)_d \rvert$ be a bounded subset. By \Cref{lemma:bounded_movement}, for any $x \in \mathfrak{S}$ and parabolic reduction $E_P$ of the $G$-bundle underlying $x$ we have
\[
\| w^*_{\delta, x} - \phi_{V}(d_{E_P})\|_{b} = \| w^*_{\delta,x} - w^{*}_{0,x}\|_{b} \leq c \lvert\delta\rvert
\]
\Cref{lemma: distance vs distance to cone} then implies that if $w^\ast_{\delta,x} \in \sigma_{E_P} \subset \Span(\sigma_{E_P})$, then
\[\phi_{V}(d_{E_P}) \in - c_1 c \lvert\delta\rvert \cdot \rho_P^{\vee} + \overline{\sigma}_{E_{P}}\]
It follows from \Cref{lemma: boundedness of filtrations vs cones} that the set $\mathfrak{S}_{par, \delta}$ of \Cref{C:hn_boundedness_criterion} is bounded.

Let $\mathfrak{S}' \subset |\cM^G_n(X)_d|$ be the subset of points whose underlying $G$ bundle is the associated graded of a $P$-bundle in $\mathfrak{S}_{par,\delta}$.\endnote{By the associated graded $G$ bundle for a parabolic bundle $E_P$, we mean the $G$ bundles induced by the homomorphism of groups $P \to P/U \to G$, where $U$ is the unipotent radical and $P/U \to G$ is the unique inclusion that is the identity on $T \subset P$.} Because $\cM^G_n(X)_d \to \Bun_G(C)$ is finite type, the subset $\mathfrak{S}'$ is bounded. Finally, we let $\mathfrak{S}'' = \mathfrak{S}' \cup \cH \subset |\cM^G_n(X)_d|$, where $\cH$ is the $\mu_{\delta}^0$-semistable locus. $\mathfrak{S}''$ is then bounded by \Cref{lemma: boundedness of constant term semsitable stack}.

To complete the proof of HN boundedness, we show that for any filtration $f$ of $x \in \mathfrak{S}$ with $\mu_{\delta}(f)>0$, there is another filtration $f'$ of $x$ with $\mu_{\delta}(f') \geq \mu_{\delta}(f)$ and $f'(0) \in \mathfrak{S}''$. Indeed, if $x$ is $\mu_{\delta}^0$-unstable, this follows immediately from \Cref{C:hn_boundedness_criterion}. On the other hand, if $x \in \mathfrak{S}$ is $\mu_{\delta}$-unstable but belongs to $\cH$, the for any $\mu_{\delta}$-destabilizing filtration $f$, the constant term $\mu^0_{\delta}(f)=0$. Since $\mu_{\delta}^0$ is $\Theta^2$-monotone (\Cref{thm: monotonicity of gauged maps}), we can apply \cite[Lemma 5.2.12]{hl_instability} to conclude that the associated graded $f(0)$ also lands in the $\mu_{\delta}^0$-semistable locus $\cH \subset \mathfrak{S}''$.
\end{proof}

\subsection{HN-boundedness for \texorpdfstring{$\cM_{n}^{G}(X)$}{MGn(X)}: general case}
In this subsection, we consider $\cM_{n}^{G}(X)$ for a projective-over-affine $G$-scheme $X$ over $S$ with affinization $A_{X}$, as in our setup in \Cref{subsection: notation}.

\begin{lem} \label{lemma: embedding projective over affine}
Suppose that $S$ is affine and $A_X /\!/ \, G = S$. There exist an $S$-projective $G$-scheme $R$, a closed $G$-equivariant embedding $X \hookrightarrow R \times A_X$, and positive integers $q,m >0$ such that
\begin{enumerate}
    \item There is a $S$-very-ample $L_{R,\amp} \in \Pic^G(R)$ such that $L_{\amp} \in \Pic^{G}(X)_\bQ$ is the restriction of $\frac{1}{mq}L_{R,\amp}$.
    \item There is an $S$-semi-ample $L_{R,\gen} \in \Pic^{G}(R)$ such that $L_{\gen}$ is the restriction of $\frac{1}{q}L_{R,\gen}$.
    \item There is an $S$-semi-ample $L^{R, \mrk}\in \Pic^{G^n}(R^n)$ such that $L^{\mrk}$ is the restriction of $\frac{1}{q}L^{R, \mrk}$.
    \item There is an $S$-projective $G$-scheme $Q$ equipped with an $S$-ample $G$-equivariant line bundle $\cO_Q(1)$ and a $G$-equivariant open immersion $A_X \subset Q$ such that $A_X$ is the $\cO_Q(1)$-semistable locus. Moreover, there is a $G$-equivariant isomorphism $\cO_Q(1)|_{A_X} \cong \cO_{A_X}$, where the left-hand side is equipped with the trivial $G$-equivariant structure.
\end{enumerate}
\end{lem}
\begin{proof}
For (1)-(3): By \cite[\href{https://stacks.math.columbia.edu/tag/0FVC}{Tag 0FVC}]{stacks-project}, there exists a positive integer $m>0$ such that the line bundle $F:=L_{\amp}^{\otimes m} \otimes L_{\gen}^{-1} \otimes (\bigotimes_{i=1}^n L_i^{\mrk})^{-1}$ is $S$-very-ample. For each of the semi-ample line bundles $F$, $L_{\gen}$, $L_i^{\mrk}$, there is a morphism to the projectivization of a $G$-equivariant coherent sheaf over $S$ such that the line bundle is the pullback of $\frac{1}{q} \cO(1)$, where $q>0$ is the same sufficiently divisible integer for each bundle. $R$ is the product of these projectivizations of coherent sheaves. The argument for (4) is similar to the proof \cite[Prop 5.11]{AHLH}.\endnote{Here is the argument in more detail: We denote by $\pi: X/G \to BG$ the structure morphism. Since $L_{\gen}$ is $S$-semi-ample, there exists some positive integer $l$ such that the line bundle $L_{\gen}^{\otimes l}$ is globally generated on $X$, in other words $\pi^* \pi_*L_{\gen} \to _{\gen}$ is surjective. Since $L_{\gen}$ is coherent and $\pi_*L_{\gen}$ is the colimit of its coherent subsheaves \cite[Pro. 15.4]{lmb-champs_algebriques}, there exits a coherent $G$-equivariant subsheaf $\cW_{\gen} \subset \pi_*L_{\gen}$ such that the induced morphism $\pi^*\cW_{\gen} \to L_{\gen}$ remains surjective. This defines a morphism from $X$ into $\bP(\cW_{\gen}) := \Proj_{S}(\Sym^{\bullet}(\cW_{\gen}))$ such that the pullback of the universal line bundle $\cO_{\bP(\cW_{\gen})}(1)$ is $L_{\gen}^{\otimes l}$. A similar reasoning yields for each $L_i^{\mrk}$ a $G$-equivariant coherent sheaf $\cW_i$ over $S$ and a morphism $X \to \bP(\cW_i)$ such that the pullback of $\cO_{\bP(\cW_i)}(1)$ is $(L_i^{\mrk})^{\otimes l_i}$ for some positive integer $l_i$. After replacing $l, l_1, l_2, \ldots, l_n$ with their common multiple $q$, we can assume that all of these integers are equal. We set $R_1 := \bP(\cW_{\gen}) \times \prod_{i=1}^n \bP(\cW_i)$, and consider $X \to R_1$. There are $S$-semi-ample $G$-equivariant line bundles $L_{R_1,\gen}$ and $L_i^{R_1, \mrk}$ that pullback to the $q^{th}$ power of the corresponding line bundles on $X$.

By \cite[\href{https://stacks.math.columbia.edu/tag/0FVC}{Tag 0FVC}]{stacks-project}, there exists a positive integer $m>0$ such that the line bundle $F:=L_{\amp}^{\otimes m} \otimes L_{\gen}^{-1} \otimes (\bigotimes_{i=1}^n L_i^{\mrk})^{-1}$ is $S$-very-ample. The $q^{th}$-multiple $F^{\otimes q}$ is also $S$-very-ample, and so the morphism $X \to \bP(\pi_*(F^{\otimes q}))$ is a (locally closed) immersion of schemes. Using the fact that $\pi_*(F)$ is the colimit of its $G$-equivariant subsheaves and the proof of \cite[\href{https://stacks.math.columbia.edu/tag/02NP}{Tag 02NP}]{stacks-project}, we see that there exists a coherent $G$-equivariant subsheaf $\cW \subset \pi_*(F)$ such that the morphism $X \to \bP(\cW)$ remains well-defined and an immersion. Set $R := R_1 \times \bP(\cW)$, and take the product morphism $X \to R \times A_X$, which is an immersion. Since $X \to A_X$ is proper and $R$ is separated, the morphism $X \to R \times A_X$ is also proper, and so it must be a closed immersion. By construction, the $G$-equivariant $S$-very-ample line bundle $L_{R, \amp} := \cO_{\bP(\cW_{\gen})}(1) \boxtimes (\boxtimes_{i=1}^n \cO_{\bP(\cW_i)}(1)) \boxtimes \cO_{\bP(\cW)}(1)$ restricts to $L_{\amp}^{\otimes mq}$. Similary the $S$-semiample line bundle $L_{R,\gen} : = L_{R_1,\gen} \boxtimes \cO_{\bP(\cW)}$ and $L_i^{\mrk,R} : = L_i^{R_1, \mrk} \boxtimes \cO_{\bP(\cW)}$ restrict to the $q^{th}$ powers of their corresponding line bundles on $X$. By setting $L^{R, \mrk}: = \boxtimes_{i=1}^n L_i^{R, \mrk}$, we have concluded the verification of (1), (2) and (3).

We are left to find the compactification $A_X \subset Q$ as in (4). Since $A_X$ is of finite type over $S$, we can use \cite[Pro. 15.4]{lmb-champs_algebriques} to find a $G$-equivariant coherent subsheaf $\cK \subset A_X$ over $S$ such that the induced morphism $A_X \to \Spec(\Sym^{\bullet}(\cK))$ is a closed immersion. The scheme $\Spec(\Sym^{\bullet}(\cK))$ admits a $G$-equivariant open embedding into the $S$-projective scheme $\bP(\cO_{S} \oplus \cK)$, where we equip $\cO_{S}$ with the trivial $G$-equivariant structure. The complement of $\Spec(\Sym^{\bullet}(\cK)) \subset \bP(\cO_{S} \oplus \cK)$ is cut out by a $G$-invariant section $s \in H^0(\cO_{\bP(\cO_{S} \oplus \cK)}(1))$, which shows that $\cO_{\bP(\cO_{S} \oplus \cK)}(1)$ restricts to the trivial $G$-equivariant line bundle on $\Spec(\Sym^{\bullet}(\cK))$. Set $Q$ to be the scheme theoretic closure of the morphism $A_X \to \Spec(\Sym^{\bullet}(\cK)) \to \bP(\cO_{S} \oplus \cK)$, and set $\cO_{Q}(1)$ to be the restriction of $\cO_{\bP(\cO_{S} \oplus \cK)}(1)$ to $Q$. Since $A_X \subset \Spec(\Sym^{\bullet}(\cK))$, the restriction of $\cO_{Q}(1)$ to $A_X$ is the trivial $G$-equivariant line bundle. By a similar argument as in the proof \cite[Prop 5.11]{AHLH} (note that reducedness in that proof is not needed, since $A_X$ is scheme theoretically dense inside $Q$), it follows that $A_X \subset Q$ is the semistable locus of $\cO_{Q}(1)$.}
\end{proof}

\begin{lem} \label{lemma: behavior of gauged maps under immersions}
If $X \to Y$ is a closed (resp. open) immersion of projective-over-affine $G$-schemes over $S$, then the induced map $\cM_n^G(X) \to \cM_n^G(Y)$ is a closed (resp. open) immersion.
\end{lem}
\begin{proof}
Let $T$ be a Noetherian $S$-scheme, and choose $T \to \cM_{n}^{G}(Y)$ corresponding to a tuple $(E, \; u: \widetilde{C} \to E(Y), \; \marked)$. The fiber product $\cM_{n}^{G}(X) \times_{\cM_{n}^{G}(Y)} T$ is the subfunctor of $\Hom(\dash, T)$ consisting of morphisms $T' \to T$ such that the base-change $u_{T'}: \widetilde{C}_{T'} \to E(Y)$ factors through $E(X)$. This is represented by a closed (resp. open) subscheme of $T$.\endnote{For closed immersions: let $Z\subset \widetilde{C}$ denote the closed subscheme defined by $Z:= \widetilde{C} \times_{E(Y)} E(X)$. Then $\cM_{n}^{G}(X) \times_{\cM_{n}^{G}(Y)} T$ parametrizes $T' \to T$ such that the immersion $Z_{T'} \hookrightarrow \widetilde{C}_{T'}$ is an isomorphism. If we denote by $\varphi: \cI_{Z} \hookrightarrow \cO_{\widetilde{C}}$ the ideal sheaf of $Z$ inside $\widetilde{C}$, then $T' \to T$ belongs to $\cM_{n}^{G}(X) \times_{\cM_{n}^{G}(Y)} T$ if and only if $\varphi|_{\widetilde{C}_{T'}}: \cI|_{\widetilde{C}_{T'}} \to \cO_{\widetilde{C}_{T'}}$ is the zero morphism. Let $\Hom(\cI, \cO_{\widetilde{C}})$ denote the functor that sends a $T$-scheme $T'$ to the group of homomorphisms $\Hom_{\cO_{\widetilde{C}_{T'}}}(\cI|_{\widetilde{C}_{T'}}, \cO_{\widetilde{C}_{T'}})$, as in \cite[\href{https://stacks.math.columbia.edu/tag/08JS}{Tag 08JS}]{stacks-project}. The zero morphism and $\varphi$ define two sections $0, \varphi: T \to \Hom(\cI, \cO_{\widetilde{C}})$, and we have seen that $\cM_{n}^{G}(X) \times_{\cM_{n}^{G}(Y)} T$ is the equalizer of these two sections. Since $\widetilde{C} \to T$ is proper and flat, we can apply \cite[\href{https://stacks.math.columbia.edu/tag/08K6}{Tag 08K6}]{stacks-project} to see that $\Hom(\cI, \cO_{\widetilde{C}})$ is represented by a separated scheme over $T$, and hence the equalizer is a closed subscheme of $T$.

For open immersions: let $Z \hookrightarrow E(Y)$ denote the closed complement of $E(X)$ equipped with the reduced subscheme structure. Define $W \vcentcolon = \widetilde{C} \times_{E(Y)} Z \hookrightarrow \widetilde{C}$. Since $\widetilde{C} \to T$ is proper, the set-theoretic image $\im(W)$ of the composition $W \hookrightarrow \widetilde{C} \to T$ is closed. A morphism $f: T' \to T$ belongs to $\cM_{n}^{G}(X) \times_{\cM_{n}^{G}(Y)} T$ if and only if the image $f(T') \subset T$ belongs to the complement $T \setminus \im(W)$. Hence $\cM_{n}^{G}(X) \times_{\cM_{n}^{G}(Y)} T \hookrightarrow T$ is represented by the open immersion $T \setminus \im(W) \hookrightarrow T$.}
\end{proof}
For the remainder of the subsection, we set $Y = R \times Q$ and $U = R \times A_X$ as in \Cref{lemma: embedding projective over affine}. For $\beta,\gamma \in \bR$, we define $\cO_{Y}(\beta, \gamma) \vcentcolon = L_{R,\gen}^{\otimes \beta} \boxtimes \cO_{Q}(\gamma) \in \Pic_G(Y)_\bR$, and we let $\cO_U(\beta)$ denote the restriction of this $\bR$-line bundle to $U$. We denote by $\mu_{\delta, \gamma}$ the numerical invariant on $\cM_{n}^{G}(Y)$ given by
\[ \mu_{\delta,\gamma} = \frac{-\wt(\cD(V)) - \delta_{\mrk} \wt(\ev^*(L^{Y,\mrk})) +\delta_{\gen} \ell_{\gen} -  \epsilon \wt(\cL_{Cor}(\epsilon))}{\sqrt{b}}.\]
where $\ell_{\gen}$ is formed using the line bundle $\cO_{Y}(1, \gamma)$, and we use the equivariant bundles $L^{Y, \mrk}:= L^{R, \mrk} \boxtimes \cO_{Q^n}$ and $L_{Y, \amp} := L_{R, \amp} \boxtimes \cO_{Q}(1)$ to define the other parts of the numerical invariant. $\mu_{\delta, \gamma}$ is strictly $\Theta$-monotone by \Cref{thm: monotonicity of gauged maps} and satisfies HN-boundedness by \Cref{prop: hn boundedness gauged maps projective case}, because $Y$ is $S$-projective. It follows that every $\mu_{\delta, \gamma}$-unstable point $x=(E,u, \marked) \in \cM_{n}^{G}(Y)$ has a unique (up to scaling) HN-filtration $HN_{x, \gamma}$ for the numerical invariant $\mu_{\delta, \gamma}$.\endnote{Indeed, \Cref{thm: theta stability paper theorem}(1) implies that the numerical invariant $\mu_{\delta, \gamma}$ defines a weak $\Theta$-stratification of $\cM_{n}^{G}(Y)$.}

Note that, by Lemma \ref{lemma: behavior of gauged maps under immersions} applied to the open immersion $U \to Y$, we have that $\cM_{n}^{G}(U)$ is an open substack of $\cM_{n}^{G}(Y)$.
\begin{prop} \label{prop:HN-boundedness}
Suppose that we are in \Cref{context:boundedness_simplification}; we shall also denote by $d$ the image of the fixed class $d \in H_2(X/G)$ inside the groups $H_2(U/G)$ and $H_2(Y/G)$. Let $\mathfrak{S} \subset |\cM_{n}^{G}(U)_d|$ be a bounded set of points, and choose $\delta \in \bR^{2}_{\geq 0}$. Regarding the elements of $\mathfrak{S}$ as points of $\cM_{n}^{G}(Y)$, there is a constant $\gamma_{\mathfrak{S}, \delta}$ such that for all $\gamma>\gamma_{\mathfrak{S}, \delta}$, the HN-filtration $HN_{x, \gamma}$ of any point $x \in \mathfrak{S}$ is independent of $\gamma$, lies in $\cM_{n}^{G}(U)$, and is the HN-filtration in $\cM_{n}^{G}(U)$ for the numerical invariant $\mu_{\delta}$. In particular, HN-boundedness holds for the numerical invariant $\mu_{\delta}$ on $\cM_{n}^{G}(U)$.
\end{prop}

We shall prove this proposition after introducing a lemma.
\begin{lem} \label{lemma:generically_semistable_filtrations}
Fix a field $k$ over $\cY$. Let $f : \Theta_k \to \cM_{n}^{G}(Y)$ be a filtration such that $f(1) \in \cM_{n}^{G}(U)$. The following are equivalent:
\begin{enumerate}
\item $f$ lies in $\cM_{n}^{G}(U)$.
\item The filtration at the generic point $f_{\gen} : \Theta_{k(C_k)} \to Y /G$ lies in the open substack $U/G$.
\item $\wt(f_{\gen}^\ast(\cO_{Y}(0,1))|_{0}) =0$. (By construction, this weight is always $\geq 0$ because $f(1)$ lands in the $\cO_{Y}(0,1)$-semistable locus $U$.)
\end{enumerate}
\end{lem}
\begin{proof}
The equivalence $(2) \Leftrightarrow (3)$ follows from \cite[Lemma 6.15]{AHLH}.\endnote{We offer an alternative proof using ideas from GIT. We restrict $u: \widetilde{C} \to E(Y)$ and the $G$-bundle $E$ to $\bA^1_{k(C_{k})} \subset \widetilde{C}, C_{\bA^1_{k}}$ in order to obtain $f_{\gen}: \Theta_{k(C_{k})} \to Y/G$. By assumption $f_{\gen}(1)$ lands in the $\cO_{Y}(0,1)$-semistable locus $U$. We want to show that the associated graded point $f_{\gen}(0)$ lands in $U$ if and only if $\wt(f_{\gen}^\ast(\cO_{Y}(0,1))|_{0}) =0$. We can focus on the composition with the projection into the second coordinate \[g_{\gen}: \Theta_{k(C_k)} \to Y_{k(C_k)}/G_{k(C_k)} \to Q_{k(C_{k})}/G_{k(C_k)}\]
Note that $\wt(g_{\gen}^\ast(\cO_{Y}(0,1))|_{0})= \wt(f_{\gen}^\ast(\cO_{Y}(0,1))|_{0})$ by definition.
We shall use the Hilbert-Mumford criterion as in \cite[Chpt. 2]{mumford-git}. The filtration $g_{\gen}$ amounts to the choice of a one-parameter subgroup $\lambda: \bG_{m} \to G_{k(C_{k})}$, and $g_{\gen}(0) = \lim_{t\to 0} \lambda(t) \cdot g_{\gen}(1)$. By our convention, $\bG_m$ acts on $\bA^1_{k(C_{k})}$ with weight $-1$, and so the weight $\wt(g_{\gen}^\ast(\cO_{Y}(0,1))|_{0})$ is $\mu^{\cO_{Y}(0,1)}(\lambda, g_{\gen}(0))$ using the notation in \cite[Chpt. 2]{mumford-git}. Suppose first that $\wt(g_{\gen}^\ast(\cO_{Y}(0,1))|_{0}) > 0$. The point $g_{\gen}(0) \in Q_{k(C_{k})}$ admits a $\mathbb{G}_m$ action arising from $\lambda$. Note that $\mu^{\cO_{Y}(0,1)}(\lambda^{-1}, g_{\gen}(0)) = - \mu^{\cO_{Y}(0,1)}(\lambda, g_{\gen}(0))$, and so the Mumford weight for $\lambda^{-1}$ will be negative. By the Hilbert-Mumford criterion, this means that $g_{\gen}(0)$ is unstable.

Conversely, suppose that $\wt(g_{\gen}^\ast(\cO_{Y}(0,1))|_{0}) = 0$. Assume for the sake of contradiction that $g_{\gen}(0)$ is unstable. Then there exists a one-parameter subgroup $\eta: \bG_m \to G_{k(C_{k})}$ such that the $\mu^{\cO_{Y}(0,1)}(\lambda, g_{\gen}(0))<0$. The images of $\lambda, \eta$ belong to a common apartment in the flag complex of $G_{k(C_{k})}$, so we can assume without loss of generality that they are one-parameter subgroups of a common maximal split torus $T$ inside $G_{k(C_{k})}$. The concrete description of the Hilbert-Mumford weight in \cite[Chpt. 2, S3]{mumford-git} yields a finite subset of characters $\{\chi_i\}_{i \in I} \subset X^{*}(T)$ associated to the point $g_{\gen}(1)$ (they correspond to the weights of the nonzero homogeneneous coordinates for a linearlized projective embedding). Our assumption $\wt(g_{\gen}^\ast(\cO_{Y}(0,1))|_{0}) = 0$ translates into $\max_{i\in I}\{ -\chi_i(\lambda)\} = 0$. Let $J \subset I$ denote the subset of indexes $j \in I$ such that $\chi_i (\lambda) = 0$. The assumption $\mu^{\cO_{Y}(0,1)}(\lambda, g_{\gen}(0))<0$ translates into $\max_{j \in J} \{ -\chi_j(\eta)\} <0$. Consider cocharacters of the form $\eta -N\lambda \in X_{*}(T)$ for $N \in \mathbb{Z}$. For any $i\not\in J$, we have $\chi_i(\lambda) <0$, and so $\chi_i(\eta -N\lambda)>0$ for $N\gg 0$. On the other hand, for $j \in J$ we have $\chi_j(\eta -N \lambda) = \chi_j(\eta) >0$. We conclude that $\max_{i \in I} \{ -\chi_i(\eta -N\lambda)\} <0$ for $N\gg0$. This translates into $\mu^{\cO_{Y}(0,1)}(\eta -N\lambda, \, g_{\gen}(1))<0$, thus contradicting the semistability of $g_{\gen}(1)$.} Also, the implication $(1) \Rightarrow (2)$ is clear, so we shall assume $(2)$ holds and show $(1)$.

The filtration $\Theta_{k} \to \cM_{n}^{G}(Y)$ corresponds to the data of a $\mathbb{G}_m$-equivariant $\bA^1_{k}$-point $(E, u: \widetilde{C} \to E(Y), \marked)$. Consider the composition with the projection
\[ g: \widetilde{C} \xrightarrow{u} E(Y) = E(R) \times E(Q) \to E(Q).\]
We want to prove that $g$ factors through the open subscheme $E(A_X) \subset E(Q)$. %
Since $g|_{\widetilde{C}_{\bA^1 \setminus 0}}$ factors through the affine $C_{\bA^1}$-scheme $E(A_X)$, we must have that $g|_{\widetilde{C}_{\bA^1 \setminus 0}}$ factors through a section $s : C_{\bA^1 \setminus 0} \to E(A_X)$. The condition (2) implies that $s$ can be extended to a section defined 
away from finitely many closed points in the $0$-fiber. Since $C_{\bA^1}$ is a smooth surface, Hartogs's theorem implies that $s$ extends to a section $\widetilde{s}: C_{\bA^1} \to E(A_X)$. By construction, $g|_{\widetilde{C}_{\bA^1 \setminus 0}}$ factors through the closed subscheme $\im(\widetilde{s}) \hookrightarrow E(Q)$, which is contained in $E(A_X)$. Since $\widetilde{C}$ is $\bA^1_{k}$-flat, the open $\widetilde{C}_{\bA^1 \setminus 0} \subset \widetilde{C}$ is scheme-theoretically dense, and therefore $g: \widetilde{C} \to E(Q)$ also factors through $\im(\widetilde{s}) \subset E(A_X)$.
\end{proof}

\begin{proof}[Proof of \Cref{prop:HN-boundedness}]
Suppose that $\delta_{\gen} =0$. Then the term with $\ell_{\gen}$ is irrelevant, and we can choose the trivial $G$-equivariant line bundle $\cO_{X}$ to define $\ell_{\gen}$. In that case $\ell_{\gen} = 0$, and so $\mu_{(\delta_{\mrk},\delta_{\gen})} = \mu_{(\delta_{\mrk}, 1)}$. Therefore, we can assume without loss of generality that $\delta_{\gen} \neq 0$.

Consider a $\mu_{\delta,\gamma}$-unstable point $x = (E,u,\marked) \in \mathfrak{S}$, and let $f = HN_{x,\gamma}$. By \Cref{lemma:generically_semistable_filtrations}, $f$ lies in $\cM^G_n(U)$ if and only if $\ell_{\aux}(f):= -\wt(f^\ast_{\gen}(\cO_Y(0,1))|_0) \geq 0$. Therefore, our goal is to find some universal $\gamma_{\mathfrak{S}, \delta}$ such that if $\ell_{\aux}(f)<0$, then $\gamma \leq \gamma_{\mathfrak{S},\delta}$.

Let $E_P$ be the parabolic reduction of $E$ such that $HN_{x,\gamma} \in \sigma_{E_P} \subset |\DF(\cM_n^G(Y),x)|$. By \Cref{lemma:mumford_weight}, there is a finite set $\Gamma \subset X^\ast(T)_\bQ$, independent of $x,\delta,$ and $\gamma$, and $\Sigma_{x,\mrk},\Sigma_{x,\gen},\Sigma_{x,\aux} \subset \Gamma$ such that the constant term $\mu_{\delta,\gamma}^0 = L(w) / \lVert w \rVert_b$ on $\sigma_{E_P}$, where
\begin{equation}\label{E:two_term_linear_function}
L(w) := (\phi_V(d_{E_P}),w)_b + \delta_{\mrk} \min \limits_{\ell \in \Sigma_{x,\mrk}} \ell(w) + \overbrace{\delta_{\gen} \min \limits_{\ell \in \Sigma_{x,\gen}} \ell(w) + \delta_{\gen} \gamma \min \limits_{\ell \in \Sigma_{x,\aux}} \ell(w)}^{\delta_{\gen} \ell_{\gen}}.
\end{equation}
We have split $\ell_{\gen}$ into a first term corresponding to $\cO_Y(1,0)$ and a second term equal to $-\wt(f^\ast_{\gen}(\cO_Y(0,\gamma))|_0)$. A point $w \in \sigma_{E_P}$ corresponds to a filtration in $\cM_n^G(U)$ if and only if $\ell(w) \geq 0$ for all $\ell \in\Sigma_{x,\aux}$.

The function $L(w)$ on $\Span(\sigma_{E_P})$ is piecewise linear, so we may choose a decomposition $\overline{\sigma}_{E_P} = \bigcup_{i \in I} \cC_i$ as a disjoint union of relatively open convex polyhedral cones $\cC_i$ such that
\[
L(w) = (w, \phi_V(d_{E_P}) + \delta_{\mrk} \ell_{i,\mrk}^\dagger + \delta_{\gen} \ell_{i,\gen}^\dagger + \delta_{\gen} \gamma \ell_{i,\aux}^\dagger)_b \text{ for all } w \in \cC_i,
\]
where $\ell_{i,\mrk} \in \Sigma_{x,\mrk}$, $\ell_{i,\gen} \in \Sigma_{x,\gen}$, and $\ell_{i,\aux} \in \Sigma_{x,\aux}$ are weights realizing the respective minima in \eqref{E:two_term_linear_function} on the cone $\cC_i$. We can choose a single decomposition of this kind that works for all $\delta$ and $\gamma$. Finally, for each $i \in I$, choose a finite (possibly empty) subset $S_i \subset \cC_i$ with the following properties:\endnote{Here is a way to choose such a finite set $S_i$. Fix a finite covering of $X_{*}(T)_{\bR} = \bigcup_j \cD_j$, where each $\cD_j$ is a closed $b$-positive cone (i.e. $(w,v)_b\geq 0$ for all $w,v \in \cD_j$). Then for each $\cC_i$ we get a decomposition into $b$-positive cones 
\[ \cC_i \cap \{w \, \lvert \, \ell_{\aux}(w)<0\} = \bigcup_j (\cC_i \cap \{w \,\lvert \,\ell_{\aux}(w)<0\}\cap \cD_j)\]
We can let $S_i$ consists of a choice of a unit vector in each nonempty $(\cC_i \cap \{w \lvert \ell_{\aux}(w)<0\}\cap \cD_j)$.}
\begin{enumerate}
\item $\forall s \in S_i$, $\ell_{i,\aux}(s) < 0$ and $\|s\|_b=1$, and
\item for any $w \in \cC_i$ with $\ell_{i,\aux}(w)<0$, there is some $s \in S_i$ such that $(w,s)_b \geq 0$.
\end{enumerate}

Because our numerical invariant takes values in $\bR[\epsilon]$, we split the analysis into two cases:

\medskip
\noindent\textit{Case 1: The constant term $\mu_{\delta,\gamma}^0(HN_{x,\gamma})>0$
:}
\medskip

In this case, $HN_{x,\gamma}$ maximizes $\mu_{\delta,\gamma}^0$. By \Cref{L:optimization}, $HN_{x,\gamma}$ corresponds to a positive multiple of the unique unconstrained maximizer $w^\ast$ of the strictly concave function $L(w) - \frac{1}{2} \lVert w \rVert^2_b$ on $W_i := \Span(\cC_i)$, where $\cC_i \subset \sigma_{E_P}$ is the unique cone containing $HN_{x,\gamma}$. \Cref{L:optimization}(4) implies that $w^\ast = \proj^\perp_{W_i}(\phi_V(d_{E_P})+ \delta_{\mrk} \ell_{i,\mrk}^\dagger + \delta_{\gen} \ell_{i,\gen}^\dagger + \delta_{\gen} \gamma \ell_{i,\aux}^\dagger)$, or equivalently
\[
\phi_V(d_{E_P}) \in -(\delta_{\mrk} \ell_{i,\mrk}^\dagger + \delta_{\gen} \ell_{i,\gen}^\dagger + \delta_{\gen} \gamma \ell_{i,\aux}^\dagger) + w^\ast + W_i^\perp.
\]

Now assume $\ell(w^\ast)<0$ for some $\ell \in \Sigma_{x,\aux}$, and hence $\ell_{i,\aux}(w^\ast)<0$. We can choose an $s \in S_i \subset \cC_i$ such that $(s,w^\ast)_b \geq 0$. For this $s$, we have
\begin{equation} \label{E:first_inequality}
(s,\phi_V(d_{E_P}))_b \geq -(\delta_{\mrk} \ell_{i,\mrk}(s) + \delta_{\gen} \ell_{i,\gen}(s) + \delta_{\gen} \gamma \ell_{i,\aux}(s)).
\end{equation}

On the other hand, we know that $\phi_{V}(d_{E_P}) \in  c_{\mathfrak{S}} \cdot \rho_P^{\vee} -\overline{\tau}_{E_{P}}$ by \Cref{lemma: boundedness of filtrations vs cones}, so $\phi_{V}(d_{E_P}) = d_E + c_{\mathfrak{S}} \cdot \rho_P^{\vee} + u_{-}$, where $d_E$ is the image of $d$ under the homomorphism $H_2(U/G) \to H_2(BG)$ and $u_{-} \in -\overline{\tau}_{E_{P}} \cap X_{*}(T')_{\bR}$. Since $(\overline{\sigma}_{E_{P}}, \overline{\tau}_{E_{P}}\cap X_{*}(T')_{\bR})_b \subset \bR_{\geq0}$ (by \Cref{notn:root data torus}), we can use $s \in \cC_{i} \subset \overline{\sigma}_{E_{P}}$ to get
\begin{equation} \label{E:second_inequality}
(s, \phi_{V}(d_{E_P}))_b = (s, d_E)_b + c_{\mathfrak{S}} (s, \rho_P^{\vee})_b + (s, u_{-})_b \leq \|d_E\|_b + c_{\mathfrak{S}} (s, \rho_P^{\vee})_b.
\end{equation}

Combining inequalities \eqref{E:first_inequality} and \eqref{E:second_inequality}, and using the fact that $\ell_{i,\aux}(s)<0$ to solve for $\gamma$, we obtain
\[
\gamma \leq \frac{\|d_E\|_b + c_{\mathfrak{S}} (s, \rho_P^{\vee})_b + \delta_{\mrk} \ell_{i,\mrk}(s) + \delta_{\gen} \ell_{i,\gen}(s)}{- \delta_{\gen} \ell_{i,\aux}(s)}.
\]
As the point $x = (E,u,\marked) \in \cM_n^G(U)_d$ and parabolic reduction $E_P$ of $E$ vary, only finitely many cones $\sigma_{E_P} \subset X_*(T)_{\bR}$, points $\rho_P^{\vee}$, decompositions $\sigma_{E_P} = \bigcup_{i \in I} \cC_i$, and choices $\ell_{i,\mrk},\ell_{i,\gen},\ell_{i,\aux} \in \Gamma$ can possibly arise, and for each one we can choose the subsets $S_i \subset \cC_i$ independently of $E_P$, $\delta$, and $\gamma$. We conclude that there are constants $h_1, h_2, h_3$ independent of $\mathfrak{S}$ or $\delta$ such that
\begin{equation} \label{eqn: bound for gamma}
\gamma \leq \frac{h_1 \|d_E\|_b + h_2 |c_{\mathfrak{S}}| + h_3\lvert\delta\rvert}{\delta_{\gen}}.
\end{equation}

\medskip
\noindent \textit{Case 2:  $\mu_{\delta,\gamma}^0(HN_{x,\gamma})=0$:}
\medskip

The $w \in \sigma_{E_P}$ corresponding to $HN_{x,\gamma}$ maximizes $\mu_{\delta, \gamma}^0$ on the cone $\cC_{i}$ containing it, so the linear numerator of $\mu_{\delta,\gamma}^0|_{\cC_{i}}$ must be identically $0$. It follows that 
\[ \phi_{V}(d_{E_P}) \in -\delta_{\mrk} \ell_{i, \mrk}^{\dagger} - \delta_{\gen} \ell_{i,\gen}^{\dagger} - \gamma \delta_{\gen} \ell_{i, \aux}^{\dagger} + W_i^{\perp}.\]
Assume that $\ell_{i, \aux}(w)<0$. Choose $s \in S_i \subset \cC_i$ with $(s, w)_b\geq0$, as in Case 1. By using $(s, W_i^{\perp})_b =0$, we get
\begin{equation} \label{eqn: first equation second case}
    (s, \phi_{V}(d_{E_P}))_b = -\delta_{\mrk} (s, \ell_{i, \mrk}^{\dagger})_b - \delta_{\gen} (s, \ell_{i,\gen}^{\dagger})_b- \gamma \delta_{\gen} (s, \ell_{i, \aux}^{\dagger})_b
\end{equation}
On the other hand, the same reasoning as in Case 1 yields
\begin{equation} \label{eqn: second equation second case}
    (s, \phi_{V}(d_{E_P}))_b \leq \|d_E\|_{b} + c_{\mathfrak{S}} (s, \rho_P^{\vee})_b
\end{equation}
Combining equations \eqref{eqn: first equation second case} and \eqref{eqn: second equation second case} yields the same bound for $\gamma$ as in Case 1.
\end{proof}

\begin{proof}[Proof of \Cref{thm: HN-boundedness nonconnected}]
By \Cref{prop: reduction to connected groups}, we can assume that we are in \Cref{context:boundedness_simplification}. By Lemma \ref{lemma: behavior of gauged maps under immersions} applied to the closed immersion $X \hookrightarrow R \times A_X$, HN-boundedness for $\cM_{n}^{G}(X)$ follows from HN-boundedness for $\cM_{n}^G(U)$, which was proven in \Cref{prop:HN-boundedness}.
\end{proof}


\section{The \texorpdfstring{$\Theta$}{Theta}-stratification}
\begin{thm} \label{thm: theta stratification in general}
For any $\delta \in (\mathbb{R}_{\geq 0})^2$, the numerical invariant $\mu_{\delta}$ defines a weak $\Theta$-stratification on $\cM_{n}^{G}(X)$. Moreover, if the base scheme $S$ is defined over $\bQ$, then $\mu_{\delta}$ defines a $\Theta$-stratification.
\end{thm}
\begin{proof}
This follows from Theorem \ref{thm: theta stability paper theorem} (1), because $\mu_{\delta}$ is strictly $\Theta$-monotone (Theorem \ref{thm: monotonicity of gauged maps}) and satisfies HN-boundedness (\Cref{thm: HN-boundedness nonconnected}).
\end{proof}

In this section, we will study further properties of this $\Theta$-stratification. For the first two subsections, we will restrict to the following context.
\begin{context} \label{C:boundedness_section_context} \quad
\begin{itemize}
\item The base $S$ is affine connected and Noetherian, and $\mathcal{Y}=S$. $G$ is a connected split reductive group over $S$. $X$ is the total space of a trivial vector bundle on $S$, equipped with a linear action of $G$. 
\item We restrict our attention to $\cM_{0}^{G}(X) = \MpC(X/G)$ without marked points.
\item We assume that $C \to S$ admits a section $p \in C(S)$.
\item We fix a line bundle $\cL$ on $\MpC(X/G)$ that is numerically equivalent to $\cD(V)$ for some almost faithful $G$-representation $V$.
\item We fix a split maximal torus $T \subset G$ with corresponding Weyl group $W$. We will use the notation $N:= X_*(T)$ and $M:= X^*(T)$. 
\item We fix a rational quadratic $W$-invariant norm $b$ on $N_{\mathbb{Q}}$, which gives a norm on graded points of $\MpC(X/G)$.
\item We will assume that the $G$-equivariant rational line bundle $L_{\gen} \in \Pic^{G}(X)_{\mathbb{Q}}$ is of the form $\cO_X(\chi) = \cO_{X} \otimes_{\cO_S} (-\chi)$ (see \Cref{subsection: notation}) for some $\chi\in  \Hom(G, (\mathbb{G}_{m})_S) \otimes \mathbb{Q} = M^W_{\mathbb{Q}}$. We set $\delta =1$.
\item We will take $L_{\amp}$ to be the trivial line bundle on $X/G$.
\end{itemize}
\end{context}

\begin{notn}
For any $\lambda \in N_{\mathbb{Q}}$, we denote by $P_{\lambda} \supset T$ the associated parabolic subgroup. We denote by $L_{\lambda} \supset T$ the corresponding Levi subgroup, and we write $W_{\lambda}$ for the Weyl group of $L_{\lambda}$.
\end{notn}

\begin{defn} \label{defn: l-chi}
Assume we are in Context \ref{C:boundedness_section_context}. For any rational character $\chi \in M^W_{\mathbb{Q}}$, we let
\[\cL(\chi):= \cL \otimes (\ev \circ p)^*(\cO_X(\chi)) \in \Pic(\MpC(X/G))_{\mathbb{Q}}\]
where $\cO_X \otimes \chi \in \Pic(X/G)_{\mathbb{Q}}$ is the trivial line bundle on $X/G$ twisted by the rational character $\chi$, and the composition
\[ \ev \circ p: \MpC(X/G) \xrightarrow{p \times \id} C\times \MpC(X/G) \xrightarrow{\ev} X/G\]
is induced by the section $p$ and the universal evaluation morphism $\ev: C\times \MpC(X/G) \to X/G$.
\end{defn}

In \Cref{defn: numerical invariant}, we defined a general form of the numerical invariant, but in \Cref{C:boundedness_section_context} it is given by the simple formula $\mu= -\wt(\cL(\chi)) / \sqrt{b}$. In particular, it is valued in $\bR$ rather than $\bR[\epsilon]$.

In the last subsection, we will return to our more general setup and prove a boundedness statement for strata with bounded numerical invariant.

\begin{remark}
All the arguments in this section and \Cref{section: index formulas for affine targets} hold with minor modifications in the more general case where we allow $\cL$ to be a line bundle on $\MpC(X/G)$ such that a positive power $\cL^{\otimes m}$ is numerically equivalent to $\cD(V)$.
\end{remark}

\subsection{Constraints on semistable maps for affine \texorpdfstring{$X$}{X} and split \texorpdfstring{$G$}{G}}

In this subsection, we work in \Cref{C:boundedness_section_context}. In particular, we have $L_{\gen} = \cO_X(\chi)$ for some $\chi\in M_{\bQ}^{W}$. A morphism $B\bG_m \to X/G$, corresponding to a cocharacter $\lambda$ of $G$ and a point of $X^{\lambda=0}$, is destabilizing with respect to this $\cO_X(\chi)$ if $\langle \lambda, \chi \rangle > 0$. We call the semistable locus $X^{\mathrm{ss}}(\chi)/G\subset X/G$ of $\cO_{X}(\chi)$ the $\chi$-semistable locus. Note that if $\phi : C \to X/G$ is a map of degree $d \in H_2(BG)$, then $\phi^\ast(L_{\gen})$ has degree $-\langle d,\chi \rangle$.

For any representation $V$ of $G$, it is possible to choose a norm on cocharacters $b$ that agrees with $\ch_2(V) = (\cdot,\cdot)_V$ on the orthogonal complement of its kernel, i.e. $\phi_V^2 = \phi_V$ (with $\phi_V$ as in \Cref{defn: ch_2(V) and adjoint}). In this case $1-\phi_V$ is the orthogonal projection onto $\ker(\phi_V)$.

\begin{prop} \label{prop: constraints degree semistable locus}
Let $\cL$ be a line bundle on $\MpC(X/G)$ numerically equivalent to $\cD(V)$ for some almost faithful representation $V$, let $b$ be a norm on cocharacters such that $\phi_V^2 = \phi_V$ as above, let $d \in H_2(BG)$ and $\chi \in M_\bQ^W$.
\begin{enumerate}
\item Suppose the GIT stratification of $X/G$ defined by the norm on cocharacters $b$ and the line bundle $\mathcal{O}_{X}(\phi_V(d)^{\dagger})$ is a $\Theta$-stratification. Then the $\cL(\chi)$-semistable locus of $\MpC(X/G)_d$ is empty if either
\[
\left\{ \begin{array}{l} \ker(\|\dash\|_V) \subset \ker(\chi) \text{ and } \|d\|_V > \|\chi\|_V,\text{ or}\\ \ker (\|\dash\|_V) \not \subset \ker(\chi) \text{ and } \|d\|_V >0,\end{array} \right.
\]
where $\|\chi\|_V$ refers to the norm of $\chi$ regarded as linear form on $N_\bR / \ker(\|\dash\|_V)$ equipped with the inner product induced by $(-,-)_V$.\\
\item Suppose the GIT stratification of $X/G$ defined by the norm on cocharacters $b$ and the line bundle $\mathcal{O}_{X}(\chi)$ is a $\Theta$-stratification. There is a constant $m_{\chi}>0$, depending only on $\chi / \|\chi\|_b$, such that when $\|\chi\|_b > \|(1-\phi_V)(d)\|_b/m_\chi + \|d\|_b/m_\chi^2$ the following holds: If $\xi : C \to X/G$ is an $\cL(\chi)$-semistable map of degree $d$, then $\xi$ maps the generic point to the $\chi$-semistable locus of $X/G$, and $\langle d, \chi' \rangle \leq 0$ for any $\chi' \in M_\bQ^W$ in the cone of elements such that $X^{\rm ss}(\chi')=X^{\rm ss}(\chi)$.
\end{enumerate}
\end{prop}

\begin{proof}
We can restrict to a geometric fiber of $S$, and hence assume that we are working over an algebraically closed field. Let $\psi \in M_{\mathbb{Q}}^W$ be a character such that the GIT stratification induced by $b$ and $\mathcal{O}_X(\psi)$ is a $\Theta$-stratification. Choose a map $\xi : C \to X/G$ of degree $d$.

Assume that $\xi$ maps the generic point of $C$ to the $\psi$-unstable locus, and let $\lambda \in X_\ast(T)$ be the maximally $\cO_{X}(\psi)$-destabilizing cocharacter for this point. Then, at the generic point, $\xi$ admits a unique lift along the morphism $X^{\lambda \geq 0} / P_\lambda \to X/G$, because of the assumption that the GIT stratification is a $\Theta$-stratification. Because $X$ is affine, this map is representable by projective schemes, so $\xi$ lifts uniquely to a map $C \to X^{\lambda \geq 0} / P_\lambda$. Let $\xi' : C \to X^{\lambda =0}/L_\lambda$ denote the composition of this unique lift with the projection $X^{\lambda \geq 0} / P_\lambda \to X^{\lambda=0}/L_\lambda$, and let $d' \in X_*(L_{\lambda})_{\mathbb{Q}} = X_\ast(T)_\bQ^{W_\lambda}$ denote the degree of this morphism.

By construction, $\xi'$ maps the generic point of $C$ to the center of the stratum in $X/G$, which is the semistable locus of $X^{\lambda=0}$ for the action of $L_\lambda$ with respect to the shifted character $\psi' := \psi - \frac{\langle \lambda, \psi \rangle}{\|\lambda\|^2_{b}} \lambda^\dagger$. This means that, after a positive scaling of $\psi'$, there is an invariant section of $\cO_{X^{\lambda=0}}(\psi')$ that is nonvanishing at the generic point of $\xi'$, so
\begin{gather*}
\deg(\xi'^\ast(\cO_{X^{\lambda=0}}(\psi'))) = -\langle d', \psi - \frac{\langle \lambda, \psi \rangle}{\|\lambda\|^2_{b}} \lambda^\dagger \rangle \geq 0 \\
\Rightarrow (d,\psi^\dagger)_b\leq (\hat{\lambda},\psi^\dagger)_b\cdot (\hat{\lambda},d')_b,
\end{gather*}
where $\hat{\lambda}:= \lambda / \|\lambda\|_{b}$, and we have used the fact that $(d,\psi^\dagger)_b=(d',\psi^\dagger)_b$ because $\psi^\dagger$ is $W$-invariant and the projection of $d'$ onto $X_\ast(T)^W_\bQ$ is $d$.

Note that $(\hat{\lambda}, \psi^\dagger)_b > 0$ by construction. If $(d,\psi^\dagger)_b \geq 0$, then the estimate $(\hat{\lambda},\psi^\dagger)_b \leq \|\psi^\dagger\|_b$ gives the inequality $(d,\hat{\psi}^\dagger)_b \leq (\hat{\lambda},d')_b$, where $\hat{\psi}^\dagger = \psi^\dagger / \|\psi^\dagger\|_b$. On the other hand, if $(d,\psi^\dagger)_b \leq 0$, then we can let $m_\psi>0$ denote the minimum value of $(\hat{\lambda},\hat{\psi}^\dagger)_b$ among the finite list of maximally destabilizing cocharacters in $X/G$, and we still have $(d,\hat{\psi}^\dagger)_b \leq m_\psi (\hat{\lambda},d')_b$. Note that the $t\psi$-stratification of $X/G$ does not depend on $t>0$, and thus $m_{t\psi} = m_\psi$ for any $t>0$.

If we regard the lift $C \to X^{\lambda \geq 0}/P_\lambda$ as a filtration of the original point $\xi$ in $\MpC(X/G)$, the value of the numerical invariant $-\wt(\cL(\chi)) / \sqrt{b}$ is given by $(\hat{\lambda}, \phi_V(d')+\chi^\dagger)_b$. Using the hypothesis that $\ker(\phi_V) \subset X_\ast(T)_\bQ^W$, and $d'-d$ lies in the orthogonal complement of this space, one has $(1-\phi_V)(d')=(1-\phi_V)(d)$. It follows that $(\hat{\lambda}, \phi_V(d')+\chi^\dagger)_b = (\hat{\lambda}, d' - (1-\phi_V)(d)+\chi^\dagger)_b$. The discussion above then shows that
\begin{equation} \label{E:destabilizing_inequalities}
(\hat{\lambda},\phi_V(d')+\chi^\dagger)_b \geq (\hat{\lambda},\chi^\dagger)_b - \|(1-\phi_V)(d)\|_b + \left\{ \begin{array}{ll} \frac{1}{m_\psi}(d,\hat{\psi}^\dagger)_b, & \text{ if } (d,\psi^\dagger) \leq 0 \\  (d,\hat{\psi}^\dagger)_b, & \text{ if } (d,\psi^\dagger) \geq 0\end{array} \right.
\end{equation}
We will prove claims (1) and (2) by showing that under the hypotheses of the proposition, this number is forced to be $>0$, hence any such $\xi \in \MpC(X/G)$ that maps the generic point of $C$ to the $\psi$-unstable locus is destabilized by this $\lambda$.

\medskip
\noindent{\textit{Proof of (1):}}
\medskip

First, we show that the $\cL(\chi)$-semistable locus is empty if $\|d\|_V>\|\chi\|_b+\|(1 - \phi_V)(d)\|_b$. We let $\psi = \phi_V(d)^\dagger$, which satisfies $(d,\hat{\psi}^\dagger)_b = \|\phi_V(d)\|_b=\|d\|_V$. Then \eqref{E:destabilizing_inequalities} combined with the estimate $|(\hat{\lambda},\chi^\dagger)_b|\leq \|\chi^\dagger\|_b = \|\chi\|_b$ give
\[
(\hat{\lambda},\phi_V(d')+\chi^\dagger)_b \geq \|d\|_V-\|\chi\|_b-\|(1-\phi_V)(d)\|_b,
\]
which is positive by the hypotheses of (1). If follows that any semistable map $\xi : C \to X/G$ must map the generic point of $C$ to the $\psi$-semistable locus of $X/G$. By the same argument as in the case of $X^{\lambda=0}/L_\lambda$ above, this implies that $\deg(\xi^\ast(\cO_X(\psi)))=-\langle d, \psi \rangle = -\|d\|_V^2 \geq 0$. This can not happen, because $\|d\|_V>0$.

To conclude (1), note that neither $\|d\|_V$ nor the $\cL(\chi)$-semistable locus depends on $b$. We are therefore free to rescale $b$ on $\ker(V)$ to be arbitrarily small, in which case $\|\chi\|_b+\|(1 - \phi_V)(d)\|_b$ becomes arbitrarily close to $\|\proj^{\perp}_{\ker(\phi_V)^{\perp}}(\chi)\|_b$.

\medskip
\noindent{\textit{Proof of (2):}}
\medskip

Here we take $\psi = \chi$. In this case, because $\lambda$ is a maximal destabilizer for $\chi$, we have the inequality $(\hat{\lambda},\chi^\dagger)_b\geq m_\chi \|\chi\|_b $. The inequality \eqref{E:destabilizing_inequalities} implies that
\[
(\hat{\lambda},\phi_V(d')+\chi^\dagger)_b \geq m_\chi \|\chi\|_b  - \|(1-\phi_V)(d)\|_b - \frac{1}{m_\chi} \|d\|_b.
\]
It follows that under the hypotheses of (2), this $\lambda$ is destabilizing, so any semistable $\xi$ must map the generic point of $C$ to the $\chi$-semistable locus of $X/G$. The same reasoning as above implies that if such a map exists, then for any $\chi'$ such that $X^{\rm ss}(\chi)=X^{\rm ss}(\chi')$, we have $-\langle d, \chi' \rangle \geq 0$.

\end{proof}


\begin{remark}
The assumption that the GIT stratification is a $\Theta$-stratification (as opposed to a weak $\Theta$-stratification) is automatic if the characteristic of $k$ is $0$ \cite[Cor. 2.1.9]{hl_instability}, or more generally if the stabilizers of $G$ acting on $X$ are smooth (cf. \cite[Thm. 2.3]{ramanan-instability-flag} or by a slight modification of the proof of \cite[Cor. 2.1.9]{hl_instability}).
\end{remark}

\begin{rem}
When $G=T$, the proof of the previous theorem can be sharpened, because $d'=d$ always. A similar analysis shows that if $\psi = \phi_V(d)+\chi^\dagger$, then $\lambda$ is always destabilizing. It follows that any $\cL(\chi)$-semistable map $C \to X/T$ maps the generic point to the $(\phi_V(d)^\dagger+\chi)$-semistable locus of $X/T$, and for a semistable map to exist one must have $\|d\|_V^2 + \langle d,\chi \rangle \leq 0$.
\end{rem}

\subsection{Centers of the unstable strata for affine \texorpdfstring{$X$}{X} and split connected \texorpdfstring{$G$}{G}} \label{subsection: centers of unstable strata}
In this subsection, we will work again in \Cref{C:boundedness_section_context}. We fix once and for all Borel subgroup $B \supset T$ containing the split maximal torus in $G$. We shall use the notation from \Cref{subsection: hn boundedness for bung}.
\begin{lem} \label{lemma: components of grad for affine target}
Under the assumptions of \Cref{C:boundedness_section_context}, we have
\[
\Grad(\MpC(X/G)) \cong \bigsqcup_{\lambda} \MpC(X^{\lambda=0}/L_{\lambda}),
\]
where the disjoint union is over dominant cocharacters $\lambda \in X_\ast(T)$ (i.e. $\langle \lambda, \alpha \rangle \geq 0$ for every positive root $\alpha$), $L_\lambda$ denotes the centralizer of $\lambda$, and $X^{\lambda=0} \subset X$ denotes the closed affine subscheme fixed by $\lambda$.
\end{lem}
\begin{proof}
We have the following chain of canonical equivalences
\begin{align*}
    \Grad(\MpC(X/G)) & = \Map_{\mathcal{Y}}((\mathbb{G}_m)_{\mathcal{Y}}, \Map_{\mathcal{Y}}(C, (X/G) \times \mathcal{Y})) \\
    & = \Map_{\mathcal{Y}}((\mathbb{G}_m)_{\mathcal{Y}} \times_{\mathcal{Y}} C, (X/G) \times \mathcal{Y})\\
    & = \Map_{\mathcal{Y}}(C, \Map_{\mathcal{Y}}((\mathbb{G}_m)_{\mathcal{Y}}, (X/G) \times \mathcal{Y}))\\
    & = \Map_{\mathcal{Y}}(C, \mathcal{Y}\times \Map_S((\mathbb{G}_m)_S, X/G))
\end{align*} 
\cite[Thm. 1.4.8]{hl_instability} applies to yield the following (notice that, even though the theorem in \cite{hl_instability} is stated over a field, the proof outlined by reducing to the case of $\GL_N$ works for any split connected reductive group over a Noetherian scheme $S$)
\[ \Map_S((\mathbb{G}_m)_S, X/G) = \Grad(X/G) \cong \bigsqcup_{\lambda} X^{\lambda =0}/L_{\lambda}\]
We conclude that
\[ \Grad(\MpC(X/G)) \cong \Map_{\mathcal{Y}}\left(C, \bigsqcup_{\lambda} (X^{\lambda =0}/L_{\lambda})\times \mathcal{Y}\right) = \bigsqcup_{\lambda} \MpC(X^{\lambda=0}/L_{\lambda}),\]
where we are using that $C \to \mathcal{Y}$ is proper has geometrically connected fibers for the second equality.
\end{proof}

The ``recognition theorem'' \cite[Sect. 4.2]{hl_instability} implies that a graded point of $\cM$ is the associated graded point of a HN filtration if and only if it is semistable in a graded sense, meaning that the value of the numerical invariant is maximized by the canonical filtration coming from the given grading. We will assume that the value of the numerical invariant is positive on this canonical filtration, as that is the only case that is relevant to the stratification of $\MpC(X/G)$.

We spell out graded semistability more explicitly in our situation. A filtration $f : \Theta_k \to \MpC(X^{\lambda=0}/L_{\lambda})$ provides the following data:
\begin{enumerate}
    \item The point $f(1)$ has a degree $d \in H_2(BL_{\lambda}) \subset N_\bQ^{W_{\lambda}}$, where $W_\lambda$ denotes the Weyl group of $L_\lambda$.
    \item The associated graded $\gr(f)$ has a weight type $\lambda' \in N$, and
    \item $\gr(f)$ also has degree $d' \in H_2(B(L_{\lambda} \cap L_{\lambda'})) \subset N_\bQ^{W_\lambda \cap W_\lambda'}$ corresponding to the degree of the underlying $L_{\lambda} \cap L_{\lambda'}$-bundles (when viewed as an element of $\MpC(X^{\lambda,\lambda'=0}/(L_\lambda \cap L_{\lambda'})) \subset \Grad(\MpC(X^{\lambda=0}/L_{\lambda}))$).
\end{enumerate}
Then we can write the value of the numerical invariant on $\MpC(X^{\lambda = 0}/L_{\lambda})$as
\begin{equation} \label{E:numerical_invariant_simple_formula}
\mu(f) = \frac{-\wt(\cL(\chi))}{\sqrt{b}} = \frac{(\lambda',d')_{V}+\langle \lambda',\chi \rangle}{\|\lambda'\|_b},
\end{equation}
where $(-,-)_V$ denotes the bilinear form on $N$ corresponding to $\ch_2(V)$ (see \Cref{notation: beginning proof hn-boundedness bun_g}). We're assuming that $G$ acts almost faithfully on $V$, i.e., $\ch_2(V)$ is positive definite on the semisimple part of $N_{\bQ}$.
\begin{prop}\label{P:graded_semistability}
A point of $\Grad(\MpC(X/G))$ is graded-semistable if and only if the underlying point of $\MpC(X^{\lambda=0}/L_{\lambda})_d$ is semistable with respect to the line bundle $\cL(\chi_{\lambda, d}')\in \Pic(\MpC(X^{\lambda=0}/L_{\lambda}))_{\mathbb{Q}}$, where we view $V$ as a representation of $L_\lambda$, and we use the shifted character
\[
\chi'_{\lambda,d} := \chi - \mu(f_{\rm canon}) \left(\frac{1}{\|\lambda\|_b} \lambda, -\right)_b \in M_{\bQ}^{W_{\lambda}},
\]
where $\mu(f_{\rm canon}) = \frac{(\lambda,d)_V+\langle \lambda,\chi\rangle}{\|\lambda\|_b}$ and $d \in H_2(BG) \subset N_{\bQ}^W$ is the degree of the $G$-bundle underlying the point.
\end{prop}

\begin{proof}
Since $\lambda$ acts centrally on $\MpC(X^{\lambda=0}/L_{\lambda})$, for any given filtration $f$ one gets a family of filtrations $f_t$ which is the same parabolic reduction but with weight vector $\lambda'+t\lambda$. In particular, $\gr(f_t)$ has the same degree $d'$ for all $t$. As $t \to \infty$, this filtration converges to the canonical filtration in $\iDeg(\MpC(X^{\lambda=0}/L_{\lambda}),f(1))$. Graded semistability is equivalent to the assertion that $\mu(f) = \mu(f_0) \leq \lim_{t \to \infty} \mu(f_t)$ for any filtration $f$ in $\MpC(X^{\lambda=0}/L_{\lambda})$. Let us compute the function explicitly:
\[
\mu(f_t) = \frac{(\lambda+\frac{\lambda'}{t},d')_{V}+\langle \lambda+\frac{\lambda'}{t},\chi \rangle}{\|\lambda+\frac{\lambda'}{t}\|_b}.
\]
The function $\mu(f_t)$ is smooth in $t$ and strictly quasi-concave along the interval in $t$ on which $\mu(f_t)>0$. Therefore, a point is graded semistable if and only if for any filtration $f$,
\[
\frac{d \mu(f_t)}{dt} \geq 0 \text{ for all } t \gg 0.
\]
We can compute this derivative explicitly as
\[
\frac{d \mu(f_t)}{dt} = \frac{\left[ (\lambda,d')_V+\langle \lambda,\chi \rangle \right] \cdot \left[(\lambda,\lambda')_b+\frac{1}{t} \|\lambda'\|_b^2 \right]- \left[(\lambda',d')_V+\langle\lambda',\chi\rangle\right] \cdot (\lambda+\frac{\lambda'}{t},\lambda)_b}{\|\lambda+\frac{\lambda'}{t}\|_b^3}
\]
The key observation is that $d'$ is not arbitrary -- any character of $L_\lambda$ must take the same value on $d'$ and $d$. In particular, because $\lambda$ is central in $L_\lambda$ and $(-,-)_V$ is $W$-equivariant, $(\lambda,d')_V = (\lambda,d)_V$. A direct translation of the condition $d\mu(f_t)/dt \geq 0$ for $t \gg 0$ states that
\begin{equation} \label{E:graded_semistable}
(\lambda',d')_V+\langle\lambda',\chi\rangle - \overbrace{\frac{(\lambda,d)_V+\langle \lambda,\chi\rangle}{\|\lambda\|_b}}^{\mu(f_{\rm canon})} \left(\frac{1}{\|\lambda\|_b} \lambda,\lambda' \right)_b \leq 0,
\end{equation}
and when equality holds one has
\[
\left((\lambda',d')_V+\langle\lambda',\chi\rangle \right) (\lambda',\lambda)_b - \left((\lambda,d)_V+\langle \lambda,\chi \rangle\right) \|\lambda'\|^2_b \leq 0.
\]
However, if \eqref{E:graded_semistable} is an equality, then the second inequality is equivalent to the Cauchy-Schwarz inequality $(\lambda',\lambda)_b^2 \leq \|\lambda\|_b^2 \cdot \|\lambda'\|_b^2$ and thus holds automatically. This shows that a given point in $\MpC(X^{\lambda=0}/L_{\lambda})$ is graded-semistable if and only if for any filtration with weight type $\lambda'$ and degree $d'$, the inequality \eqref{E:graded_semistable} holds. This in turn is just semistability with respect to $\cL(\chi'_{\lambda,d})$.
\end{proof}

\subsubsection{Relationship between the degree and the weight type}

For the classical HN stratification of $\Bun_{\GL_n}(C)$, the canonical weights of the Harder-Narasimhan filtration are proportional to the slopes of the associated graded bundles. We generalize this relationship to $\MpC(X/G)$.

\begin{notn} \label{notn: pseudoinverse}
Let $V$ be an almost faithful representation, and recall the endomorphism $\phi_V: N_{\mathbb{Q}} \to N_{\mathbb{Q}}$ defined in \Cref{defn: ch_2(V) and adjoint}. As $\phi_V$ is self-adjoint with respect to $b$, we have an orthogonal decomposition $N_\bQ = \ker(\phi_V) \oplus \im(\phi_V)$, and $\phi_V$ restricts to an automorphism of $\im(V)$. The \emph{Moore-Penrose pseudoinverse} $\phi_V^+$ is the endomorphism of $N_\bQ$ that acts by $0$ on $\ker(\phi_V)$ and by the inverse of $\phi_V$ on $\im(\phi_V)$.

We note that $\phi^+_V \phi_V$ is the orthogonal projection $\proj^{\perp}_{(\ker(\phi_V)^{\perp})}$ onto $\ker(\phi_V)^{\perp}$.
\end{notn}

\begin{lem} \label{L:relationship_topological_weight}
Consider a connected component of $\Grad(\MpC(X/G))$, corresponding to a choice of weight cocharacter $\lambda \in N$ and a degree $d \in N_\bQ^{W_\lambda}$, and use the numerical invariant \eqref{E:numerical_invariant_simple_formula}. Assume that $\mu >0$, where $\mu := \left((\lambda,d)_V+\langle\lambda,\chi\rangle \right) / \left\| \lambda \right\|_b > 0$, and let
\[
Y_{\lambda} = \left\{ \left. w \in  N_\bQ^{W_\lambda} \right| \langle w, \beta \rangle = 0 \text{ for all weights } \beta \text{ of } X \text{ s.t. } \langle \lambda, \beta \rangle = 0 \right\}.
\]
If the graded semistable locus is non-empty, then the following hold
\begin{enumerate}
\item $\lambda = \gamma \cdot \proj^{\perp}_{Y_{\lambda}}(\phi_{V}(d) + \chi^{\dagger})$, for some $\gamma>0$.

\item $\mu = \|\proj^{\perp}_{Y_\lambda}(\phi_V(d)+\chi^\dagger)\|_b$.
 
\item Suppose that the GIT stratification of $X/G$ defined by the norm on cocharacters $b$ and the line bundle $\mathcal{O}_{X}(\phi_V(d)^{\dagger})$ is a $\Theta$-stratification. Then we have either $\chi^\dagger - \mu \hat{\lambda} \notin \im(\phi_V)$ and $\|d\|_V=0$, or $\chi^\dagger - \mu \hat{\lambda} \in \im(\phi_V)$ and
    \[
    \|d\|_V \leq \| (\phi_V^+)^{1/2} (\chi^\dagger - \mu \hat{\lambda}) \|_b,
    \]
    where $\phi^+_V$ is the pseudoinverse defined in \Cref{notn: pseudoinverse}. In particular, $\|d\|_V \leq (\|\chi\|_b+\mu)/\sqrt{a}$, where $a$ is the smallest non-zero eigenvalue of $\phi_V$.
\end{enumerate}

\end{lem}
\begin{proof}
By \Cref{P:graded_semistability}, graded semistability is equivalent to semistability on $\MpC(X^{\lambda=0}/L_{\lambda})_d$ with respect to the line bundle $\cL(\chi')$, where $\chi' := \chi - \mu (\hat{\lambda},\dash)_b$. $Y_\lambda$ parameterizes rational cocharacters in the maximal central torus $Z(L_\lambda)$ that act trivially on $X^{\lambda=0}$, and these cocharacters act centrally on any point of $\MpC(X^{\lambda=0}/L_{\lambda})$. The weight $-\wt(\cL(\chi'))$ with respect to this generic stabilizer group is given by the linear functional on $Y_\lambda$
\begin{equation} \label{E:relationship_topological_weight}
(-,d)_V + \langle -,\chi \rangle - \mu (-,\frac{1}{\|\lambda\|_b} \lambda)_b = (-,\phi_V(d)+\chi^\dagger)_b-\mu(-,\hat{\lambda})_b.
\end{equation}
So if this weight is non-zero for any cocharacter of $Y_\lambda$, this cocharacter will destabilize every point of $\MpC(X^{\lambda=0}/L_{\lambda})_d$. The claim (1) follows from setting this linear functional $=0$ on $Y_{\lambda}$ and observing that $\lambda \in Y_\lambda$. The claim (2) follows from the claim (1) after evaluating the functional \eqref{E:relationship_topological_weight} on $\hat{\lambda}$ and setting it equal to $0$.

The condition (3) is just a translation of the bound in \Cref{prop: constraints degree semistable locus}(1). Using the orthogonal decomposition $N_\bR = \ker(\phi_V) \oplus \im(\phi_V)$, the condition $\ker(\|\dash\|_V) \subset \ker(\chi')$ is equivalent to the condition that $(\chi')^\dagger = \chi^\dagger - \mu \hat{\lambda} \in \im(\phi_V)$. Assuming this condition holds, we observe that $\|\chi'\|_V$ agrees with the norm of the restriction of the character $\chi'$ to $\im(\phi_V)$, with respect to the inner product induced by $\|\dash\|_V$ on $\im(\phi_V)$, and this is equivalent to $\|(\phi_V^+)^{1/2}(\chi^\dagger-\mu \hat{\lambda})\|_b$. Finally, the Cauchy-Schwarz and triangle inequalities imply that $\|(\phi_V^+)^{1/2}(\chi^\dagger-\mu \hat{\lambda})\|_b \leq \|\phi_V^+\|^{1/2}_b (\|\chi\|_b + \mu)$, and $\|\phi_V^+\|_b$ is the reciprocal of the smallest non-zero eigenvalue of $\phi_V$.

\end{proof}

\subsubsection{Explicit indexing for the $\Theta$-stratification}
We have seen that a $\Theta$-stratum in $\MpC(X/G)$ with respect to the numerical invariant $-\wt(\mathcal{L}(\chi)) / \sqrt{b}$ corresponds to a pair $(d, \lambda) \in N_{\mathbb{Q}}^{\oplus 2}$, where $\lambda$ is a cocharacter determining a parabolic subgroup $P_{\lambda}$, and $d$ corresponds to the degree of the associated graded viewed as an element of $\MpC(X^{\lambda=0}/L_{\lambda})$. Not all pairs $(d, \lambda)$ arise this way; there is a relationship between $d$ and $\lambda$, as demonstrated in \Cref{L:relationship_topological_weight}. We shall describe the relation concretely in terms of a fan on $N_{\mathbb{Q}}$.

For the following definition, we use $ \overline{\sigma}_{P} = X_{*}(Z)_{\bR} + \sum_{ j \in I_{P}} \mathbb{R}^{\geq 0}\omega^{\vee}_j$ for any parabolic subgroup $B \subset P \subset G$ as in \Cref{notation: beginning proof hn-boundedness bun_g}.
\begin{defn}
Given the linear representation $X$ of $G$, we denote by $\Sigma_X$ the fan of all cones $\sigma\subset N_{\mathbb{Q}}$ of the form 
\[\sigma = \left\{ w \, | \, \beta (w) \geq 0 \, \text{ for all } \beta \in S_+ \right\} \cap \left\{ w \leq 0 \, | \, \beta(w) \leq 0\, \text{ for all } \beta \in S_- \right\} \cap \overline{\sigma}_P,\]
where $S_+,S_- \subset M$ are (possibly empty) subsets such that $S_+ \cup S_-$ is the set of $T$-weights of $X$ and $P\subset G$ is a parabolic subgroup containing $B$.
\end{defn}

It follows from definition that the fan $\Sigma_X$ has finitely many cones.

\begin{defn} \label{defn: chi indexing datum}
Let $(d,\lambda)$ be a pair in $N_{\mathbb{Q}}^{\oplus 2}$. We say that $(d, \lambda)$ is a $\chi$-active indexing datum if the following are satisfied.
\begin{enumerate}
    \item $\langle d, \psi\rangle \in \mathbb{Z}$ for all $\psi: P_{\lambda} \to \mathbb{G}_m$ (viewed as a character in $M$ by restricting to $T$).
    \item $\lambda \in \bar{\sigma}_B$, and we have
    \[ \lambda = \proj^{\perp}_{\Span(\sigma)}(\phi_V(d) + \chi^{\dagger})\]
    where $\sigma$ is the minimal cone in $\Sigma_X$ containing $\lambda$.
    \item If $\chi^\dagger - \lambda \notin \im(\phi_V)$, then $\|d\|_V = 0$, and if $\chi^\dagger-\lambda \in \im(\phi_V)$, then $\|d\|_V \leq \|(\phi_V^+)^{1/2}(\chi^\dagger-\lambda)\|_b$.
\end{enumerate}
\end{defn}

\begin{cor} \label{cor: theta strata indexing}
Suppose that $(d, \lambda) \in N_{\mathbb{Q}}^{\oplus 2}$ corresponds to a nonempty $\Theta$-stratum for the numerical invariant $-\wt(\mathcal{L}(\chi)) / \sqrt{b}$, and assume that the GIT stratification induced by $\mathcal{}O_X(\phi_V(d)^{\dagger})$ on $X/G$ is a $\Theta$-stratification. Then $\nu = (d, \lambda)$ is a $\chi$-active indexing datum.
\end{cor}
\begin{proof}
This is a direct consequence of \Cref{L:relationship_topological_weight}.
\end{proof}

\begin{lem} \label{lemma: coxeter lattice for indexing data}
If $(d, \lambda)$ is a $\chi$-active indexing datum, then $d \in \frac{1}{h(W_{P_{\lambda}})}\cdot N^{W_{P_\lambda}}$, where $h(W_{P_{\lambda}})$ is the order of a Coxeter element. In particular, if $H = {\rm lcm}\{h(\Span(\sigma_{E_P})) | P \subset G\}$, then for all $\chi$-active indexing data $\nu = (d, \lambda)$ the degree $d$ lies in the lattice $\frac{1}{H} N \subset N_\bQ$.
\end{lem}
\begin{proof}
If $\varphi \in W_{P_{\lambda}}$ is a Coxeter element, then $M^{W_{P_{\lambda}}} = M^\varphi$ \cite[30-1 Prop.A (i), pg. 311]{kane-reflection-groups}, so if $h$ is the order of $\varphi$, then $\sum_{i=1}^h \varphi^i m \in M^{W_{P_{\lambda}}}$ for any $m \in M$. Therefore, if $d$ satisfies (1) in \Cref{defn: chi indexing datum}, then $\langle d, \sum_{i=1}^{h} \varphi^i m \rangle = \langle \sum_{i=1}^{h} \varphi^i d, m \rangle = h \langle d, m \rangle \in \bZ$ for any $m \in M$.
\end{proof}

For any $\chi$-active indexing datum $\nu = (d, \lambda)$, we set $\mu(\nu):= \|\lambda\|_b$. This corresponds to the value of numerical invariant for the corresponding $\Theta$-stratum. We call $\proj^{\perp}_{\ker(\phi)}(d)$ the kernel part of the indexing datum $\nu$.
\begin{lem} \label{lemma: finiteness of indexing data with bounded numerical invariatn}
Fix a number $\gamma>0$ and an element $d_{\ker} \in \ker(\phi_V)$. Then the set of $\chi$-active indexing data $\nu= (d, \lambda)$ with $\mu(\nu)\leq \gamma$ and kernel part $d_{\ker}$ is finite.
\end{lem}
\begin{proof}
For a fixed value of $d$, the possible choices for $\lambda$ depend on the finitely many cones in the fan $\Sigma_X$. Therefore, it suffices to show that there are finitely many such $d$ fitting into a $\chi$-active pair $\nu = (d, \lambda)$ with central part $d_{\ker}$ and numerical invariant $\mu(\nu)\leq \gamma$. By Lemma \ref{lemma: coxeter lattice for indexing data}, such $d$ lie inside the lattice $\frac{1}{H} M \subset M_{\mathbb{Q}}$. Therefore, it suffices to show that $d$ lies in a bounded subset of $N_\bR$. Because we are fixing the projection of $d$ to $\ker(\phi_V)$, and $\|\dash\|_V$ is positive definite on $\ker(\phi_V)^\perp$, it suffices to bound $\|d\|_V$ above. Such bound follows from condition (3) of \Cref{defn: chi indexing datum} and the Cauchy-Schwarz and triangle inequalities.

\end{proof}
\begin{remark} \label{remark: finitely many strata for bounded numerical invariant}
If the representation $V$ is infinitesimally faithful (i.e. the kernel of $G \to \GL(V)$ is finite), then the only possible kernel part is $0$. Therefore, \Cref{remark: finitely many strata for bounded numerical invariant} implies that for infinitesimally faithful $V$ there are finitely many strata in $\MpC(X/G)$ strata with numerical invariant $\leq \gamma$.
\end{remark}

\begin{remark}
One can alternatively index the strata by pairs $(d,\sigma)$, where $d \in N_\bQ$ and $\sigma \in \Sigma_X$. In this case one defines $\lambda = \proj^\perp_{\Span(\sigma)}(\phi_V(d)+\chi^\dagger)$, and one says that $(d,\sigma)$ is $\chi$-active if $\lambda$ lies in the relative interior of $\sigma$, and the conditions (1) and (3) of \Cref{defn: chi indexing datum} hold for the pair $(d,\lambda)$. The alternative indexing emphasizes that there are finitely many indices for each value of $d$.
\end{remark}

\subsection{Boundedness properties of the \texorpdfstring{$\Theta$}{Theta}-stratification and moduli spaces}
In this subsection we return to our more general context as in \Cref{subsection: gauged maps introduction}. We shall investigate the boundedness properties of the $\Theta$-stratification defined in Theorem \ref{thm: theta stratification in general}. For simplicity, we restrict to the case when $S$ is a $\mathbb{Q}$-scheme (i.e. we work in characteristic $0$). In particular, all weak $\Theta$-stratifications on stacks over $S$ will automatically be $\Theta$-stratifications.

\begin{prop} \label{prop: boundedness of the semistable locus affine case}
Suppose that $S$ is a $\bQ$-scheme, $X \to S$ is affine, and $G=G_0$.
Consider the numerical invariant $\mu = \left(-\wt(\cD(V)) + \ell_{gen}\right)/\sqrt{b}$. Then for every connected component $\cX \subset \Bun_{G}(C)$, the intersection $\Forget^{-1}(\cX) \cap \MpC( X/G)^{\mu \dash ss}$ is quasi-compact, where $\MpC( X/G)^{\mu \dash ss}$ denotes the open substack of $\mu$-semistable points.
\end{prop}
\begin{proof}
Boundedness was proven by Schmitt \cite[Prop. 2.8.3.1]{schmitt-decorated-principal-bundles} in \Cref{C:boundedness_section_context} under the assumption that $\mathcal{Y} = S = \Spec(k)$ for some algebraically closed field $k$ of characteristic $0$. We will reduce to this case.

The morphism $\Forget: \MpC(X/G) \to \Bun_{G}(C)$ is affine and of finite type. Therefore, it suffices to show that the image $\Forget(\MpC( X/G)^{\mu \dash ss}) \cap \cX$ is bounded. By \Cref{lemma: boundedness of HN stratification relative case}, it is enough to prove that there exists some $\gamma'$ such that $\Forget(\MpC( X/G)^{\mu \dash ss}) \cap \cX$ is contained in $\Bun_{G}(C)^{\mu \leq \gamma'}$.

After replacing $X$ with the total space of the $\bG_m$-bundle corresponding to $L_{\gen}$, replacing $G$ with $G \times\bG_m$, and replacing $b$ with $b + \eta^2$, where $\eta$ is the tautological character of $\bG_m$, we can assume that $L_{\gen} = \cO_X(\chi)$ for some character $\chi$ of $G$.\endnote{Indeed, if $\pi : E \to X$ is the principal $\bG_m$-bundle corresponding to $L_{\gen}$, then $\pi^{\ast} L_{\gen} \cong \cO_E(\eta)$, where $\eta$ is the tautological character of $\bG_m$. Because $\bG_m$ acts freely on $E$, the norm on graded points $b$ is still a norm on graded points of $E/{G\times \bG_m}$. This does not strictly fit into our framework, though, because $b$ is not a norm on graded points of $BG \times B\bG_m$. $b+ \eta^2$ is, however, a norm on graded points of $BG \times B\bG_m$, and it induces the same norm on graded points of $E/(G\times \bG_m)$ as $b$, so it does not change the numerical invariant.} After passing to Noetherian smooth covers of $\mathcal{Y}$ and $S$, we can assume that $\mathcal{Y}=S$ is connected and affine, that $C \to S$ admits a section, and that $G$ is split over $S$. This means that $G$ is the base change of a split reductive group $G_{\mathbb{Q}}$ over $\mathbb{Q}$. Finally, one can choose an equivariant closed embedding $X \hookrightarrow S \times_{\Spec(\bQ)} M$ for some finite dimensional linear representation $M$ of $G_{\mathbb{Q}}$.\endnote{Regarding $\cO_X$ as a $G_\bQ$ representation, one can write $\cO_X$ as a filtered union $\bigcup_\alpha M_\alpha$ of finite dimensional $G_\bQ$-representations. Because $\cO_X$ is globally generated and finite type as an $\cO_S$-algebra, the canonical map $\cO_S \otimes_\bQ \Sym(M_\alpha) \to \cO_X$ is surjective for some $\alpha$. This induces the desired closed immersion.} This induces a closed immersion $\MpC(X/G) \hookrightarrow \MpC(S \times M/G)$. Since we have $\MpC(X/G)^{\mu \dash ss} = \MpC( S \times M/G)^{\mu \dash ss} \cap \MpC( X/G)$, it suffices to prove the claim with $S \times M$ in place of $X$.


We can check $\Forget(\MpC(S \times M/G)^{\mu \dash ss}) \cap \cX \subset \Bun_{G}(C)^{\mu \leq \gamma'}$ on geometric fibers of $\cY$, as long as the bound $\gamma'$ is the same on every fiber. This follows from Schmitt's boundedness result \cite[Prop. 2.8.3.1]{schmitt-decorated-principal-bundles}. In the proof, the bound constructed only depends on the group $G_{\bQ}$, the representation $M$, the degree of $G$-bundles on the component $\cX$, the character $\chi$ and the genus $g$ of the curve $C$.

\end{proof}

\begin{cor} \label{cor: boundedness of strate affine case}
Suppose that $S$ is a $\bQ$-scheme, and $X \to S$ is affine. Consider $\mu = -\wt(\cD(V)) / \sqrt{b}$. Then for every connected component $\cX \subset \Bun_{G}(C)$ and every $\gamma \in \bR$, the intersection $\Forget^{-1}(\cX) \cap \MpC( X/G)^{\mu \leq \gamma}$ is quasi-compact.
\end{cor}
\begin{proof}
By a similar reasoning as in the proof of \Cref{prop: reduction to connected groups}, we can reduce to proving the claim in \Cref{C:boundedness_section_context}, and in this case $\chi$ is the trivial character. Then $\cX \subset \Bun_{G}(C)_{d}$ for some $d \in H_2(BG)$, and we will complete the proof by showing that $\MpC(X/G)_d^{\mu \leq \gamma} := \MpC(X/G)^{\mu \leq \gamma} \cap \Forget^{-1}(\Bun_G(C)_d)$ is bounded.

By \Cref{cor: theta strata indexing}, $\MpC(X/G)_d^{\mu \leq \gamma}$ is a disjoint union of locally closed $\Theta$-strata $\cS_\nu$ indexed by the set of $\chi$-active indexing data $\nu = (d', \lambda)$ such that $\mu(\nu) \leq \gamma$ and $d'$ has central part $d \in N_\bQ^W$. Because $d_{\ker}:= \proj^{\perp}_{\ker(\phi_V)}(d') = \proj^{\perp}_{\ker(\phi_V)}(d)$ is the same for every $\nu$, \Cref{lemma: finiteness of indexing data with bounded numerical invariatn} implies that there are finitely many such $\nu$. For each $\nu$, the morphism from the stratum to its center $\gr: \cS_\nu \to \cZ_\nu^{\rm ss} \subset \Grad(\MpC(X/G))$ is quasi-compact \cite[Lem.~1.3.8]{hl_instability}, so it suffices to show that each of the finite collection of stacks $\cZ_\nu^{\rm ss}$ is bounded. Finally, \Cref{P:graded_semistability} identifies the center $\cZ_\nu^{\rm ss}$ with the semistable locus in $\MpC( X^{\lambda = 0}/L_{\lambda})$ for the numerical invariant $\mu'$ formed using a character $\chi'_{\lambda, d}\in M^{W_{\lambda}}_{\mathbb{Q}}$, so \Cref{prop: boundedness of the semistable locus affine case} implies that $\cZ_\nu^{\rm{ss}}$ is bounded.

\end{proof}
\begin{prop} \label{prop: boundedness of strata general case}
Suppose that $S$ is a $\bQ$-scheme. For every $\gamma \in \mathbb{R}$, every connected component $\cX \subset \Bun_{G}(C)$, and every $d \in H_2(X/G)$, the open substack $\cM_{n}^{G}(X)_d^{\mu_{\delta}\leq \gamma} \cap \Forget^{-1}(\cX)$ is quasi-compact. In particular, there are finitely many $\Theta$-strata inside $\Forget^{-1}(\cX)\cap \cM_{n}^{G}(X)_d$ having bounded numerical invariant, and each $\Theta$-stratum is quasi-compact.
\end{prop}
\begin{proof}
By a similar argument as in the proof of \Cref{prop: reduction to connected groups}, we can assume that we are in \Cref{context:boundedness_simplification}. The forgetful morphism factors as  $\Forget: \cM_{n}^{G}(X)_d \xrightarrow{p} \MpC( A_X/G) \xrightarrow{F} \Bun_{G}(C)$. By the proof of \Cref{prop: formal line bundle is ample on p fibers} the morphism $p$ is quasi-compact, and hence it suffices to show that $p(\cM_{n}^{G}(X)^{\mu_{\delta}\leq \gamma}_d) \cap F^{-1}(\cX)$ is quasi-compact. By \Cref{cor: boundedness of strate affine case}, we are reduced to showing that there exists some $\gamma'$ such that $p(\cM_{n}^{G}(X)^{\mu_{\delta}\leq \gamma}) \subset \MpC( A_X/G)^{\mu \leq \gamma'}$, where we equip $\MpC(A_X/G)$ with the numerical invariant $\mu = -\wt(\cD(V))/\sqrt{b}$.


Choose $x \in \cM_{n}^{G}(X)^{\mu_{\delta}\leq \gamma}(k)$ for some algebraically closed field $k$. By \Cref{T:degeneration_fan_DM_morphism} the proper morphism $p$ induces a homeomorphism $|\DF(\cM_{n}^{G}(X),x)| \xrightarrow{\sim} |\DF(\MpC(A_X/G), p(x))|$. Under this identification, the constant term $\mu_{\delta}^0$ is related to the numerical invariant $\mu$ on $\MpC(A_X/G)$ by
\begin{equation} \label{equation: numerical invariant for gauged vs affine2}
\mu_{\delta}^0 = \mu + \frac{ -\delta_{\mrk} \wt(\cL_{\mrk}) +\delta_{\gen} \ell_{\gen}}{\sqrt{b}} 
\end{equation}
The affine morphism $\MpC(A_X/G) \to \Bun_{G}(C)$ naturally identifies $|\DF(\MpC(A_X/G), p(x))|$ with a subset of $|\DF(\Bun_G(C),E)|$, where $E$ is the underlying $G$-bundle of $p(x)$. The $\mu$-maximizing HN filtration of $p(x)$ yields a point $w$ in some cone $\sigma_{E_P} \subset |\DF(\Bun_G(C),E)|$. We need to show that $\mu(w) \leq \gamma'$, for some $\gamma'$ that does not depend on $x$.

By \Cref{lemma:mumford_weight}, applied to the $G$-projective $S$-scheme $R \times Q$ of \Cref{lemma: embedding projective over affine},\endnote{More precisely, we have a locally closed immersion $X \hookrightarrow R \times Q$, which induces a locally closed immersion $\cM_{n}^{G}(X) \hookrightarrow \cM_{n}^G(R \times Q)$ by \ref{lemma: behavior of gauged maps under immersions}. It follows that $|\DF(\cM_{n}^{G}(X),x)| \to |\DF(\cM_{n}^{G}(R \times Q),x)|$ is a closed embedding \cite[Prop.~3.2.12]{hl_instability}. The functionals $\ell_\mrk$ and $\ell_\gen$ do not necessarily extend to $\DF(\cM_{n}^{G}(R \times Q),x)$, but for some integer $q>0$, $q \ell_{\mrk}$ and $q \ell_\gen$ extend by \Cref{lemma: embedding projective over affine}. \Cref{lemma:mumford_weight} then describes the extended functionals $q \ell_\mrk$ and $q \ell_\gen$.} there is a finite subset $\Gamma \subset X^\ast(T)_{\mathbb{R}}$, which does not depend on $x$, and two subsets $\Sigma_{x,\gen}, \Sigma_{x,\mrk} \subset \Gamma$ such that \eqref{equation: numerical invariant for gauged vs affine2} can be rewritten as
\[ \mu^0_{\delta}(w) = \mu(w) + \frac{ \delta_{\mrk} \cdot \min_{\ell \in \Sigma_{p,\mrk}} \ell(w) +\delta_{\gen} \cdot \min_{\ell \in \Sigma_{p,\gen}} \ell(w)}{\|w\|_b}
\]
Defining $h:=\max_{\ell \in \Gamma} \|\ell\|_b$, we have $\min_{\ell \in \Sigma_{x,\gen}} \ell(w) \geq - h \cdot \|w\|_b$, and similarly for $\Sigma_{x,\mrk}$. Therefore, we can bound the expression above by
\begin{equation} \label{equation: last inequality boundedness of strata general2}
\mu_{\delta}^0(w) \geq \mu(w) + \frac{ -\delta_{\mrk} \cdot h \|w\|_b -\delta_{\gen} \cdot h \|w\|_b}{\|w\|_b} \geq \mu(w) -  |\delta|h
    \end{equation}
The assumption $x \in \cM_{n}^{G}(X)^{\mu_{\delta}\leq \gamma}(k)$ implies that $\mu_{\delta}^0(w) \leq \gamma$. Therefore, \eqref{equation: last inequality boundedness of strata general2} implies that $\mu(w) \leq \gamma' := \gamma + |\delta|h$.
\end{proof}

\begin{thm} \label{thm: moduli space in the case of centrally contractible}
Suppose that $S$ is a $\bQ$-scheme. The semistable locus $\cM_{n}^{G}(X)^{\mu_{\delta}\dash \mathrm{ss}}$ admits a relative good moduli space that is separated and locally of finite type over $\cY$. Each connected component of the moduli space is proper over $(A_{X} /\!/G) \times_{S} \cY$.

Moreover, under the additional assumptions of \Cref{C:boundedness_section_context}, the connected components of the center $\cZ_{\nu}^{{\rm ss}}$ of each stratum of the $\Theta$-stratification defined by $\mu$ also admit good moduli spaces that are proper over $(A_{X} /\!/G) \times_{S} \cY \cong X /\!/G$.
\end{thm}
\begin{proof}
The first part of the statement is a consequence of \Cref{thm: theta stability paper theorem}. Indeed, the numerical invariant $\mu_{\delta}$ is strictly $\Theta$-monotone and strictly $S$-monotone by \Cref{thm: monotonicity of gauged maps}, and it satisfies HN-boundedness by Proposition \ref{prop:HN-boundedness}. Moreover, an application of \Cref{prop: boundedness of strata general case} with $\gamma = 0$ shows that $\cM_{n}^{G}(X)^{\mu_{\delta}\dash \mathrm{ss}}$ is a disjoint union of bounded open and closed substacks. By \Cref{prop: valuative criterion over affine GIT quotient}, the existence part of the valuative criterion is satisfied for $\cM_{n}^{G}(X)\to (A_{X} /\!/G) \times_{S} \cY$.

For the second part, we recall that the connected components of the centers of strata are indexed by $\chi$-active indexing data $\nu = (d', \lambda)$ as in \Cref{defn: chi indexing datum}. The component corresponding to such $\nu$ is of the form $\MpC(X^{\lambda =0}/L_{\lambda})_{d'}^{\mu' \dash {\rm ss}}$, where $\mu'$ is the shifted numerical invariant induced by the line bundle $\cL(\chi_{\lambda, d'}')$ as described in \Cref{P:graded_semistability}. In particular, it admits a good moduli space that is proper over $X^{\lambda =0}/\!/L_{\lambda}$ by the first part of the theorem. It therefore suffices to show that the canonical morphism $X^{\lambda=0}/\!/L_{\lambda} \to X /\!/G$ is finite. This morphism is affine, so it suffices to show that the pushforward of $\cO_{X^{\lambda=0}}$ along $X^{\lambda=0}/L_\lambda \to X/\!/G$ is coherent. To see this, observe that the morphism $X^{\lambda\geq 0}/P_{\lambda} \to X/\!/G$ factors through the morphism $X^{\lambda \geq 0} / P_\lambda \to X^{\lambda=0}/(L_\lambda/\im(\lambda))$ and the pushforward of $\cO_{X^{\lambda\geq 0}}$ along the latter is $\cO_{X^{\lambda=0}}$. It therefore suffices to show that the pushforward of $\cO_{X^{\lambda \geq 0}}$ along $X^{\lambda \geq 0}/P_\lambda \to X/\!/G$ is coherent, and this follows from the fact that $X^{\lambda \geq 0}/P_\lambda \to X/G$ is proper and $X/G \to X/\!/G$ is a good moduli space.
\end{proof}

\begin{remark}
We note that the results about the existence of the moduli space and the properness of the generalized Hitchin morphism were proven by Schmitt under the assumptions of \Cref{context:boundedness_simplification} when $S = \Spec(k)$ for an algebraically closed field $k$ \cite{schmitt-decorated-principal-bundles}. Schmitt's GIT construction also yields the projectivity of the moduli space over the Hitchin base. Our results yield an alternative construction of the moduli space as a relative algebraic space that avoids GIT and works in the generality of an arbitrary bases $\cY$ and $S$. One cannot a priori deduce the projectivity of the generalized Hitchin morphism from our methods. 

However, we believe that it is possible to use Schmitt's results to conclude the projectivity. Under the assumptions of \Cref{context:boundedness_simplification}, a power of the line bundle $\cL(\chi)$ descends to the good moduli space. In order to show that this line bundle is relatively ample over the generalized Hitchin base, it would suffice to check ampleness fiberwise over $S$. We believe that proving ampleness of this specific line bundle over $S= \Spec(k)$ should be possible by tracing Schmitt's construction of the moduli space.
\end{remark}

\section{Index formulas for affine targets} \label{section: index formulas for affine targets}

\begin{context} \label{C:index_section_context}
For most of this section, we will work under the assumptions and notation in \Cref{C:boundedness_section_context}. We will furthermore assume the following.
\begin{itemize}
    \item $S = \Spec(k)$ for a field $k$ of characteristic $0$.
    \item There is a line bundle $\sqrt{K}$ on $C$ such that its corresponding rational $K$-theory class is $\sqrt{K_C}$ as in \Cref{subsection: notation}.
    \item Let $\mathbf{T}$ denote the virtual representation $X \oplus \mathfrak{g}[1] \in K_0(BG)$. For some arguments we will need to arrange for the following condition to be satisfied:
\begin{equation} \label{E:positivity_condition}
c_{X/V}:= \|\phi^+_V\|_b \max_{ \lambda \in N} \left( \lVert \phi_{\mathbf{T}^{\lambda<0}} \rVert_b \right) < 1
\end{equation}
where $\lVert - \rVert_b$ denotes the operator norm, $\phi^+_V$ is the pseudoinverse defined in \Cref{notn: pseudoinverse}, and $\mathbf{T}^{\lambda<0} \in K_0(BT) \cong \bZ[M]$ is the summand whose weights pair negatively with $\lambda$.\endnote{\eqref{E:positivity_condition} is a $\max$ instead of a $\sup$ because the class $\mathbf{T}^{\lambda<0}$ takes on only finitely many values as $\lambda$ varies, and thus there are finitely many such endomorphisms $\phi_{\mathbf{T}^{\lambda<0}}$.} \eqref{E:positivity_condition} always holds after replacing $V$ with $V^{\oplus m}$ and $\cL$ with $\cL^m$ for a large enough integer $m>0$. If one also replaces $\chi$ with $m\chi$, the semistable locus and the stratification on $\MpC(X/G)$ remain unchanged.
\end{itemize}
\end{context}

For any $\chi \in M_\bQ^W$, we let $\MpC(X/G)^{\chi\dash\rm{ss}}$ denote the semistable locus with respect to $\cL(\chi) \in \Pic(\MpC(X/G))_\bQ$. Our goal is to compute the index of certain tautological $K$-theory classes on $\MpC(X/G)^{\chi\dash\rm{ss}}$.

Let $\ev : C \times \MpC(X/G) \to X/G$ be the universal evaluation map and let $\pi : C \times \MpC(X/G) \to \MpC(X/G)$ be the projection. 
\begin{defn}[Atiyah-Bott complexes] \label{defn: atiyah-bott complexes} We define for any $a \in \Perf(C)$ and $\xi \in \Perf(X/G)$
\begin{equation} \label{E:atiyah_bott_class}
E_a^\ast(\xi) := R\pi_\ast(a \otimes \ev^\ast(\xi)) \in \Perf(\MpC(X/G))
\end{equation}
\end{defn}

For instance, the line bundle associated in \Cref{defn: line bundle almost faithful reps} to a $G$-representation $V$ can be identified with
\[
\cD(V) \cong \det(E_{\sqrt{K}}^\ast(\cO_X \otimes V))^\dual.
\]

Throughout this section, we will let $\MpC(X/G)$ denote the \emph{derived} mapping stack $\uMap(C,X/G)$, which is the canonical derived enhancement of the stack studied earlier in this paper (see \cite[Sect. 2.2.6.3]{HAG2}, \cite{halpernleistner2019mapping}). The classes $E^\ast_a(\xi)$ extend to this derived enhancement using the same definitions. We state our main theorem with respect to this derived stack, and explain below how to translate it to a statement on the underlying classical stack $\MpC(X/G)^{\rm cl}$ by introducing a virtual structure sheaf.
\begin{thm} \label{T:main_index}
In \Cref{C:index_section_context}, let $\chi \in M^W$ and $d \in H_2(BG)$. Let $\MpC(X)^{\chi-\rm{ss}}_d \subset \MpC(X)_d$ denote the semistable locus with respect to the line bundle $\cL(\chi)$. Then for any $a \in \Perf(C)$ and $\xi \in \Perf(X/G)$, the following hold:
\begin{enumerate}
\item For all $m>0$, $R\Gamma(\MpC(X/G)_d, \cL(\chi)^m \otimes E_a^\ast(\xi))$ and $R\Gamma(\MpC(X/G)^{\chi-\rm{ss}}_d, \cL(\chi)^m \otimes E_a^\ast(\xi))$ have bounded homology groups that are coherent as $\cO_X^G$-modules;
\item For $m \gg 0$ large enough, the following restriction map is an isomorphism
\begin{equation} \label{E:quantization_reduction}
R\Gamma(\MpC(X/G)_d, \cL(\chi)^m \otimes E_a^\ast(\xi)) \to R\Gamma(\MpC(X/G)^{\chi-\rm{ss}}_d, \cL(\chi)^m \otimes E_a^\ast(\xi)); \text{ and}
\end{equation}
\item If $\xi$ is also equivariant for the action of $\bG_m$ on $X$ via scaling, i.e., $\xi \in \Perf(X/(G \times \bG_m))$, then both sides of \eqref{E:quantization_reduction} have natural auxiliary gradings whose weight space have finite dimension over $k$ (see \Cref{L:graded_sections}).
\end{enumerate}
\end{thm}

\subsubsection{Elaborations on the main theorem} Suppose that $\pi_1(G)$ is free. Then, in the context of (3) above, the graded Euler characteristic (i.e., Hilbert polynomial) of the left-hand side of \eqref{E:quantization_reduction} is computed exactly using the Teleman-Woodward index formula in \Cref{T:index_formula}. Also, \Cref{P:recursive_index_formula} strengthens the statement (2) above to an exact formula for the relationship between the graded Euler characteristic of both sides of \eqref{E:quantization_reduction} that applies for any $m$. So, the results of this section give a method for computing the graded Euler characteristic of the right-hand side of \eqref{E:quantization_reduction} for any $m>0$.

\subsection{Admissibility of tautological classes}

\begin{defn}
    A derived algebraic stack $\cM$ over a locally Noetherian base derived algebraic space $S$ is said to be \emph{quasi-smooth} if it is locally almost of finite presentation and its relative cotangent complex is perfect with Tor amplitude in $[-1,1]$.
\end{defn}

We begin by recalling some general results regarding a quasi-smooth derived algebraic stack with a $\Theta$-stratification, and then we apply these to the stack $\MpC(X/G)$.

\subsubsection{Good moduli spaces, finiteness of cohomology, and the derived structure sheaf}

If $\cM$ is a quasi-smooth algebraic derived stack and $\pi : \cM \to Y$ is a morphism to a derived scheme $Y$, we will say that $\cM$ admits a good moduli space that is proper over $Y$ if its underlying classical stack $\cM^{\rm cl} \hookrightarrow \cM$ does.

If $\cM$ is quasi-compact, then quasi-smoothness implies that the homology sheaves $H^i(\cO_\cM) = 0$ for $|i|\gg 0$.\endnote{The claim on vanishing of $H^i(\cO_\cM)$ is smooth-local, and by \cite[Cor.~2.1.6]{arinkingaitsgory}, such a stack has a smooth cover by an affine complete intersection, i.e., a stack of the form $V(f_1,\ldots,f_m) \subset \bA^n$ for some polynomials $f_1,\ldots,f_m$. In this case $\cO_{\cM}$ is presented as an $\cO_{\bA^n}$-module by a Koszul complex, which has bounded homology.} Each of these sheaves has a canonical structure of an $\cO_{\cM^{\rm cl}}$-module, which we also denote $H^i(\cO_\cM)$. We define $\cO_\cM^{\rm vir} := \bigoplus_i H^i(\cO_{\cM})[i] \in \Dqc(\cM^{\rm cl})$, which has bounded coherent homology. Using the derived projection formula for the closed immersion $j : \cM^{\rm cl} \hookrightarrow \cM$, one sees that $R(\cM \to Y)_\ast(F) \cong R(\cM \to Y)_*(\cO_\cM \otimes F)$ has a finite filtration\endnote{By a filtration of an object $E_0$ in a stable $\infty$-category, we just mean a diagram of the form $E_n \to E_{n-1} \to E_{n-2} \to \cdots \to E_0$. The $i^{th}$ associated graded object is $\cofib(E_{i+1} \to E_i)$. The filtration we have in mind here is the filtration coming from taking right truncations of $\cO_{\cM}$, \[\cdots \to \tau^{\leq -2}(\cO_\cM) \otimes F \to \tau^{\leq -1}(\cO_\cM) \otimes F \to E_0 = \cO_\cM \otimes F.\] The associated graded complexes are $\cofib(\tau^{\leq -i -1}(\cO_\cM) \to \tau^{\leq -i}(\cO_\cM)) \otimes F \cong H^i(\cO_\cM)[i] \otimes F$.} whose associated graded complex is
\[
R(\cM \to Y)_\ast(j_\ast(\cO_{\cM}^{\rm vir}) \otimes F) \cong R(\cM^{\rm cl} \to Y)_*( j^\ast(F) \otimes \cO_{\cM}^{\rm vir}).
\]
Because $\cM^{\rm cl} \to Y$ factors uniquely through $Y^{\rm cl}$, the right-hand side is canonically the pushforward of a complex on $Y^{\rm cl}$, which has bounded coherent homology if $\cM^{\rm cl}$ admits a good moduli space that is proper over $Y$ \cite[Thm. 4.16(x)]{alper-good-moduli}. Therefore, if $\cM$ admits a good moduli space that is proper over $Y$ and $F \in \Perf(\cM)$, then $R\pi_\ast(\cM,F)$ has bounded coherent homology over $Y$, and
\[
\left[R(\cM \to Y)_*(F)\right] = \left[R(\cM^{\rm cl} \to Y)_*(F|_{\cM^{\rm cl}} \otimes \cO_{\cM}^{\rm vir})\right] \; \; \;  \text{ in } K_0(Y).
\]

\subsubsection{The non-abelian localization formula}

\begin{context} \label{context: nonabelian localization formula}
    Let $\cM$ be a quasi-smooth algebraic derived stack locally of finite type and with affine diagonal over $k$, and let $\cM = \bigcup_{\alpha \geq 0} \cM_{\leq \alpha}$ be a $\Theta$-stratification of $\cM$ indexed by some totally ordered set $I$. We will let $\pi : \cM \to Y$ be a morphism to a derived scheme. For the following, we assume that $\cM$ is a \emph{local quotient derived stack}, by which we mean a derived stack that is a filtered union of stacks that can be expressed as a quotient of a derived scheme by an affine algebraic group.
\end{context} 
Recall from \cite[Lem.~1.2.3]{HL-D-equivalence} that the data of a $\Theta$-stratification of $\cM$ is equivalent to the data of a $\Theta$-stratification of the underlying classical stack $\cM^{\rm cl}$ -- the modular interpretation of the strata as open substacks of $\Filt(\cM)$ equips them with a canonical derived structure.

Let us denote the semistable locus $\cM^{\rm ss}:= \cM_{\leq 0}$, let us denote the unstable strata $\cS_\alpha \hookrightarrow \cM_{\leq \alpha} \subset \cM$, and let $\cZ_\alpha^{\rm ss} \to \cS_\alpha \to \cM$ denote the centers of the strata \cite[Def.~2.1.6]{hl_instability}. By definition, each $\cZ_\alpha^{\rm ss}$ is an open substack of $\Grad(\cM)$ \cite[Lem. 1.3.4]{HL-D-equivalence}. We will also use the following notation:
\begin{enumerate}
    \item Let $\bN_{\alpha} := \bL_{\cZ_\alpha^{\rm ss} / \cM}^\dual$ be the linear dual of the relative cotangent complex of the morphism $\cZ_\alpha^{\rm ss} \to \cM$. It is a perfect complex because $\cM$ is quasi-smooth. The central $\bG_m$-action on $\cZ_\alpha$ splits $\bN_\alpha$ as a direct sum $\bN_\alpha \cong \bN_\alpha^+ \oplus \bN_\alpha^-$, and we define
    \begin{equation} \label{E:euler_complex_general}
    \bE_\alpha := \Sym((\bN_\alpha^+)^\dual) \otimes \Sym(\bN_\alpha^-) \otimes \det(\bN_\alpha^-)[\rank(\bN_\alpha^-)],
    \end{equation}
    where the tensor product is in $\Dqc(\cZ_\alpha^{\rm ss})$. (See \cite[\S 1.0.6 \& \S 2.1]{verlinde} for further detail.) 
    \item Let $\eta_\alpha : \cZ_\alpha \to \bZ$ be the locally constant function assigning each point to the $\bG_m$-weight of the invertible sheaf $\det(\bN_\alpha^-)^\dual$.
\end{enumerate}
The complexes $\bE_\alpha \in \Dqc(\cZ_\alpha^{\rm ss})$ play the role of a reciprocal of the Euler class of the virtual normal bundle of the maps $\cZ_\alpha \to \cM$. Their key property is that under the weight decomposition for the central action of $\bG_m$ on $\cZ_\alpha$, for any $w \in \bZ$ the summand $\bE^w_\alpha$ of $\bE_\alpha$ in weight $w$ is perfect, and $\bE_\alpha^w = 0$ for $w > \eta_\alpha$ on each component of $\cZ_\alpha$.

\begin{prop}[Quantization commutes with reduction] \cite[Prop.~2.2]{verlinde}, \cite[Prop.~2.1.3]{HL-D-equivalence} \label{P:quantization commutes with reduction}
In \Cref{context: nonabelian localization formula}, let $F \in \Perf(\cM)$ be a complex such that $(F|_{\cZ_\alpha^{\rm ss}})^w \cong 0$ for $w \geq \eta_\alpha$ on each component of $\cZ_\alpha^{\rm ss}$. Then the restriction map $R\Gamma(\cM,F) \to R\Gamma(\cM^{\rm ss}, F)$ is an isomorphism.
\end{prop}

\begin{rem}
In the proof of \Cref{P:non-abelian localization} below, we will see that in fact $R(\cM \to Y)_\ast(F) \to R(\cM^{\rm ss} \to Y)_\ast(F)$ is an isomorphism under the hypotheses of \Cref{P:quantization commutes with reduction}.
\end{rem}

\begin{defn}\cite[Sect.~2.1]{verlinde} \label{defn: almost admissible}
We say that a perfect complex $F \in \Perf(\cM)$ is \emph{almost admissible} if for all but finitely many $\alpha \in I$, we have $(F|_{\cZ_{\alpha}^{\rm ss}})^w \cong 0$ for all $w \geq \eta_\alpha$. Given an invertible sheaf $\cL$ on $\cM$, we say that $F$ is \emph{$\cL$-admissible} if $\cL \otimes F^{\otimes m}$ is almost admissible for all $m > 0$.
\end{defn}

\begin{lem}\cite[Lem.~2.3]{verlinde} \label{L:admissible_complex_properties}
The category of $\cL$-admissible complexes is closed under shifts, cones, direct summands, tensor products, and symmetric powers.
\end{lem}

\begin{prop}[Virtual non-abelian localization] \label{P:non-abelian localization}
Let $\cM$ be an algebraic derived stack as above, and assume it is a local quotient derived stack. Let $F \in \Perf(\cM)$ be almost admissible, and assume that $\cM^{\rm ss}$ and every center $\cZ_\alpha^{\rm ss}$ admit proper good moduli spaces over $Y$. Then $R(\cM \to Y)_*(F) \in \Dqc(Y)$ has bounded coherent homology, and one has
\begin{equation}\label{E:non-abelian-localization}
\left[R(\cM \to Y)_*(F)\right] = \left[R(\cM^{\rm ss} \to Y)_*(F)\right] + \sum_{\alpha > 0} \left[R(\cZ_\alpha^{\rm ss} \to Y)_*(F|_{\cZ_\alpha^{\rm ss}} \otimes \bE_\alpha) \right] \in K_0(\DCoh(Y)).
\end{equation}
Furthermore, if $\pi$ is equivariant with respect to the action of an algebraic group $H$ on both $\cM$ and $Y$, the $\Theta$-stratification is $H$-equivariant (See \Cref{appendix:equivariant}), and $F \in \Perf(\cM/H)$, then the formula \eqref{E:non-abelian-localization} holds in $K_0(\DCoh(Y/H))$ as well.
\end{prop}
\begin{proof}
This is a slight generalization of \cite[Thm.~2.1]{verlinde}, but for the convenience of the reader we sketch the proof.

Let $\cU \subset \cM$ be an open union of finitely many $\Theta$-strata such that $F|_{\cZ_\alpha^{\rm{ss}}}$ has highest weight $<\eta_\alpha$ for the center of every stratum in the complement of $\cU$. We claim that the restriction map $R(\cM \to Y)_*(F) \to R(\cU \to Y)_*(F)$ is an isomorphism in $\Dqc(Y)$. The $\Theta$-stratification of $\cM$ induces one on the preimage in $\cM$ of every affine open subscheme of $Y$, and the claim is local on $Y$, so we may assume $Y$ is affine. In this case the claim follows from \Cref{P:quantization commutes with reduction}

We may therefore replace $\cM$ with $\cU$ and assume for the remainder of the proof that $\cM$ is a global quotient derived stack and the stratification is \emph{finite}.

The local cohomology filtration of $(R\cM \to Y)_*(F)$ has associated graded pieces $R(\cM^{\rm ss} \to Y)_*(F)$ and $R(\cM \to Y)_*(R\underline{\Gamma}_{\cS_{\alpha}} F)$ for each unstable stratum $\cS_\alpha$.\endnote{This follows from the general properties of local cohomology: If you have a finite stratification $\cX = \cS_0 \cup \cS_1 \cup \cdots \cup \cS_n$, where $\cS_0 \cup \cS_1 \cup \cdots \cup \cS_i$ is open for any $i$ and contains $\cS_i$ as closed substack. Then $\cO_\cX$ has a finite filtration in $\Dqc(\cX)$ whose associated graded object is $(j_0)_\ast(\cO_{\cS_0}) \oplus (j_1)_\ast(R\Gamma_{\cS_1} \cO_X) \oplus \cdots \oplus (j_n)_\ast(R\Gamma_{\cS_n} \cO_X)$, where $j_i$ denotes the open immersion $j_i: \cS_0 \cup \cS_1 \cup \ldots \cS_i \hookrightarrow \cM$. For a single stratum this comes from the defining exact triangle $R\Gamma_{\cS_1} \cO_\cX \to \cO_\cX \to (j_0)\ast \cO_{\cS_0}$, and in general it can be constructed inductively. Tensoring with $F$ and applying $R(\cM \to Y)_\ast$ gives the desired filtration of $R(\cM \to Y)_\ast(F)$.} $R(\cM^{\rm ss}\to Y)_*(F)$ has bounded coherent cohomology over $Y$ by the hypothesis that $\cM^{\rm ss}$ admits a good moduli space that is proper over $Y$. So it suffices to show that $R(\cM \to Y)_*(R\underline{\Gamma}_{\cS_\alpha}F)$ has bounded coherent cohomology.

The finiteness of $R(\cM \to Y)_*(R\underline{\Gamma}_{\cS_\alpha}F)$ is shown in \cite[Thm.~2.1]{verlinde} in the special case where $Y = \Spec(k)$, and the proof extends to our setting. Specifically, it is argued there that $R\Gamma(\cM,R\underline{\Gamma}_{\cS_\alpha}F)$ admits a convergent (as a colimit) ascending filtration that is a posteriori finite, and each of the graded pieces admit convergent (as a limit) descending filtrations that are also a posteriori finite. The associated graded of these filtrations is identified with $R(\cZ_\alpha^{\rm ss} \to Y)_*(F|_{\cZ_\alpha^{\rm ss}} \otimes \bE_\alpha)$. The filtrations involved in this argument extend naturally to filtrations of $R(\cM \to Y)_*(R\underline{\Gamma}_{\cS_\alpha} F)$,\endnote{The idea is that the first filtration of $R\Gamma(\cM,R\underline{\Gamma}_{\cS_\alpha}F)$ constructed in \cite[Thm.~2.1]{verlinde} is the result of applying $R\Gamma(\cM,-)$ to an ascending filtration of the complex $R\underline{\Gamma}_{\cS_{\alpha}}(F)$ itself, with associated graded pieces $j_\ast(\Sym^n(\bL^\dual_{\cS_\alpha/\cM}[1]) \otimes j^!(F))$, where $j : \cS_\alpha \to \cM$ is the inclusion. One can simply apply $R(\cM \to Y)_\ast(-)$ instead of $R\Gamma(\cM,-)$.

The second filtration is a convergent descending filtration of $R\Gamma(\cS_\alpha,G)$ on any object $G \in \Dqc(\cS_\alpha)$ whose associated graded pieces are $R\Gamma(\cZ_\alpha,\sigma^\ast(\Sym^m(\bL_{\cS_\alpha/\cZ_\alpha})\otimes G))$, where $\sigma : \cZ_\alpha \to \cS_\alpha$ is the canonical section of the projection $\cS_\alpha \to \cZ_\alpha$ induced by the projection $\Theta \to B\bG_m$. This filtration on global sections is obtained from a canonical filtration of the pushforward $(\cS_\alpha \to \cZ_\alpha)_\ast(G)$. So in order to get the desired filtration, one must show that the pushforward of $G$ along the composition $\cS_\alpha \to \cM \to Y$ agrees with the pushforward along the composition $\cS_\alpha \to \cZ_\alpha \to Y$. 

In fact, something stronger is true: the map $\cS_\alpha \to \cZ_\alpha$ is universal for maps to derived schemes, in the sense that it induces an equivalence of $\infty$-groupoids $\Map(\cZ_\alpha,X) \to \Map(\cS_\alpha,X)$ for any derived scheme $X$. By descent, it suffices to prove this after base change along a smooth cover $\cZ':=\Spec(A) \times B\bG_m \to \cZ_\alpha$, where the $B\bG_m$ maps to the canonical central $\bG_m$ in the inertia of $\cZ_\alpha$. Let $\cS' := \cS_\alpha \times_{\cZ_\alpha} \cZ'$.

We claim that both of the following maps are equivalences
\[
\Map(\Spec(A),X) \to \Map(\Spec(A) \times B\bG_m,X) \to \Map(\cS',X)
\]
for any derived scheme $X$. To see this, one first shows that for any morphism $f : \cS' \to X$ and any open derived subscheme $U \subset X$, $f^{-1}(U)$ is the preimage under $\cS' \to \Spec(A)$ of a unique open subset of $\Spec(A)$. Likewise for maps from $\cZ'$. This is equivalent to the same claim on underlying classical stacks, and it is known for classical $\Theta$-strata. Using this fact, one can choose a covering of $X$ by affines and reduce to the case where $X$ is affine. In this case, the claim amounts to the isomorphisms $A \cong R\Gamma(\cZ',\cO_{\cZ'}) \cong R\Gamma(\cS',\cO_{\cS'})$. The first isomorphism is clear, and for the second one can use the same universal filtration for $R\Gamma(\cS',\cO_{\cS'})$ discussed above. All of the associated graded pieces vanish except $R\Gamma(\cZ',\sigma^\ast(\cO_{\cS'})) \cong R\Gamma(\cZ',\cO_{\cZ'})$.} and the finiteness of these filtrations can be argued Zariski locally over $Y$, where it follows from the case where $Y=\Spec(k)$. So $R(\cM \to Y)_*(R\underline{\Gamma}_{\cS_\alpha} F)$ has a finite filtration whose associated graded is $R(\cZ_\alpha^{\rm ss} \to Y)_*(F|_{\cZ_\alpha^{\rm ss}} \otimes \bE_\alpha) \in \Dqc(Y)$.

The weight summands of $F|_{\cZ_\alpha^{\rm ss}} \otimes \bE_\alpha$ are perfect, and only finitely many contribute to the pushforward to $Y$. Therefore, $R(\cM \to Y)_*(\underline{\Gamma}_{\cS_\alpha}(F))$ has bounded coherent cohomology by the hypothesis that $\cZ_\alpha^{\rm ss}$ admits a good moduli space that is proper over $Y$, and $\left[R(\cM \to Y)_*( R\underline{\Gamma}_{\cS_\alpha}(F))\right] =  \left[R(\cZ_\alpha^{\rm ss} \to Y)_*(F|_{\cZ_\alpha^{\rm ss}} \otimes \bE_\alpha)\right]$ in $K_0(\DCoh(Y))$.

If $\pi : \cM \to Y$ is equivariant for the action of an algebraic group $H$, the $\Theta$-stratification is $H$-equivariant, and $F \in \Perf(\cM/H)$, then the local cohomology filtration of $R(\cM \to Y)_*(F)$, as well as the two filtrations on $R(\cM \to Y)_*(R\underline{\Gamma}_{\cS_\alpha} F)$ are $H$-equivariant, so the formula \eqref{E:non-abelian-localization} holds in $K_0(\DCoh(Y/H))$.

\end{proof}

\subsubsection{\texorpdfstring{$\cL$}{L}-admissibility of tautological classes}

We turn to the stack $\MpC(X/G)$ in \Cref{C:index_section_context}. It is well-known that $\Bun_G(C)$ is a smooth local quotient stack, and the morphism $\cM_C(X/G) \to \Bun_G(C)$ is derived affine, so $\cM_C(X/G)$ is a local quotient derived stack. The cotangent complex of the derived mapping stack $\MpC(X/G) := \Map(C,X/G)$ is described by \cite{halpernleistner2019mapping}
\[
\bL_{\MpC(X/G)} \cong R\pi_\ast(\ev^\ast(\bT_{X/G}))^\dual \cong E_{\cO_C}^\ast(\bT_{X/G})^\dual,
\]
where $\bT_{X/G}$, the linear dual of the cotangent complex $\bL_{X/G}$, is the two-term complex on $X/G$ in degrees $[-1,0]$ given by the infinitesimal action homomorphism $\cO_X \otimes \mathfrak{g} \to \cO_X \otimes X$, and $\pi: C \to \Spec(k)$ is the structure morphism. Because $\pi$ is proper of relative dimension $1$, this shows that $\MpC(X/G)$ is quasi-smooth.

For a coweight $\lambda \in N$, $\lambda$ acts trivially on $X^{\lambda=0}$, so we can decompose any quasi-coherent complex on $X^{\lambda=0}/L_\lambda$ as a direct sum of summands of constant $\lambda$-weight. We let $\bT_\lambda^+$ and $\bT_\lambda^-$ denote the summands of $\bT_{X/G}|_{X^{\lambda=0}/L_\lambda}$ with positive and negative $\lambda$-weights, respectively. We consider the complex on $\MpC(X^{\lambda=0}/L_\lambda)$,
\begin{equation} \label{E:euler_complex_specific}
\bE_\lambda := \Sym(E_{\cO_C}^\ast(\bT_\lambda^+)^\dual) \otimes \Sym(E_{\cO_C}^\ast(\bT_\lambda^-)) \otimes \det(E_{\cO_C}^\ast(\bT_\lambda^-))[\rank(E_{\cO_C}^\ast(\bT_\lambda^-))].
\end{equation}
By \Cref{lemma: components of grad for affine target}, $\Grad(\MpC(X/G))$ is a disjoint union of stacks of the form $\MpC(X^{\lambda=0}/L_\lambda)$ where $\lambda$ ranges over all anti-dominant coweights $\lambda \in N$, and the complexes $\bE_\lambda$ correspond to the complexes $\bE_\alpha$ of \eqref{E:euler_complex_general} when restricted to the graded semistable locus.

\begin{prop} \label{P:admissible complexes}
In \Cref{C:index_section_context}, fix $d_{\ker} \in \ker(\phi_V) \subset H_2(BG)$. Then, for any $a \in \Perf(C)$ and $\xi \in \Perf(X/G)$, the complex $E_a^\ast(\xi)$ is $\cL(\chi)$-admissible for the $\Theta$-stratification of \[\bigsqcup_{d \in d_{\ker} + \im(\phi_V)} \MpC(X/G)_d\] defined by $\cL(\chi)$ and the norm on graded points $b$.

\end{prop}

\begin{proof}
For any linear representation $U$ of $G$, there is an associated vector bundle $U \otimes \mathcal{O}_{X/G}$ on $X/G$ obtained by pulling back under $X/G \to BG$. As $U$ runs over all isomorphism classes of representations, the set $\{U \otimes \mathcal{O}_{X/G}\}$ split generates the triangulated category $\Perf(X/G)$\endnote{By \cite[Cor. 9.2(v)]{hall-rydh-perfect} the triangulated category $D(\QCoh(X/G))$ is compactly generated. Hence by \cite[Thm. 2.1]{neeman-bousfield} it suffices to show that for every object $A$ in $D(\QCoh(X/G))$ there exists a nonzero homomorphism $f: U \otimes \cO_{X/G}[i] \to A$ from a shift of the vector bundle corresponding to a representation $U$ of $G$. Since the morphism $X/G \to BG$ is affine, by adjunction it suffices to check that every object in $D(\QCoh(BG))$ admits a nonzero homomorphism from a shift of a linear representation $U$. This is immediate from the description of $D(\QCoh(BG))$ as the derived category of algebraic representations for the linearly reductive group $G$.}. By \Cref{L:admissible_complex_properties}, in order to prove admissibility it suffices to consider $E_a^\ast(\xi)$ where $\xi = U \otimes \cO_{X/G}$, so without loss of generality we reduce to this case.

Consider a graded-semistable point of $\Grad(\MpC(X/G)_d)$, corresponding to a cocharacter $\lambda \in N$ and a point $u \in \MpC(X^{\lambda=0}/L_\lambda)_{d'}$ for some degree $d' \in N_\bQ^{W_\lambda}$. For the corresponding HN stratum indexed by $(d', \lambda)$, the vanishing condition from \Cref{defn: almost admissible} (i.e. $(\cL(\chi) \otimes E_a^\ast(\xi)^{\otimes m}|_{\cZ_{\alpha}^{\rm ss}})^w \cong 0$ for $w \geq \eta_\alpha$)  can be written as
\begin{equation} \label{E:weight constraint}
\wt_{\lambda} (\cD(V)_u) - \langle \lambda,\chi \rangle + \maxwt_{\lambda}((E_a^\ast(\xi)_u)^{\otimes m}) < \wt_{\lambda} (\det(\bL_{\MpC(X/G),u}^{\lambda >0})),
\end{equation}
where $\maxwt_\lambda$ denotes the maximum non-zero $\lambda$-weight space in homology. The right-hand-side is the negative of the highest $\lambda$-weight of $\bE_\lambda$, and using the fact that $[\bT_{X/G}] = [\cO_X \otimes \mathbf{T}] \in K_0(X/G)$ it evaluates to $- \left(\lambda, d' \right)_{\mathbf{T}^{\lambda<0}} + (g-1) \cdot \langle \lambda, \det(\mathbf{T}^{\lambda<0}) \rangle$. One can also compute $\wt_{\lambda} (\cD(V)_u) = -(\lambda, d')_V$, and $\maxwt_\lambda((E_a^{\ast}(\xi)_u)^{\otimes m}) \leq \max \left\{ \langle \lambda, m\theta \rangle \, \left\lvert \, U_\theta \neq 0 \right. \right\}$ where the maximum runs over all $T$-weights $\theta$ of $U$ \endnote{Let $E_{\lambda}$ denote the underlying $L_{\lambda}$-bundle of the graded point $u$. By cohomology and base-change $E^*_a(\xi)_u^{\otimes m} = R\Gamma(a \otimes E_{\lambda}(U))^{\otimes m}$. The one-parameter subgroup $\lambda$ induces a grading $U = \bigoplus_{i \in \mathbb{Z}} U^{\lambda =i}$ of the $L_{\lambda}$-representation $U$. Here $\lambda$ acts with weight $i$ on $U^{\lambda =i}$, and $U^{\lambda =i}$ is nonzero if and only if there exists a $T$-weight $\theta$ of $U$ such that $\langle \lambda, \theta \rangle =i$.

By taking associated bundles, we get a decomposition $E_{\lambda}(U) = \bigoplus_{i \in \mathbb{Z}} E_{\lambda}(U^{\lambda =i})$ by $\lambda$-weights. Therefore we get a direct sum decomposition 
\[E^*_a(\xi)_u^{\otimes m} = \bigoplus_{i_1, i_2, \ldots , i_n \in  \in \mathbb{Z}} \left(\bigotimes_{j=1}^n R\Gamma_*(a \otimes E_{\lambda}(U^{\lambda =i_j})) \right)\]
$\lambda$ acts with weight $\sum_{j=1}^n i_j$ on the homologies of the graded piece $\bigotimes_{j=1}^n R\Gamma_*(a \otimes E_{\lambda}(U^{\lambda =i_j}))$. Furthermore, this graded piece is zero unless there are $T$-weights $\theta_1, \theta_2, \ldots, \theta_n$ of $U$ such that $\langle \lambda, \theta_j \rangle \leq i_j$. In particular we have $\sum_{j=1}^n i_j \leq m \cdot \max \left\{ \langle \lambda, \theta \rangle \, \left\lvert \, U_\theta \neq 0 \right. \right\}$ whenever the corresponding graded piece $\bigotimes_{j=1}^n R\Gamma_*(a \otimes E_{\lambda}(U^{\lambda =i_j}))$ is nonzero, and therefore the $\lambda$-weights of the homologies of $E^*_a(\xi)_u^{\otimes m}$ are at most $ m \cdot \max \left\{ \langle \lambda, \theta \rangle \, \left\lvert \, U_\theta \neq 0 \right. \right\}$.}. Therefore, \eqref{E:weight constraint} is implied by the inequality
\begin{equation} \label{eq: admissibility eq 2}
-\wt_{\lambda} (\cD(V)_u) + \langle \lambda, \chi \rangle > \max \left\{ \langle \lambda, m\theta + (1-g)\det(\mathbf{T}^{\lambda<0}) + \phi_{\mathbf{T}^{\lambda <0}}(d')^\dagger \rangle \, \left\lvert \, U_\theta \neq 0 \right. \right\} .
\end{equation}
Upon dividing both sides by $\|\lambda\|_b$ and replacing $\lambda$ with $\hat{\lambda} := \lambda / \|\lambda\|_b$ we see that the contribution of the stratum vanishes unless
\begin{equation} \label{eq: admissibility eq 3}
\mu \leq \max \limits_{\substack{ \theta \text{ with} \\  U_{\theta} \neq 0}} \left\langle \hat{\lambda}, m\theta + (1-g)\det(\mathbf{T}^{\lambda<0}) + \phi_{\mathbf{T}^{\lambda <0}}(d')^\dagger \right\rangle.
\end{equation}
Writing $d' = d_{\ker} + (d'-d_{\ker})$ and applying the triangle inequality, we can replace \eqref{eq: admissibility eq 3} with
\begin{equation} \label{eq: admissibility eq 4}
\mu \leq \max\limits_{\substack{ \theta \text{ with} \\  U_{\theta} \neq 0}} \left( \left\|m\theta + (1-g)\det(\mathbf{T}^{\lambda<0}) + \phi_{\mathbf{T}^{\lambda <0}}(d_{\ker})^\dagger \right\|_b \right) + \|\phi_{\mathbf{T}^{\lambda <0}}\|_b \cdot \|d'-d_{\ker}\|_b.
\end{equation}
Because $d'-d_{\ker} \in \im(\phi_V)$ and $d_{\ker} \in \ker(\phi_V)$, we have  $\|d'-d_{\ker}\|_b = \| (\phi_V^+)^{1/2} \phi_V^{1/2}(d')\|_b$, so
\[
\begin{array}{rll}
    \|d'-d_{\ker}\|_b & \leq \| \phi_V^+ \|_b^{1/2} \|d'\|_V & \text{ by Cauchy-Schwarz} \\
     & \leq \| \phi_V^+ \|_b (\|\chi\|_b + \mu) & \text{ by \Cref{L:relationship_topological_weight}(3)}
\end{array}
\]
Combining this with \eqref{eq: admissibility eq 4} and the assumption \eqref{E:positivity_condition}, we conclude that any stratum with nontrivial contribution satisfies:
\begin{equation} \label{eq: admissibility eq 5}
(1-c_{X/V})\mu \leq \max\limits_{\substack{ \theta \text{ with} \\  U_{\theta} \neq 0}} \left( \left\|m\theta + (1-g)\det(\mathbf{T}^{\lambda<0}) + \phi_{\mathbf{T}^{\lambda <0}}(d_{\ker})^\dagger \right\|_b \right) + c_{X/V}\|\chi\|_b.
\end{equation}
If $d \in d_{\ker}+\im(\phi_V)$, then any stratum in $\MpC(X^{\lambda=0}/L_\lambda)_d$ is indexed by a $d' \in d_{\ker}+\im(\phi_V)$ as well, because $d'-d \in (N_\bQ^W)^\perp \subset \im(\phi_V)$. \Cref{cor: theta strata indexing} and \Cref{lemma: finiteness of indexing data with bounded numerical invariatn} then imply the claim, because the right-hand-side of \eqref{eq: admissibility eq 5} admits an upper bound that holds simultaneously for all $\lambda$ and $d' \in d_{\ker}+\im(\phi_V)$, and only finitely many strata indexed by such pairs $(d',\lambda)$ can satisfy an upper bound on $\mu$.
\end{proof}

\begin{rem}
Because $\cL(\chi) = \cL \otimes E_{\cO_p}^\ast(\cO_X(\chi))$ for any $\chi$, \Cref{P:admissible complexes} and \Cref{L:admissible_complex_properties} imply that the category of $\cL(\chi)$-admissible complexes with respect to the $\cL(\chi')$-stratification of $\MpC(X/G)$ is independent of the choice of $\chi, \chi' \in M_\bQ^W$.
\end{rem}

\subsubsection{Graded global sections} \label{SSS:graded_global_sections}

Because $X$ is a linear representation of $G$, we can define the following algebraic derived stack
\[
\MpC^{\rm gr}(X/G): T \mapsto \left\{
\begin{array}{c} \text{invertible sheaf }L\text{ on }T, \text{ and} \\ 
G \text{ bundle } E \text{ on } C_T, \text{ and} \\
\text{section of } L \otimes E(X) \text{ over }C_T
\end{array} \right\}.
\]
The forgetful map assigning a $T$-point to the invertible sheaf $L$ on $T$ defines a morphism
\[
\MpC^{\rm gr}(X/G) \to B\bG_m.
\]
The fiber over the trivial point is $\MpC(X/G)$, which defines a $\bG_m$-action on $\MpC(X/G)$ and identifies $\MpC^{\rm gr}(X/G) \cong \MpC(X/G) / \bG_m$. There is also a forgetful morphism
\[
\MpC^{\rm gr}(X/G) \to \Bun_G(C)
\]
taking a $T$-point to its underlying $G$-bundle on $C_T$. The classical fiber of this morphism over a point $E \in \Bun_G(C)(k)$ is the quotient stack $\Gamma(C_k,E(X))) / \bG_m$, where $\bG_m$ acts on this affine space linearly with weight $1$.

\begin{lem} \label{L:explicit_computation_graded_maps}
There is an equivalence of algebraic derived stacks
\[
\MpC^{\rm gr}(X/G) \cong \Spec_{\Bun_G(C) \times B\bG_m}\left(\Sym(E^\ast_K(X^\ast[1])) \right),
\]
where on the right we regard $X$ as a linear representation of $G \times \bG_m$ that is concentrated in weight $1$ with respect to $\bG_m$, and we denote by $X^*$ its dual.
\end{lem}
\begin{proof}
By definition a $T$-point of $\MpC^{\rm gr}(X/G)$ is a $G$-bundle $E:C_T \to BG$, an invertible sheaf $L$ on $T$, and a map of $\cO_T$-module complexes $\cO_T \to L \otimes p_\ast(E^\ast(X))$, where $p : C_T \to T$ is the projection. This is naturally equivalent to a map $L^\dual \otimes p_\ast(E^\ast(X))^\dual \to \cO_T$ in $\Dqc(T)$, which by relative Serre duality (\cite[\href{https://stacks.math.columbia.edu/tag/0A9Q}{Tag 0A9Q}]{stacks-project} \cite{neeman-traces-residues}) is naturally equivalent to a map $L^\dual \otimes p_\ast(E^\ast(X^\ast) \otimes K[1]) \to \cO_T$ in $\Dqc(T)$, or equivalently a map of $dg$-algebras $\Sym_{T}(L^\dual \otimes p_\ast(E^\ast(X^\ast) \otimes K[1])) \to \cO_T$. The latter is precisely a $T$-point of the functor on the right-hand-side in the lemma.
\end{proof}

The evaluation morphism $\ev : \MpC(X/G) \times C \to X/G$ is $\bG_m$-equivariant, and therefore descends to a morphism
\[
\ev : \MpC^{\rm gr}(X/G) \times C \to X / (G\times \bG_m).
\]
Given $\xi \in \Perf(X/(G\times \bG_m))$, we can define the complex
\[
E_a^{\ast,\rm{gr}}(\xi) := \pi_\ast(a \otimes \ev^\ast(\xi)) \in \Perf(\MpC^{\rm gr}(X/G)),
\]
which is a natural $\bG_m$-equivariant structure on $E_a^\ast(\xi)$. Similarly, we equip $\cL = \cD(V)$ with a canonical $\bG_m$-equivariant structure, by letting $\bG_m$ act with weight $0$ on $V$, i.e.,
\[
\cL \cong \det\left(E_{\sqrt{K}}^{\ast,\rm{gr}}(\cO_X \otimes V)\right)^\dual.
\]
We note that the $\Theta$-stratification of $\MpC(X/G)$ descends to a $\Theta$-stratification of $\MpC^{\rm gr}(X/G)$ as well by \Cref{L:all_strata_equivariant}.

For $m \in\bZ$, we let $\cO\langle m \rangle$ denote the pullback to $\MpC^{\rm gr}(X/G)$ of the invertible sheaf on $B\bG_m$ corresponding to the character of weight $-m$, and let $E\langle m \rangle := E \otimes \cO\langle m \rangle$ for any $E \in \Dqc(\MpC^{\rm gr}(X/G))$. The pushforward of $\cL \otimes E_a^{\ast,\rm{gr}}(\xi)$ along the morphism $\MpC^{\rm gr}(X/G) \to B\bG_m$ results in a graded complex whose weight $m$ direct summand is $R\Gamma(\MpC^{\rm gr}, \cL \otimes E_a^{\ast,\rm{gr}}(\xi)\langle m\rangle)$.
\begin{lem} \label{L:graded_sections}
If $V$ is infinitesimally faithful (i.e., the kernel of $G \to \GL(V)$ is finite), then the canonical map of complexes of vector spaces gives an isomorphism
\[
\bigoplus_{m \in \bZ} R\Gamma\left(\MpC^{\rm gr}(X/G), \cL \otimes E_a^{\ast,\rm{gr}}(\xi)\langle m\rangle \right) \to R\Gamma \left(\MpC(X/G),\cL \otimes E_a^{\ast}(\xi) \right).
\]
The same holds for global sections over $\MpC^{\rm gr}(X/G)_d$ for a fixed $d \in N^W_{\bQ}$ when $V$ is almost faithful.
\end{lem}
\begin{proof}
For an arbitrary complex $F$, the complex of vector spaces underlying $R\Gamma(\MpC^{\rm gr}(X/G),F) = \bigoplus_{m \in \bZ} R\Gamma\left(\MpC^{\rm gr}(X/G), F\langle m\rangle \right)$ does not necessarily agree with $R\Gamma(\MpC(X/G),F)$, because $\MpC(X/G)$ is not quasi-compact.\endnote{For instance, if $Y = \bZ \times \Spec(k)$ is an infinite disjoint union of copies of $\Spec(k)$ indexed by $\bZ$, equipped with a trivial $\bG_m$-action, and $F$ is the equivariant sheaf on $Y$ that is $\cO \langle w \rangle$ on the component $\{[w]\}$, then $R(Y/\bG_m \to B\bG_m)_\ast(F)$ corresponds to the graded vector space $\bigoplus_{w \in \bZ} k\langle w \rangle$. The underlying vector space is not isomorphic to $R\Gamma(Y,F) \cong \prod_{w \in \bZ} k\langle w \rangle$.} However, as we have mentioned, the $\Theta$-stratification of $\MpC(X/G)$ descends to $\MpC^{\rm gr}(X/G)$, and the complex $F = \cL \otimes E_a^{\ast,\rm{gr}}(\xi)$ is admissible with respect to this stratification by \Cref{P:admissible complexes}. It follows that there is an open finite union of strata $\cU^{\rm gr} \subset \MpC^{\rm gr}(X/G)$ such that $R\Gamma(\MpC^{\rm gr}(X/G),F) \cong R\Gamma(\cU^{\rm gr},F) \in \Dqc(B\bG_m)$. The preimage $\cU := \cU^{\rm gr} \times_{B\bG_m} \Spec(k) \subset \MpC(X/G)$ is also an open finite union of strata such that $R\Gamma(\MpC(X/G),F) \cong R\Gamma(\cU,F) \in \Dqc(\Spec(k))$. The claim now follows from the fact that, by \Cref{cor: boundedness of strate affine case}, $\cU$ is quasi-compact.
\end{proof}

\subsubsection{A generalization of the Hitchin base}

The canonical morphism $X/G \to X/\!/G = \Spec(\cO_X^G)$ induces a morphism $p : \MpC(X/G) \to \MpC(X/\!/G)$. Here $Y:=  \MpC(X/\!/G)$ is a derived affine scheme whose underlying classical scheme is $X/\!/G$, because any map $C_T \to X/\!/G$ for over a classical base $T$ factors uniquely through a map $T \to X/\!/G$, but typically $H_i(\cO_Y) \neq 0$ for $i>0$. The grading on $\cO_X^G$ induces a $\bG_m$-action on $X/\!/G$. This allows one to define a $\bG_m$-action on $Y$, where
\[
Y/\bG_m : T \mapsto \left\{ \begin{array}{c} 
\bG_m \text{-bundle } L \text{ on } T  \; \; \text{and}\\
\text{section of } (C_T \times_T L) \times^{\bG_m} (X/\!/G) \text{ over }C_T
\end{array} \right\}.
\]
The morphism $p : \MpC(X/G) \to Y$ descends to the quotient $p : \MpC^{\rm gr}(X/G) \to Y/\bG_m$.

\subsubsection{Proof of \Cref{T:main_index}}

We have now collected most of the inputs for the proof. The complexes $\cL(\chi) \otimes E_a^\ast(\xi)$ are almost admissible for the $\cL(\chi)$-stratification of $\MpC(X/G)$ by \Cref{P:admissible complexes}. \Cref{P:non-abelian localization} applies to the morphism $p : \MpC(X/G) \to Y$ described above by \Cref{thm: moduli space in the case of centrally contractible}, hence $R(\MpC(X/G) \to Y)_*(\cL(\chi) \otimes E_a^\ast(\xi)) \in \Dqc(Y)$ has bounded coherent homology. $Y$ is affine with underlying classical scheme $\Spec(\cO_X^G)$, which implies claim (1).

By definition, $\cL(\chi)|_{\cZ_\alpha^{\rm ss}}$ has strictly negative weight for the canonical $\bG_m$-action on the center of any stratum. \Cref{P:admissible complexes} implies that already for $m=1$, the highest weight of $\cL(\chi)^m \otimes E_a^\ast(\xi)|_{\cZ_\alpha^{\rm ss}}$ is less than $\eta_\alpha$ for all but finitely many unstable strata, so for sufficiently large $m$ this will hold for all $\alpha$. The claim (2) then follows from \Cref{P:quantization commutes with reduction}.

When $\xi$ is the restriction of a class in $\Perf(X/(G\times \bG_m))$, \Cref{L:graded_sections} establishes the auxiliary grading on $R\Gamma(\MpC(X/G),\cL(\chi) \otimes E_a^\ast(\xi))$. In fact, this complex canonically has the structure of a graded $dg-\cO_Y$-module via pushforward of $\cL(\chi) \otimes E_a^{\ast,\rm{gr}}(\xi)$ along the morphism $p : \MpC^{\rm gr}(X/G) \to Y/\bG_m$. The fact that the weight pieces are finite dimensional follows from (1), because $\cO_X^G$ is strictly positively graded, so any coherent graded module has finite dimensional weight pieces.

\qed


\subsection{A recursive formula for gauged Gromov-Witten invariants}

We continue to work in \Cref{C:index_section_context}. \Cref{T:main_index}(3) allows us to make the following:


\begin{defn}[Gauged Gromov-Witten invariants]
For $\chi \in M_\bQ^W$, $d \in (M^W)^\dual$, and an $\cL$-admissible complex $F \in \Perf(\MpC^{\rm gr}(X/G))$, we introduce the following notation for the graded index both over the stack and the $\chi$-semistable locus
\begin{align*}
\GGW_d(X/G,F) &:= \chi^{\rm gr}(\MpC^{\rm gr}(X/G)_d,\cL \otimes F) \\
\GGW^{\chi}_d(X/G,F) &:= \chi^{\rm gr}(\MpC^{\rm gr}(X/G)_d^{\chi\rm{-ss}},\cL \otimes F),
\end{align*}
where $\chi^{\rm gr}(E) = \sum_{m} \chi(E\langle m \rangle) t^m \in \bZ(\!(t)\!)$ denotes the graded Euler characteristic.
\end{defn}

\Cref{T:main_index}(2) says that $\GGW_d(X/G,F) = \GGW^\chi_d(X/G,F)$ under certain ``positivity'' assumptions. In this section we establish a more precise relationship between $\GGW_d(X/G,F)$ and $\GGW_d^\chi(X/G,F)$ that holds more generally.

Given $F \in \Perf(\MpC^{\rm gr}(X/G))$ almost admissible, we will abuse notation by letting
\[
F \otimes \bE_\lambda \in \Dqc(\MpC^{\rm gr}(X^{\lambda=0}/L_\lambda))
\]
denote the restriction of $F$ along the canonical morphism $\MpC^{\rm gr}(X^{\lambda=0}/L_\lambda) \to \MpC^{\rm gr}(X/G)$ tensored with the quasi-coherent complex \eqref{E:euler_complex_specific}. Although $F \otimes \bE_\lambda$ is quasi-coherent, each $\lambda$-weight summand is perfect and $\cL$-admissible by \Cref{P:admissible complexes} and \Cref{L:admissible_complex_properties} \endnote{Note that for any $\lambda$ we still have $c_{X/V} <1$ after replacing $X$ with $X^{\lambda=0}$ and $G$ with $L_{\lambda}$}, and the weight pieces vanish for $\lambda$-weight sufficiently positive. It follows that the $\lambda$-weight $0$ piece of $F \otimes \bE_{\lambda}$ is $\cL$-admissible, and so $\GGW_d(X^{\lambda=0}/L_\lambda,F \otimes \bE_\lambda)$ and $\GGW^\chi_d(X^{\lambda=0}/L_\lambda,F \otimes \bE_\lambda)$ are well-defined.

\begin{prop} \label{P:recursive_index_formula}
In \Cref{C:index_section_context}, for any $\cL$-admissible $F \in \Perf(\MpC(X/G))$, $d \in H_2(BG)$, and $\chi \in M^W_\bQ$, one has
\begin{equation}\label{E:recursive_formula}
I_d^\chi(X/G,F) = I_{d}(X/G,F) - \sum_{\substack{ (d', \lambda) \text{ $\chi$-active} \\\text{with $\proj^{\perp}_{N^W_{\mathbb{Q}}}(d')=d$}}} I_{d'}^{\chi- \lambda^{\dagger}}(X^{\lambda=0}/L_{\lambda},F \otimes \bE_{\lambda}),
\end{equation}
Furthermore, if $F = E_a^{\ast,\rm{gr}}(\xi)$ for some $\xi \in \Perf(X/G\times\bG_m)$, then the only non-vanishing contributions to \eqref{E:recursive_formula} come from terms with $(d', \lambda)$ satisfying
\begin{equation} \label{eq: admissibility prop inequality}
     \|d'\|_V \leq \frac{\|\phi_{V}^{+}\|^{1/2}_b }{1-c_{X/V}} \left( \max\limits_{\substack{ \theta \text{ with} \\  H^\ast(\xi|_0)_{\theta} \neq 0}}\left( \left\|\theta +\chi + (1-g)\det(\mathbf{T}^{\lambda<0}) + \phi_{\mathbf{T}^{\lambda<0}}(d_{\ker})^{\dagger}\right\|_b \right) + \|\chi\|_b \right)
 \end{equation}
where $\xi|_0$ is the restriction to the origin $0 \in X(k)$, $H^*(\xi|_0)_{\theta}$ denotes the $\theta$-weight space of the graded vector space $H^\ast(\xi|_0)$ for any given $\theta \in M$, and $d_{\ker}$ is the $b$-orthogonal projection of $d$ onto $\ker(\phi_V)$.
\end{prop}

\begin{proof}
The formula \eqref{E:recursive_formula} is just a translation of the non-abelian localization formula \eqref{E:non-abelian-localization} into this context. For any indexing datum $(d',\lambda)$ that appears, the projection of $d'$ onto $N^W_{\mathbb{Q}}$ is $d$, since we are considering the strata of the stack $\MpC(X/G)_d$ with fixed $d \in H_2(BG)$. In particular, because $\ker(\phi_V) \subset N^W_\bQ$, all $d'$ that appear have the same projection $d_{\ker} \in \ker(\phi_V)$ as well.

For the explicit estimate, we first choose a quasi-isomorphism between $\xi$ and a unique minimal $G \times \bG_m$-equivariant complex $(\cO_X \otimes U_\bullet, \partial)$, meaning $\partial|_0 = 0$.\endnote{Equivariant coherent sheaves on $X$ are equivalent to $G$-equivariant $\cO_X$ modules. The objects $\cO_X \otimes U$ for $U \in \Rep(G)$ are projective, because $G$ is linearly reductive, and they generate the category in the sense that any object admits a surjection from a direct sum of objects of this kind. Therefore, the usual way of constructing a projective resolution of a bounded complex of $\cO_X$-modules works just as well in the category of $G$-equivariant $\cO_X$-modules, resulting in a presentation of the form $(\cO_X \otimes U_\bullet,\partial)$ for any bounded above complex $\xi$ of the $G$-equivaraint $\cO_X$-modules. Furthermore, if the complex $\xi$ is also invariant for the scaling action of $\bG_m$ on $X$, i.e., it is a homogeneous $G$-equivariant complex of \emph{graded} equivariant $\cO_X$-modules, then the usual algorithm for constructing a presentation $(\cO_X \otimes U,\partial)$ such that $\partial_0=0$, as discussed in \cite[1B]{eisenbud-geometry-syzygies} for instance, carries over in the $G$-equivariant setting as well.} Note that $U_i \cong H^i(\xi|_0)$. We now apply the bound \eqref{eq: admissibility eq 5} in the proof of \Cref{P:admissible complexes} to the terms of the complex $\cO_X \otimes U_\bullet(-\chi)$, since we are computing the index of $\cL \otimes E^*_a(\xi) = \cL(\chi) \otimes E^*_a(\xi(-\chi))$. It says that any stratum that contributes to the cohomology must satisfy
\[
(1-c_{X/V}) \mu \leq \max\limits_{\substack{ \theta \text{ with} \\  H^\ast(\xi|_0)_{\theta} \neq 0}}\left( \left\|\theta +\chi + (1-g)\det(\mathbf{T}^{\lambda<0}) + \phi_{\mathbf{T}^{\lambda<0}}(d_{\ker})^{\dagger}\right\|_b \right) + c_{X/V}\|\chi\|_b .
\]
To deduce \eqref{eq: admissibility prop inequality}, we add $(1-c_{X/V} )\|\chi\|_b$ to both sides of this inequality and combine it with the inequality $\|d'\|_V \leq \|\phi_V^+\|_b^{1/2}(\|\chi\|_b + \mu)$ from \Cref{L:relationship_topological_weight}(3).

\end{proof}

\begin{rem}
By \Cref{lemma: coxeter lattice for indexing data}, any $\chi$-active indexing data $(d', \lambda)$ satisfies $d' \in \frac{1}{H} N$. We define the radius of $V$ to be
\[\mathrm{r}_V := \frac{1}{\|\phi_V^{+}\|_b^{1/2}}\left(\min_{\substack{w \in \frac{1}{H}N \\ \|w\|_V \neq 0}} \|w\|_V\right)\]
By \Cref{P:recursive_index_formula}, if we have
\[ \mathrm{r}_V >\frac{1}{1-c_{X/V}} \left(\max\limits_{\substack{ \theta \text{ with} \\  H^\ast(\xi|_0)_{\theta} \neq 0}}\left( \left\|\theta +\chi + (1-g)\det(\mathbf{T}^{\lambda<0}) + \phi_{\mathbf{T}^{\lambda<0}}(\proj_{\ker(\phi_V)}^{\perp}(d))^{\dagger}\right\|_b \right) + \|\chi\|_b\right)\]
then a $\chi$-active indexing datum ($d', \lambda)$ contributes to the formula \eqref{E:recursive_formula} for $F = E_a^{\ast,\rm{gr}}(\xi)$ only if $d' =d=\proj^{\perp}_{\ker(\phi_V)}(d)$. If in addition $\phi_V(d)+\chi^{\dagger} =0$, then this forces $I_{d}(X/G,E_a^{\ast,\rm{gr}}(\xi)) = I_d^\chi(X/G,E_a^{\ast,\rm{gr}}(\xi))$.\endnote{Indeed, recall that any $\chi$-active datum $(d', \lambda)$ indexing a stratum satisfies the relation $\lambda = \proj^{\perp}_{\Span(\sigma)}(\phi(d')+\chi^{\dagger})$ (\Cref{defn: chi indexing datum}(2)). If $d = d'$ and $\phi(d) + \chi^{\dagger} =0$, it follows that $\lambda =0$, and hence the stratum must have been the semistable locus.}
\end{rem}

\subsubsection{The algorithm for computing gauged Gromov-Witten invariants}
Let $Z_X \subset Z$ denote the kernel of the representation $Z \hookrightarrow G \to \GL(X)$ in the maximal central torus $Z \subset G$. If the linear functional $(d,-)_V + \chi$ does not restrict identically to $0$ on the subspace of cocharacters $X_*(Z_X)_{\mathbb{R}} \subset X_*(T)_{\mathbb{R}}$, then the semistable locus $\MpC(X/G)_d^{\chi\dash {\rm ss}}$ is empty\endnote{Suppose that the restriction $(d,-)_V + \chi$ to $X_*(Z_X)_{\mathbb{R}}$ is not identically zero. This means that there exists some central cocharacter $\lambda $ in $Z_X$ such that $(\lambda, d)_V + \chi(\lambda) >0$. Since $\lambda$ acts trivially on $X$, for any point in $\MpC(X/G)$ there is a filtration defined by $\lambda$ and the trivial parabolic reduction to $G= P_{\lambda}$. The assumption that $(\lambda, d)_V + \chi(\lambda) >0$ implies that this filtration is always destabilizing, and therefore the semistable locus $\MpC(X/G)_d^{\chi\dash {\rm ss}}$ is empty.}, and so in particular the index satisfies $I^{\chi}_d(X/G,F)=0$ for all $F= E_a^{\ast,\rm{gr}}(\xi)$.

Therefore, for the purposes of computing indexes, we are mostly interested in the case when $(d,-)_V + \chi$ restricts identically to $0$ on $X_*(Z_X)_{\mathbb{R}}$. In this case the right-hand side of \eqref{E:recursive_formula} involves the following terms:
\begin{itemize}
    \item The index $I_d(X/G, F)$ over the whole stack, which we will compute explicitly using the Teleman-Woodward index formula (see \Cref{T:index_formula} below).
    \item Indexes of the form $I_{d'}^{\chi- \lambda^{\dagger}}(X^{\lambda=0}/L_{\lambda},F \otimes \bE_{\lambda})$ where either $L_{\lambda} \subsetneq G$ is a strictly smaller parabolic subgroup, or $L_{\lambda}=G$ and $X^{\lambda =0} \subsetneq X$ is a strictly smaller subrepresentation.\endnote{Suppose for the sake of contradiction that an unstable stratum indexed by $(d', \lambda)$ satisfies $L_{\lambda}=G$ and $X^{\lambda =0} = X$. The condition $L_{\lambda} =G$ implies that $\lambda$ is a rational cocharacter of $Z \subset G$, and furthermore the condition $X^{\lambda=0} = X$ implies that $\lambda$ is in $X_*(Z_X)_{\mathbb{R}}$. By the assumption that $(d',\lambda)$ parametrizes an unstable stratum, it follows that $(d,\lambda)_V + \chi(\lambda) >0$, which shows that $(d,-)_V + \chi$ restricts to a nonzero linear functional on $X_*(Z_X)_{\mathbb{R}}$.}
\end{itemize}

This in principle allows one to recursively compute $I_d^\chi(X/G,F)$ for any $d$ and $\chi$.



\subsubsection{Recollection of the Teleman-Woodward index formula}

The Teleman-Woodward formula \cite{teleman-woodward} computes the index of certain tautological $K$-theory classes on $\Bun_G(C)$. Let us establish some terminology:

\begin{itemize}
\item As before, $N$ denotes the cocharacter lattice of a maximal torus $T \subset G$. Analytically, we identify $T \cong N_\bC / N(1)$, where $N(1) = 2\pi i N \subset N_\bC := N \otimes \bC$, and we denote the quotient map $e^{(-)} : N_\bC \to T$. Likewise, the $\bZ$-dual of $N$ is $M = \Hom(N(1),\bZ(1)) = \Hom(T,\bC/\bZ(1))$. In particular, for $\alpha \in M \subset M_\bC$, the formula $e^\alpha$, which is a-priori only a map $N_\bC \to \bC^\ast$ actually descends to $T = N_\bC/N(1)$. \\

\item We fix a Borel subgroup $B \supset T$, which yields a subset of positive roots $\Phi^{+} \in M$. We denote by $\Delta: N_{\bC} \to \bC$ the Weyl denominator given by
\[ \Delta(\xi) = \prod_{\alpha \in \Phi^{+}} \left( e^{-\alpha(\xi)/2} - e^{\alpha(\xi)/2}\right)\]

\item An invertible sheaf $\cL \in \Pic(\Bun_G(C))$ has a \emph{level} $h \in H^4(BG;\bQ)$ if $c_1(\cL) = c_1(\cD(U))$ for some virtual representation $U \in \Rep(G)_\bQ$ with ${\rm ch}_2(U) = -h$. The level $h$ is said to be positive if it lies in the cone spanned by $-{\rm ch}_2(U)$ for all infinitesimally faithful representations $U$. $h$ is positive if and only if it is identified with a positive definite bilinear form under the correspondence $H^4(BG;\bR) \cong \Sym^2(N(1)_\bR)^W$ that identifies $\ch_2(U)$ with the trace form $(\eta,\xi)_U := {\rm Tr}_U(\eta \xi)$. We let $c = -\frac{1}{2} {\rm Tr}_{\mathfrak{g}}$.
\end{itemize}

\begin{defn}
A line bundle $\cL \in \Pic(\Bun_G(C))$ with positive level $h\in H^4(BG;\bQ)$ is TW-admissible if $h' := h+c$ is a negative definite integer valued symmetric bilinear form on $N$.
\end{defn} 

Contraction with $h'$ defines a homomorphism $N_\bC \to M_\bC$ taking $\xi \mapsto (\xi,-)_{h+c}$. Note that this map takes $N(1)$ to $M(1)$ and thus descends to an isogeny $\chi' : T \to T^\dual$. Let $F_{\rho} \subset N_\bC$ denote the set of solutions of the equation
\[
(\xi,-)_{h+c} = 2 \pi i \rho \mod M(1),
\]
where $\rho \in M$ is the half sum of positive roots. Note that the group $\widetilde{W} := N(1) \rtimes W$ acts on $F_{\rho}$, and we let $F^{\rm reg}_\rho \subset F_\rho$ denote the points whose conjugacy class in the compact form of $G$ is regular, i.e., where $W$ acts freely. Note that $F_{\rho}/N(1) \subset T$ is the preimage of $(2\pi i \rho  \mod M(1) )\in T^{\vee}$ under the isogeny $\chi'$, and so in particular $F_{\rho}/N(1)$ is finite. It follows that there are finitely many $\widetilde{W}$-orbits in $F_{\rho}$ and in $F_{\rho}^{reg}$.

For $U \in \Rep(G)$, let $\ch_U : T=N_\bC/N(1) \to \bC$ be the trace function. It will be convenient to regard $\ch_U$ as a function on the Lie algebra $N_\bC$ itself, and to regard the differential $d\ch_U$ as a holomorphic map $N_\bC \to M_\bC$. If $U$ has weights $\beta_1,\ldots,\beta_l \in M$, then $\ch_U(\xi) = \sum_{i=1}^l e^{\beta_i(\xi)}$, and
\[
d\ch_U(\xi) = \sum_{i=1}^l e^{\beta_i(\xi)} \beta_i \in M_\bC.
\]
We let $H_U$ denote the Hessian of the map $\ch_U : N_\bC \to \bC$, so $H_U(\xi)$ is the bilinear form $\ch_U(e^\xi \mu \nu)$ for $\mu,\nu \in N_\bC$. We use $h'$ to convert this to an endomorphism $H_U(\xi)^\dagger : N_\bC \to N_\bC$. If $\beta_1^{\dagger},\ldots,\beta_l^{\dagger} \in N_{\mathbb{Q}}$ are the duals of the weights of $U$ with respect to $h+c$, then
\[
H_U(\xi)^\dagger = \sum_i e^{\beta_i(\xi)} \beta_i^{\dagger} \beta_i.
\]

Both $H_U(\xi)^\dagger$ and $d\ch_U(\xi)$ depend additively on the character of $U$, so for any $\xi \in N_\bC$, the map $U \mapsto d\ch_U(\xi)$ extends to a $\bZ[\![\mathbf{t}]\!]$-linear homomorphism $K_0(BG)[\![\mathbf{t}]\!] \to M_\bC[\![\mathbf{t}]\!]$, where $\mathbf{t}$ denotes some (possibly infinite) collection $t_1,t_2,\ldots$ of formal variables. Likewise, we define $H_{U_{\mathbf{t}}}(\xi)^\dagger \in \mathrm{End}(N_\bC)[\![\mathbf{t}]\!]$ and $E_K^\ast(U_{\mathbf{t}}) \in K_0(\Bun_G(C))[\![\mathbf{t}]\!]$ for $U_{\mathbf{t}} \in K_0(BG)[\![\mathbf{t}]\!]$ by linear extension.

\begin{thm}[\cite{teleman-woodward}]\label{T:teleman_woodward_index}\endnote{Rather than considering an equation involving a formal deformation of a morphism of lattices $N_\bC \to M_\bC$ taken modulo $M(1)$, as we do here, the paper \cite{teleman-woodward} expresses the sum over solutions of the exponentiated formal map of tori $T \to T^\dual$. The two formulations are equivalent.

In addition, the main result \cite[Thm.~2.18]{teleman-woodward} is only stated for the case when there is a single formal variable and $W_{\mathbf{t}} = t W$. However, the extension to possibly infinitely many indeterminates and arbitrary $W_{\mathbf{t}}$ is straightforward, and already appears in the application of the formula in \cite[Sect.~6]{teleman-woodward}.}  Suppose that $\pi_1(G)$ is free. Let $\cL$ be a TW-admissible line bundle on $\Bun_{G}(C)$ with positive level $h$ as above, let $U \in K_0(BG)$, and let $W_{\mathbf{t}} \in K_0(BG)[\![\mathbf{t}]\!]$ with vanishing constant term. Let $F^{\rm reg}_{\rho,\mathbf{t}} \subset N_\bC[\![\mathbf{t}]\!]$ denote be the set of solutions $\xi_{\mathbf{t}}$ of the equation
\[
(\xi_{\mathbf{t}},-)_{h+c} + d\ch_{W_{\mathbf{t}}}(\xi) = 2 \pi i \rho \mod M(1),
\]
that lie in $F^{\rm reg}_\rho$ when $0 = \mathbf{t}=(t_1,t_2,\ldots)$. If we denote
\[
\theta_{\mathbf{t}}(\xi) = \det(1+H_{W_{\mathbf{t}}}(\xi)^\dagger)^{-1} \frac{\Delta(\xi)^2}{|F_\rho/N(1)|},
\]
then
\[
\chi\left(\Bun_G(C),\cL \otimes E_{\cO_p}^\ast(U) \otimes \exp(E_{\sqrt{K}}^\ast(W_{\mathbf{t}}))\right) = \sum_{\xi_{\mathbf{t}} \in F^{\rm reg}_{\rho,\mathbf{t}}/\widetilde{W}} \ch_U(\xi_{\mathbf{t}}) \theta_{\mathbf{t}}(\xi_{\mathbf{t}})^{1-g}.
\]
\end{thm}

\subsubsection{Index formula on \texorpdfstring{$\MpC(X/G)$}{M_C(X/G)}}

Our goal is to use \Cref{T:teleman_woodward_index} to compute the index $I_d(X/G,F)$ for Atiyah-Bott classes $F \in \Perf(\MpC^{\rm gr}(X/G))$. Instead of individual $d$, we will compute the generating function
\[
I(X/G,F) := \sum_{d \in H_2(BG)} z^d I_d(X/G,F) = \sum_{d \in H_2(BG)} z^d \chi^{\rm gr}(\MpC^{\rm gr}(X/G)_d,\cL \otimes F).
\]

We will also use the equivariant parameter $t$ for $\bG_m$, which we regard as a formal parameter for the purposes of computing the index.

By hypothesis, contraction with the level $h$ of $\cL$ induces a morphism of lattices $N \to M$ given by, which induces an inclusion of finite index
\[
H_2(BG) = (M^W)^\dual \subset (N^W)^\dual, \; \;  \; \; \; d \mapsto h(d,-).
\]
Note that $(N^W)^\dual$ is the space of characters for the maximal central torus $Z \subset G$, whereas $H_2(BG)$ is the group of characters for the torus dual to the abelianization $G_{\rm ab}^\dual$. The inclusion $H_2(BG) \subset (N^W)^\dual$ yields an inclusion of corresponding group rings $\cO_{G_{\rm ab}^\dual} \subset \cO_{Z}$.

For $\mu \in (N^W)^\dual$, we will let $z^\mu$ denote the corresponding monomial in the group ring $\bC[(N^W)^\dual] = \cO_{Z}$. For $\beta \in M$, we will also regard $z^\beta \in \cO_{Z}$ under the canonical restriction map $M \to (N^W)^\dual$. For any $U \in \Rep(G)$, we will use the notation $U_\mu \in \Rep(G)$ to denote the summand of $U$ of $Z$-weight $\mu \in (N^W)^\dual$.

\begin{thm}\label{T:index_formula}
Suppose that $\pi_1(G)$ is free and that $\cL$ is TW-admissible. Let $U,U' \in \Rep(G)$ be regarded as representations of $G \times \bG_m$ on which the second factor acts trivially. Let $\beta_1,\ldots,\beta_n \in M$ be the weights of $X^\ast$. Then
\begin{equation} \label{E:index_formula}
I(X/G, E^{\ast,\rm{gr}}_{\cO_p}(U') \otimes \exp(s E_{\sqrt{K}}^{\ast,\rm{gr}}(U))) = \sum_{\xi_{z,s,t} \in F^{\rm reg}_{z,s,t} / \widetilde{W}} \sum_\mu z^{-\mu} \ch_{U'_\mu}(\xi_{z,s,t}) \theta_{z,s,t}(\xi_{z,s,t})^{1-g}.
\end{equation}
Here $F_{z,s,t}^{\rm reg} \subset (N_\bC \otimes \cO_{Z})[\![s,t]\!]$ is the set of solutions $\xi_{z,s,t}$ of the equation
\[
(\xi,-)_{h+c}+ \sum_{\mu \in (N^W)^\dual} s z^{-\mu} d\ch_{U_\mu}(\xi) + \sum_{i=1}^n  \ln(1-tz^{-\beta_i} e^{\beta_i(\xi)}) \beta_i = 2\pi i \rho \mod M(1)
\]
that are regular at $0=s=t$. Also, we let $\beta_i^{\dagger} \in N_{\mathbb{Q}}$ be the duals of the elements $\beta_i$ under the bilinear form $h+c$, and define
\[
\theta_{z,s,t}(\xi) = \det\left(1+\sum_{\mu \in (N^W)^\dual} s z^{-\mu} H_{U_\mu}(\xi)^\dagger - \sum_{i=1}^n \frac{tz^{-\beta_i} e^{\beta_i(\xi)}}{1-tz^{-\beta_i}e^{\beta_i(\xi)}}  \beta_i^\dagger \beta_i\right)^{-1} \frac{\Delta(\xi)^2}{|F_\rho / N(1)|}.
\]
\end{thm}

\begin{rem}
The right-hand-side of \eqref{E:index_formula} a priori takes values in $\cO_{Z}[\![s,t]\!]$, but part of the content of the theorem is that a posteriori the right-hand side is valued in $\cO_{G_{\rm ab}^\dual}[\![s,t]\!] \subset \cO_{Z}[\![s,t]\!]$, i.e., the coefficient of $z^\mu$ vanishes for any $\mu \in (N^W)^\dual$ that is not of the form $h(d,-)$ for some $d \in H_2(BG)$.
\end{rem}

\begin{proof}[Proof of \Cref{T:index_formula}]
To ease notation, let $F := E_{\cO_p}^{\ast,\rm{gr}}(U') \otimes \exp(s E_{\sqrt{K}}^{\ast,\rm{gr}}(U))$ which regard as a formal sum (with parameter $s$) of complexes on $\MpC^{\rm gr}(X/G)$, which are in fact pulled back from $\Bun_{G}(C)$ because $U$ and $U'$ have $\bG_m$-weight $0$. By pushing forward along the map $\MpC^{\rm gr}(X/G) \to \Bun_G(C) \times B\bG_m$, one identifies $I_d(X/G,F)$ with
\[
\chi(\Bun_{G}(C)_d, F \otimes \cL \otimes \lambda_{-t}(E_K^\ast(X^\ast))),
\]

There is a canonical embedding $Z \times \Bun_G(C) \hookrightarrow I_{\Bun_G(C)}$ of group schemes over $\Bun_G(C)$, where $I_{\Bun_G(C)}$ denotes the inertia group. Using this, one can decompose any $E \in \Dqc(\Bun_G(C))$ uniquely as a direct sum $E = \bigoplus_{\mu \in (N^W)^\dual} E_\mu$, where $E_\mu$ has constant $Z$-weight $\mu$. On the component of $\Bun_G(C)$ indexed by $d \in H_2(BG) \cong (M^W)^\dual$, the $Z$-weight of $\cL$ is $h(d,-) \in (N^W)^\dual$ \cite[Lem.~3.1]{verlinde}. Therefore, for any $E$,
\[
R\Gamma(\Bun_G(C), \cL \otimes E_{\mu}) = \left\{ \begin{array}{l} R\Gamma(\Bun_{G}(C)_d, \cL \otimes E), \text{ if }\mu = -h(d,-) \\ 0,\text{ if }\mu \notin H_2(BG) \subset (N^W)^\dual \end{array} \right. .
\]
It follows that
\[
I(X/G,F) = \chi\left(\Bun_G(C), \sum_{\mu \in (N^W)^\dual} z^{-\mu} \cL \otimes \left[F \otimes \lambda_{-t}(E_K^\ast(X^\ast)) \right]_\mu \right),
\]
where $z^{\mu}$ is treated as a formal variable. This sum makes sense because only finitely many $z^\mu$ terms are nonvanishing modulo $t^m$ for any $m$. Finally, the formal sum $\sum z^{-\mu} [F \otimes \lambda_{-t}(E_K^\ast(X^\ast))]_\mu$ in $(K_0(\Bun_G(C))\otimes \cO_{Z})[\![s,t]\!]$ can be rewritten as
\[
(\sum_\mu z^{-\mu} E_{\cO_p}^\ast(U'_\mu)) \otimes \bigotimes_\mu \exp(sz^{-\mu} E_{\sqrt{K}}^\ast(U_\mu)) \otimes \bigotimes_\mu \lambda_{-t z^{-\mu}}(E_K^\ast(X^\ast_\mu)).
\]

In algebraic $K$-theory modulo algebraic equivalence (see \cite[Sect.~2.1.1]{verlinde}), one has $[K] \sim [\sqrt{K}]+(g-1)[\cO_p]$, so
\[
\lambda_{-tz^{\mu}}(E^\ast_{K}(X^\ast_\mu)) \sim \lambda_{-tz^{\mu}}(E_{\sqrt{K}}^\ast(X^\ast_\mu)) \otimes \lambda_{-tz^{\mu}}(E_{\cO_p}^\ast(X^\ast_\mu))^{\otimes(g-1)}.
\]
Expressing the class $\lambda_{-tz^\mu}$ in terms of Adams operations $\psi^p$ as in \cite[Sect.~6]{teleman-woodward}, and using the fact that $\psi^p(E_{\sqrt{K}}^\ast(-)) = \frac{1}{p} E_{\sqrt{K}}^\ast(\psi^p(-))$, one obtains
\[
\lambda_{-tz^{-\mu}}(E_{\sqrt{K}}^\ast(X^\ast_\mu)) = \exp \left[-\sum_{p>0} \frac{(tz^{-\mu})^p}{p^2} E_{\sqrt{K}}^\ast(\psi^p(X^\ast_\mu)) \right].
\]
Now if $X^\ast$ has character $\sum_i e^{\beta_i}$, then $\psi^p(X^\ast)$ has character $\sum_i e^{p\beta_i}$, $d\ch_{\psi^p(X^\ast)}(\xi) = p \sum_i e^{p \beta_i(\xi)} \beta_i$, and $H_{\psi^p(X^\ast)}(\xi)^\dagger = \sum_i e^{p \beta_i(\xi)} p^2 \beta_i^\dual \beta_i$. The proposition now results from directly applying \Cref{T:teleman_woodward_index}, treating each $(tz^{-\mu})^p$ as a different formal variable.
\end{proof}

\subsubsection{Index of a general Atiyah-Bott class}

For a general element $a \in K_0(C)$, one has an equality in algebraic $K$-theory modulo algebraic equivalence
\[
a \sim \rank(a) [\sqrt{K}] + (\deg(a)+1-g)[\cO_p].
\]
This decomposition gives
\[
E_a^\ast(U) \sim \rank(a)  E_{\sqrt{K}}^\ast(U) + (\deg(a)+1-g) E_{\cO_p}^\ast(U) \in K_0(\MpC^{gr}(X/G))
\]
Therefore it suffices to compute $I(X/G,E_a^\ast(U))$ for $a = \cO_p$ and $a = \sqrt{K}$. When $a = \cO_p$, this is obtained from \eqref{E:index_formula} by setting $U=U'$ and $s=0$. When $a = \sqrt{K}$, it is obtained by setting $U'$ to be the trivial representation of dimension $1$, taking the derivative of \eqref{E:index_formula} with respect to $s$, then setting $s=0$.\endnote{ $\chi(\MpC^{gr}(X/G),E_{\sqrt{K}}^\ast(U) \otimes \cL)$ is given by
\[
\sum_{\xi_t \in F^{reg}_{t}/\widetilde{W}} \theta_t(\xi_t)^{1-g} \left(\sum_{\alpha>0} \frac{1+e^{\alpha(\xi_t)}}{1-e^{\alpha(\xi_t)}} \alpha(\sum_j n_j \mu_j(\xi_t) \mu_j^\dual) -  \Tr \left( (1-\sum_{i=1}^N \frac{t e^{\beta_i(\xi)}}{1-te^{\beta_i(\xi)}}  \beta_i^\dual \beta_i)^{-1} ( \sum_j n_j \mu_j^\dual \mu_j ) \right) \right),
\]
where $F^{reg}_t$ are the regular solutions to $(\xi,-)_{h+c} + \sum_i \ln(1-te^{\beta_i(\xi)})\beta_i = 2\pi i \rho \mod M(1)$, and 
\[
\theta_{t}(\xi) = \det\left(1 - \sum_{i=1}^N \frac{te^{\beta_i(\xi)}}{1-te^{\beta_i(\xi)}}  \beta_i^\dagger \beta_i\right)^{-1} \frac{\Delta(\xi)^2}{|F_\rho / N(1)|}.
\]
}


\section{The case of orbifold curves} \label{section: orbifold curves}
In this section, we work over an algebraically closed field $k$ of characteristic $0$. We let $S = \mathcal{Y} = \Spec(k)$. We replace $C$ with a proper smooth orbifold curve $\mathcal{C}$, and show that the same results of monotonicity of HN boundedness hold for appropriate numerical invariants on $\Map(\mathcal{C}, X/G)$, where $X$ is affine (\Cref{prop: montonicity for orbifolds} and \Cref{prop: hn boundedness for orbifolds}).

\subsection{Smooth orbifold curves}
\begin{defn}
An algebraic stack $\mathcal{C}$ over $\Spec(k)$ is called an orbifold curve if it is a separated, connected, DM stack of finite type and of dimension 1 over $\Spec(k)$ with generically trivial stabilizers.
\end{defn}
If $\cC$ is a smooth orbifold curve, then its coarse moduli space $C$ is a smooth quasiprojective curve over $k$, and by the proof of \cite[Thm. 4.2.1(1)]{abramovich-graber-vistoli-gw-dm-stacks} there are reduced effective Cartier divisors $D_1, D_2, \ldots, D_r \subset C$ and positive integers $n_1, n_2, \ldots, n_r$ such that $\cC$ is isomorphic to the fiber product of root stacks
\[C[n_1, n_2, \ldots, n_r] := \sqrt[n_1]{D_1/C} \times_{C} \sqrt[n_2]{D_2/C} \times_{C} \ldots \sqrt[n_r]{D_r/C}\]
as in \cite[Section 4.2]{abramovich-graber-vistoli-gw-dm-stacks}.\endnote{By the Keel-Mori theorem \cite[\href{https://stacks.math.columbia.edu/tag/0DUT}{Tag 0DUT}]{stacks-project}, the stack $\mathcal{C}$ has a coarse moduli space $f: \mathcal{C} \to C$, which by the assumptions is separated, of finite type and of dimension 1 over $\Spec(k)$. By applying Nagata compactification \cite{clo-nagata} and \cite[Cor. B.2]{smyth-compactifications}, the coarse space $C$ is seen to be a scheme.

Suppose that $\mathcal{C}$ is a smooth orbifold curve. Then its coarse space $C$ is normal, and hence it is a smooth curve. By the triviality of the generic stabilizer, there are only finitely many closed points $p_1, p_2, \ldots, p_r \in \mathcal{C}(k)$ that have nontrivial stabilizers. The local structure theorem for smooth tame stacks in \cite[Thm. 1.2]{alper-hall-rydh-1} and the assumption that the generic stabilizers are trivial shows that the stabilizer of each $p_i$ must act faithfully on the one-dimensional tangent space at $p_i$, and so it is of the form $\mu_{n_i}$ for some positive integer $n_i$. The image $f(p_i) \in C(k)$ of each point in the smooth coarse space $C$ corresponds to an effective Cartier divisor $D_i$. We denote by 
\[C[n_1, n_2, \ldots, n_r] = \sqrt[n_1]{D_1/C} \times_{C} \sqrt[n_2]{D_2/C} \times_{C} \ldots \sqrt[n_r]{D_r/C}\]
the fiber product of root stacks, as in \cite[Section 4.2]{abramovich-graber-vistoli-gw-dm-stacks}. The proof of \cite[Thm. 4.2.1(1)]{abramovich-graber-vistoli-gw-dm-stacks} applies verbatim to show that there is a canonical isomorphism $\mathcal{C} \cong C[n_1, n_2, \ldots, n_r]$.}

\begin{prop} \label{prop: global presentation orbifold curves}
Let $\mathcal{C}$ be a smooth orbifold curve. Then, there exists a smooth curve $\widetilde{C}$ over $k$ equipped with an action of a finite constant cyclic abelian group $\Gamma= \mu_N$ for some $N$, and a proper tame moduli space morphism of stacks $\widetilde{C}/\Gamma \to \mathcal{C}$.
\end{prop}
\begin{proof}
By the previous discussion, we have $\mathcal{C} = C[n_1, \ldots, n_r]$ for some smooth curve $C$ over $k$. By taking the unique smooth projective model of $C$, we can assume without loss of generality that $C$ is projective. Let $N$ denote the product $N = \prod_i n_i$. Then the canonical morphism $C[N, N, \ldots, N] \to C[n_1, n_2, \ldots, n_r]$ is a proper tame moduli space morphism. By replacing $\mathcal{C}$ with $C[N, N, \ldots, N]$, we can assume that $\mathcal{C}$ is a root stack $\mathcal{C} = \sqrt[N]{D/C}$, where $D = \sum_i D_i$. For any other reduced effective divisor $D'$ disjoint from $D$, we have a proper tame moduli space morphism $\sqrt[N]{(D+D')/C} \to \sqrt[N]{D/C} = \mathcal{C}$. Since $D$ is ample, we can find such reduced effective divisor $D'$ that is linearly equivalent to $(Nh-1)D$ for some $h\gg 0$. After replacing $D$ with $D+D' \equiv NhD$, we can assume that the line bundle $\mathcal{O}(D)$ corresponding to the Cartier divisor $D$ has an $N^{th}$ root. A choice of such root yields an $N^{th}$-cover $\widetilde{C} \to C$ ramified along the divisor $D$. By construction $\widetilde{C}$ is equipped with an action of $\Gamma = \mu_N$, and the choice of $N^{th}$-root of the line bundle $\mathcal{O}(D)$ induces an isomorphism $\sqrt[N]{D/C} \cong \widetilde{C}/\Gamma$, thus concluding the proof.
\end{proof}

\subsection{Mapping stacks from orbifold curves}
We keep our assumption that $\mathcal{Y} = S = \Spec(k)$ for an algebraically closed field $k$ of characteristic $0$. We fix a reductive group $G$ over $k$, and we shall assume that $X= A_X$ is an affine scheme of finite type over $k$ equipped with a $G$-action.
\begin{defn}
Let $\mathcal{C}$ be a smooth proper orbifold curve. We denote by $\mathcal{M}_{\mathcal{C}}(X/G)$ the mapping stack $\Map_k(\mathcal{C}, X/G)$. We define $\mathcal{M}_{\mathcal{C}}(X/G)^{rep} \subset \mathcal{M}_{\mathcal{C}}(X/G)$ to be substack consisting of morphisms $\mathcal{C} \to X/G$ that are representable.
\end{defn}

\begin{lem} \label{lemma: algebraicity of orbifold mapping stack}
The stack $\mathcal{M}_{\mathcal{C}}(X/G)$ is an algebraic stack with affine diagonal and locally of finite type over $k$. The subfunctor $\mathcal{M}_{\mathcal{C}}(X/G)^{rep} \subset \mathcal{M}_{\mathcal{C}}(X/G)$ is represented by an open and closed substack of $\mathcal{M}_{\mathcal{C}}(X/G)$.
\end{lem}
\begin{proof}
The algebraicity and properties of $\mathcal{M}_{\mathcal{C}}(X/G)$ follow from \cite[Thm.~1.2]{hall-rydh-tannakahom} or \cite[Thm.~5.1.1]{halpernleistner2019mapping}. Furthermore, $X/G \cong Y/\GL_N$ for some affine $\GL_N$-scheme $Y$, so it suffices to prove the claim for $G=\GL_N$.\endnote{Since $G$ is an affine algebraic group over the field $k$, there exists a closed embedding $G \hookrightarrow \GL_N$ for some $N$. The quotient $Y:= \GL_N \times^G X$ is an affine scheme of finite type over $k$ equipped with an action of $\GL_N$, and we have an isomorphism of stacks $X/G \cong Y/\GL_N$.}

Let $T$ be a $k$-scheme and choose a $T$-point $T \to \mathcal{M}_{\mathcal{C}}(X/\GL_N)$, corresponding to a pair $(E,s)$, where $E$ is a vector bundle of rank $N$ on $\mathcal{C}_T$ (thought of as a $\GL_N$-bundle) and $s: \mathcal{C}_T \to E(X)$ is a section of the associated fiber bundle $E(X) \to \mathcal{C}_T$ with fiber $X$. Let $\mathcal{I}_{\mathcal{C}_T} \to \mathcal{C}_T$ denote the inertia stack of $\mathcal{C}_T$, which is a finite unramified relatively representable group stack over $\mathcal{C}_T$. There is a homomorphism of relatively affine $\mathcal{C}_T$-group stacks $\varphi: \mathcal{I}_{\mathcal{C}_T} \to \uAut(E)$, where $\uAut(E) \to \mathcal{C}_T$ is the reductive group stack classifying automorphisms of $E$. The fiber product we are interested in $F:=T \times_{\mathcal{M}_{\mathcal{C}}(X/\GL_N)} \mathcal{M}_{\mathcal{C}}(X/\GL_N)^{rep}$ classifies morphisms $S \to T$ such that the base-change $\varphi_S$ is a monomorphism. The inertia stack $\mathcal{I}_{\mathcal{C}}$ has nontrivial fibers over finitely many closed points $p_1, p_2,\ldots, p_r \in \mathcal{C}(k)$. At each $p_i$, the stabilizer $(\mathcal{I}_{\mathcal{C}})_{p_i}$ is a finite constant group $\Gamma_i$. By pulling back the homomorphism $\varphi$ to the base-change $(p_i)_T: T \to \mathcal{C}_T$, we get $r$ different homomorphisms $\varphi_i: (\Gamma_i)_T \to \Aut(E|_{p_i})$, where the automorphism group of the restriction $\Aut(E|_{p_i})$ is a reductive group scheme over $T$. Then the fiber product $F$ is the subfunctor of $T$ classifying morphisms $S \to T$ such that for all $i$ the base-change $(\varphi_i)_S: (\Gamma_i)_S \to \Aut(E|_{p_i})_S$ is a monomorphism.  Since $(\Gamma_i)_T$ is a linearly reductive constant group scheme, $\varphi_i$ induces an isotypic decomposition $E|_{p_i} = \bigoplus_{\chi} (E|_{p_i})_{\chi}$ where $\chi$ ranges over the irreducible representations of $\Gamma_i$. The property of a base-change of $\varphi_i$ being a monomorphism amounts to determining which of the pullbacks of the isotypic vector subbundles $(E|_{p_i})_{\chi}$ have rank $0$. This property is locally constant in families, and so $F$ is represented by an open and closed subscheme of $T$.
\end{proof}

Let $\widetilde{C}/\Gamma \to \mathcal{C}$ be a proper tame moduli space morphism as in \Cref{prop: global presentation orbifold curves}. There is an induced morphism $\mathcal{M}_{\mathcal{C}}(X/G) \to \mathcal{M}_{\widetilde{C}/\Gamma}(X/G)$. We shall use the following lemma to reduce the study of $\mathcal{M}_{\mathcal{C}}(X/G)$ to the case when $\mathcal{C}$ is a quotient stack $\widetilde{C}/\Gamma$.
\begin{lem} \label{lemma: maps from orbifold stack vs maps from quotient stack}
The morphism $\mathcal{M}_{\mathcal{C}}(X/G) \to \mathcal{M}_{\widetilde{C}/\Gamma}(X/G)$ is an open and closed immersion.
\end{lem}
\begin{proof}
As in the proof of \Cref{lemma: algebraicity of orbifold mapping stack}, we can assume without loss of generality that $G$ is the general linear group $\GL_N$ for some $N$.

We can factor the morphism $\mathcal{M}_{\mathcal{C}}(X/\GL_N) \to \mathcal{M}_{\widetilde{C}/\Gamma}(X/\GL_N)$ as follows
\[ \mathcal{M}_{\mathcal{C}}(X/\GL_N) \xrightarrow{f} \Bun_{\GL_N}(\mathcal{C}) \times_{\Bun_{\GL_N}(\widetilde{C}/\Gamma)} \mathcal{M}_{\widetilde{C}/\Gamma}(X/\GL_N) \xrightarrow{g} \mathcal{M}_{\widetilde{C}/\Gamma}(X/\GL_N)\]
We prove that $f$ is an isomorphism and $g$ is an open and closed immersion., thus concluding the proof of the lemma.

\medskip

\noindent \textit{Proof that $f$ is an isomorphism.}
Let $T$ be a Noetherian $k$-scheme, and choose a $T$-point $T \to \Bun_{\GL_N}(\mathcal{C}) \times_{\Bun_{\GL_N}(\widetilde{C}/\Gamma)} \mathcal{M}_{\widetilde{C}/\Gamma}(X/\GL_N)$. This corresponds to a pair $(F, s)$, where $F$ is a vector bundle of rank $N$ on $\mathcal{C}_T$ (thought of as a $\GL_N$-bundle) with pullback $F|_{(\widetilde{C}/\Gamma)_T}$ and $s: (\widetilde{C}/\Gamma)_T \to (F|_{(\widetilde{C}/\Gamma)_T})(X)$ is a section of the corresponding fiber bundle $(F|_{(\widetilde{C}/\Gamma)_T})(X) \to (\widetilde{C}/\Gamma)_T$ with fiber $X$. There is a Cartesian diagram
\[\begin{tikzcd}
 (F|_{(\widetilde{C}/\Gamma)_T})(X)\ar[d] \ar[r] & (\widetilde{C}/\Gamma)_T \ar[d]\\ F(X) \ar[r] & \mathcal{C}_T
\end{tikzcd}\]
The fiber product $Z: = T \times_{(\Bun_{\GL_N}(\mathcal{C}) \times_{\Bun_{\GL_N}(\widetilde{C}/\Gamma)} \mathcal{M}_{\widetilde{C}/\Gamma}(X/\GL_N))} \mathcal{M}_{\mathcal{C}}(X/\GL_N)$ classifies sections $\overline{s}: \mathcal{C}_T \to F(X)$ such that the corresponding morphism $(\widetilde{C}/\Gamma)_T \to (F|_{(\widetilde{C}/\Gamma)_T})(X)$ induced by the universal property of the Cartesian diagram is $s$. We need to show that there exists a unique such $\overline{s}$. By \'etale descent, we can check this after base-changing to an \'etale cover $U \to \mathcal{C}_T$. Then the base-change $(\widetilde{C}/\Gamma)_U \to U$ is a good moduli space by \cite[Prop. 3.10(vii)]{alper-good-moduli}, and so existence and uniqueness follows from the universal property of maps from the good moduli space $(\widetilde{C}/\Gamma)_U \to U$ into the algebraic space $F(X)$ \cite[Thm. 6.6]{alper-good-moduli}.

\medskip

\noindent \textit{Proof that $g$ is an open and closed immersion.}
Since $g: \mathcal{M}_{\widetilde{C}/\Gamma}(X/\GL_N) \to \mathcal{M}_{\widetilde{C}/\Gamma}(X/\GL_N)$ is a base-change of $\Bun_{\GL_N}(\mathcal{C}) \to \Bun_{\GL_N}(\widetilde{C}/\Gamma)$, it suffices to show that this latter morphism is an open and closed immersion. Let $T$ be a $k$-scheme and choose a $T$-point $T \to \Bun_{\GL_N}(\widetilde{C}/\Gamma)$, corresponding to a vector bundle $E$ of rank $N$ on $(\widetilde{C}/\Gamma)_T$. The fiber product $T \times_{\Bun_{\GL_N}(\widetilde{C}/\Gamma)} \Bun_{\GL_N}(\mathcal{C})$ classifies pairs $(F, \psi)$ of a vector bundle $F$ on $\mathcal{C}_T$ along with an isomorphism $\psi: F|_{(\widetilde{C}/\Gamma)_T} \xrightarrow{\sim} E$. Since $(\widetilde{C}/\Gamma)_T \to \mathcal{C}_T$ is a good moduli space morphism, there exists at most one such pair $(\mathcal{F}, \psi)$ up to unique isomorphism, and it exists if and only if the relative inertia stack $\mathcal{I}_{(\widetilde{C}/\Gamma)_T / \mathcal{C}_T}$ acts trivially on the fibers of $E$ \cite[Thm. 10.3]{alper-good-moduli}. The relative inertia stack $\mathcal{I}_{(\widetilde{C}/\Gamma) / \mathcal{C}}$ is unramified and has nontrivial fibers over finitely many closed points $p_1, p_2,\ldots, p_r \in (\widetilde{C}/\Gamma)(k)$, and at each $p_i$, the stabilizer $(\mathcal{I}_{(\widetilde{C}/\Gamma) / \mathcal{C}})_{p_i}$ is a subgroup $\Gamma_i \subset \Gamma$. There is a homomorphism of relatively affine $(\widetilde{C}/\Gamma)_T$-group stacks $\varphi: \mathcal{I}_{(\widetilde{C}/\Gamma)_T / \mathcal{C}_T} \to \uAut(E)$, where $\uAut(E) \to (\widetilde{C}/\Gamma)_T$ is the reductive group stack classifying automorphisms of $E$. Now the same argument as in \Cref{lemma: algebraicity of orbifold mapping stack} shows that $T \times_{\Bun_{\GL_N}(\widetilde{C}/\Gamma)} \Bun_{\GL_N}(\mathcal{C})$ (which is the locus where $\varphi$ is a monommorphism) is represented by an open and closed subscheme of $T$. 

\end{proof}

The morphism $\widetilde{C} \to \widetilde{C}/\Gamma$ induces a map $\mathcal{M}_{\widetilde{C}/\Gamma}(X/G) \to \mathcal{M}_{\widetilde{C}}(X/G)$. The action of $\Gamma$ on $\widetilde{C}$ induces an action of $\Gamma$ on the stack $\mathcal{M}_{\widetilde{C}}(X/G)$ given by precomposition. By \'etale descent, the morphism $\mathcal{M}_{\widetilde{C}/\Gamma}(X/G) \to \mathcal{M}_{\widetilde{C}}(X/G)$ witnesses $\mathcal{M}_{\widetilde{C}/\Gamma}(X/G)$ as the fixed point stack $\mathcal{M}_{\widetilde{C}/\Gamma}(X/G)^{\Gamma}$ as in \cite{romagny-actions-stacks}.
\begin{prop} \label{prop: fixed point morphism for orbifolds}
The morphism $\mathcal{M}_{\widetilde{C}/\Gamma}(X/G) \to \mathcal{M}_{\widetilde{C}}(X/G)$ is affine and of finite type.
\end{prop}
\begin{proof}
This is a direct consequence of \cite[Pro. 3.7]{romagny-actions-stacks}, since $\mathcal{M}_{\widetilde{C}}(X/G)$ has affine diagonal.

We also present a hands-on proof, since we will need it for our construction of the affine grassmannian. We can factor the morphism $\mathcal{M}_{\widetilde{C}/\Gamma}(X/G) \to \mathcal{M}_{\widetilde{C}}(X/G)$ as
\[ \mathcal{M}_{\widetilde{C}/\Gamma}(X/G) \xrightarrow{f} \Bun_{\GL_N}(\widetilde{C}/\Gamma) \times_{\Bun_{\GL_N}(\mathcal{C})} \mathcal{M}_{\widetilde{C}/\Gamma}(X/G) \xrightarrow{g}  \mathcal{M}_{\widetilde{C}}(X/G) \]
We show that $f$ is a closed immersion and $g$ is affine and of finite type, thus concluding the proof.

\medskip

\noindent \textit{$f$ is a closed immersion}. Let $T$ be a scheme over $k$ and choose a point $T \to \Bun_{\GL_N}(\widetilde{C}/\Gamma) \times_{\Bun_{\GL_N}(\mathcal{C})} \mathcal{M}_{\widetilde{C}/\Gamma}(X/G)$ corresponding to a pair $(E, s)$ where $E$ is a rank $N$ $\Gamma$-equivariant vector bundle on $\widetilde{C}_T$ (thought of as a $\GL_N$-bundle on $(\widetilde{C}/\Gamma)_T$), and $s$ is a section $s: \widetilde{C}_T \to E(X)$ of the associated fiber bundle $E(X) \to \widetilde{C}_T$ with fiber $X$. The fiber product $Z:= T \times_{\Bun_{\GL_N}(\widetilde{C}/\Gamma) \times_{\Bun_{\GL_N}(\mathcal{C})} \mathcal{M}_{\widetilde{C}/\Gamma}(X/G)} \mathcal{M}_{\widetilde{C}/\Gamma}(X/G)$ is the subfunctor of $T$ consisting of all morphisms $B \to T$ such that the base change $s_B: \widetilde{C}_B \to E_B(X)$ is $\Gamma$-equivariant. In this case, being equivariant for the finite constant group scheme $\Gamma$ amounts to the equality of finitely many sections $\gamma \cdot s: \widetilde{C}_T \to E(X)$ for $\gamma \in \Gamma$ obtained by using the action of $\Gamma$ on $\widetilde{C}_T$ and $E$. Since the functor of sections $\Sect_{\widetilde{C}}(E(X))$ is represented by a separated (in fact relatively affine) scheme over $T$ \cite[Thm. 2.3 (i)]{hall-rydh-hilbert-quot}, the equality of these sections is a closed condition, and so the fiber product $Z \to T$ is a closed immersion, as desired.

\medskip

\noindent \textit{$g$ is affine and of finite type} (cf. the proof in \cite[Prop. 5.4.2.4]{lurie-gaitsgory-tamagawabook}). Since $g: \Bun_{\GL_N}(\widetilde{C}/\Gamma) \times_{\Bun_{\GL_N}(\mathcal{C})} \mathcal{M}_{\widetilde{C}/\Gamma}(X/G) \to \mathcal{M}_{\widetilde{C}}(X/G)$ is a base-change of $\Bun_{G}(\widetilde{C}/\Gamma) \to \Bun_{G}(\widetilde{C})$, it suffices to show that this latter morphism is affine and of finite type. Let $T$ be a scheme over $k$, and choose a $T$-point $T \to \Bun_{G}(\widetilde{C})$ corresponding to a $G$-bundle $E$ on $\widetilde{C}_T$. The fiber product $T \times_{\Bun_{G}(\widetilde{C})} \Bun_{G}(\widetilde{C}/\Gamma)$ classifies $\Gamma$-equivariant structures on $E$. For each $\gamma \in \Gamma$, let $\gamma^*E$ denote the pullback of $E$ under the action of $\gamma$ on $\widetilde{C}_T$. Let $\prod_{\gamma} \uIso(\gamma^*E,E) \to \widetilde{C}_T$ be the product of the corresponding relatively affine schemes of isomorphims over $\widetilde{C}_T$. For any $\Gamma$-equivariant structure on $E$, the corresponding cocycle yields a section $\psi:\widetilde{C}_T \to  \prod_{\gamma} \uIso(\gamma^*E,E) $. The condition of $\psi$ being a cocycle amounts to $\psi$ factoring through a closed subscheme $Z \subset \prod_{\gamma} \uIso(\gamma^*E,E)$. In other words, the set of $\Gamma$-equivariant structures on $E$ is in bijection with the set of sections $\psi: \widetilde{C}_T \to Z$ of the affine morphism of finite type $Z \to \widetilde{C}_T$. The functor of such sections is represented by a relatively affine scheme \cite[Thm. 2.3 (i)]{hall-rydh-hilbert-quot} of finite type \cite[Thm. 1.1]{olsson-homstacks} over $T$, as desired.
\end{proof}

\begin{notn} \label{notn: numerical invariants on orbifold mapping stacks}
Let $\mathcal{C}$ be a proper smooth orbifold curve. Choose $\widetilde{C}/\Gamma \to \mathcal{C}$ as in \Cref{prop: global presentation orbifold curves}. Let $V$ be an almost faithful representation of $G$, and let $b$ be a rational quadratic norm on graded points of $BG$. Let $\mu = -\wt(\cD(V))/\sqrt{b}$ denote the corresponding numerical invariant on the stack $\mathcal{M}_{\widetilde{C}}(X/G)$. We also denote by $\mu = -\wt(\cD(V))/\sqrt{b}$ the pullback of the numerical invariant $\mu$ under the representable morphism $\mathcal{M}_{\widetilde{C}/\Gamma}(X/G) \to \mathcal{M}_{\widetilde{C}}(X/G)$, and also write $\mu$ for the corresponding numerical invariant on the open and closed substack $\mathcal{M}_{\mathcal{C}}(X/G) \subset \mathcal{M}_{\widetilde{C}/\Gamma}(X/G)$.
\end{notn}
The numerical invariant $\mu$ on $\mathcal{M}_{\mathcal{C}}(X/G)$ defined above depends not only on the choice of $V$ and $b$, but also on the choice of ``cover'' $\widetilde{C}/\Gamma \to \mathcal{C}$.

\subsection{Monotonicity and HN boundedness}
We keep the same assumptions as in the previous subsection.
\begin{prop} \label{prop: montonicity for orbifolds}
Let $\mu$ be a numerical invariant on $\mathcal{M}_{\mathcal{C}}(X/G)$ defined as in \Cref{notn: numerical invariants on orbifold mapping stacks}. Then $\mu$ is strictly $\Theta$-monotone and strictly $S$-monotone.
\end{prop}

The proof of \Cref{prop: montonicity for orbifolds} is analogous to the proof of monotonicity in \Cref{prop: monotonicity of mapping stack general group}. As in \Cref{subsection: affine grassmannian}, we shall introduce a version of affine grassmannians for $\mathcal{M}_{\widetilde{C}/\Gamma}(X/G)$ in the special case when $G= \GL_N \times H$, where $H$ is an algebraic group over $k$ fitting into a short exact sequence
\[ 1 \to T \to H \to F \to 1\]
with $T$ a torus and $F$ finite \'etale. 

\begin{defn}
The stack of rational maps $\mathcal{M}_{\widetilde{C}/\Gamma}^{rat}(X/(\GL_N\times H))$ is the pseudofunctor that sends any $S$-scheme $T \in \Aff_{S}$ to the groupoid of triples $(\mathcal{E} \times \cH, D, s)$, where
\begin{enumerate}
    \item $D$ is a $T$-flat $\Gamma$-stable effective Cartier divisor on $\widetilde{C}_T$.
    \item $\mathcal{E} \times \cH$ is $\GL_N \times H$-bundle on $\widetilde{C}_T$ consisting of a vector bundle $\mathcal{E}$ of rank $N$ and an $H$-bundle $\cH$. 
    \item Both the $H$-bundle $\cH$ and the restriction $\mathcal{E}|_{\widetilde{C}_T\setminus D}$ of the vector bundle are equipped with $\Gamma$-equivariant structures.
    \item $s: \widetilde{C}_{T} \setminus D \to (\mathcal{E} \times \cH)(X)$ is an $\Gamma$-equivariant section of $(\mathcal{E} \times \cH)(X) \to \widetilde{C}_{T}$ defined away from $D$.
\end{enumerate}
A morphism from $(\mathcal{E}_1\times \cH_1, D_1, s_1)$ into $(\mathcal{E}_2\times \cH_2, D_2, s_2)$ is the data of equality $D_1 = D_2$ and a pair $(\phi, \psi)$ where
\begin{enumerate}
\item $\psi$ is a $\Gamma$-equivariant isomorphism of $H$-bundles $\psi: \cH_2 \xrightarrow{\sim} \cH_1$.
\item $\phi$ is a $\Gamma$-equivariant isomorphism of vector bundles $\phi: \mathcal{E}_1|_{\widetilde{C}_{T} \setminus D_1} \to \mathcal{E}_2|_{\widetilde{C}_{T} \setminus D_2}$ that is compatible with the sections $s_1$ and $s_2$.
\end{enumerate}
\end{defn}
We denote by $\Div^{\Gamma}(\widetilde{C})$ the scheme classifying $\Gamma$-stable relative Cartier divisors in $\widetilde{C}$. There is a natural morphism
\[
q: \mathcal{M}_{\widetilde{C}/\Gamma}(X/(\GL_N \times H)) \times \Div^{\Gamma}(\widetilde{C}) \to \mathcal{M}_{\widetilde{C}/\Gamma}^{rat}(X/(\GL_N\times H)), \; \; \; (\mathcal{E} \times \cH, s, D) \mapsto (\mathcal{E}\times \cH, D, s|_{\widetilde{C} \setminus D})
\]

\begin{defn}
Let $T$ be an $S$-scheme of finite type, and let $m: T \to \mathcal{M}_{\widetilde{C}/\Gamma}^{rat}(X/(\GL_N\times H))$ be a morphism represented by a triple $(\mathcal{E}\times \cH, D, s)$. We define the generalized BD Grassmanian ${}^\Gamma\Gr_{X/(\GL_N \times H)}^{(\mathcal{E}\times \cH, D, s)}$ to be the fiber product $(\mathcal{M}_{\widetilde{C}/\Gamma}(X/(\GL_N \times H))) \times \Div^\Gamma(D)) \times_{\mathcal{M}_{\widetilde{C}/\Gamma}^{rat}(X/(\GL_N\times H))} T$.
\end{defn}
 By definition, ${}^\Gamma\Gr_{X/(\GL_N \times H)}^{(\mathcal{E}\times \cH, D, s)}$ is the pseudofunctor that sends $Q \in \Aff_{T}$ to the set of isomorphism classes of triples $(\mathcal{F}, \phi, s)$ where
\begin{enumerate}
    \item $\mathcal{F}$ is a $\Gamma$-equivariant vector bundle of rank $N$ on $\widetilde{C}_{Q}$.
    \item $\phi$ is a $\Gamma$-equivariant isomorphism $\phi: \mathcal{E}|_{\widetilde{C}_{Q} \setminus D_{Q}}\xrightarrow{\sim} \mathcal{F}|_{\widetilde{C}_{Q} \setminus D_{Q}}$.
    \item $\widetilde{s}: \widetilde{C}_{Q} \to (\mathcal{F} \times \cH|_{\widetilde{C}_{Q}})(X)$ is an extension of the section
    \[ \widetilde{C}_{Q} \setminus D_{Q} \xrightarrow{s_{Q}} (\mathcal{E} \times \cH)|_{\widetilde{C}_{Q} \setminus D_{Q}}(X) \xrightarrow{\sim} (\mathcal{F} \times \cH)|_{\widetilde{C}_{Q} \setminus D_{Q}}(X)\]
with the last isomorphism induced by the identification $\phi: \cE|_{\widetilde{C}_{Q} \setminus D_{Q}} \xrightarrow{\sim} \cF|_{\widetilde{C}_{Q} \setminus D_{Q}}$ \endnote{Note that such extension $\widetilde{s}$ will automatically be $\Gamma$-equivariant by the density of $\widetilde{C}_Q \setminus D_Q \subset \widetilde{C}_Q$.}.
\end{enumerate}
There is a projection morphism ${}^\Gamma\Gr_{X/(\GL_N \times H)}^{(\mathcal{E}\times \cH, D, s)} \to \mathcal{M}_{\widetilde{C}/\Gamma}(X/(\GL_N \times H))$ that sends $(\mathcal{F}, \phi,\widetilde{s})$ to $(\mathcal{F}\times (\cH|_{C_{Q}}), \widetilde{s})$.
There is also a forgetful morphism $\mathcal{M}_{\widetilde{C}/\Gamma}^{rat}(X/(\GL_N\times H)) \to \mathcal{M}_{\widetilde{C}}^{rat}(X/(\GL_N\times H))$ that forgets the $\Gamma$-equivariance. This induces a forgetful morphism of affine grassmannians ${}^\Gamma\Gr_{X/(\GL_N \times H)}^{(\mathcal{E}\times \cH, D, s)} \to \Gr_{X/(\GL_N \times H)}^{(\mathcal{E}\times \cH, D, s)} $ that forgets $\Gamma$-equivariant structures. These morphisms fit into a commutative diagram
\[
\begin{tikzcd}
 {}^\Gamma\Gr_{X/(\GL_N \times H)}^{(\mathcal{E}\times \cH, D, s)} \ar[d] \ar[r] & \mathcal{M}_{\widetilde{C}/\Gamma}(X/(\GL_N \times H)) \ar[d]\\ \Gr_{X/(\GL_N \times H)}^{(\mathcal{E}\times \cH, D, s)} \ar[r] & \mathcal{M}_{\widetilde{C}}(X/(\GL_N \times H))
\end{tikzcd}
\]
\begin{lem} \label{lemma: forgetful morphism of equivariant affine grassmannian is closed}
For any $m: T \to \mathcal{M}_{\widetilde{C}/\Gamma}^{rat}(X/(\GL_N\times H))$ represented by a triple $(\mathcal{E}\times \cH, D, s)$, the forgetful morphism ${}^\Gamma\Gr_{X/(\GL_N \times H)}^{(\mathcal{E}\times \cH, D, s)} \to \Gr_{X/(\GL_N \times H)}^{(\mathcal{E}\times \cH, D, s)} $ is represented by closed immersions.
\end{lem}
\begin{proof}
Let $Q$ be a $T$-scheme, and choose a point $Q \to \Gr_{X/(\GL_N \times H)}^{(\mathcal{E}\times \cH, D, s)}$ represented by a triple $(\mathcal{F}, \phi, \widetilde{s})$. The isomorphism $\phi: \mathcal{E}|_{\widetilde{C}_Q \setminus D_Q} \xrightarrow{\sim} \mathcal{F}|_{\widetilde{C}_Q \setminus D_Q}$ induces a $\Gamma$-equivariant structure on $\mathcal{F}|_{\widetilde{C}_Q \setminus D_Q}$. By definition, the fiber product $Q \times_{\Gr_{X/(\GL_N \times H)}^{(\mathcal{E}\times \cH, D, s)}} {}^\Gamma\Gr_{X/(\GL_N \times H)}^{(\mathcal{E}\times \cH, D, s)}$ classifies extensions to $\widetilde{C}_Q$ of the $\Gamma$-equivariant structure on $\mathcal{F}|_{\widetilde{C}_Q \setminus D_Q}$. We have seen in the proof of \Cref{prop: fixed point morphism for orbifolds} that such equivariant structures correspond to sections $\widetilde{C}_Q \to Z$ of an affine morphism of finite type $Z \to \widetilde{C}_Q$. The equivariant structure on $\mathcal{F}|_{\widetilde{C}_Q \setminus D_Q}$ yields a section $\psi: \widetilde{C}_Q \setminus D_Q \to Z$ defined over $\widetilde{C}_Q \setminus D_Q$, and the fiber product we are interested is the functor classifying extensions $\widetilde{\psi}$ of the section $\psi$ to the whole $\widetilde{C}_Q$. By \cite[Lemma 4.12]{rho-sheaves-paper} this functor is represented by a closed subscheme of $Q$.
\end{proof}
\begin{proof}[Proof of \Cref{prop: montonicity for orbifolds}]
Since $\mathcal{M}_{\mathcal{C}}(X/G)$ is an open and closed substack of $\mathcal{M}_{\widetilde{C}/\Gamma}(X/G)$, it suffices to show strict monotonicity for the corresponding numerical invariant on $\mathcal{M}_{\widetilde{C}/\Gamma}(X/G)$, and so we can assume that $\mathcal{C}$ is a quotient stack of the form $\widetilde{C}/\Gamma$. The same reductions as in the proof of \Cref{prop: monotonicity of mapping stack general group} allow us to reduce to the case when $G$ is of the form $\GL_N \times H$ as above. Indeed, the proof of the main lemma used in the reduction (\Cref{lemma: finiteness of map to isogenous group}) goes through verbatim, using \Cref{lemma: bun_g isogeny orbifold case} instead of \Cref{lemma: finiteness stack of bundles under isogeny}.

In the case $G=\GL_N \times H$ the proof of \Cref{prop: monotonicity of mapping stack general group} goes through word for word, by working with the curve $\widetilde{C}$, noting that the Cartier divisor from the rational filling condition can always be arranged to be $\Gamma$-stable by taking the sum of translations $\sum_{\gamma \in \Gamma} \gamma\cdot D$, and using the fact that the corresponding affine grassmannian for $\mathcal{M}_{\widetilde{C}/\Gamma}(X/\GL_N\times H)$ is closed inside the affine grassmannian for $\mathcal{M}_{\widetilde{C}}(X/\GL_N \times H)$ (\Cref{lemma: forgetful morphism of equivariant affine grassmannian is closed}).
\end{proof}

\begin{prop} \label{prop: hn boundedness for orbifolds}
Any numerical invariant $\mu$ on $\mathcal{M}_{\mathcal{C}}(X/G)$ as in \Cref{notn: numerical invariants on orbifold mapping stacks} satisfies the HN boundedness property.
\end{prop}
\begin{proof}
Since $\mathcal{M}_{\mathcal{C}}(X/G)$ is an open and closed substack of $\mathcal{M}_{\widetilde{C}/\Gamma}(X/G)$, it suffices to show HN boundedness for the corresponding numerical invariant $\mu$ on $\mathcal{M}_{\widetilde{C}/\Gamma}(X/G)$. The affine morphism of finite type $\varphi: \mathcal{M}_{\widetilde{C}/\Gamma}(X/G) \to \mathcal{M}_{\widetilde{C}}(X/G)$ in \Cref{prop: fixed point morphism for orbifolds} exhibits $\mathcal{M}_{\widetilde{C}/\Gamma}(X/G)$ as the stack of $\Gamma$-fixed points of $\mathcal{M}_{\widetilde{C}}(X/G)$. As in the end of the proof of \Cref{prop: reduction to connected groups}, \'etale descent implies the set of isomorphism classes of filtrations of any fixed geometric point $x \in \mathcal{M}_{\widetilde{C}/\Gamma}(X/G)$ corresponds to the subset of isomorphism classes of filtrations of $\varphi(x)$ fixed by $\Gamma$. We can then apply the same argument as in the end of the proof of \Cref{prop: reduction to connected groups} to deduce the HN boundedness statement for $\mathcal{M}_{\widetilde{C}/\Gamma}(X/G)$ from the HN boundedness for $\mathcal{M}_{\widetilde{C}}(X/G)$ proven in \Cref{thm: HN-boundedness nonconnected}.
\end{proof}

\subsection{\texorpdfstring{$\Theta$}{Theta}-stratification and moduli space for orbifold curves}
Let $\mathcal{C}$ be a smooth proper orbifold curve over $k$ along with the choice of a proper tame moduli space morphism $\widetilde{C}/\Gamma \to \mathcal{C}$ as in \Cref{prop: global presentation orbifold curves}.
\begin{thm} \label{thm: theta-stratification orbifolds}
Let $\mu$ be a numerical invariant on $\mathcal{M}_{\mathcal{C}}(X/G)$ defined as in \Cref{notn: numerical invariants on orbifold mapping stacks}. Then $\mu$ defines a $\Theta$-stratification on $\mathcal{M}_{\mathcal{C}}(X/G)$.
\end{thm}
\begin{proof}
This is a consequence of \Cref{thm: theta stability paper theorem}, because $\mu$ is strictly-$\Theta$-monotone (\Cref{prop: montonicity for orbifolds}) and satisfies HN boundedness (\Cref{prop: hn boundedness for orbifolds}).
\end{proof}

\begin{defn}
For any $d \in H_2(BG)$, we define $\mathcal{M}_{\mathcal{C}}(X/G)_d\subset \mathcal{M}_{\mathcal{C}}(X/G)$ to be the open and closed substack given by the preimage $\varphi^{-1}(\mathcal{M}_{\widetilde{C}}(X/G)_d)$ under the morphism $\varphi: \mathcal{M}_{\mathcal{C}}(X/G) \to \mathcal{M}_{\widetilde{C}/\Gamma}(X/G) \to \mathcal{M}_{\widetilde{C}}(X/G)$ (cf. \Cref{lemma: maps from orbifold stack vs maps from quotient stack} and \Cref{prop: fixed point morphism for orbifolds}).
\end{defn}

\begin{prop} \label{prop: boundedness strata orbifolds}
Let $\mu$ be a numerical invariant on $\mathcal{M}_{\mathcal{C}}(X/G)$ defined as in \Cref{notn: numerical invariants on orbifold mapping stacks}. For any $d \in H_2(BG)$ and $\gamma \in \mathbb{R}$, the open substack $\mathcal{M}_{\mathcal{C}}(X/G)_d^{\mu \leq \gamma}$ of points whose HN filtration has numerical invariant $\leq \gamma$ is of finite type over $k$.
\end{prop}
\begin{proof}
The stack $\mathcal{M}_{\mathcal{C}}(X/G)_d^{\mu \leq \gamma}$ is locally of finite type over $k$; we need to show that it is quasi-compact. Since the composition $\varphi: \mathcal{M}_{\mathcal{C}}(X/G)_d \to \mathcal{M}_{\widetilde{C}/\Gamma}(X/G)_d \to \mathcal{M}_{\widetilde{C}}(X/G)_d$ is affine and of finite type, and since $\mathcal{M}_{\widetilde{C}}(X/G)_d^{\mu \leq  \gamma}$ is quasi-compact (see the proof of \Cref{cor: boundedness of strate affine case}), it suffices to show that $\mathcal{M}_{\mathcal{C}}(X/G)^{\mu \leq \gamma}_d$ is contained in the preimage $\varphi^{-1}(\mathcal{M}_{\widetilde{C}}(X/G)_d^{\mu \leq  \gamma})$. Suppose for the sake of contradiction that $p$ is a geometric point of $\mathcal{M}_{\mathcal{C}}(X/G)^{\mu \leq \gamma}_d$ such that $\varphi(p)$ is not in $\mathcal{M}_{\widetilde{C}}(X/G)_d^{\mu \leq \gamma}$. This means that $\mu(f) > \gamma$ for an HN filtration $f$ of $p$. By the uniqueness of the HN filtration and the description of filtrations of $p$ as $\Gamma$-fixed points in the degeneration fan of $\varphi(p)$ (see the proof of \Cref{prop: hn boundedness for orbifolds}), there is a filtration $\widetilde{f}$ of $p$ such that $\varphi \circ \widetilde{f} \cong f$ as filtrations of $\varphi(p)$. In this case we have $\mu(\widetilde{f}) = \mu(f) >\gamma$, contradicting the assumption that $p$ was in $\mathcal{M}_{\mathcal{C}}(X/G)_d^{\mu \leq \gamma}$.
\end{proof}

\begin{thm}
Let $\mu$ be a numerical invariant on $\mathcal{M}_{\mathcal{C}}(X/G)$ defined as in \Cref{notn: numerical invariants on orbifold mapping stacks}. For any $d \in H_2(BG)$, the $\mu$-semistable locus $\mathcal{M}_{\mathcal{C}}(X/G)^{\mu \dash \rm{ss}}_d$ admits a separated good moduli space of finite type over $k$.
\end{thm}
\begin{proof}
This follows from \Cref{thm: theta stability paper theorem}, since $\mu$ is strictly $\Theta$-monotone and strictly $S$-monotone (\Cref{prop: montonicity for orbifolds}), it satisfies HN boundedness (\Cref{prop: hn boundedness for orbifolds}), and the open and closed substack $\mathcal{M}_{\mathcal{C}}(X/G)^{\mu \dash \rm{ss}}_d \subset \mathcal{M}_{\mathcal{C}}(X/G)^{\mu \dash \rm{ss}}$ is of finite type over $k$ (\Cref{prop: boundedness strata orbifolds}).
\end{proof}
\appendix
\section{Some results on \texorpdfstring{$\Bun_{G}(X)$}{BunG(C)}}
In this appendix, we establish some results on $\Bun_{G}(X)$, where $X \to S$ is a proper flat morphism and $G$ is a geometrically reductive group scheme over $S$. We will make use of the following lemma.
\begin{lem} \label{lemma: sections of finite morphism finite}
Let $T$ be a Noetherian scheme and let $X \to T$ be a flat proper algebraic stack with geometrically normal fibers over $T$. For any schematic morphism $Y \to X$, let $\Gamma_{X}(Y)$ denote the functor 
\[ (\Sch/T)^{\op{op}} \to \Set, \; \quad \; (T' \to T) \mapsto \{ \text{ sections of the morphism } Y_{T'} \to X_{T'} \; \} \]
If $Y \to X$ is finite, then $\Gamma_{X}(Y)$ is represented by a finite scheme over $T$.
\end{lem}
\begin{proof}
The functor $\Gamma_{X}(Y)$ is represented by a scheme that is affine \cite[Thm. 2.3 (i)]{hall-rydh-hilbert-quot} and of finite type (\cite[Thm. 1.1]{olsson-homstacks}+\cite[Thm. B]{rydh-noetherian-approximation}) over $S$. It remains to show that $\Gamma_{X}(Y)$ is $S$-proper; for this we use the existence part of the valuative criterion for properness, since we already know that the scheme $\Gamma_{X}(Y)$ is separated. Let $R$ be a discrete valuation ring over $S$, and choose a $\Frac(R)$-point in $\Gamma_{X}(Y)$ corresponding to a section $s: X_{\Frac(R)} \to Y_{R}$ defined over the open dense generic fiber $X_{\Frac(R)} \subset X_{R}$. Since $X_{R}$ is flat and has geometrically normal fibers over the discrete valuation ring, $X_{R}$ is itself normal \cite[Prop. 6.14.1]{egaivII}. Let $Z \subset Y_{R}$ be the scheme theoretic image of $s: X_{\Frac(R)} \to Y_{R}$. The composition $\psi: Z \to Y_{R} \to X_{R}$ is a finite birational morphism into the normal stack $X_{R}$, and so it is an isomorphism. The inverse morphism yields a section $X_{R} \xrightarrow{\psi^{-1}} Z \to Y_{R}$ extending the generic section $s$, as desired.
\end{proof}

\begin{lem} \label{lemma: stack of h-bundles is theta red and s complete}
Let $H$ be a flat group scheme over $S$ fitting into an extension sequence
\[ 1 \to T \to H \to F \to 1\]
where $T$ is a torus and $F$ is a finite \'etale group scheme over $S$. Let $X \to S$ be a smooth proper algebraic stack over $S$. Then, $\Bun_{H}(X)$ is $\Theta$-reductive and S-complete.
\end{lem}
\begin{proof}
For the purposes of checking $\Theta$-reductivity and S-completeness, we can assume that $S$ is the spectrum of a complete discrete valuation ring $R$. Since we are allowed to pass to finite extensions of $R$, we can assume that $T \cong \bG_m^r$ for some nonnegative integer $r$. Let $\XX$ denote either $\Theta_{R}$ or $\overline{ST}_{R}$. Choose a morphism $\XX \setminus \closedpt \to \Bun_{H}(X)$ corresponding to a $\mathbb{G}_m$-equivariant $H$-bundle $\cH$ on $C_{Y_{\XX} \setminus \closedpt}$, Since the diagonal of $\Bun_{H}(X)$ is affine, it suffices to show that we can extend $\cH$ to a bundle on $C_{Y_{\XX}}$.

Let $\cF$ denote the $\mathbb{G}_m$-equivariant $F$-bundle on $Y_{\XX} \setminus \closedpt$ determined by the composition $\XX \setminus \closedpt \to \Bun_{H}(X) \to \Bun_{F}(X)$. Lemma \ref{lemma: sections of finite morphism finite} implies that the stack $\Bun_{F}(X)$ has finite diagonal.\endnote{Indeed, for any Noetherian $S$-scheme $T$ and any two $F$-bundles $\cF_1, \cF_2$ on $X_{T}$, there is a finite \'etale stack $\c{I}so(\cF_1, \cF_2)$ over $X_{T}$ that parametrizes isomorphisms between $\cF_1$ and $\cF_2$. Then the functor of isomorphisms $\Iso(\cF_1, \cF_2)$ in $\Bun_{F}(X)$ is parametrized by the functor $\Gamma_{X_{T}}(\c{I}so(\cF_1, \cF_2))$, which is represented by a finite $T$-scheme by Lemma  \ref{lemma: sections of finite morphism finite}.} By \cite[Prop. 3.22, 3.44]{alper-good-moduli}, it follows that $\Bun_{F}(X)$ is $\Theta$-reductive and $S$-complete. Therefore, there is a unique extension of $\cF$ to a $\mathbb{G}_m$-equivariant $F$-bundle $\widetilde{\cF}$ on $Y_{\XX}$. The total space $\widetilde{X} \to X_{Y_{\XX}}$ of $\widetilde{\cF}$ is a finite $\mathbb{G}_m$-equivariant \'etale cover of $X$. We have natural morphisms $\Bun_{H}(X_{Y_{\XX}}) \to \Bun_{H}(\widetilde{X})$ and $\Bun_{F}(X_{Y_{\XX}}) \to \Bun_{F}(\widetilde{X})$ given by pulling back bundles, and both of these morphisms are schematic and affine by \'etale descent (cf. the proof in \cite[Prop. 5.4.2.4]{lurie-gaitsgory-tamagawabook}). The composition $Y_{\XX} \to \Bun_{F}(X_{Y_{\XX}}) \to \Bun_{F}(\widetilde{X})$ is a canonically trivialized $F$-bundle on $\widetilde{X}$. Therefore, the fiber product $\Bun_{H}(\widetilde{X}) \times_{\Bun_{F}(\widetilde{X})} Y_{\XX}$ is canonically isomorphic to the stack $\Bun_{T}(\widetilde{X}) \cong Pic(\widetilde{X})^r$. Hence, we can think of the section $Y_{\XX} \setminus \closedpt \to \Bun_{F}(X_{Y_{\XX}}) \to \Bun_{F}(\widetilde{X})$ induced by $\cH$ as a tuple of line bundles on the open substack $\widetilde{X}_{Y_{\XX} \setminus \closedpt}$ with complement of codimension $2$ in the regular stack $\widetilde{X}$. Pushing forward each line bundle in the tuple yields an extension to a section $s: Y_{\XX} \to  Pic(\widetilde{X})^r \cong \Bun_{F}(X_{Y_{\XX}}) \to \Bun_{F}(\widetilde{X})$. Consider the Cartesian diagram
\[
\begin{tikzcd}
  Z \ar[rr] \ar[d] & & \Bun_{H}(X_{Y_{\XX}}) \ar[d]\\ Y_{\XX} \ar[r, "s"] &
  \Bun_{H}(X_{Y_{\XX}}) \times_{\Bun_{F}(X)} Y_{\XX} \ar[r] & \Bun_{H}(\widetilde{X}) \times_{\Bun_{F}(\widetilde{X})} \Bun_{F}(X_{Y_{\XX}})
\end{tikzcd}
\]
The bundle $\cH$ yields a section $Y_{\XX} \setminus \closedpt \to Z$ of the morphism $Z \to Y_{\XX}$. To conclude the proof, we shall show that $Z \to Y_{\XX}$ is affine, which by Hartogs's theorem implies that the section extends to $Y_{\XX}$. Since $\Bun_{F}(X) \to \Bun_{F}(\widetilde{X})$ is affine, we also know that $\Bun_{H}(\widetilde{X}) \times_{\Bun_{F}(\widetilde{X})} \Bun_{F}(X_{Y_{\XX}}) \to \Bun_{H}(\widetilde{X})$ is affine. Using cancellation for the composition
\[ \Bun_{H}(X) \to \Bun_{H}(\widetilde{X}) \times_{\Bun_{F}(\widetilde{X})} \Bun_{F}(X_{Y_{\XX}}) \to \Bun_{H}(\widetilde{X}) \]
 and the affineness of $\Bun_{H}(X) \to \Bun_{H}(\widetilde{X})$ we conclude that $\Bun_{H}(X_{Y_{\XX}}) \to \Bun_{H}(\widetilde{X}) \times_{\Bun_{F}(\widetilde{X})} \Bun_{F}(X_{Y_{\XX}})$ is affine. Therefore, the pullback $Z \to Y_{\XX}$ is also affine. 
\end{proof}
Our next goal is to study the behavior of the stack $\Bun_{G}(X)$ under central isogenies. We will use the following lemma.
\begin{lem} \label{lemma: splitting smooth morphisms into connected components}
Let $X$ be a flat proper tame stack with geometrically normal fibers over a Noetherian integral domain $A$. Then, there exists a finite \'etale extension of integral domains $A \subset B$ such that $X_{B}$ is a disjoint union of open and closed substacks with geometrically integral $B$-fibers.
\end{lem}
\begin{proof}
By passing to the tame moduli space of $X$ (which is flat proper and has geometrically normal fibers over $A$), we can assume without loss of generality that $X$ is a scheme.

We induct on the number $n$ of geometric connected components of the generic fiber $X_{\Frac(A)}$. If $n=1$, then the fibers of the morphism $X \to \Spec(A)$ are geometrically connected \cite[\href{https://stacks.math.columbia.edu/tag/0E0N}{Tag 0E0N}]{stacks-project}. Equivalently, the $A$-fibers are geometrically integral, since we already know that they are geometrically normal. This concludes the base case.

We proceed with the induction step. Since $X$ is proper, $H^0(\cO_{X})$ is finite over $A$. Write the Noetherian ring $H^0(\cO_{X})$ as a finite product $H^0(\cO_{X}) = \prod_i D_i$, where each $D_i$ has connected spectrum. After writing $X$ as a union of open and closed subschemes $X = \sqcup_i X \times_{\Spec(H^0(\cO_{X}))} \Spec(D_i)$, we can reduce to the case when $H^0(\cO_{X})$ has no idempotents. By \cite[Chpt.8, Prop. 8.5.16]{fga-explained}, the algebra $C \vcentcolon = H^0(\cO_X)$ is finite \'etale over $A$ and its formation commutes with base-change over $A$. By $A$-flatness, $C \hookrightarrow C \otimes_{A} \Frac(A)$ is injective. Since the right-hand side is finite over $\Frac(A)$, reduced and connected, we conclude that $C \otimes_{A} \Frac(A)$ is a field. Hence, $C$ must be an integral domain. 

By \cite[\href{https://stacks.math.columbia.edu/tag/03H0}{Tag 03H0}]{stacks-project}, the fibers of the morphism $X \to \Spec(C)$ are geometrically connected. If $C = A$, then we must have $n=1$, and we are in the base case. Suppose otherwise that $C \neq A$. Write $C \otimes_{A} C$ as a finite product $C \otimes_{A} C = \prod_i D_i$ of rings with no nontrivial idempotents. This decomposition is not trivial; there are at least two factors because the relative diagonal $\Spec(C) \to \Spec(C \otimes_{A} C)$ of the finite \'etale morphism $\Spec(C) \to \Spec(A)$ yields a proper nontrivial open and closed subscheme of $\Spec(C \otimes_{A} C)$. For each $i$, we set $X_i \vcentcolon = X_{C} \times_{\Spec(C \otimes_{A} C)} \Spec(D_i)$. This way $X_{C}$ can be expressed as a disjoint union of open and closed subschemes $\sqcup_i X_i$ flat and proper over $C$. In particular, each $X_i$ has nonempty generic fiber $(X_i)_{\Frac(C)}$. Since the number of geometric connected components $n$ of $X_{\Frac(C)}$ is the sum of the number of geometric connected components of all $(X_i)_{\Frac(C)}$, each $(X_i)_{\Frac(C)}$ must have strictly less than $n$ geometric connected components. We conclude by induction.
\end{proof}

\begin{lem} \label{lemma: finiteness stack of bundles under isogeny}
Let $G \to \overline{G}$ be an isogeny of (not necessarily connected) geometrically reductive group schemes over $S$, with kernel $K \hookrightarrow G$ finite and of multiplicative type. Let $X$ be a scheme that is flat, proper and has geometrically normal fibers over $S$. Then, the morphism of stacks $\Bun_{G}(X) \to \Bun_{\overline{G}}(X)$ is quasi-finite and has finite diagonal. If in addition $X \to S$ is smooth, then $\Bun_{G}(X) \to \Bun_{\overline{G}}(X)$ is proper.
\end{lem}
\begin{proof}
Since the claim is \'etale local on the base $S$, we can reduce to the case when $K$ is diagonalizable and both $G$ and $\overline{G}$ are split. This means in particular that $K \cong \prod_{i} \mu_{n_i}$ for some tuple of positive integers $(n_i)$. We shall show all required properties for the morphism $\Bun_{G}(X) \to \Bun_{\overline{G}}(X)$.\\

\noindent $\bullet$ \textit{Finite relative diagonal:} Choose a morphism $T \to \Bun_{G}(X) \times_{\Bun_{\overline{G}}(X)} \Bun_{G}(X)$ represented by a pair of $G$-bundles $\cP_1, \cP_2$ along with an isomorphism of the associated $\overline{G}$-bundles $\overline{\cP_1} \cong \overline{\cP_2}$. Consider the induced morphism of $X_{T}$-schemes $\mathcal{I}so(\cP_1, \cP_2) \to \mathcal{A}ut(\overline{\cP_1})$. The fiber $\mathcal{I}so_{\Bun_{G}(X)/\Bun_{\overline{G}}(X)}(\cP_1, \cP_2)$ of this morphism over the identity section in $\mathcal{A}ut(\overline{\cP_1})$ is \'etale locally isomorphic to $K \times X_{T}$, and hence it is finite and flat over $X_{T}$. The isomorphisms $\cP_1 \xrightarrow{\sim} \cP_2$ relative to $\Bun_{\overline{G}}(X)$ are parametrized by the functor of sections $\Gamma_{X_{T}}(\mathcal{I}so_{\Bun_{G}(X)/\Bun_{\overline{G}}(X)}(\cP_1, \cP_2))$, which by Lemma \ref{lemma: sections of finite morphism finite} is finite over $T$.\\

\noindent $\bullet$ \textit{Quasicompactness:} Let $A$ be a Noetherian integral domain over $S$. Choose a morphism $\Spec(A) \to \Bun_{\overline{G}}(X)$ represented by a $\overline{G}$-bundle $\overline{P}$ on $X_{A}$. It suffices to show that $\Bun_{G}(X) \times_{\Bun_{G}(X)} \Spec(A)$ is quasi-compact. We can view $\overline{P}$ as a morphism $X_{A} \to B\overline{G}$. The fiber product $\cG \vcentcolon = BG \times_{B\overline{G}} X_{A}$ is an \'etale gerbe over $X_{A}$ banded by the group scheme $K^{\overline{P}}$ obtained from $K$ by twisting with the $\overline{G}$-bundle $\overline{P}$ via the conjugation action. Using the terminology of \cite{olsson-boundedness}, we have $\Bun_{G}(X) \times_{\Bun_{G}(X)} \Spec(A) = \Sect(\cG/X_{A})$. By \cite[Thm. B]{rydh-noetherian-approximation}, there exists a scheme $Z$ with a finite surjective morphism $Z \to \cG$. Applying the same reasoning as in \cite[6.4]{olsson-boundedness}, we obtain a Noetherian integral domain $A'$ with a finite type, dominant, generically finite morphism $\Spec(A') \to \Spec(A)$ and a smooth proper $A'$-scheme $W$ equipped with a surjection $W \to X_{A'}$ and a section $W \to \cG$ over $X_{A}$. By Noetherian induction on $A$, we can reduce to showing that $\Sect(\cG_{A'}/X_{A'})$ is quasi-compact. By the proof of \cite[Prop. B.5]{hall-rydh-tannakahom} (with $Z = X_{A'}$, $\cX = BG$ and $\cY = B\overline{G}$), it is sufficient to show that $\Sect(\cG_{A'}/W)$ is quasi-compact. Hence, we have reduced to the case when $\cG \to X_{A}$ admits a section, meaning that the gerbe $\cG$ is trivial. Then the relevant fiber product $\Bun_{G}(X) \times_{\Bun_{G}(X)} \Spec(A)$ is isomorphic to the stack of $K^{\overline{P}}$-bundles $\Bun_{K^{\overline{P}}}(X)$.

We shall actually show that $\Bun_{K^{\overline{P}}}(X)$ is quasi-finite, and that $\Bun_{K^{\overline{P}}}(X)$ is proper when $X \to S$ is smooth. Consider the nonabelian \'etale cohomology class in $H^1_{et}(X_{A}, \Aut(K))$ obtained from the $\overline{G}$-bundle $\overline{\cP}$ via the conjugation action of $\overline{G}$ on $K$. We can equivalently think of this as a cohomology class $H^1_{et}(X_{A}, \Aut(\prod_{i} \bZ/n_i\bZ))$ with values in the automorphisms of the constant Cartier dual group $\prod_{i} \bZ/n_i\bZ$ of the diagonalizable group $K$. The total space of the corresponding $\Aut(\prod_{i} \bZ /N_i \bZ)$-torsor yields a finite \'etale cover $\widetilde{X} \to X_{A}$ such that $K^{\cP} \times_{X_{A}} \widetilde{X}$ is isomorphic to the diagonalizable group scheme $K \times \widetilde{X}$. By \'etale descent (c.f. the proof of \cite[Prop. 5.4.2.4]{lurie-gaitsgory-tamagawabook}) combined with Lemma \ref{lemma: sections of finite morphism finite}, the pullback morphism $\Bun_{K^{\overline{\cP}}}(X_{A}) \to \Bun_{K}(\widetilde{X})$ is finite. Hence, we are reduced to showing that $\Bun_{K}(\widetilde{X}) \cong \prod_{i} \Bun_{\mu_i}(\widetilde{X})$ is quasi-finite and proper. By Lemma \ref{lemma: splitting smooth morphisms into connected components}, after replacing $\Spec(A)$ with a finite surjective cover $\Spec(B) \to \Spec(A)$ and replacing $\widetilde{X}$ with its connected components, we can assume that $\widetilde{X}$ has geometrically integral fibers.

For each $i$, $\Bun_{\mu_i}(\widetilde{X})$ is the $n_i$-torsion substack of the Picard stack $\Bun_{\bG_m}(\widetilde{X})$ (cf. the proof of \cite[Lemma 2.2.1]{biswas-hoffmann-line-bundles}), meaning that it fits into the Cartesian diagram
\[
\begin{tikzcd}
  \Bun_{\mu_{n_i}}(\widetilde{X}) \ar[d] \ar[r] & \Spec(A) \ar[d, "\triv"] \\
  \Bun_{\bG_m}(\widetilde{X}) \ar[r, "\; \; \; \; (-)^{\otimes n_i}"] & \Bun_{\bG_m}(\widetilde{X})
\end{tikzcd}
\]
where $\triv$ denote the trivial $\bG_m$-bundle on $\widetilde{X}$. Let $\Pic_{\widetilde{X}/\Spec(A)}$ denote the Picard group algebraic space with $\tau$-component $\Pic^{\tau}_{\widetilde{X}/\Spec(A)}$ (as in \cite[pg. 233]{blr-neron}). The natural morphism $\Bun_{\mathbb{G}_m}(\widetilde{X}) \to \Pic_{\widetilde{X}/\Spec(A)}$ is a $\mathbb{G}_m$-gerbe.\endnote{Indeed, this can be checked ffpf-locally on the base $\Spec(A)$, so after passing to a flat cover of $\Spec(A)$ we can assume that $\widetilde{X} \to \Spec(A)$ admits a section $s$. In this case the fact that $\Bun_{\mathbb{G}_m}(\widetilde{X}) \to \Pic_{\widetilde{X}/\Spec(A)}$ is a $\mathbb{G}_m$-gerbe follows from the concrete description of $\Pic_{\widetilde{X}/\Spec(A)}$ as a the functor of $\mathbb{G}_m$-bundles $\cL$ on $\widetilde{X}$ equipped with a trivialization of the pullback $s^{*}\cL$ \cite[Lemma 2.9]{kleiman-picarscheme}.} A similar reasoning, using the Cartesian diagram for $\Bun_{\mu_{n_i}}(\widetilde{X})$ above, implies that $\Bun_{\mu_i}(\widetilde{X})$ is a $\mu_{n_i}$-gerbe over the $n_i$-torsion subgroup $\Pic_{\widetilde{X}/\Spec(A)} = \Pic^{\tau}_{\widetilde{X}/\Spec(A)}[n_i]$. By \cite[8.4, Thm. 3 + Thm. 4 (c)]{blr-neron}, the group algebraic space $\Pic^{\tau}_{\widetilde{X}/ \Spec(A)}$ is of finite type over $A$, and moreover it is proper if $X \to S$ is smooth. Therefore, the closed $n_i$-torsion subgroup $\Pic^{\tau}_{\widetilde{X}/\Spec(A)}[n_i]$ is of finite type, and it is proper if $X \to S$ is smooth. In particular, the same holds for the $\mu_i$-gerbe $\Bun_{\mu_i}(\widetilde{X}) \to \Pic^{\tau}_{\widetilde{X}/\Spec(A)}[n_i]$.

We are left to show that $\Pic^{\tau}_{\widetilde{X}/\Spec(A)}[n_i]$ is quasi-finite. We can check this on geometric fibers, so we can assume that $A$ is an algebraically closed field $k$. In this case \cite[Lemma 4.2]{artin-algebraization} implies that the algebraic space $\Pic^{\tau}_{\widetilde{X}/\Spec(k)}$ is a group scheme. The argument in \cite[Thm. 5.4]{kleiman-picarscheme} then shows that $\Pic^{\tau}_{\widetilde{X}/\Spec(k)}$ is proper over $k$. To check quasi-finiteness of the $n_i$-torsion, we can pass to the reduced subgroup scheme $\left(\Pic^{\tau}_{\widetilde{X}/\Spec(k)}\right)_{red}$, which is smooth. Using the short exact sequence induced by the quotient of the neutral component, we see that $\left(\Pic^{\tau}_{\widetilde{X}/\Spec(k)}\right)_{red}$ is an extension of an abelian variety by a finite \'etale group. So the finiteness of the $n_i$-torsion readily follows from the finiteness of $n_i$-torsion for abelian varieties \cite[\href{https://stacks.math.columbia.edu/tag/03RP}{Tag 03RP (5)}]{stacks-project}.\\

\noindent $\bullet$ \textit{Discrete fibers:} Let $k$ be an algebraically closed field over $S$, and let $\Spec(k) \to \Bun_{\overline{G}}(C)$ be a point represented by a $\overline{G}$-bundle $\overline{\cP}$ on $\Bun_{\overline{G}}(X_{k})$. If the fiber of $\Bun_{G}(X_{k}) \to \Bun_{\overline{G}}(X_{k})$ over $\cP$ is nonempty, then it is isomorphic to $\Bun_{K^{\overline{\cP}}}(X_{k})$, where $K^{\cP}$ is the multiplicative $X_{k}$-scheme obtained by twisting $K\times X_{k}$ by $\overline{\cP}$ via the conjugation action of $\overline{G}$ on $K$. We have seen in the course of the proof of quasi-compactness that $\Bun_{K^{\overline{P}}}(X)$ is quasi-finite, so it follows that the fiber over the $k$-point $\overline{P}$ is discrete.\\

\noindent $\bullet$ \textit{Properness:} Assume that $X \to S$ is smooth. We use the valuative criterion for properness. Let $R$ be a complete discrete valuation ring over $S$. Assume that we are given a $\overline{G}$-bundle $\overline{\cP}$ on $X_{R}$ and a lift of $\overline{\cP}|_{X_{\Frac(R)}}$ to a $G$-bundle. We shall extend the lift to a $G$-bundle on $X_{R}$ after perhaps passing to a finite extension of $R$. The bundle $\overline{\cP}$ yields a cohomology class in $H^1_{et}(X_{R}, \Aut(K \times X_{R}))$ via the conjugation action $\overline{G} \to \Aut(K)$. Again, we view this cohomology class as an element $H^1_{et}(X_{R}, \Aut(\prod_{i} \bZ/n_i\bZ))$ with values in the automorphisms of the constant Cartier dual group $\prod_{i} \bZ/n_i\bZ$ of the diagonalizable group $K$. Let $\widetilde{X} \to X_{R}$ be a finite \'etale cover that splits this cohomology class, so that $K^{\overline{\cP}} \times_{X_{R}} \widetilde{X} \cong K \times \widetilde{X}$. Consider the following diagram of Cartesian squares
\[
\begin{tikzcd}
  \Bun_{G}(X_{R}) \times_{\Bun_{\overline{G}}(X_{R})} \Spec(R) \ar[d] \ar[r] & \Bun_{G}(\widetilde{X}) \times_{\Bun_{\overline{G}}(\widetilde{X})} \Spec(R) \ar[d] \ar[r] & \Spec(R) \ar[d, "\overline{\cP}"] \\ \Bun_{G}(X_{R}) \ar[r] & \Bun_{G}(\widetilde{X}) \times_{\Bun_{\overline{G}}(\widetilde{X})} \Bun_{\overline{G}}(X_{R}) \ar[d] \ar[r] & \Bun_{\overline{G}}(X_{R}) \ar[d] \\   & \Bun_{G}(\widetilde{X}) \ar[r] & \Bun_{\overline{G}}(\widetilde{X})
\end{tikzcd}
\]
We conclude the proof of the lemma by showing that both of the top horizontal morphisms are proper.

\medskip
\noindent(1) \textit{$ \Bun_{G}(\widetilde{X}) \times_{\Bun_{\overline{G}}(\widetilde{X})} \Spec(R) \to \Spec(R)$ is proper:} By Lemma \ref{lemma: splitting smooth morphisms into connected components}, after passing to a finite discrete valuation extension of $R$ and replacing $\widetilde{X}$ with its connected components, we can assume that $\widetilde{X}$ has geometrically integral $R$-fibers, and hence is integral. There is an obstruction $\ob \in H^2_{et}(\widetilde{X}, K^{\overline{\cP}}) = H^2_{et}(\widetilde{X}, K) $ to lifting the $\overline{G}$-bundle $\overline{\cP}|_{\widetilde{X}}$ to a $G$-bundle. If $\ob =0$, then the fiber $\Bun_{G}(\widetilde{X}) \times_{\Bun_{\overline{G}}(X)} \Spec(R)$ is isomorphic to $\Bun_{K}(\widetilde{X})$. During the course of the proof of quasi-compactness, we have seen that $\Bun_{K}(\widetilde{X})$ is proper over $R$. Hence, it suffices to show that $\ob$ vanishes, after maybe replacing $R$ with a finite flat extension (since properness can be checked flat-locally). For some positive integer $r$, we can realize the diagonalizable group $K = \prod_i \mu_{n_i}$ as the kernel of a morphism of group schemes
\[ 1 \to K \to \bG_m^r \xrightarrow{(n_i)} \bG_m \to 1.  \]
Consider the corresponding Kummer sequence
\[ H^1_{et}(\widetilde{X}, \bG_m^r) \xrightarrow{(n_i)} H^1_{et}(\widetilde{X}, \bG_m^r) \to H^2_{et}(\widetilde{X}, K) \xrightarrow{\varphi} H^2_{et}(\widetilde{X}, \bG_m^r)\]
$\ob|_{\widetilde{C}_{\Frac(R)}}$ vanishes, because the bundle $\overline{\cP}|_{\widetilde{X}_{\Frac(R)}}$ admits a lift to $G$. Hence, the restriction of the image $\varphi(\ob)|_{\widetilde{X}_{\Frac(R)}}$ vanishes. But by \cite[Cor. 1.8]{grothendieck-brauer} (or \cite[III Ex. 2.22, pg. 107]{milne-etalecohomology}) the restriction morphism $H^2_{et}(\widetilde{X}, \bG_m^r) \to H^2_{et}(\widetilde{X}_{\Frac(R)}, \bG_m^r)$ is injective. We conclude that $\varphi(\ob)$ vanishes, and hence $\ob$ comes from a class in $H^1_{et}(\widetilde{X}, \bG_m^r)$, corresponding to a tuple $(\cL_i)_{i=1}^r$ of line bundles. We are reduced to showing that each $\cL_i$ admits a $n_i^{th}$ root. After replacing $R$ with a finite extension, we can assume that the proper morphism $\widetilde{X} \to \Spec(R)$ admits a section. Hence, the $R$-points of the Picard algebraic space $\Pic_{\widetilde{X}/R}$ are in bijection with isomorphism classes of line bundles on $\widetilde{X}$ \cite[Thm. 2.5]{kleiman-picarscheme}. The existence of the lift $\cP$ on $\widetilde{X}|_{\Frac(R)}$ shows that the restriction $\cL_i|_{\widetilde{X}|_{\Frac(R)}}$ admits a $n_i^{th}$ root. Since $\Pic_{\widetilde{X}/R} \xrightarrow{n_i} \Pic_{\widetilde{X}/R}$ is proper (by \cite[Rem. 9.6.28]{kleiman-picarscheme} and \cite[Thm. 3]{blr-neron}), we conclude by the valuative criterion that $\cL_i$ admits a $n_i^{th}$ root, and so $\ob$ vanishes.

\medskip
\noindent (2) \textit{$\Bun_{G}(X_{R}) \times_{\Bun_{\overline{G}}(X)} \Spec(R) \to \Bun_{G}(\widetilde{X}) \times_{\Bun_{\overline{G}}(\widetilde{X})} \Spec(R)$ is finite (in particular proper):} Choose a point $T \to \Bun_{G}(\widetilde{X}) \times_{\Bun_{\overline{G}}(\widetilde{X})} \Spec(R)$ consisting of a $G$-bundle $\cP$ on $\widetilde{X}_{T}$ such that the associated $\overline{G}$-bundle is isomorphic to $\overline{\cP}|_{\widetilde{X}_{T}}$. Consider the two \'etale projections $p_1, p_2: \widetilde{X}_{T} \times_{X_{T}} \widetilde{X}_{T} \to \widetilde{X}_{T}$. By \'etale descent (cf. \cite[Prop. 5.4.2.4]{lurie-gaitsgory-tamagawabook}), the fiber product $\Bun_{G}(X) \times_{\Bun_{G}(\widetilde{X}) \times_{\Bun_{\overline{G}}(\widetilde{X})} \Spec(R)} T$ is parametrized by the functor of descent data for the cover $\widetilde{X}_{T} \to X_{T}$, which amounts to the subsets of isomorphisms $\psi \in \Iso_{\Bun_{G}(\widetilde{X}_{T} \times_{X_{T}} \widetilde{X}_{T}) / \Bun_{\overline{G}}(\widetilde{X}_{T} \times_{X_{T}} \widetilde{X}_{T})}(p_1^*(\cP), p_2^*(\cP))$ in the relative diagonal of $\Bun_{G}(\widetilde{X}_{T} \times_{X_{T}} \widetilde{X}_{T}) \to \Bun_{\overline{G}}(\widetilde{X}_{T} \times_{X_{T}} \widetilde{X}_{T})$ satisfying a cocycle condition. The same proof as in the case of $X$ shows that the relative diagonal of $\Bun_{G}(\widetilde{X}_{T} \times_{X_{T}} \widetilde{X}_{T}) \to \Bun_{\overline{G}}(\widetilde{X}_{T} \times_{X_{T}} \widetilde{X}_{T})$ is finite, and the subfunctor consisting of elements satisfying the cocycle condition is represented by a closed immersion (see \cite[Prop. 5.4.2.4]{lurie-gaitsgory-tamagawabook}).
\end{proof}

\begin{remark}
The quasi-compactness of the morphism $\Bun_{G}(X) \to \Bun_{\overline{G}}(X)$ would follow directly from \cite[Thm. 1.7]{olsson-boundedness} if $K$ was \'etale (i.e. if the $n_i$ were invertible in $S$ for all $i$).
\end{remark}

\begin{remark}
The smoothness of $X \to S$ is necessary in order to guarantee that the morphism is quasi-finite and proper in Lemma \ref{lemma: finiteness stack of bundles under isogeny}. For example, for a fixed prime $p$ we can consider the trivial isogeny $\mu_p \to 1$ with source the $p^{th}$ roots of unity, and let $X \to \Spec(k)$ to be a projective curve of geometric genus $0$ with one cusp over a field $k$ of characteristic $p$. Then the corresponding stack $\Bun_{\mu_p}(X) \to \Spec(k)$ is isomorphic to $\mathbb{A}^1_k \times B\mu_p$. 
\end{remark}

\begin{lem} \label{lemma: bun_g isogeny orbifold case}
Let $G \to \overline{G}$ be an isogeny of (not necessarily connected) geometrically reductive group schemes over $S$, with kernel $K \hookrightarrow G$ finite and of multiplicative type. Let $X$ be a quotient stack of the form $X=\widetilde{X}/\Gamma$, where $\widetilde{X}$ is a smooth proper scheme over $S$, and $\Gamma$ is a finite \'etale group scheme over $S$. Then, the morphism of stacks $\Bun_{G}(X) \to \Bun_{\overline{G}}(X)$ is quasi-finite and proper.
\end{lem}
\begin{proof}
There are morphism $\Bun_{G}(X) \to \Bun_{G}(\widetilde{X})$ and $\Bun_{\overline{G}}(X) \to \Bun_{\overline{G}}(\widetilde{X})$ induced by pulling back principal bundles under the morphism $\widetilde{X} \to X$. The morphism $\Bun_{G}(X) \to \Bun_{\overline{G}}(X)$ factors as a composition
\[ \Bun_{G}(X) \xrightarrow{f} \Bun_{G}(\widetilde{X}) \times_{\Bun_{\overline{G}}(\widetilde{X})} \Bun_{\overline{G}}(X) \xrightarrow{g} \Bun_{\overline{G}}(X)\]
The morphism $g$ is a base-change of $\Bun_{G}(\widetilde{X}) \to \Bun_{G}(\widetilde{X})$, and so it is quasi-finite and proper by \Cref{lemma: finiteness stack of bundles under isogeny} applied to the proper smooth scheme $\widetilde{X} \to S$. On the other hand, a similar reasoning as in the proof of claim (2) in the proof of properness in \Cref{lemma: finiteness stack of bundles under isogeny} shows that $f$ is finite, and so it follows that the composition $g \circ f: \Bun_{G}(X) \to \Bun_{\overline{G}}(X)$ is quasi-finite and proper.
\end{proof}
We conclude this section with a lemma on the boundedness of the Harder-Narasimhan stratification for $G$-bundles on a smooth family of curves. An analogous statement for a connected reductive group $G$ when the base $S$ is the spectrum of a field was proven in \cite[Thm. 8.2.6]{behrend-thesis}. We elaborate on how to adapt some of the arguments to apply to our more general setup. As in the body of the text, we assume that $C \to S$ is a smooth projective morphism whose fibers are geometrically connected curves. We will need the following lemma, which is a generalization of \cite[Prop. 8.5.2]{behrend-thesis} for a general base scheme.
\begin{lem} \label{lemma: boundedness connected components of bunb}
Let $B$ be an affine $S$-group scheme of the form $B = T \ltimes U$, where $T$ is a split torus over $S$ and $U$ is an iterated extension of copies of $\bG_a$ equipped with a linear $T$ action. Then the natural morphism $\Bun_{B}(C) \to \Bun_{T}(C)$ induced by $B \twoheadrightarrow T$ is of finite type. In particular, the preimage of any connected component of $\Bun_{T}(C)$ is a quasi-compact open substack of $\Bun_{B}(C)$.
\end{lem}
\begin{proof}
The same argument as in the proof of \cite[Prop. 2.5]{herrero2021quasi-compactness} applies in this situation to show that $\Bun_{B}(C) \to \Bun_{T}(C)$ is of finite type. In \cite{herrero2021quasi-compactness} the assumption on the base of characteristic $0$ is only used to conclude that $U$ is an iterated extension of $\mathbb{G}_a$ and that $G$ acts linearly, but this is already included in our assumptions above. The last sentence in the lemma follows from the fact that every connected component of $\Bun_{T}(C)$ is quasi-compact.\endnote{To see this, note that $\Bun_{T}(C)$ can be written as a product of Picard stacks $\Bun_{\bG_m}(C)$, because $T$ is split. So we can reduce to the case when $T = \bG_m$. Since $\Bun_{\bG_m}(C)$ is a $\bG_m$-gerbe over the relative Picard scheme $\Pic_{C/S}$, it suffices to show that each connected component of $\Pic_{C/S}$ is bounded. This follows from \cite[9.3 Thm. 1]{blr-neron}.}
\end{proof}
Fix an almost faithful representation $G \to \GL(V)$ and choose a rational quadratic norm $b$ on graded points of $BG$. We obtain a rational quadratic norm $b$ on $\Bun_{G}(C)$ by restricting bundles to the generic point of fibers of $C$. Equip $\Bun_{G}(C)$ with the numerical invariant $\mu = -\wt(\cD(V)) / \sqrt{b}$. By Theorem \ref{thm: theta stratification in general}, the numerical invariant $\mu$ induces a weak $\Theta$-stratification on $\Bun_{G}(C)$. For any real number $\gamma$, we will denote by $\Bun_{G}(C)^{\mu \leq \gamma}$ the open substack of $\Bun_{G}(C)$ where the Harder-Narasimhan filtration has invariant less than or equal to $\gamma$.
\begin{lem} \label{lemma: boundedness of HN stratification relative case}
For every $\gamma \geq 0$ and every connected component $\cX \subset \Bun_{G}(C)$, the open substack $\Bun_{G}(C)^{\mu \leq \gamma} \cap \cX$ is quasi-compact. In particular, there are finitely many weak $\Theta$-strata inside $\cX$ having bounded numerical invariant, and each such stratum is quasi-compact.
\end{lem}
\begin{proof}
We use the same notation and setup as in \Cref{subsection: hn boundedness for bung}. We shall prove the lemma in increasing levels of generality.

\medskip
\noindent \textit{Proof when $G$ has connected fibers.} We can assume after passing to an \'etale cover of $S$ that $G$ is split over the base and that $S$ is connected. Fix a choice of a Borel subgroup $B \subset G$ containing a split maximal torus $T$. For this special case we shall follow the strategy of proof of \cite[Prop. 8.2.6]{behrend-thesis}, and show that there is a bounded substack $\cU \subset \Bun_{B}(C)$ such that the image of $\cU$ under the natural morphism $\Bun_{B}(C) \to \Bun_{G}(C)$ contains $\Bun_{G}(C)^{\mu \leq \gamma}\cap \cX$.

Define $(-,-)_{V}$ as in \Cref{notation: beginning proof hn-boundedness bun_g}. Choose a geometric point of $\Bun_{G}(C)$, represented by a $G$-bundle on a geometric fiber of $C$. For any Borel reduction $E_{B}$, consider the cone $\sigma_{E_{B}} \subset X_{*}(T)_{\mathbb{R}}$ inside $|\DF(\Bun_{G}(C), E)|$. By \Cref{lemma: riemann-roch expression line bundle}, the numerical invariant on $\sigma_{E_{B}}$ is given by $\mu(w) = (d_{E_{B}}, w)_{V}/\|w\|_b$. Consider the set $\omega_1, \omega_2, \ldots , \omega_r$ of fundamental weights, where $r$ is the rank of $G$. There exists a positive integer $m$ such that for every fundamental weight the multiple $m\omega_j$ of $G$ arises as a character of $T$. By composing with the natural quotient homomorphism $B \twoheadrightarrow T$, we can view $m\omega_j$ as a character of $B$, and similarly for the roots $\alpha_j$. Let $g$ denote the genus of the fibers of $C$ over the connected scheme $S$. Since every fundamental weight can be written as a positive linear combination of simple roots, there exists a constant $c_g$ depending only on $g$ and $G$ such that for all $\psi \in X_{*}(T)_{\mathbb{R}}$ we have
\[ \langle \psi, \alpha_j \rangle \geq -g \; \text{for all $j$} \; \Longrightarrow \; \langle \psi, \omega_j \rangle \geq -c_g \; \text{for all $j$}\]
Set $c_{min} = \min_{1 \leq j \leq r} \left(\|\omega_j\|_V^2\right)$. Define $\cU \subset \Bun_{B}(C)$ to be the open substack parametrizing $B$-bundles $P$ whose associated $G$-bundle lies in $\cX\subset \Bun_{G}(C)$ and such that the degree $d_{E_B}$ satisfies 
\[-c_g \leq \langle d_{E_B}, \omega_j \rangle \leq \frac{1}{c_{min}}\left(\gamma \|\rho_{B}\|_b + c_g (r-1) \| \rho_B\|_V^2\right)\]
for all fundamental weights $\omega_j$. To conclude the lemma for the case with connected fibers, we show the following.
\begin{enumerate}
    \item The stack $\cU$ is bounded.
    \item $\Bun_{G}(C)^{\mu \leq \gamma}$ is contained in the image of $\cU$ under the natural morphism $\Bun_{B}(C) \to \Bun_{G}(C)$.
\end{enumerate}
For (1), note for any $B$-bundle $P$ in $\cU$ there are finitely many possible choices for the degrees of the line bundles $P(m\omega_j)$, since they must be integers and they are bounded between two constants. Hence, $\cU$ is a disjoint union of finitely many open substacks of $\Bun_{B}(C)$ where $\deg(P(m\omega_j)) =d_j$ are fixed. Since we are restricting to a connected component $\cX \subset \Bun_{G}(C)$, the degree $d_{E}$ is fixed, and so we also have that $\deg(P(\chi))= \langle d_{E}, \chi \rangle$ is fixed for any character $\chi$ of $G$. Since the fundamental weights and the characters of $G$ jointly span the character space $X^*(T)_{\mathbb{R}}$, the degrees $\deg(P(m\omega_j))$ determine the multi-degree of the associated $T$-bundle $P \times^{B} T$ under the homomorphism $B \twoheadrightarrow T$ (where for the multi-degree we are viewing $\Bun_{T}(C)$ as parametrizing tuples of line bundles on $C$ after fixing a splitting $T \cong \mathbb{G}_m^n$). Therefore, $\cU$ maps to finitely many connected component of $\Bun_{T}(C)$, and hence by Lemma \ref{lemma: boundedness connected components of bunb} $\cU$ is bounded.

For (2), it suffices to check the claim on geometric fibers. So we can assume that $S = \Spec(k)$ for an algebraically closed field $k$. Let $E$ be a $G$-bundle on $C$ lying on the stack $\Bun_{G}(C)$. We need to find a Borel reduction lying on $\cU$. By (*) in the proof of \cite[Prop. 4.3.1]{lurie-gaitsgory-tamagawabook}, the bundle $E$ admits a Borel reduction $E_{B}$ such that $\langle d_{E_{B}}, \alpha_j\rangle \geq -g$ for all $j$. This implies that $\langle d_{E_{B}},\omega_j \rangle \geq -c_g$ for all $j$. We are left to show that $\langle d_{E_{B}}, \omega_j\rangle \leq \frac{1}{c_{min}}\left(\frac{1}{m}\gamma \|\rho_{B}\|_b + c_g (r-1) \| \rho_B\|_b^2\right)$ for all $j$. Consider the point $\rho_B \in \sigma_{E_{B}}$ in the degeneration fan $|\DF(\Bun_{G}(C), E)|$. By the assumption that $E \in \Bun_{G}(C)^{\mu \leq \gamma}$, we have $\mu(\rho_{B}) \leq \gamma$. Using the formula for $\mu$, we get $\frac{(d_{E_{B}}, \, \rho_{B})_{V}}{\|\rho_{B}\|_b} \leq \gamma$. Using $d_{E_{B}} = \sum_{j=1}^r\langle d_{E_{P}}, \omega_j\rangle \omega_j^{\vee}$, we get
\[ \sum_{j=1}^r \langle d_{E_{B}}, \omega_j\rangle \cdot (\omega_j^{\vee}, \rho_{B})_{V} \leq \gamma \|\rho_{B}\|_b\]
So for each fixed $j$, we can write
\begin{equation} \label{eqn: first inequality boundedness hn bung appendix}
    \langle d_{E_{B}}, \omega_j\rangle \cdot (\omega_j^{\vee}, \rho_{B})_{V} \leq \gamma \|\rho_{B}\|_b - \sum_{i \neq j} \langle d_{E_{P}}, \omega_j\rangle (\omega_j^{\vee}, \rho_{B})_{V}
\end{equation}
By definition $\rho_{B} = \sum_{i=1}^r \omega_i^{\vee}$. Using this and the fact that $(\omega^{\vee}_i, \omega^{\vee}_j)_{V} \geq 0$ for all $i,j$ \endnote{Indeed, since $(-,-)_V$ restricts to a Weyl group invariant positive definite inner product on $X_*(T')_{\mathbb{R}}$, we know that the induced inner product on the dual space $X^{*}(T')_{\bR}$ satisfies $(\alpha_i, \alpha_j)_V \leq 0$ for all $i \neq j$. This implies $(\omega^{\vee}_i, \omega_j^{\vee})_V \geq 0$ for all $i,j$, because $\{\omega^{\vee}_i\}$ is the dual basis in $X_*(T')_{\mathbb{R}}$.}, we see that
\[ \|\omega^{\vee}_i \|_V^2 \leq (\rho_B, \omega_i^{\vee})_{V} \leq \|\rho_B\|_V^2 \; \; \text{for all $i$}\]
Moreover, we have $\langle d_{E_{B}}, \omega_i\rangle \geq -c_g$. We can use these facts to rewrite the right-hand side of the inequality in Equation \eqref{eqn: first inequality boundedness hn bung appendix}
\[\langle d_{E_{P}}, \omega_j\rangle \cdot (\omega_j^{\vee}, \rho_{B})_{V} \leq \gamma \|\rho_{B}\|_b + c_g (r-1) \| \rho_B\|_V^2\]
Using $(\rho_B, \omega^{\vee}_j)_{V} \geq \|\omega^{\vee}_j\|_V^2 \geq c_{min}$, we get our desired inequality
\[ \langle d_{E_{B}}, \omega_j) \leq \frac{1}{(\rho_B, \omega^{\vee}_j)_{V}} \left( \gamma \|\rho_{B}\|_b + c_g (r-1) \| \rho_B\|_V^2\right) \leq \frac{1}{c_{min}}\left(\gamma \|\rho_{B}\|_b + c_g (r-1) \| \rho_B\|_V^2\right)\]
This concludes the case when $G$ has connected fibers.

\medskip
\noindent \textit{Proof for arbitrary geometrically reductive $G$.} After passing to an \'etale cover of $S$, we can assume that $G$ admits a split maximal torus over $S$. Consider the short exact sequence of group schemes
\[ 1 \to G_0 \to G \to \pi_0(G) \to 1 \]
The group $\pi_0(G)$ is a finite \'etale group over $S$ by geometric reductivity. The quotient morphism $G \to \pi_0(G)$ induces a morphism of stacks $\Bun_{G}(C) \to \Bun_{\pi_0(G)}(C)$ via extension of group. By \cite[Thm. 1.7]{olsson-boundedness}, the stack $\Bun_{\pi_0(G)}(C) = \uMap_S(C, B\pi_0(G))$ is of finite type over $S$. We conclude the proof of the lemma by showing that the restriction $\varphi: \Bun_{G}(C)^{\mu \leq \gamma} \cap \cX \to \Bun_{\pi_0(G)}(C)$ is quasi-compact. Let $T$ be an affine Noetherian $S$-scheme and choose a morphism $T \to \Bun_{\pi_0(G)}(C)$ represented by a $\pi_0(G)$-bundle $\cP$ on $C_{T}$. We want to show that the fiber product $W \vcentcolon = U \times_{\Bun_{\pi_0(G)}(C)} (\Bun_{G}(C)^{\mu \leq \gamma} \cap \cX)$ is a quasi-compact stack. The torsor $\cP$ defines a finite \'etale cover $\widetilde{C} \to C_{T}$. Consider the commutative diagram
\[
\begin{tikzcd}
  \Bun_{G}(C_{T}) \ar[d, "\psi_G"] \ar[r] & \Bun_{\pi_0(G)}(C_{T}) \ar[d, "\psi_{\pi_0(G)}"] \\
  \Bun_{G}(\widetilde{C}) \ar[r]  & \Bun_{\pi_0(G)}(\widetilde{C})
\end{tikzcd}
\]
where the vertical arrows $\psi_G$ and $\psi_{\pi_0(G)}$ are defined by pulling back to the cover $\widetilde{C}$. By definition, $\psi_{G}(\cP)$ is canonically identified with the trivial $\pi_0(G)$-bundle on $\widetilde{C}$. By \'etale descent, the morphism $\psi_{G}$ is affine and of finite type (cf. \cite[Prop. 5.4.2.4]{lurie-gaitsgory-tamagawabook}). Therefore, it suffices to show that the image of $W$ in $\Bun_{G}(\widetilde{C})$ is bounded. By \Cref{lemma: splitting smooth morphisms into connected components}, after passing to a finite \'etale cover of $T$ and replacing $\widetilde{C}$ with its connected components, we can assume that $\widetilde{C}$ has geometrically connected fibers over $T$.

The image $\psi_{G}(W)\subset \psi_{G}(\cX)$ is contained in a connected component $\cX'$ of $\Bun_{G}(\widetilde{C})$. Define the numerical invariant $\widetilde{\mu} = -\wt(\cD(V)) / \sqrt{b}$ on $\Bun_{G}(\widetilde{C})$. This numerical invariant $\widetilde{\mu}$ pulls back to $D \cdot \mu$ under the map $\psi_{G}: \Bun_{G}(C_{T}) \to \Bun_{G}(\widetilde{C})$, where $\mu$ is our original numerical on invariant $\Bun_{G}(C_{T})$ and $D$ is the degree of the Galois cover $\widetilde{C} \to C_T$. For any geometric point $x$ in $\Bun_{G}(C_T)$, the unique $\widetilde{\mu}$-maximizing canonical $\Theta$-filtration $f_{max}$ of $\psi_{G}(x)$ is preserved by the covering automorphisms of the \'etale cover $\widetilde{C} \to C_T$, and hence it descends to a $\Theta$-filtration of $\overline{f}_{max}$ of $x$ with $\mu(\overline{f}_{max}) = \frac{1}{D}\widetilde{\mu}(f_{max})$\endnote{In fact it is not difficult to see that $\overline{f}_{max}$ is the canonical filtration of $x$ on $\Bun_{G}(C_T)$, by using the fact that $|\DF(\Bun_{G}(C),x)| \to |\DF(\Bun_{G}(\widetilde{C}), \psi_{G}(x)|$ is injective for the affine morphism $\psi_{G}$. We shall not need this here.}. This shows that $\psi_{G}(\Bun_{G}(C_T)^{\mu \leq \gamma})\subset \Bun_{G}(\widetilde{C})^{\widetilde{\mu} \leq D \cdot \gamma}$. Therefore, we can replace $C$ with $\widetilde{C}$, and assume without loss of generality that the $\pi_0(G)$-bundle $\cP$ representing the morphism $T \to \Bun_{\pi_0(G)}(C)$ is trivial. Then fiber product $T \times_{\Bun_{\pi_0(G)}(C)} \Bun_{G}(C)$ is naturally isomorphic to $\Bun_{G_0}(C_T)$ with its natural morphism to $\Bun_{G}(C_T)$. \Cref{lemma: sections of finite morphism finite} applied to the universal finite associated $G/G_0$-bundle on $C_T \times \Bun_{G}(C_T)$ shows that $\phi : \Bun_{G_0}(C_T) \to \Bun_{G}(C_T)$ is finite. In particular, by \cite[Prop. 3.2.12 (1)]{hl_instability}, for any geometric point $x$ in $\Bun_{G_0}(C_T)$ the morphism $\phi$ induces a homeomorphism of degeneration fans $| \DF(\Bun_{G}(C_T), x)| \to |\DF(\Bun_{G}(C_T), \phi(x)|$. Define a compatible numerical invariant $\mu$ on $\Bun_{G_0}(C_T)$ by considering the representation $G_0 \to G \to \GL(V)$. The homeomorphism of degeneration points and the compatibility of $\mu$ shows that the preimage of $\phi^{-1}(\Bun_{G}(C_T)^{\mu \leq \gamma})$ is contained in $\Bun_{G_0}(C_T)^{\mu \leq \gamma}$. By finiteness, we also have that $\phi^{-1}(\cX)$ is contained in finitely many connected components of $\Bun_{G_0}(C_T)$. Hence, quasi-compactness follows from the case of a group $G_0$ with connected fibers treated above.
\end{proof}

\section{Morphisms representable by Deligne-Mumford stacks and degeneration spaces.}

In our analysis of instability in $\cM_n^G(X)$, we needed a slight enhancement of the theory of degeneration spaces developed in \cite[Sect.~3]{hl_instability}. The main result of this appendix is the following:

\begin{thm} \label{T:degeneration_fan_DM_morphism}
Let $\cX$ and $\cY$ be algebraic stacks, locally of finite type, quasi-separated, with separated inertia and affine automorphism groups relative to a Noetherian base stack $\cB$, and let $\pi : \cX \to \cY$ be a separated morphism that is relatively representable by Deligne-Mumford stacks. Then for any algebraically closed field $k$ over $\cB$ and any point $x \in \cX(k)$, the morphism
\[
\pi_\ast : |\DF(\cX,x)_\bullet| \to |\DF(\cY,\pi(x))_\bullet|
\]
is a closed embedding whose image in any cone of $|\DF(\cY,\pi(x))_\bullet|$ is a finite union of rational simplicial sub-cones. If $\pi$ is affine, then the image of $\pi_\ast$ is locally convex, and if $\pi$ satisfies the existence part of the valuative criterion for properness, then $\pi_\ast$ is a homeomorphism.
\end{thm}

When $\pi$ is representable by algebraic spaces, then this result is \cite[Prop.~3.2.10]{hl_instability}. We will prove the theorem after establishing some initial results.

\begin{lem} \label{L:lift_Gm_homomorphism}
Let $k$ be a field, let $G$ be an affine algebraic group over $k$, and let $G \to (\bG_m)_k$ be a surjective homomorphism with finite \'etale kernel. Then for $n > 0$ a sufficiently divisible integer, there is a commutative diagram of $k$-group homomorphisms
\[
\xymatrix{ & G \ar[d] \\ (\bG_m)_k \ar[r]^{(\bullet)^n} \ar[ur] & (\bG_m)_k }.
\]
\end{lem}
\begin{proof}
Let us denote by $\overline{k}$ the separable closure of $k$. After replacing $G$ with its neutral component, we can assume without loss of generality that $G$ is connected. Since the kernel is \'etale, the sujective morphism $G \to (\mathbb{G}_m)_k$ is an \'etale cover. In particular, $G$ is a one-dimensional smooth algebraic group over $k$, since $(\mathbb{G}_m)_k$ is one-dimensional and smooth. By \cite[Thm. A.1.1]{conrad_reductive}, there exists a torus $T \subset G$ such that $T_{\overline{k}}$ is a maximal split torus of $G_{\overline{k}}$. Since $G_{\overline{k}}$ admits the nontrivial character $G_{\overline{k}} \to (\mathbb{G}_m)_{\overline{k}}$, it is not unipotent. Therefore, by \cite[Cor. 11.5(2)]{borel-linear-algebraic-book}, the linear algebraic group $G_{\overline{k}}$ admits a nontrivial torus, which must be maximal because $\dim(G_{\overline{k}}) =1$. Since all maximal tori are conjugate in $G_{\overline{k}}$ \cite[Cor. 11.3(1)]{borel-linear-algebraic-book}, we conclude that $T_{\overline{k}}$ is nontrivial, and so it must be isomorphic to $(\mathbb{G}_m)_{\overline{k}}$. The composition $T \to G \to \mathbb{G}_k$ has finite kernel, and so it must be an isogeny. This induces an inclusion of $Gal(\overline{k}/k)$-modules $X^*((\mathbb{G}_m)_{\overline{k}}) \hookrightarrow X^*(T_{\overline{k}})$. As abelian groups, both of these modules are isomorphic to $\mathbb{Z}$, and the morphism must be of the form
\[  \mathbb{Z} \to \mathbb{Z}, \; \; \; x \mapsto n \cdot x\]
for some integer $n$, which we take to be positive. Since the Galois group $Gal(\overline{k}/k)$ acts trivially on $X^*((\mathbb{G}_m)_{\overline{k}})$, we see that the image $n \cdot X^*(T_{\overline{k}}) \subset X^*(T_{\overline{k}})$ is fixed by $Gal(\overline{k}/k)$. This implies that $Gal(\overline{k}/k)$ acts trivially on $X^*(T_{\overline{k}})$, because $X^*(T_{\overline{k}})$ is torsion-free. Therefore, the torus $T$ is isomorphic to the multiplicative group $(\mathbb{G}_m)_k$. Using this isomorphism we get our desired commutative diagram
\[
\begin{tikzcd}
  T \ar[r] & G \ar[d] \\
  (\mathbb{G}_m)_k \ar[u, symbol = \xrightarrow{\sim}] \ar[r, "(\bullet)^n"]  & (\mathbb{G}_m)_k
\end{tikzcd}
\]
\end{proof}

\begin{lem} \label{L:finite_component}
Let $k$ be a field, $A$ a $\bZ^n$-graded $k$-algebra, and $\Spec(A) \to \bA^n_k$ a $\bG_m^n$-equivariant finite type quasi-finite morphism. Then there is a unique decomposition $\Spec(A) = \Spec(A') \bigsqcup \Spec(B)$ such that $\Spec(A') \to \bA^n_k$ is finite, and the image of $\Spec(B) \to \bA^n_k$ does not contain $(0,\ldots,0)$.
\end{lem}
\begin{proof}
Note that $\Spec(A)$ has finitely many orbits, so the good moduli space $\Spec(A^{\bG_m^n})$ is a finite $k$-scheme. It suffices to prove the claim for each connected component of $\Spec(A)$. So, we may assume that $\Spec(A^{\bG_m^n})$ is local artinian, or equivalently that there is a unique closed $\bG_m^n$-orbit. If $\Spec(A) \to \bA^n_k$ is finite, then $(0,\ldots,0)$ lies in the image, because every orbit in $\bA^n_k$ specializes to $(0,\ldots,0)$, and in this case the closed orbit in $\Spec(A)$ must map to $(0,\ldots,0)$. So, it will suffice to show that if the unique closed orbit in $\Spec(A)$ maps to $(0,\ldots,0)$, then $\Spec(A) \to \bA^n_k$ is finite.

It suffices to prove this for each irreducible component of $\Spec(A)$ separately, so we may assume that $A$ is an integral domain. Also, the image of the homomorphism $k[x_1,\ldots,x_n] \to A$ must be an integral $\bZ^n$-graded quotient, which therefore must be of the form $k[x_1,\ldots,x_n]/(x_{i_1},\ldots,x_{i_r})$ for some subset of the generators. We can therefore replace $k[x_1,\ldots,x_n]$ with this quotient, and assume that $k[x_1,\ldots,x_n] \to A$ is injective. Because $A$ is integral, $A^{\bG_m^n} = k'$ for some field extension $k'$, and $A = k' \cdot 1 \oplus \mathfrak{m}$, where $\mathfrak{m}$ is the ideal of the unique $\mathbb{G}_m^n$-fixed closed point. Any product of homogeneous elements of $\mathfrak{m}$ with total weight $0$ must vanish, or else it would be a unit. So the hypothesis that $A$ is integral implies that the weights of non-zero homogeneous elements of $\mathfrak{m}$,  and thus of $A$, must lie in a strictly convex cone in $\bZ^n$. This cone must contain the cone spanned by the weights of the $x_i$, because $k[x_1,\ldots,x_n] \subset A$. Therefore, we can find some cocharacter of $\bG_m^n$ such that $\mathfrak{m}$ and $(x_1,\ldots,x_n)$ are all strictly positively graded with respect to this cocharacter.

The claim now follows from a graded variant of Nakayama's lemma, which states that if $x_1,\ldots,x_n$ are given positive degree, and $M \to N$ is a map of graded $k[x_1,\ldots,x_n]$-modules such that both $M$ and $N$ have bounded-below gradings, then $M \to N$ is surjective if and only if $M \otimes_{k[x_1,\ldots,x_n]} k \to N \otimes_{k[x_1,\ldots,x_n]} k$ is surjective. In particular, a $k[x_1,\ldots,x_n]$-module $M$ with bounded-below gradings is finitely generated if and only if $M \otimes_{k[x_1,\ldots,x_n]} k$ is finite dimensional. By hypothesis $A \otimes_{k[x_1,\ldots,x_n]} k$ is a finite $k$-algebra, so this shows that $A$ is module-finite over $k[x_1,\ldots,x_n]$.
\end{proof}

\begin{lem} \label{L:coarse_moduli_fans}
Let $k$ be an algebraically closed field, let $\cX \to \Theta_k^n$ be a morphism of algebraic stacks that is relatively representable by Deligne-Mumford stacks and that is smooth locally on the source a coarse moduli space morphism, and let $x \in \cX(k)$ be a lift of $(1,\ldots,1) \in \Theta_k^n$. If there is a section $\Theta_k^n \to \cX$ that sends $(1,\ldots,1) \mapsto x$, then it is unique up to unique isomorphism, and such a section exists after base change along a morphism $\Theta_k^n \to \Theta_k^n$ given in coordinates by $(z_1,\ldots,z_n) \mapsto (z_1^{r_1},\ldots,z_n^{r_n})$ for some choice of integers $r_i>0$.
\end{lem}
\begin{proof}
For any two sections $\sigma_1,\sigma_2 : \Theta_k^n \to \cX$, an isomorphism over the open substack $(\bA^1_k \setminus \{0\})^n / \bG_m \cong \{(1,\ldots,1)\}$ extends uniquely to an isomorphism of sections over $\Theta_k^n$, because the relative diagonal $\cX \to \cX \times_{\Theta_k^n} \cX$ is finite. This implies the uniqueness of a section, if it exists.

Let $y \in \cX(k)$ be the unique $k$-point lying over $(0,\ldots,0) \in \Theta_k^n(k)$, and let $G_y$ be its automorphism group. Then because $\cX \to \Theta_k^n$ is representable by Deligne-Mumford stacks, $G_y \to (\bG_m^n)_k$ has finite \'etale kernel $F \hookrightarrow G_y$.

We first claim that after base change along a morphism $\Theta_k^n \to \Theta_k^n$ as in the statement of the lemma, we can assume that the homomorphism $\pi: G_y \to (\bG_m^n)_k$ splits. To see this, consider the fiber product $G_i := G_y \times_{(\bG_m^n)_k} (\bG_m)_k$ along the cocharacter $\bG_m \to \bG_m^n$ give by $(1,\ldots,1,z,1,\ldots,1)$ with $z$ in the $i^{th}$ coordinate. \Cref{L:lift_Gm_homomorphism} implies that after composition with the homomorphism $(\bullet)^{r_i}$ for some $r_i>0$, this cocharacter lifts to a cocharacter $\psi_i : (\bG_m)_k \to G_y$. Now for any $i,j$, the map $\psi_i(z)\psi_j(w) \psi_i(z)^{-1}\psi_j(w) : (\bG_m^2)_k \to G_y$ maps to the identity in $(\bG_m^n)_k$, and hence lands in $F$. Because $F$ is finite and $\bG_m^2$ is reduced and connected, this implies that this commutator is $1$, hence the cocharacters $\psi_i$ and $\psi_j$ pairwise commute. Together they define a homomorphism $(\bG_m^n)_k \to G_y$ whose composition with $\pi$ is $(z_1^{r_1},\ldots,z_n^{r_n})$. After pulling back under the morphism $\Theta_k^n \to \Theta_k^n$ defined by the same expression in coordinates, the homomorphism $(\bG_m^n)_k \to G_y$ just defined induces a section of $\pi$.

We think of $\cX$ as a quotient of the Deligne-Mumford stack $X:= \cX \times_{\Theta_k^n} \bA^n_k$ by the action of $(\bG_m^n)_k$. Applying \cite[Thm.~4.1]{alper-hall-rydh-1}, one can find a $\bG_m^n$-equivariant \'etale neighborhood $\Spec(A) \to X$ of $y$. By \Cref{L:finite_component}, after replacing $\Spec(A)$ with the connected component that meets the fiber over $y$, we can assume that the equivariant morphism $\Spec(A) \to \bA^n_k$ is finite.

The image of $\Spec(A) \to \bA^n_k$ is open, $\bG_m^n$-equivariant, and contains $0$, hence it is surjective. We may therefore choose a point $z \in \Spec(A)(k)$ mapping to $(1,\ldots,1) \in \bA^n_k(k)$. \cite[Cor.~1.3.4]{hl_instability} implies that there is a unique section of the morphism $\Spec(A)/(\bG_m^n)_k \to \Theta_k^n$ identifying $(1,\ldots,1)$ with $z$, and composition with the morphism $\Spec(A)/\bG_m^n \to \cX$ gives the desired section of the morphism $\cX \to \Theta_k^n$.

\end{proof}

\begin{prop} \label{P:coarse_moduli_degeneration_fan}
Let $\pi: \cX \to \cY$ be a morphism of algebraic stacks that is relatively representable by Deligne-Mumford stacks, and such that for any smooth morphism $U \to \cY$ from a scheme $U$, the base change $U \times_\cY \cX \to U$ is a coarse moduli space. Then for any algebraically closed field $k$ and any point $x \in \cX(k)$, the morphism
\[
|\DF(\cX,x)_\bullet| \to |\DF(\cY,\pi(x))_\bullet|
\]
is a homeomorphism.
\end{prop}
\begin{proof}
The argument is modeled on the proof of \cite[Prop..~3.2.10]{hl_instability}. One first shows that the morphism of fans $\DF(\cX,x)_\bullet \to \DF(\cY,\pi(x))_\bullet$ is bounded \cite[Def.~3.1.15]{hl_instability}. Given a morphism from the representable fan $h_{[n]} \to \DF(\cY,\pi(x))_\bullet$, which corresponds to a non-degenerate morphism $f : \Theta_k^n \to \cY$ and isomorphism $\pi(x) \cong f(1,\ldots,1)$, the fiber product of fans $h_{[n]} \times_{\DF(\cY,\pi(x))_\bullet} \DF(\cX,x)_\bullet$ is isomorphic to $\DF(\cZ, z)_\bullet$, where $\cZ := \cX \times_\cY \Theta_k^n$ with its marked point $z \in \cZ(k)$ induced by the isomorphism $\pi(x) \cong f(1,\ldots,1)$. The morphism $\cZ \to \Theta_k^n$ satisfies the hypotheses of \Cref{L:coarse_moduli_fans}, and the conclusion of that lemma implies that $\DF(\cZ,z)_\bullet \subset h_{[n]}$ is a sub-fan such that any cone of $h_{[n]}$ is contained in $\DF(\cZ,z)_\bullet$ after positively scaling its ray generators. In particular, $|\DF(\cZ,z)_\bullet|$ is covered by a single cone that lifts a positive multiple of the canonical $n$-cone in $h_{[n]}$. Thus $\DF(\cX,x)_\bullet \to \DF(\cY,\pi(x))_\bullet$ is bounded, and \cite[Lem.~3.1.19]{hl_instability} implies that the morphism $|\DF(\cX,x)_\bullet| \to |\DF(\cY,\pi(x))_\bullet|$ is a closed map. \Cref{L:coarse_moduli_fans} also implies that $\DF(\cX,x)_\bullet \to \DF(\cY,\pi(x))_\bullet$ is injective and $|\DF(\cX,x)_\bullet| \to |\DF(\cX,\pi(x))_\bullet|$ is surjective. The claim follows from the fact that a closed bijection is a homeomorphism.
\end{proof}

\begin{proof}[Proof of \Cref{T:degeneration_fan_DM_morphism}]
The construction of coarse moduli spaces for Deligne-Mumford stacks is functorial and commutes with smooth base change, so one can factor the morphism $\pi$ as a composition $\cX \to \cX' \to \cY$ where $\cX \to \cX'$ is a morphism such that for any morphism from a scheme $U \to \cX'$, the base change $\cX \times_{\cX'} U \to U$ is the coarse moduli space for the Deligne-Mumford stack $\cX \times_{\cX'} U$. \Cref{P:coarse_moduli_degeneration_fan} implies that $|\DF(\cX,x)_\bullet| \to |\DF(\cX',x)_\bullet|$ is a homeomorphism, so it suffices to prove the claim for the morphism $\cX' \to \cY$, which is representable by algebraic spaces. This follows from \cite[Prop.~3.2.10]{hl_instability}.
\end{proof}

\section{Equivariant \texorpdfstring{$\Theta$}{Theta}-stratifications}
\label{appendix:equivariant}

For this section, we work over a fixed Noetherian base scheme $S$. Our algebraic stacks will be assumed to be locally of finite type, quasi-separated, with separated inertia and affine automorphism groups relative to $S$. We let $G$ be a flat and finitely presented separated $S$-group scheme with affine fibers.

By an action of $G$ on a stack $\cX$, we mean a stack $\cX'$ along with a Cartesian diagram
\[
\xymatrix{\cX \ar[r] \ar[d] & \cX' \ar[d] \\
S \ar[r] & BG}.
\]
In this case we denote $\cX/G := \cX'$. By definition, a $G$-equivariant morphism $\cX \to \cY$ between two stacks with $G$-action is a morphism $\cX/G \to \cY/G$ relative to $BG$.

\begin{lem} \label{L:equivariant_appendix_1}
Let $\cX$ be an algebraic stack over $S$ as above. Then there is a Cartesian diagram
\[
\xymatrix{
\Filt(\cX)/G \ar[r]^q \ar[d] & \Filt(\cX/G) \ar[d] \\
BG \ar[r]^-{\substack{trivial\\filtration}} & \Filt(BG)
},
\]
\end{lem}
\begin{proof}
A $T$-point of $\Filt(\cX/G)$ is a morphism $\Theta_T \to \cX/G$, which amounts to a $G$-bundle $\cE \to \Theta_T$ and a $G$-equivariant map $\cE \to \cX/G$. 
 
On the other hand, for the natural action of $G$ on $\Filt(\cX)$ induced by a $G$-action on $\cX$, a $T$-point of $\Filt(\cX)/G$ is a $G$-bundle $E \to T$ along with a $G$-equivariant map $E \times_T \Theta_T \to \cX$.  This corresponds to the data of a $T$-point of $\Filt(\cX/G)$ along with an isomorphism of the associated $G$-bundle $\cE \to \Theta_T$ with $E \times_T \Theta_T$ for some $G$-bundle $E$ on $T$, as claimed in the lemma.

\end{proof}

For any algebraic stack $\cY$ as above, the map $\cY \to \Filt(\cY)$ classifying the trivial filtration is an open and closed immersion. Indeed, by \cite[Thm.~8.2]{alper-hall-rydh-2} for any $T$, $\Map(T,\cY) \subset \Map(\Theta_T,\cY)$ is the full subgroupoid of maps $f :\Theta_T \to \cY$ such that for any point $t \in T$ the homomorphism $(\bG_m)_{k(t)} \to \Aut_\cY(f_t(0))$ is trivial, and this is an open and closed condition on $T$ by \cite[Prop.~1.3.9]{hl_instability}. Applying this to $\cY=BG$, the lemma above implies that $\Filt(\cX)/G \to \Filt(\cX/G)$ is an open and closed substack.

\begin{defn}
We say that a (weak) $\Theta$-stratum $\cS \subset \Filt(\cX)$ is $G$-equivariant if it is the preimage $q^{-1}(\cS_G)$ of an open and closed substack $\cS_G \subset \Filt(\cX)/G \subset \Filt(\cX/G)$. Similarly, a (weak) $\Theta$-stratification is one in which every stratum is $G$-equivariant.
\end{defn}

Note that by faithfully flat descent, the substack $\cS_G$ is unique if it exists. Also, in an equivariant (weak) $\Theta$-stratification, every open subset $\cX_{\leq c} \subset \cX$ is the preimage of an open substack of $\cX/G$.
\begin{lem} \label{L:all_strata_equivariant}
If $G$ has connected fibers over $S$, any (weak) $\Theta$-stratification of $\cX$ is $G$-equivariant.
\end{lem}
\begin{proof}
The map that takes $\cS \subset \Filt(\cX)$ to $q(\cS) \subset \Filt(\cX)/G$ gives a bijection on connected components with inverse $\cS_G \mapsto q^{-1}(\cS_G)$, because $q$ is a $G$-bundle and hence has connected geometric fibers.
\end{proof}

\begin{lem}
The assignment that takes an open and closed substack $\cS_G \subset \Filt(\cX/G)$ to its preimage in $\Filt(\cX)$ induces a bijection between $G$-equivariant (weak) $\Theta$-stratifications of $\cX$ and (weak) $\Theta$-stratifications of $\cX/G$ such that for every HN filtration of a geometric point in $\cX/G$, the filtration in $BG$ induced by the morphism $\cX/G \to BG$ is trivializable.
\end{lem}
\begin{proof}
By \Cref{L:equivariant_appendix_1}, the condition on HN filtrations is precisely the condition that every stratum in $\cX/G$ lies in $\Filt(\cX)/G \subset \Filt(\cX/G)$. So the claim amounts to showing that $\cS_G \subset \Filt(\cX)/G \subset \Filt(\cX/G)$ is a (weak) $\Theta$-stratum in $\cX/G$ if and only if $\cS:= q^{-1}(\cS_G)$ is a (weak) $\Theta$-stratum in $\cX$. This follows from faithfully flat descent, and the commutativity of the following diagram, in which the left-most square is Cartesian
\[
\xymatrix{
\Filt(\cX) \ar[r] \ar[d]^{\ev^\cX_1} & \Filt(\cX)/G \ar@{}[r]|-{\subset} \ar[d]^{\ev^\cX_1 / G} & \Filt(\cX/G) \ar[d]^{\ev^{\cX/G}_1} \\
\cX \ar[r] & \cX/G \ar@{=}[r] & \cX/G
}.
\]
\end{proof}

\printendnotes

\bibliographystyle{alpha}
\bibliography{theta_gauged_maps.bib}

\end{document}